\documentclass{amsart}
\usepackage[top=20mm, left=20mm, right=20mm, bottom=20mm]{geometry}

\usepackage{hyperref} 

\usepackage{amsmath, amssymb, amsthm, amsfonts, mathrsfs, euscript} 

\usepackage{longtable} 
\usepackage{graphicx, float} 
\graphicspath{{img/}, {img/graph/}} 

\newcommand {\rr} {\mathbb{R}}

\newcommand {\al} {\alpha}

\newcommand {\Ga} {\Gamma}

\newcommand {\fy} {\varphi}

\newcommand {\IN} {{\subset}}

\newcommand {\mmm} {{\setminus}}

\newcommand {\0} {{\varnothing}}

\newcommand{\eA}{{\EuScript A}}

\newcommand{\eK}{{\EuScript K}}

\newcommand{\eP}{{\EuScript P}}
\newcommand{\eG}{{\EuScript G}}

\newcommand{\eS}{{\EuScript S}}
\newcommand{\mD}{{\mathcal D}}

\newtheorem{thm}{\bf Theorem}

 \newtheorem{prop}[thm]{\bf Proposition}
\newtheorem{dfn}{\bf Definition}

\title{On Bi-Lipschitz classification of fractal cubes possessing one-point intersection property.}
\author[]{A.~Tetenov$^{1,2,3}$ \and M.~Chanchieva$^{3}$ \and D.~Drozdov$^{1,2}$ \and D.~Rahmanov$^{3}$ \and V.~Safonova$^{3}$ \and I.Udin$^{3}$ \and A. Vetrova$^{3}$\vspace{0.5cm}\\
    \small{$^{1}$Mathematical Center at Akademgorodok\\
    $^{2}$Novosibirsk State University\\
    $^{3}$Gorno-Altaisk State University}}

\date{}

\begin{document}
\maketitle

\begin{abstract}

We show that there are 105 isometric classes of fractal cubes of order 3 with 7 copies  possessing one-point intersection property which are
subdivided to 5 bi-Lipschitz classes of dendrites and  7 bi-Lipschitz classes of non-dendrites.
\end{abstract}

\footnotetext{The work was supported by the Ministry of Science and Higher Education of Russia (agreement No. 075-02-2022-884)}


\section*{Introduction}

The paper presents the results of work of a project team "The topology of fractals and geometry of manifolds" during summer online Mathematical Workshop in July-August 2021, whose sessions were conducted on a base of Regional Scientific and Educational Mathematical Center of Tomsk State University. The aim of our work was to work out the methods of classification of fractal cubes up to Bi-Lipschitz equivalence.\\

{\bf Fractal cubes.} We begin with a definition of a fractal cube in the space $\rr^k$, where  $k\ge 2$.

\begin{dfn}Let $n\ge 2$ and let $\mD=\{d_1,\cdots,d_N\}\subset\{0,1,\dots,n-1\}^k$.
We call the set $\mD$  a digit set. A compact set $K\IN \rr^k$, which satisfies the equation  
\begin{equation}\label{frcube0}
nK=K+\mD. 
\end{equation}
is called a \emph{fractal $k$-cube} of order $n$ ith digit set $\mD$ and $N$ pieces.\end{dfn}

The  digit set    $\mD$  generates  a system of contraction similarities $\eS=\{S_j(x)=\frac{1}{n}(x+d_j)\}_{j=1}^N$ in $\rr^k$. The equation \eqref{frcube0} is equivalent to the  equation
\begin{equation}\label{fractalcube}
K=\bigcup\limits_{j=1}^N S_j(K)
\end{equation}
which means that  the {fractal $k$-cube} $K$ is the attractor of the system $\eS$ \cite{Hata}.
We denote by $T_\mD$ the Hutchinson operator of the system $\eS$, $T_\mD(A):=\dfrac{A+\mD}{n}$.
The images $K_j:=S_j(K)$ are called the pieces of the fractal k-cube $K$.\\
Let $P^k=[0,1]^k$ be the unit cube  in $\rr^k$.
If this does not lead to ambiguity, further we write $P$, omitting $k$.
Since $T_\mD(P)=\bigcup\limits_{j=1}^N S_j(P)\IN P $, each fractal $k$-cube $K$  is contained in $P $ and each piece $K_j$ is contained in $P_j:=S_j(P k)$. 
The faces of the unit cube are defined by
vectors $\al\in\eA$, where $\eA=\{-1,0,1\}^k$, by the formula $P_\al=P\cap (P+\al)$.

In our work we consider  the case k=3, so  whenever we mention  fractal cubes, we mean fractal 3-cubes. In  cases when  $k=2$ and $k=1$, we call a fractal k-cube $K$ a \emph{fractal square} and  a \emph{fractal segment} respectively.\\
\bigskip

{\bf One-point intersection property and its influence.} Our methods are applied to connected fractal cubes of order 3 which possess   one-point intersection property.

\begin{dfn}
Let $\eS=\{S_1,...,S_m\}$ be a system of contracting similarities in $\rr^k$ with attractor $K$. The 
system  $\eS$ possesses {\em one point intersection property,} if for any $i,j\in\{1,...,m\},i\neq j$ the intersection $K_i\cap K_j$, is either empty set or a singleton $\{p\}$.
\end{dfn}

If the system $\eS$ possesses   one point intersection property, the topological structure of the attractor $K$ is defined by the {\em bipartite intersection graph} of the system $\eS$.

\begin{dfn}

 Let $\eK=\{K_1,...,K_m\}$ and $\eP$ be the set of intersection points $\{p\}=K_i\cap K_j$. The {\em bipartite intersection graph}  $\Ga(\eS)$ is a bipartite graph  $(\eK,  \eP; E)$ with parts $\eK$ and $\eP$,  for which an edge  $\{K_i, p\}\in E$ iff  $p\in K_i$ \end{dfn}
 We call $K_i\in \eK$ {\em white vertices}  and $p\in \eP$ -- {\em black vertices}
 of the graph $\Ga$.

 The main tool for detecting the dendrites among the fractal cubes is the following theorem \cite{FProp,SSS2}
 
 \begin{thm}\label{dtest}
Let $\eS$ be a system of injective contraction maps in a complete metric space $X$   which possesses finite intersection property.     The attractor   $K$ of the system $\eS$ is a dendrite if and only if the bipartite intersection graph  of the system $\eS$ is a tree. 
\end{thm} 
 
Nevertheless, if for any of the points $p$ of the set $\eP$ the number of  $K_i\in\eK$ such that $p\in K_i$ is equal to 2, then the bipartite intersection graph is defined uniquely by usual intersection graph $\eG(\eS)$ of the system $\eS$ whose vertex set is $\eK$ and the set $E$ of the edges of  $\eG(\eS)$ is defined by $E=\{(K_i,K_j): K_i\cap K_J\neq\0\}$.

Since all fractal cubes of order 3 with 7 copies which possess one-point intersection property, have no triple points, it is sufficient to use 
the graph $\eG(\eS)$ to determine the properties of such fractal cubes.

Adapting the approach, developed for polygonal and polyhedral systems \cite{STV,STV2,FProp} to fractal cubes with one point intersection property, we  get the following statement:
\begin{prop}
 Each fractal cube with one-point intersection property is the attractor of a polyhedral system and is a continuum of bounded turning.
\end{prop}
Combining  this statement with the results of \cite{GalTet}, and taking into account that the contraction ratios of all mappings $S_i$ are equal to $1/3$, we obtain the following Theorem.
\begin{thm}\label{LipIsom}
Let $K,K'$ be fractal cubes of order 3 which possess one point intersection property. There is a bi-Lipschitz isomorphism $\fy:K\to K'$ if and only if the bipartite intersection graphs
$\Ga(\eS)$ and $\Ga(\eS')$ are  isomorphic.
\end{thm}
\bigskip

{\bf Finding the intersection graph. } For a given  fractal cube $K$ we consider its digit set $\mD=\{d_1,...,d_m\}$ along with the system of homotheties $\eS=\{S_1,...,S_m\}$, where 
$S_i(x)=(x+d_i)/3$ and Hutchinson operator
$T_\mD(A)=\dfrac{A+\mD}{3}$. 
As it was proved in \cite{FQIn}, the intersection of two pieces $K_i$ and $K_j$ may be non-empty only if the vector $\al_{ij}=d_i-d_j$ lies in the set $\eA=\{-1,0,1\}^3\mmm\{0\}$. If $\al\in\eA$ and $d_i-d_j=\al$, then  $K_i\cap K_j=\frac13 (F_\al+d_i)$, where $F_\al=K\cap (K+\al)$ depends only on $\al$. The set $F_\al$ can be represented as the intersection $K_\al\cap (K_{-\al}+\al)$, where $K_\al=K\cap P_\al$ is the $\al$-face of the fractal cube $K$. 
Therefore the fractal cube $K$ has one-point intersection property if for any $\al\in\eA$,
either (1) $\#F_\al\le 1$ or
(2) $\#F_\al> 1$ and $\{(i,j):d_i-d_j=\al\}=\0$.
Provided (2) holds, to construct the usual intersection graph for $K$ one has to take a vertex set $V=\{1,...,m\}$ and edge set
$E=\{(i,j):d_i-d_j=\al \& \#F_\al=1 \}$

{\bf The search procedure.} Using the Finder tool of the program IFSTile \cite{ifst} we created the list of all fractal cubes of order 3 with 7 components, which contained 888030 sets.

The computation tools of the IFStile package allow to filter out  only those elements of the list which are connected and which possess one-point intersection property. 

This gave us a set of 3200 elements, which were divided to 105 isometric classes of fractal cubes. The elements of these classes are equivalent under the symmetries of unit cube $P^3$, maximal number of elements of a class being equal up to 48 elements in the case when its elements have no symmetries. 

Taking a unique representative $K$ from each class, we get a family of 105 non-equivalent elements.
For each of these elements we construct its
usual intersection graph $\eG(\eS)$. If the graph is a tree, the fractal cube $K$ is a dendrite. Collecting the fractal cubes having isomorphic graphs, we get the classes of fractal cubes which are bi-Lipschitz isomorphic,
according to Theorem \ref{LipIsom}.

The results are collected in the tables below.


\section*{Tables of fractal cubes}

\begin{table}[H]
    \centering
    \caption{All possible types of intersection graphs $\eG(\eS)$ for fractal 3-cubes with 7 pieces and number $N$ of isometric types for each graph.}
    \label{tab:graphs}
    \begin{tabular}{|p{0.19\textwidth}|p{0.19\textwidth}|p{0.19\textwidth}|p{0.19\textwidth}|p{0.19\textwidth}|}
        \hline
        \multicolumn{5}{|c|}{\Large Dendrites}\\
        \hline
        \includegraphics{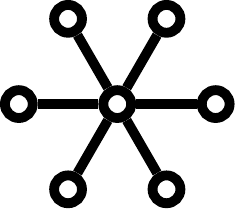}\qquad $7_{11}$ &
        \includegraphics{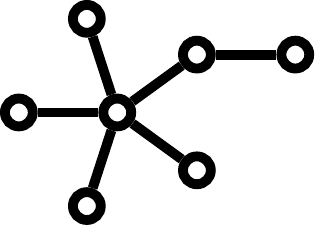} $7_{10}$&
        \includegraphics{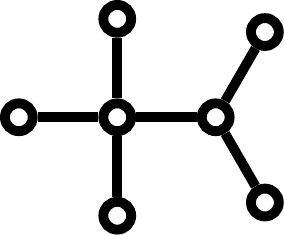}\quad $7_{9}$&
        \includegraphics{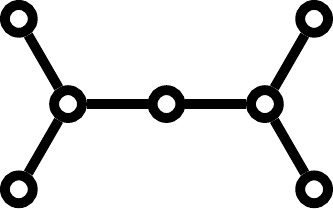} $7_{5}$&
        \includegraphics{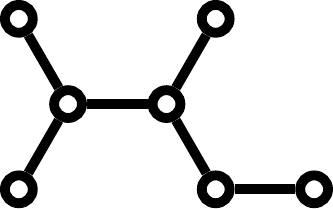} $7_{6}$\\
        \hline
       N=19 & N=3 & N=12 & N=3 &N=1 \\
        \hline
        \end{tabular}\bigskip
        
        \begin{tabular}{|p{0.24\textwidth}|p{0.24\textwidth}|p{0.24\textwidth}|p{0.24\textwidth}|}
        \hline
        \multicolumn{4}{|c|}{\Large Non dendrites}\\
        \hline
        \includegraphics{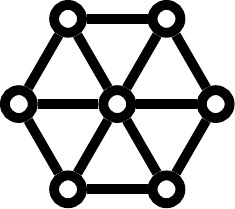} &
        \includegraphics{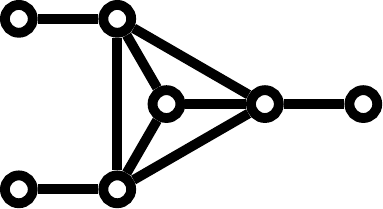} &
        \includegraphics{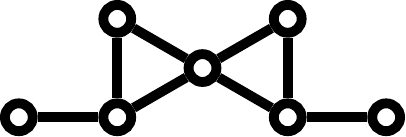} &
        \includegraphics{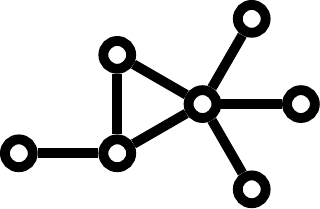} \\
        \hline
        N=3 &N=8 & N=17 & N=25  \\
        \hline
        \includegraphics{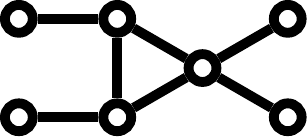} &
        \includegraphics{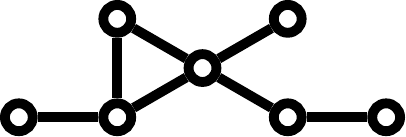} &
        \includegraphics{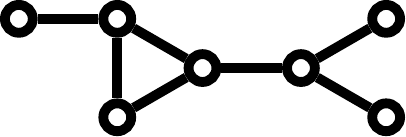} &\\
        \hline
        N=2 &  N=9 & N=3 &\\
        \hline   
    \end{tabular}
    
\end{table}
\bigskip

\begin{longtable}{|p{0.48\textwidth}|p{0.48\textwidth}|}
\caption{Dendrites corresponding to graph $7_{11}$, their sets $T(\mD)$ and digit sets}\label{tab:d1}\\
    \hline
    \multicolumn{2}{|c|}{\includegraphics{den1.pdf}} \\
    \hline
    \includegraphics[width=0.23\textwidth] {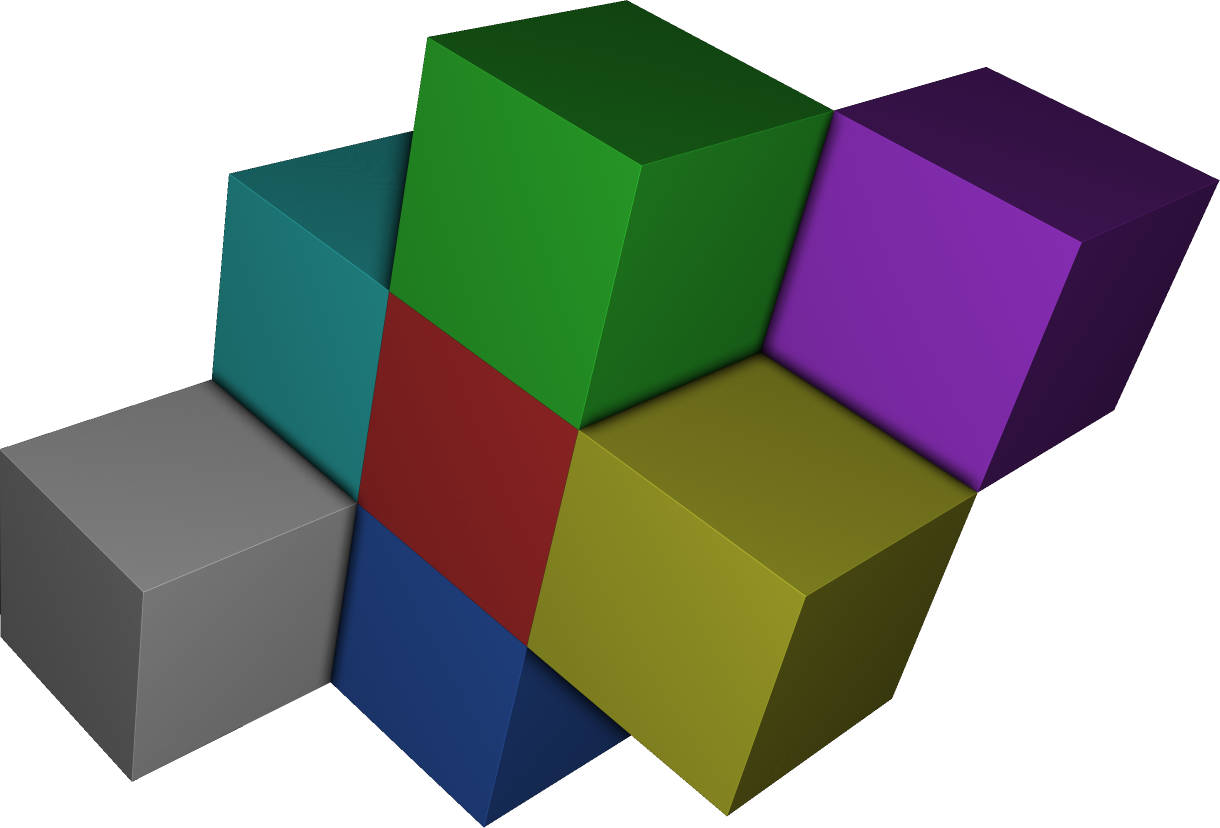}
    \includegraphics[width=0.23\textwidth] {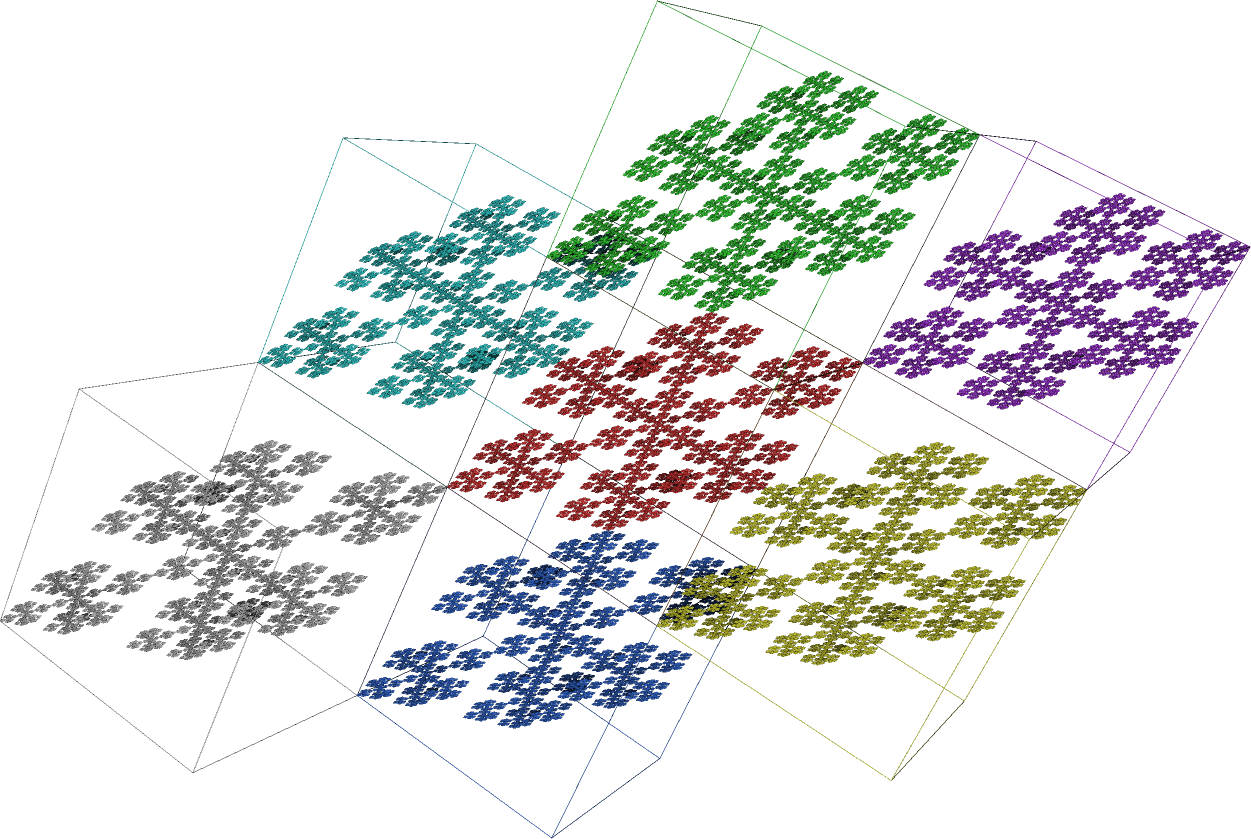} & 
    \includegraphics[width=0.23\textwidth] {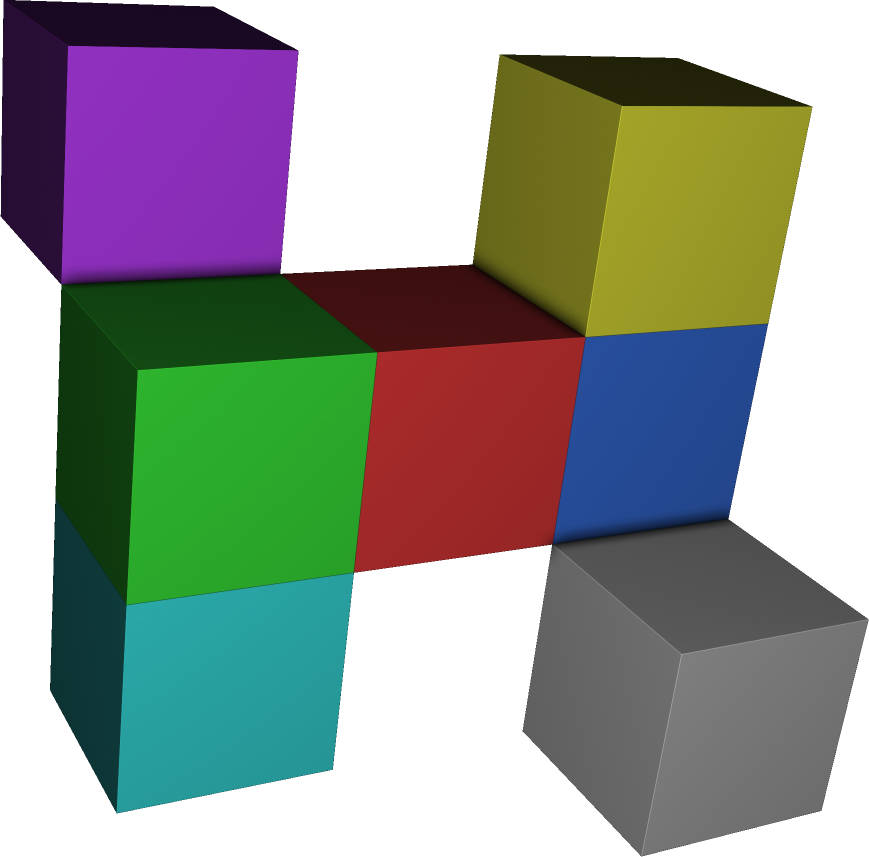} \includegraphics[width=0.23\textwidth] {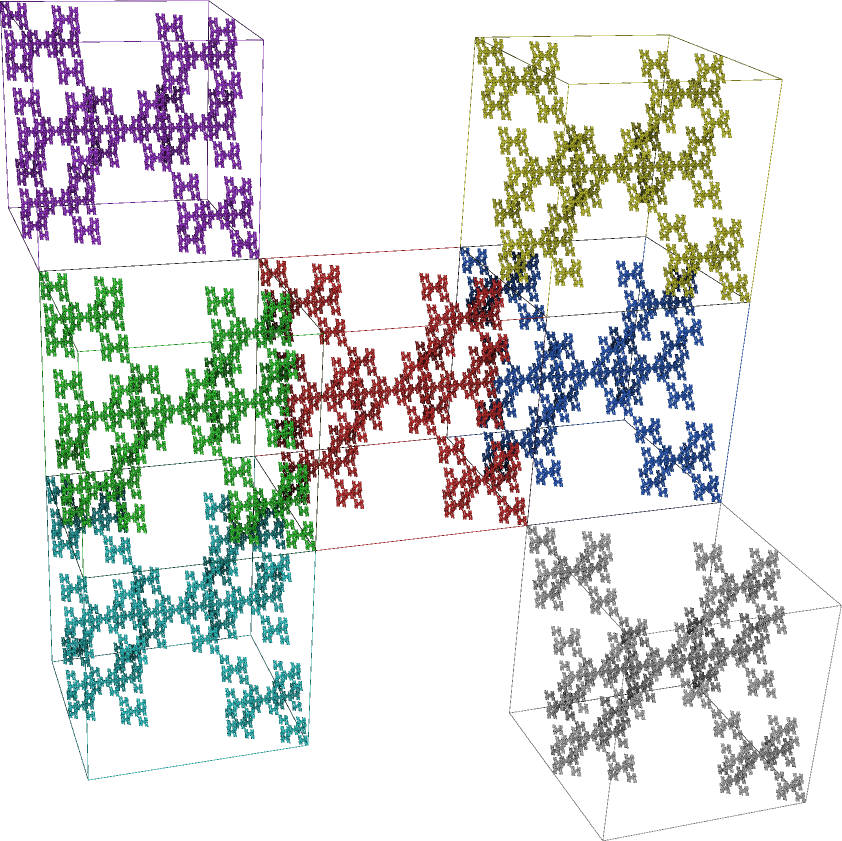}\\
    {\scriptsize$\mD=\{(0,2,0),\ (1,0,1),\ (1,1,0),\ (1,1,1),\ (1,1,2),\ (1,2,1),\ (2,0,2)\}$} &
    {\scriptsize$\mD=\{(0,2,0),\ (1,0,0),\ (1,1,0),\ (1,1,1),\ (1,1,2),\ (1,2,2),\ (2,0,2)\}$} \\
    \hline
    \includegraphics[width=0.23\textwidth] {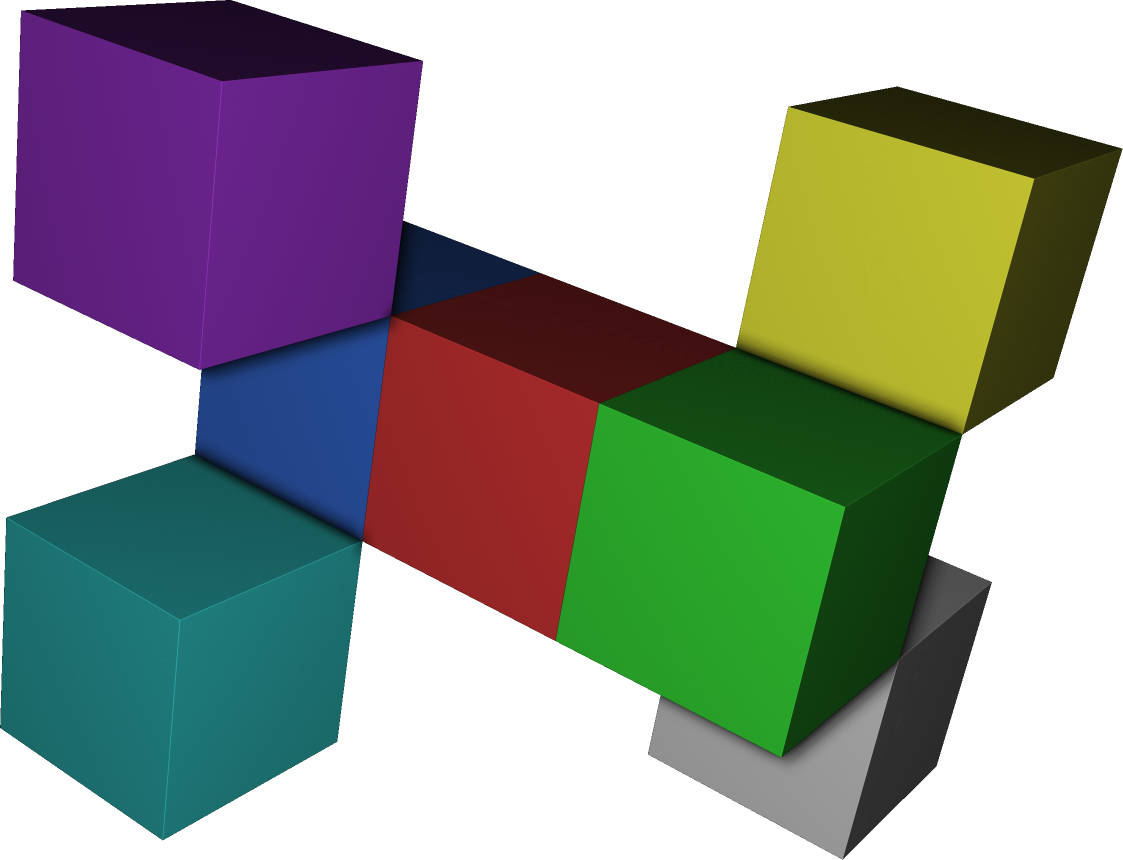}
    \includegraphics[width=0.23\textwidth] {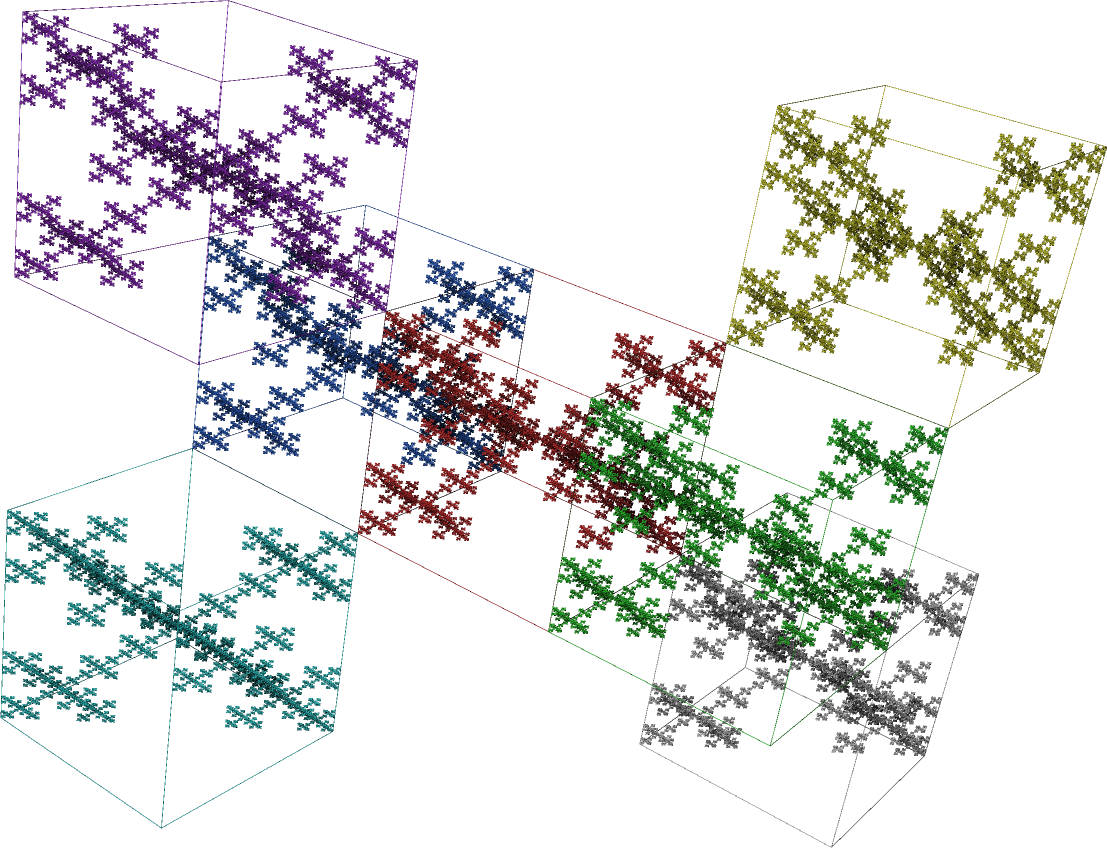} & 
    \includegraphics[width=0.23\textwidth] {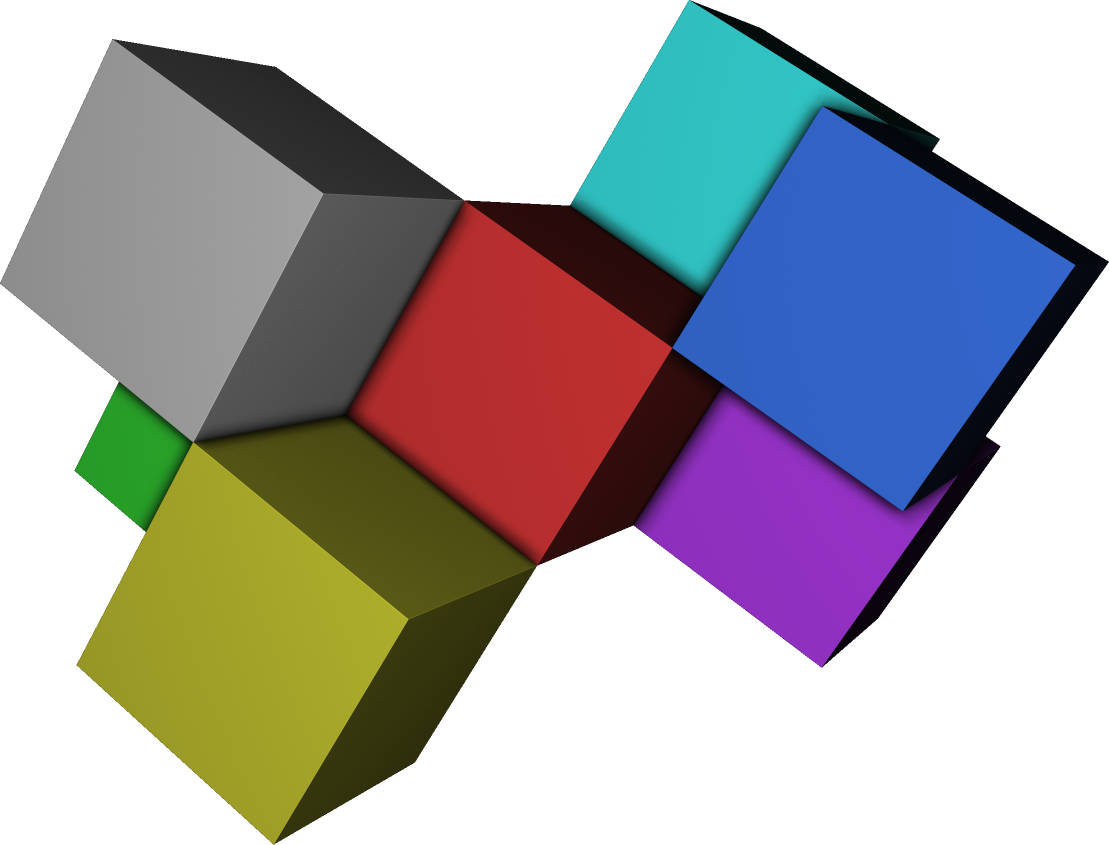} \includegraphics[width=0.23\textwidth] {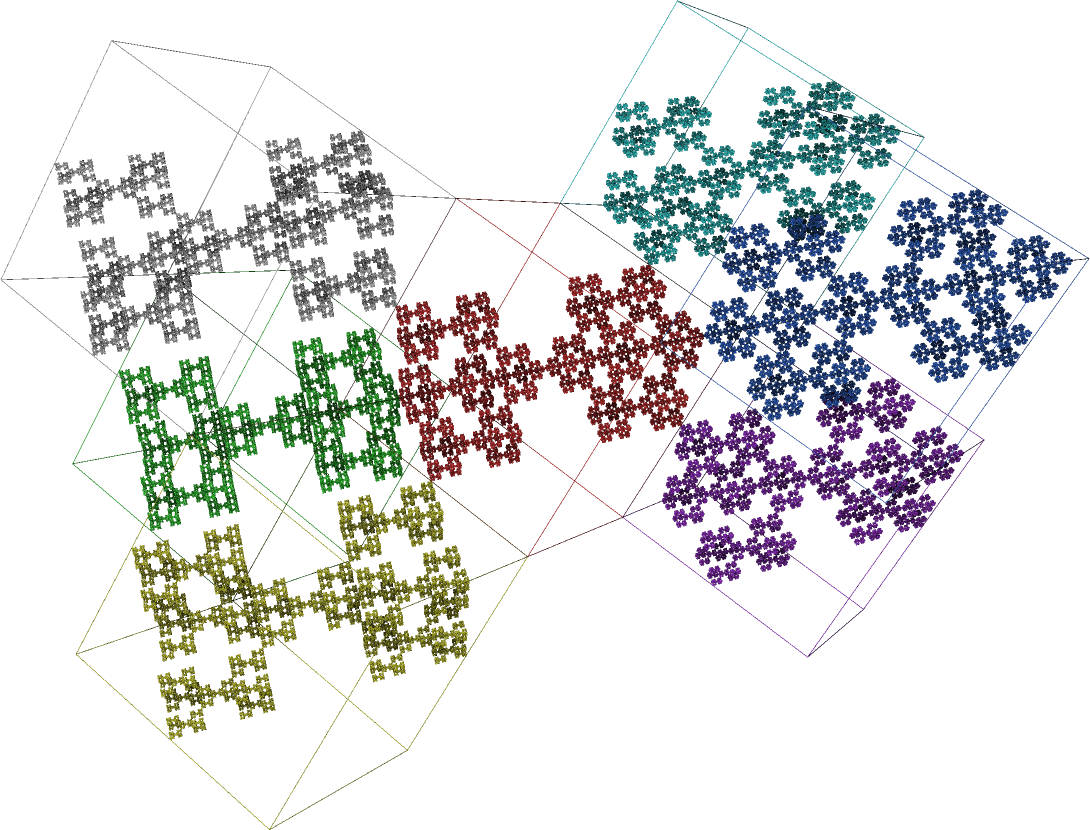}\\
    {\scriptsize$\mD=\{(0,2,0),\ (0,2,2),\ (1,0,1),\ (1,1,1),\ (1,2,1),\ (2,0,0),\ (2,0,2)\}$} &
    {\scriptsize$\mD=\{(0,1,2),\ (0,2,1),\ (1,0,0),\ (1,1,1),\ (1,2,2),\ (2,0,1),\ (2,1,0)\}$} \\
    \hline
    \includegraphics[width=0.23\textwidth] {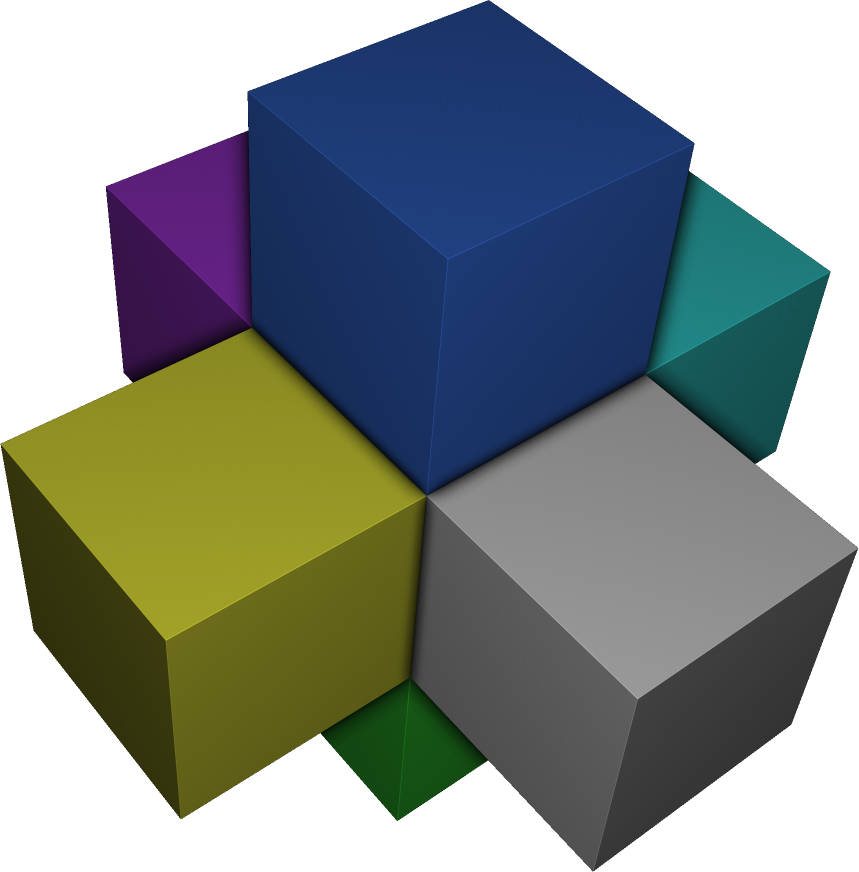}
    \includegraphics[width=0.23\textwidth] {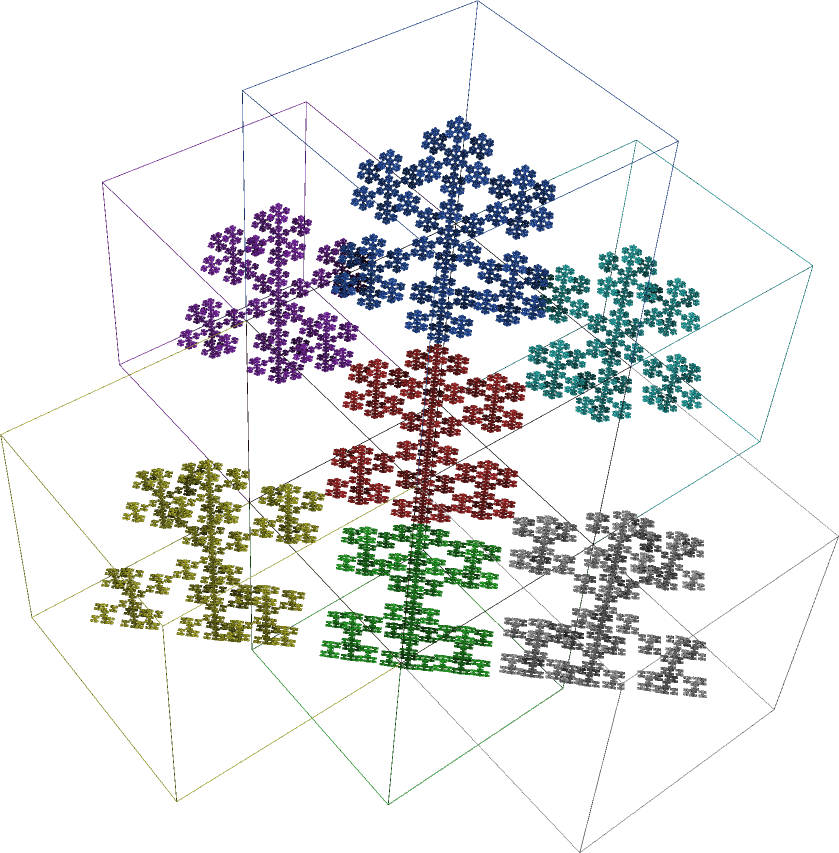} & 
    \includegraphics[width=0.23\textwidth] {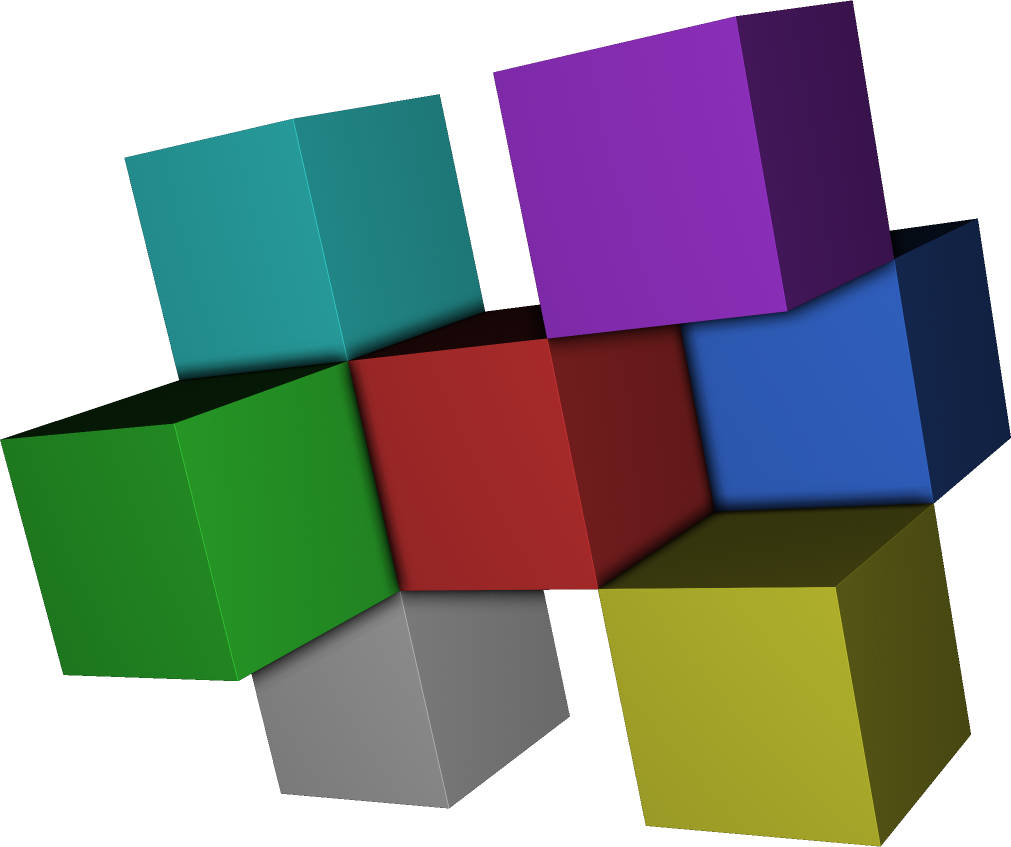} \includegraphics[width=0.23\textwidth] {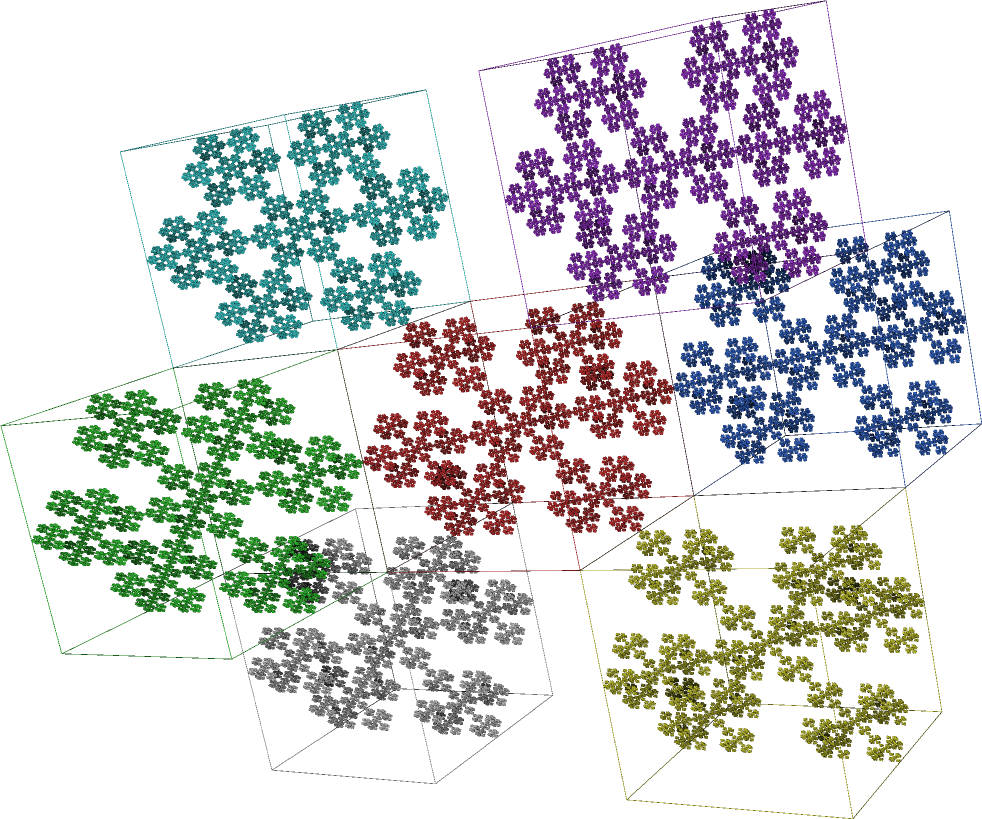}\\
    {\scriptsize$\mD=\{(0,1,1),\ (1,0,1),\ (1,1,0),\ (1,1,1),\ (1,1,2),\ (1,2,1),\ (2,1,1)\}$} &
    {\scriptsize$\mD=\{(0,1,0),\ (0,1,2),\ (1,0,2),\ (1,1,1),\ (1,2,0),\ (2,1,0),\ (2,1,2)\}$} \\
    \hline
    \includegraphics[width=0.23\textwidth] {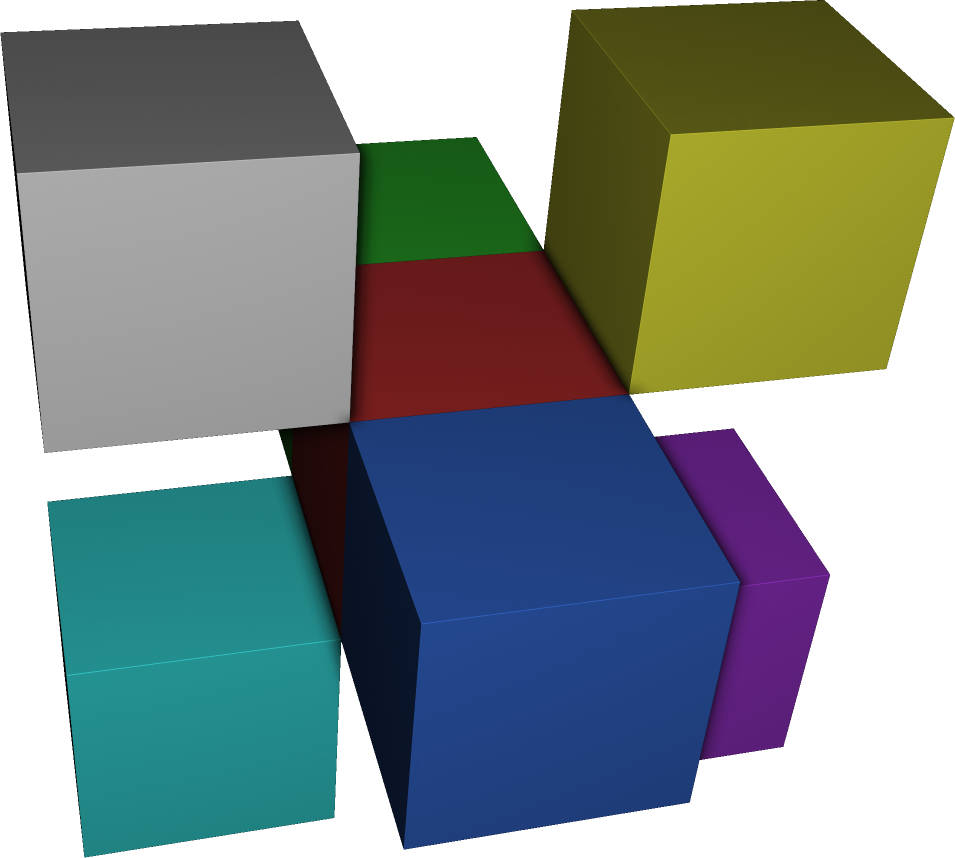}
    \includegraphics[width=0.23\textwidth] {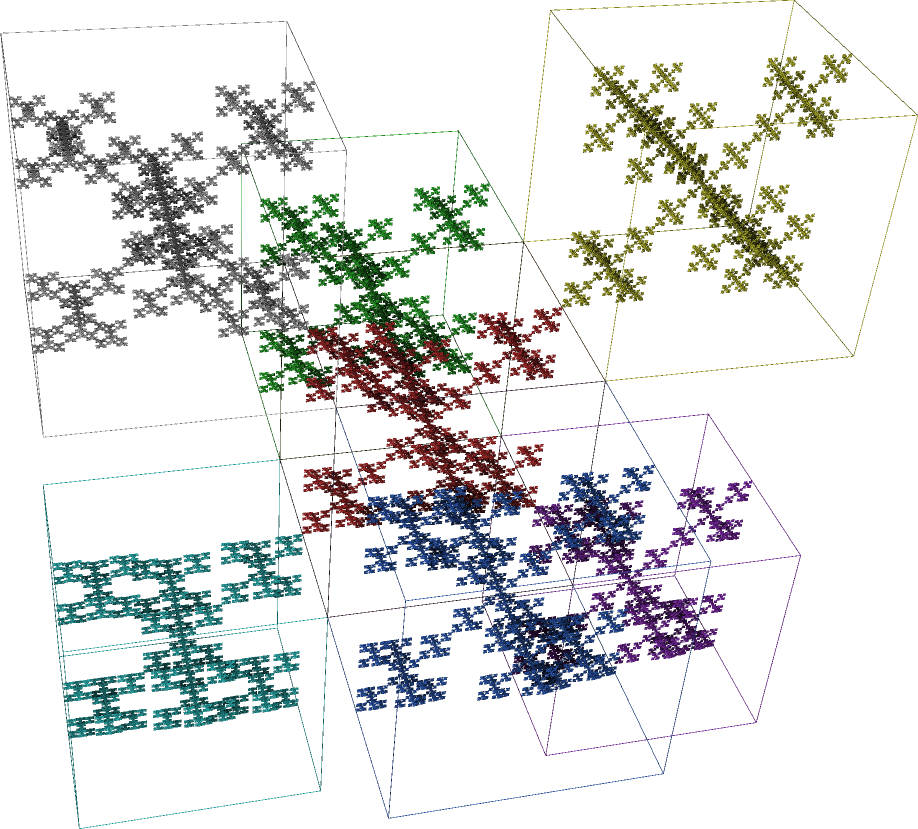} & 
    \includegraphics[width=0.23\textwidth] {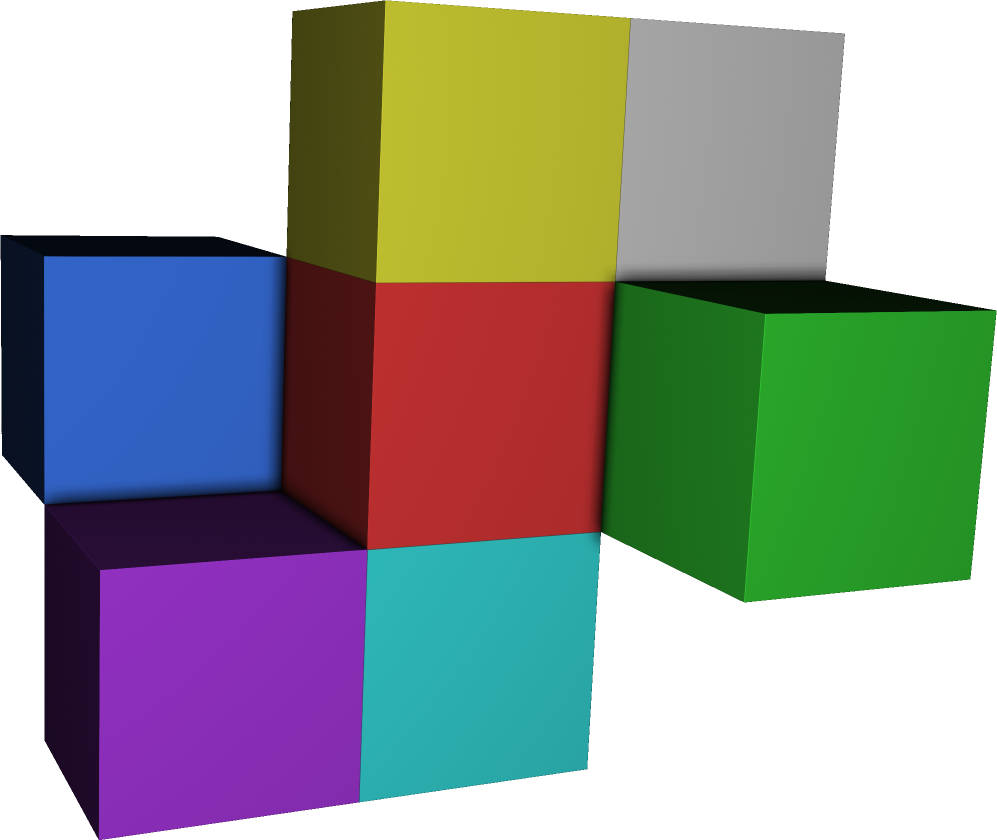} \includegraphics[width=0.23\textwidth] {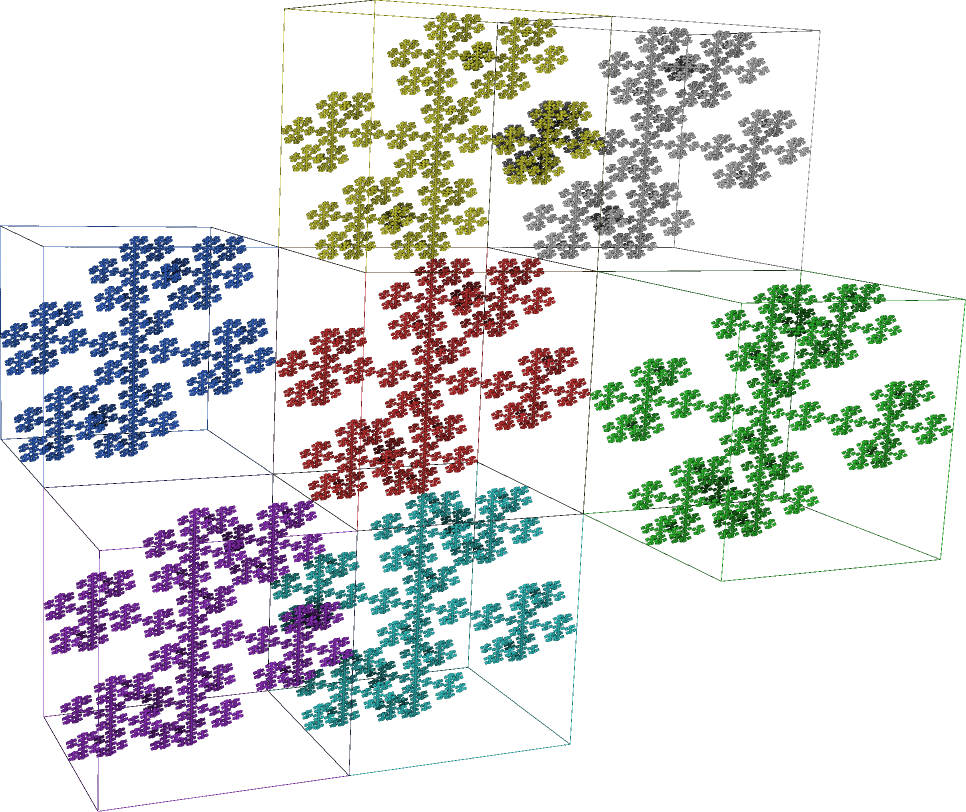}\\
    {\scriptsize$\mD=\{(0,1,0),\ (0,1,2),\ (1,0,1),\ (1,1,1),\ (1,2,1),\ (2,1,0),\ (2,1,2)\}$} &
    {\scriptsize$\mD=\{(0,1,0),\ (0,1,1),\ (1,0,2),\ (1,1,1),\ (1,2,0),\ (2,1,1),\ (2,1,2)\}$} \\
    \hline
    \includegraphics[width=0.23\textwidth] {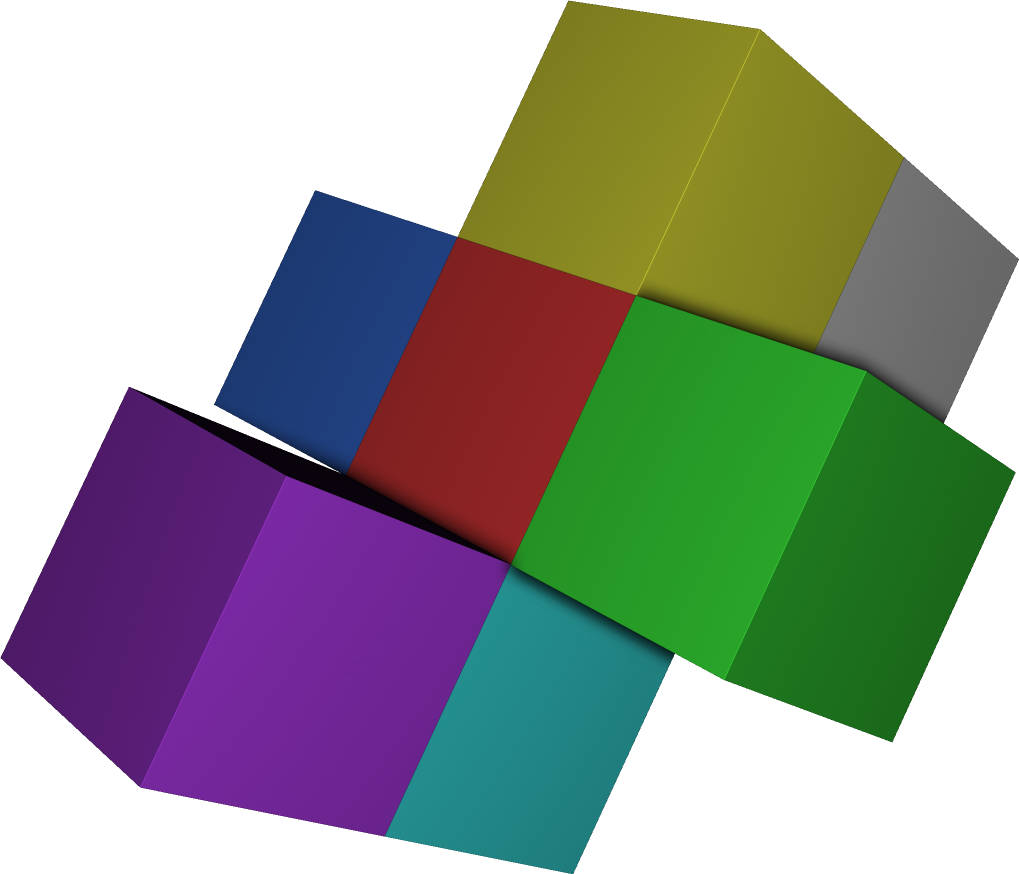}
    \includegraphics[width=0.23\textwidth] {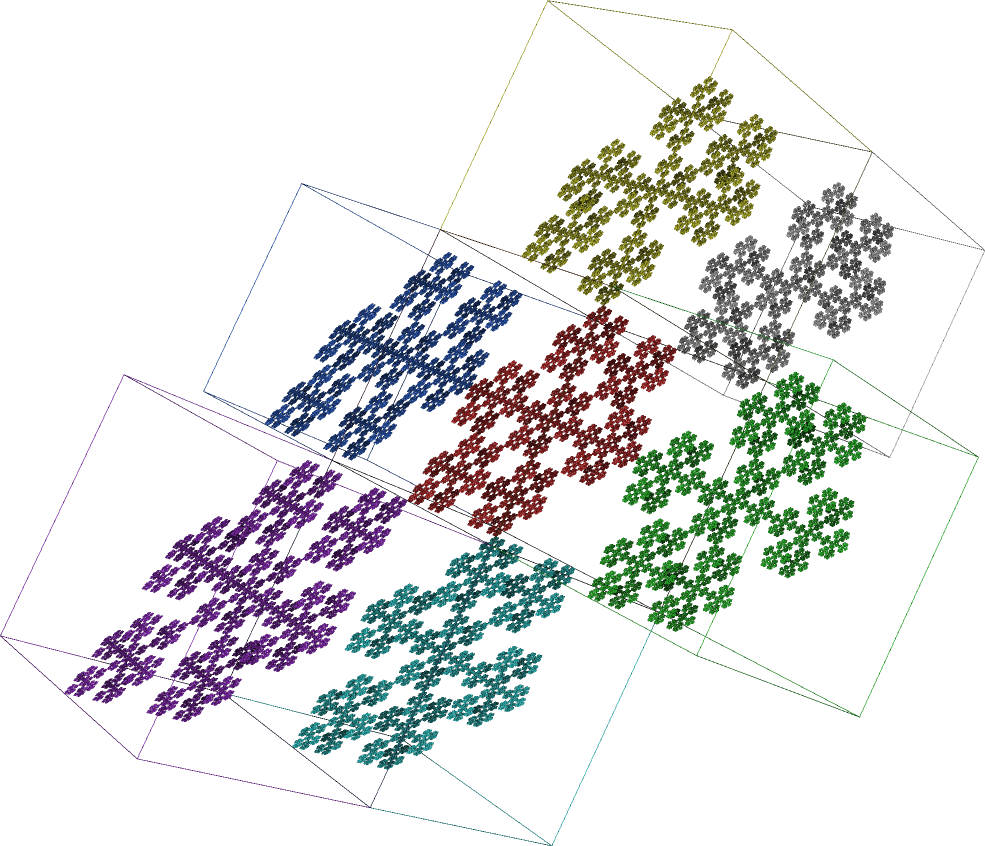} & 
    \includegraphics[width=0.23\textwidth] {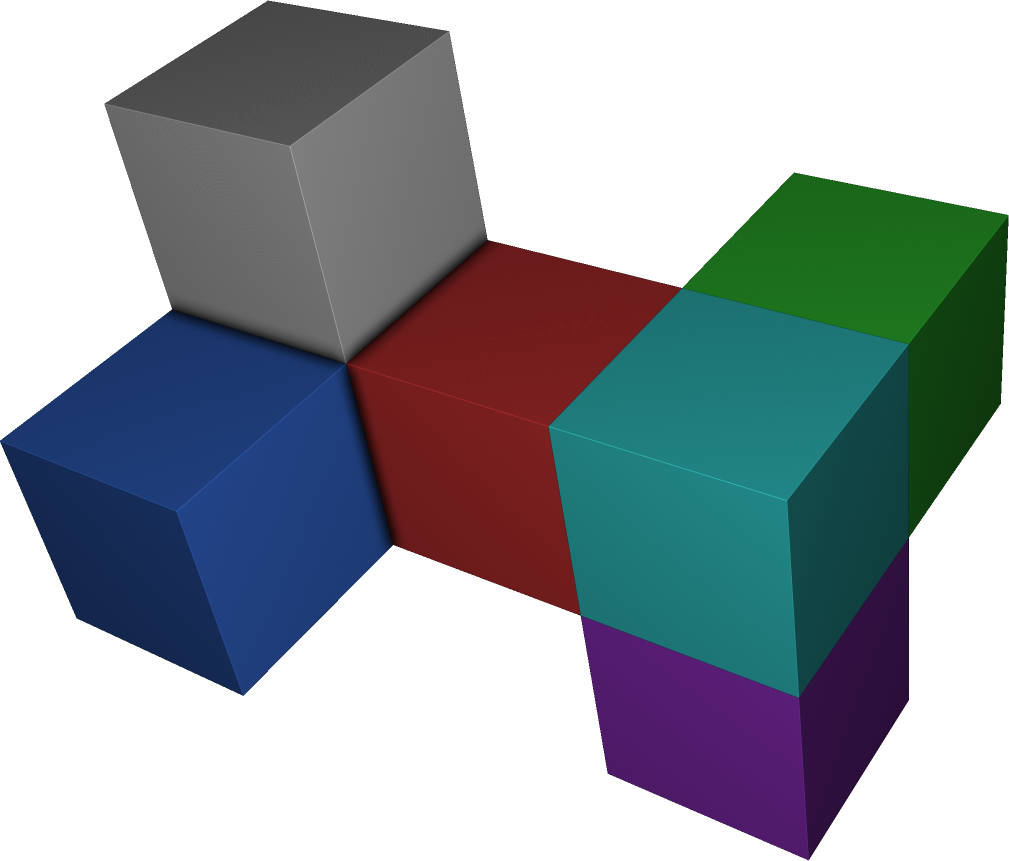} \includegraphics[width=0.23\textwidth] {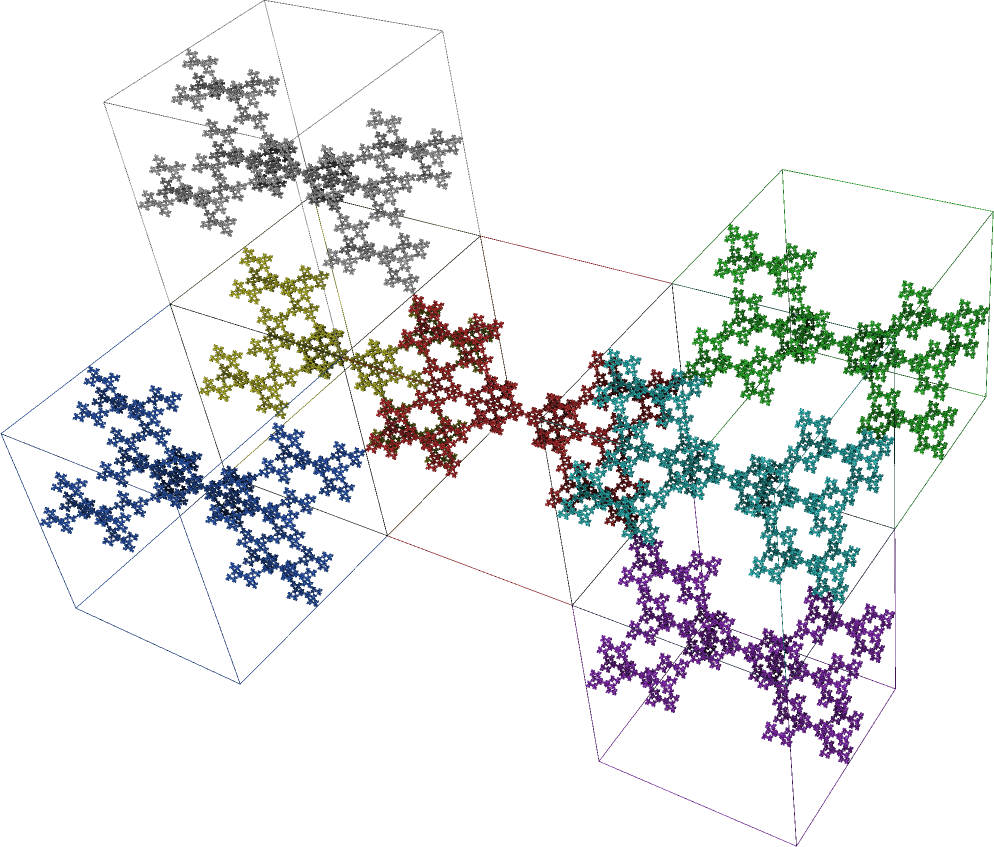}\\
    {\scriptsize$\mD=\{(0,1,0),\ (0,1,1),\ (1,0,1),\ (1,1,1),\ (1,2,1),\ (2,1,1),\ (2,1,2)\}$} &
    {\scriptsize$\mD=\{(0,1,0),\ (0,1,1),\ (0,2,1),\ (1,1,1),\ (2,0,1),\ (2,1,1),\ (2,1,2)\}$} \\
    \hline
    \includegraphics[width=0.23\textwidth] {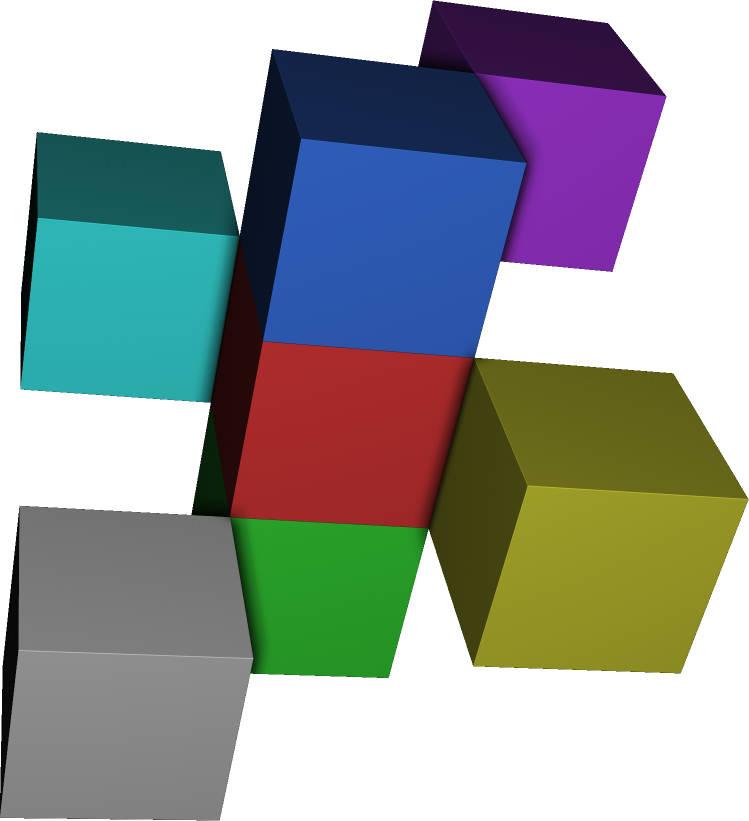}
    \includegraphics[width=0.23\textwidth] {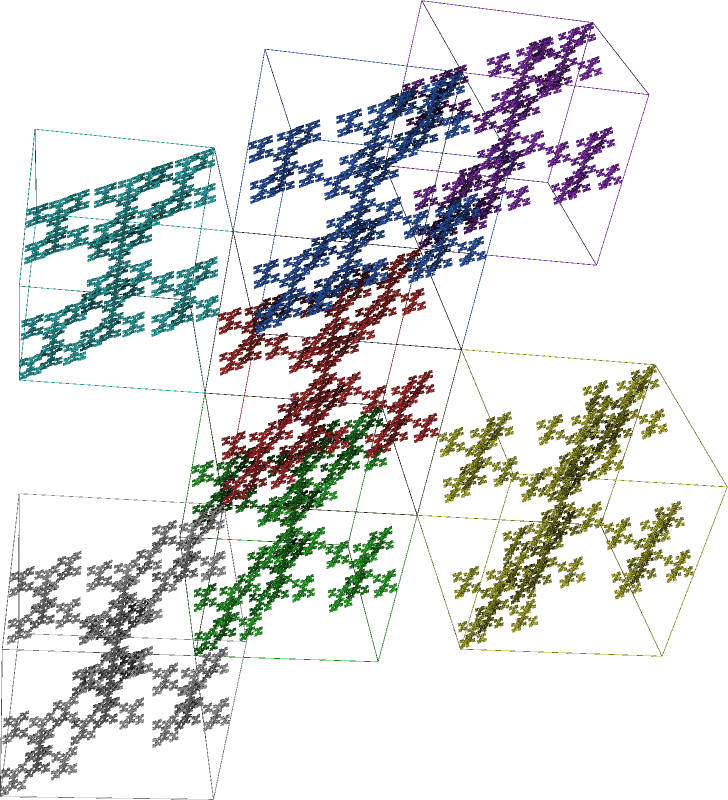} & 
    \includegraphics[width=0.23\textwidth] {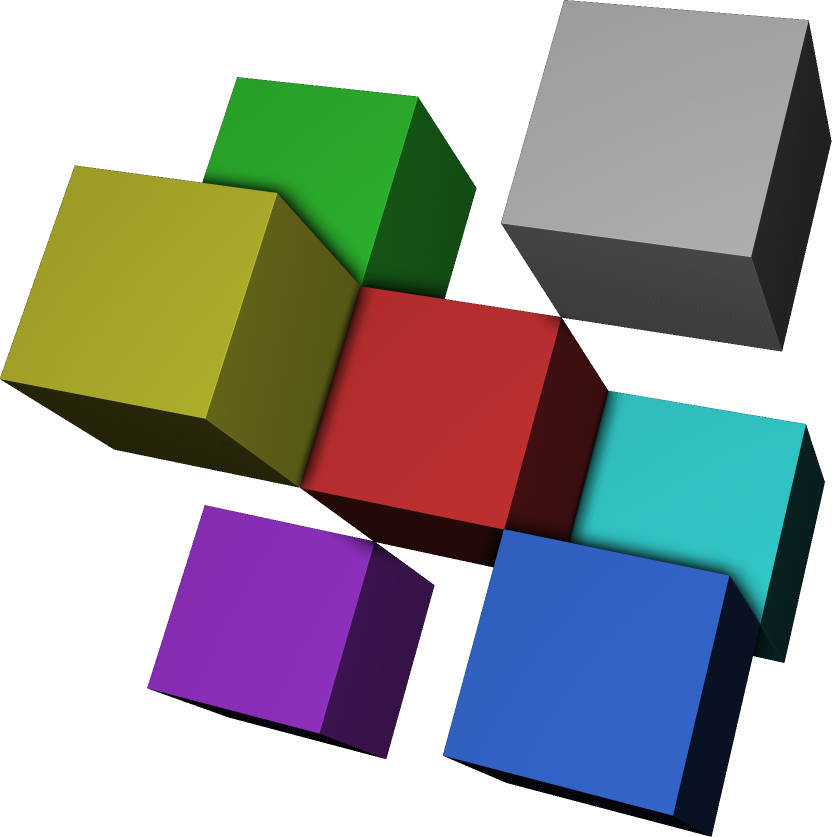} \includegraphics[width=0.23\textwidth] {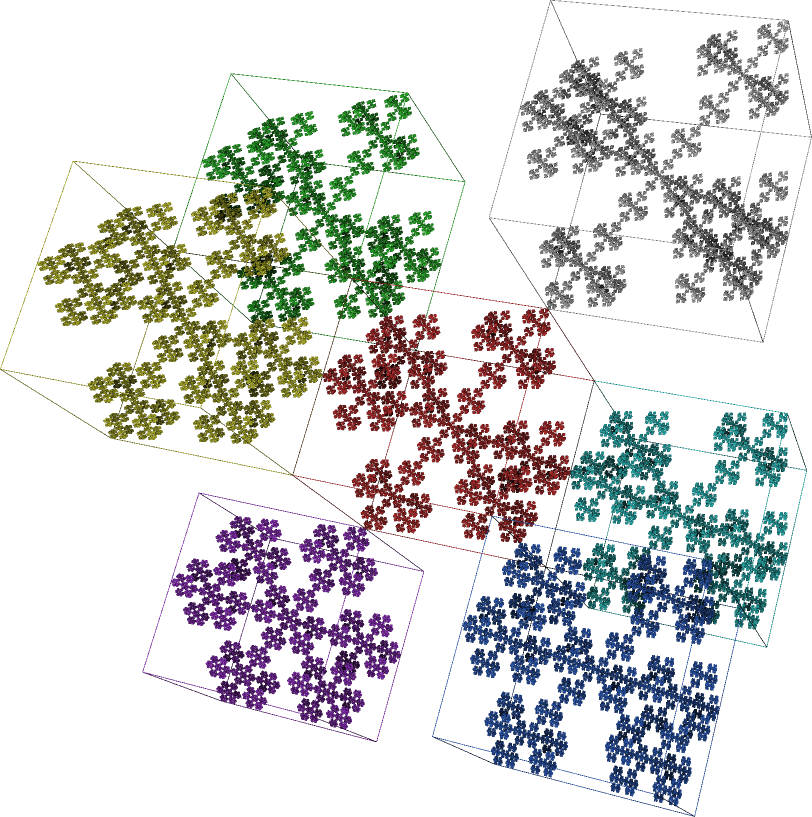}\\
    {\scriptsize$\mD=\{(0,0,2),\ (0,2,1),\ (1,1,0),\ (1,1,1),\ (1,1,2),\ (2,0,1),\ (2,2,0)\}$} &
    {\scriptsize$\mD=\{(0,0,2),\ (0,2,1),\ (1,0,0),\ (1,1,1),\ (1,2,2),\ (2,0,1),\ (2,2,0)\}$} \\
    \hline
    \includegraphics[width=0.23\textwidth] {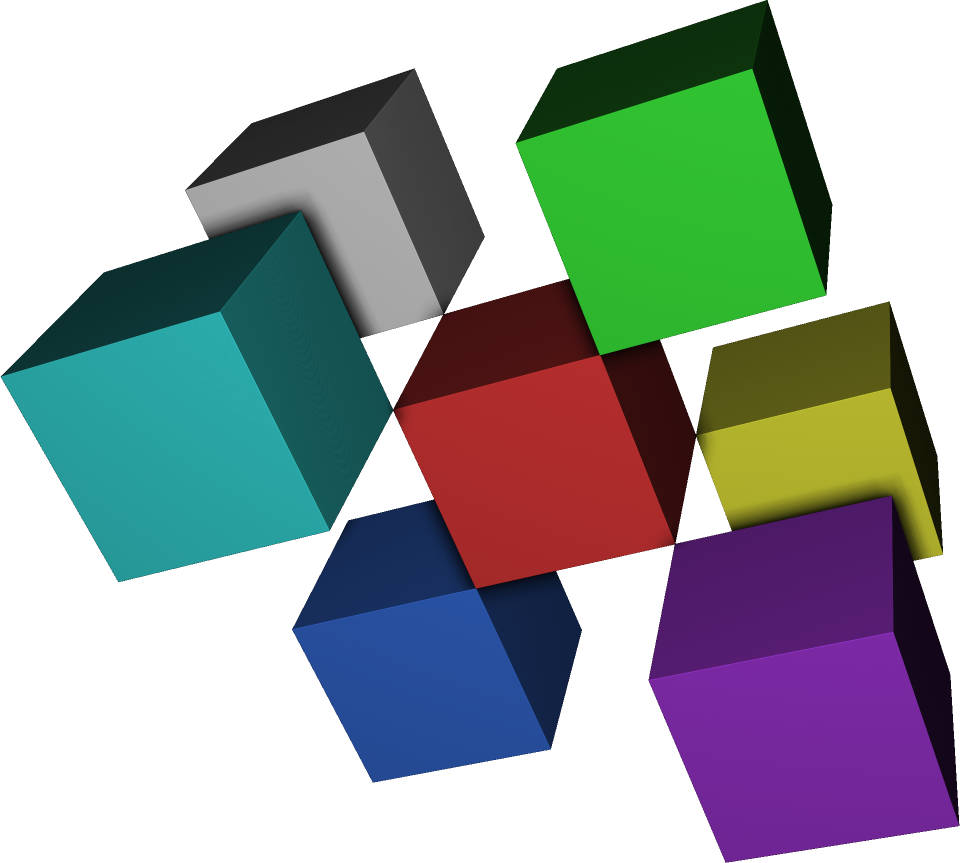}
    \includegraphics[width=0.23\textwidth] {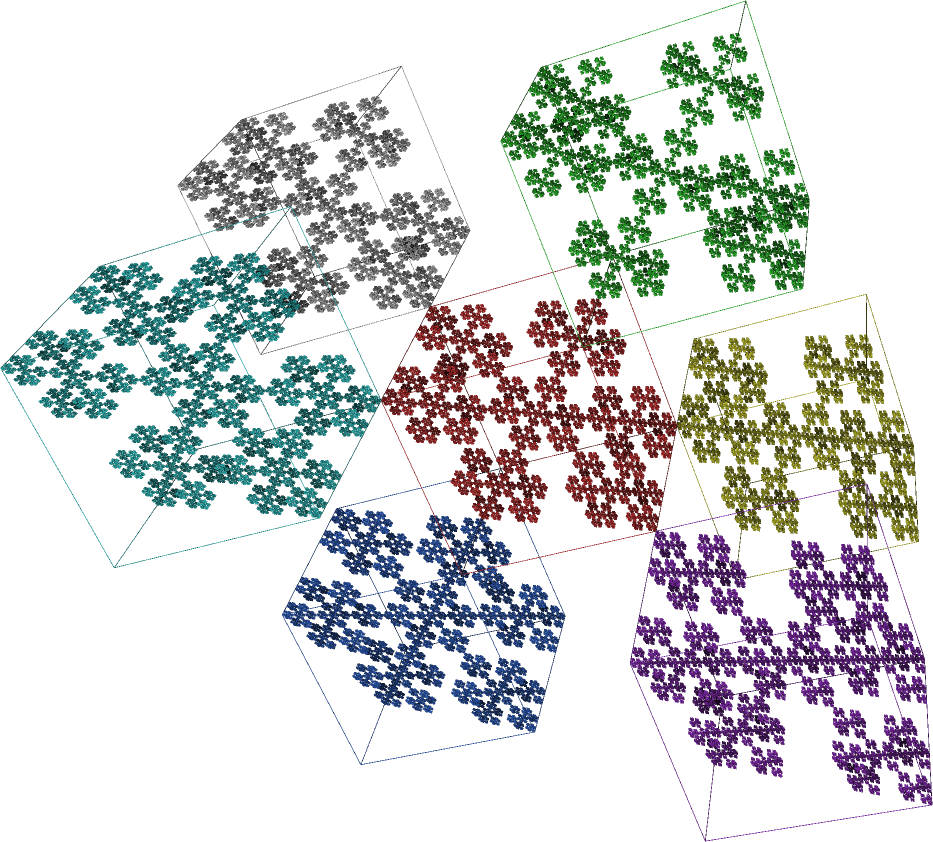} & 
    \includegraphics[width=0.23\textwidth] {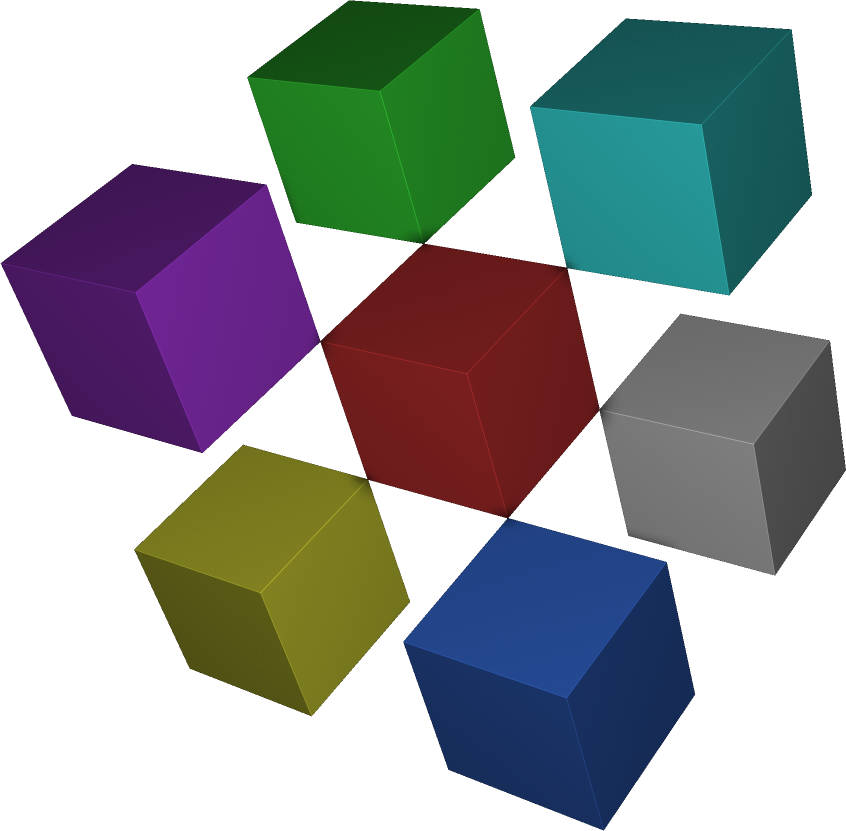} \includegraphics[width=0.23\textwidth] {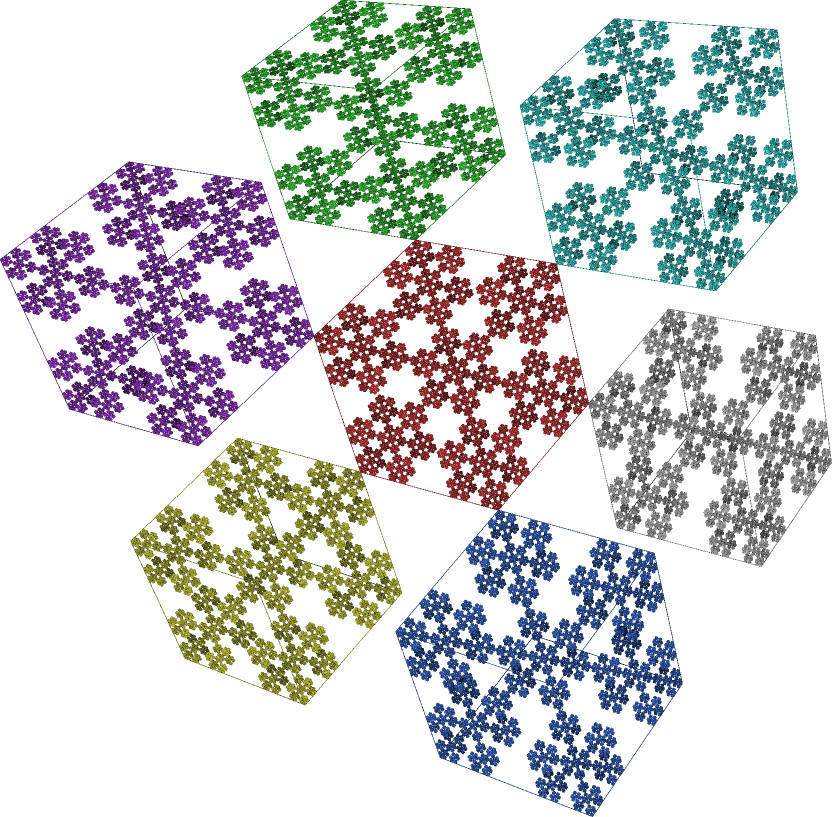}\\
    {\scriptsize$\mD=\{(0,0,2),\ (0,2,0),\ (1,0,0),\ (1,1,1),\ (1,2,2),\ (2,0,2),\ (2,2,0)\}$} &
    {\scriptsize$\mD=\{(0,0,2),\ (0,2,0),\ (0,2,2),\ (1,1,1),\ (2,0,0),\ (2,0,2),\ (2,2,0)\}$} \\
    \hline
    \includegraphics[width=0.23\textwidth] {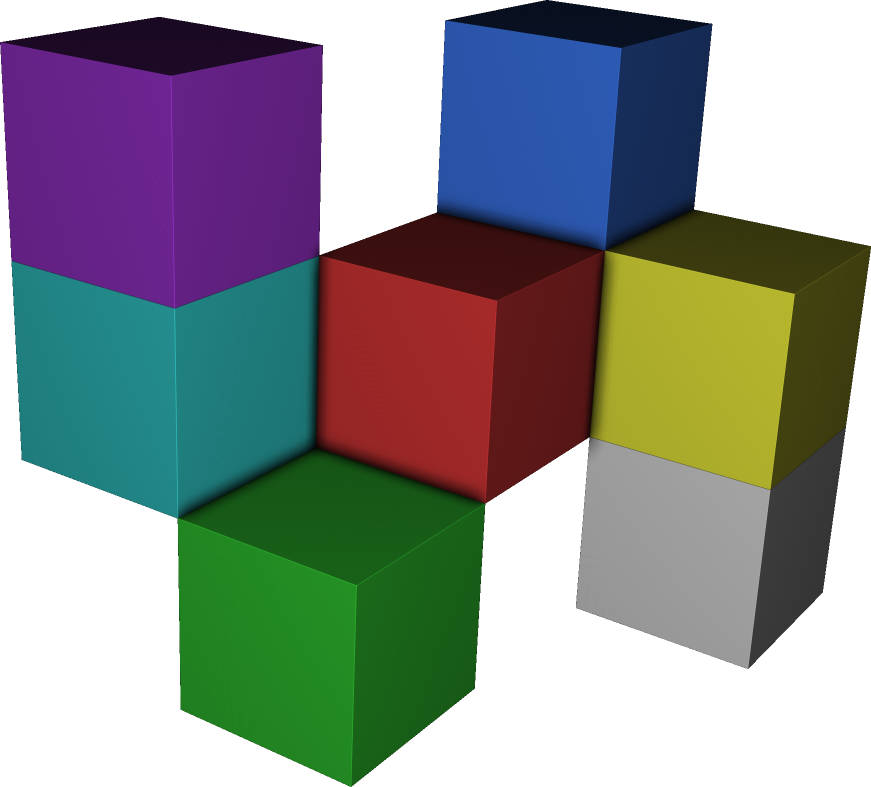}
    \includegraphics[width=0.23\textwidth] {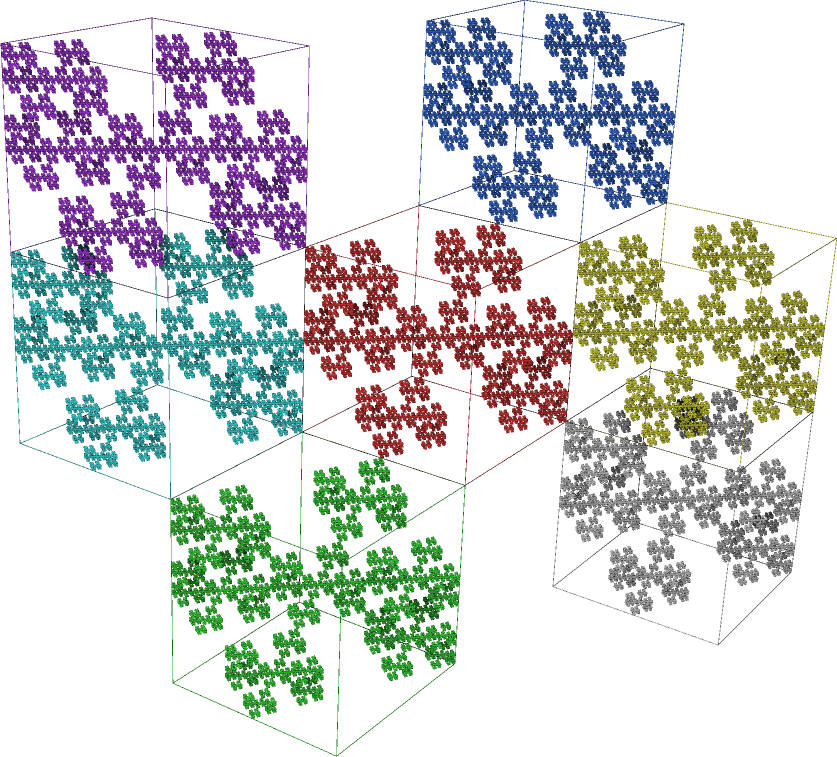} & 
    \includegraphics[width=0.23\textwidth] {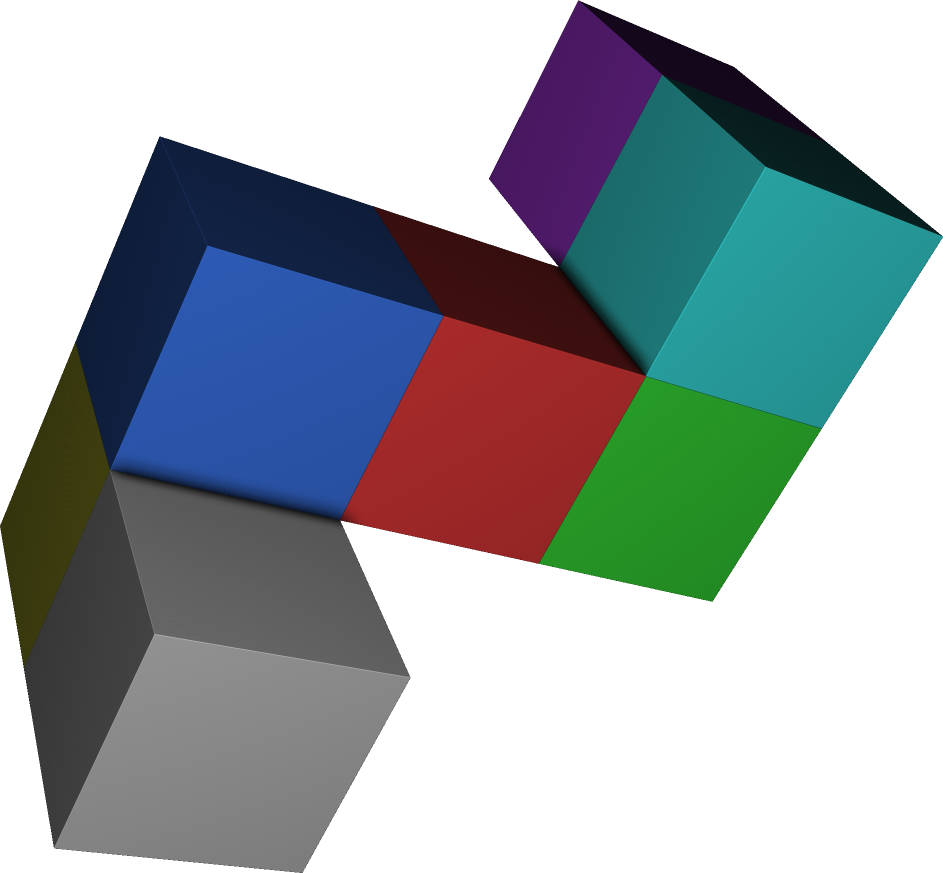} \includegraphics[width=0.23\textwidth] {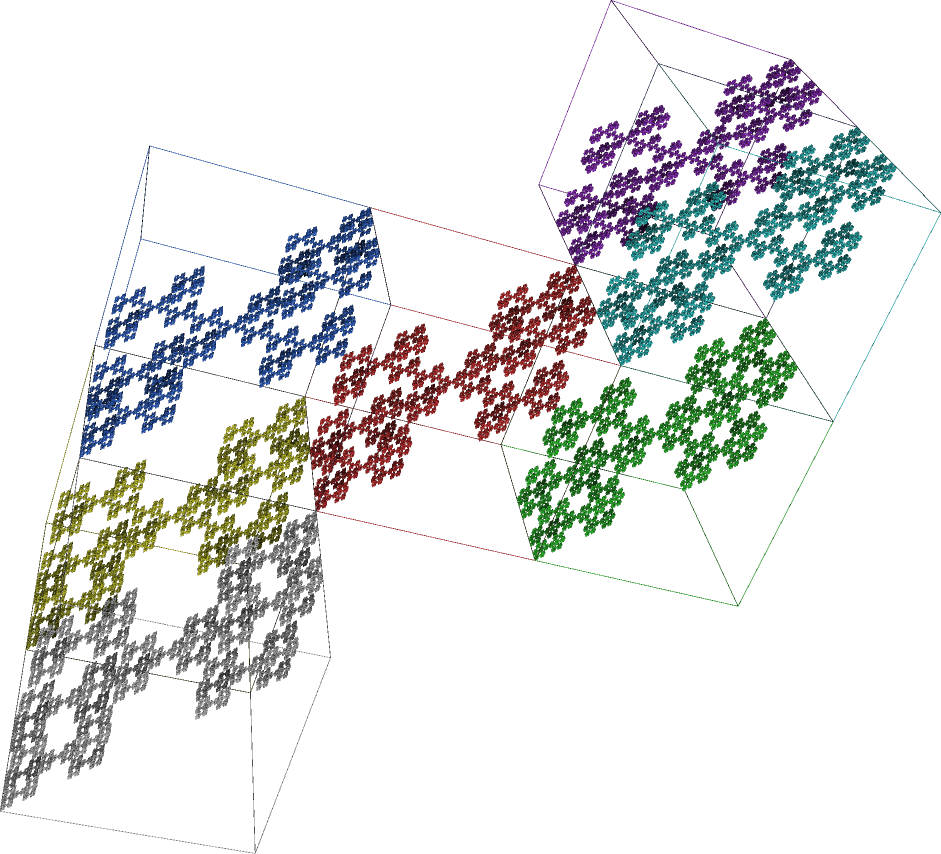}\\
    {\scriptsize$\mD=\{(0,0,0),\ (0,0,1),\ (1,0,2),\ (1,1,1),\ (1,2,0),\ (2,2,1),\ (2,2,2)\}$} &
    {\scriptsize$\mD=\{(0,0,0),\ (0,0,1),\ (1,0,1),\ (1,1,1),\ (1,2,1),\ (2,2,1),\ (2,2,2)\}$} \\
    \hline
    \includegraphics[width=0.23\textwidth] {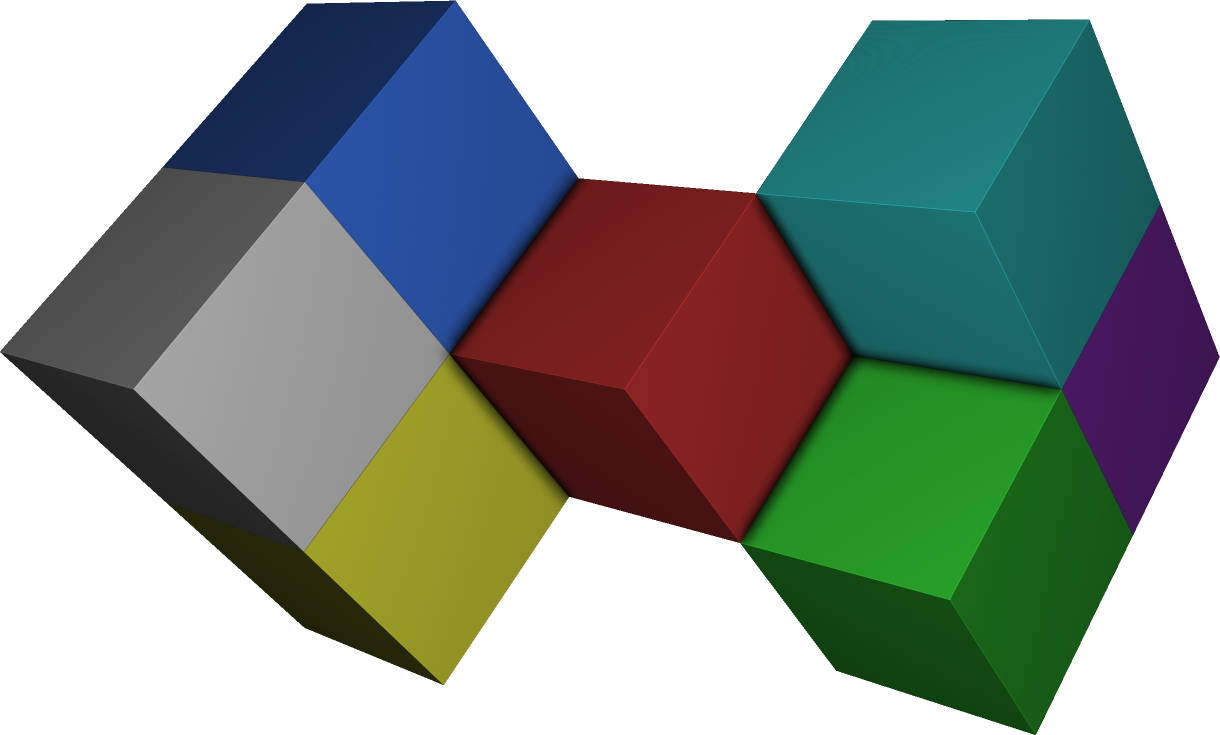}
    \includegraphics[width=0.23\textwidth] {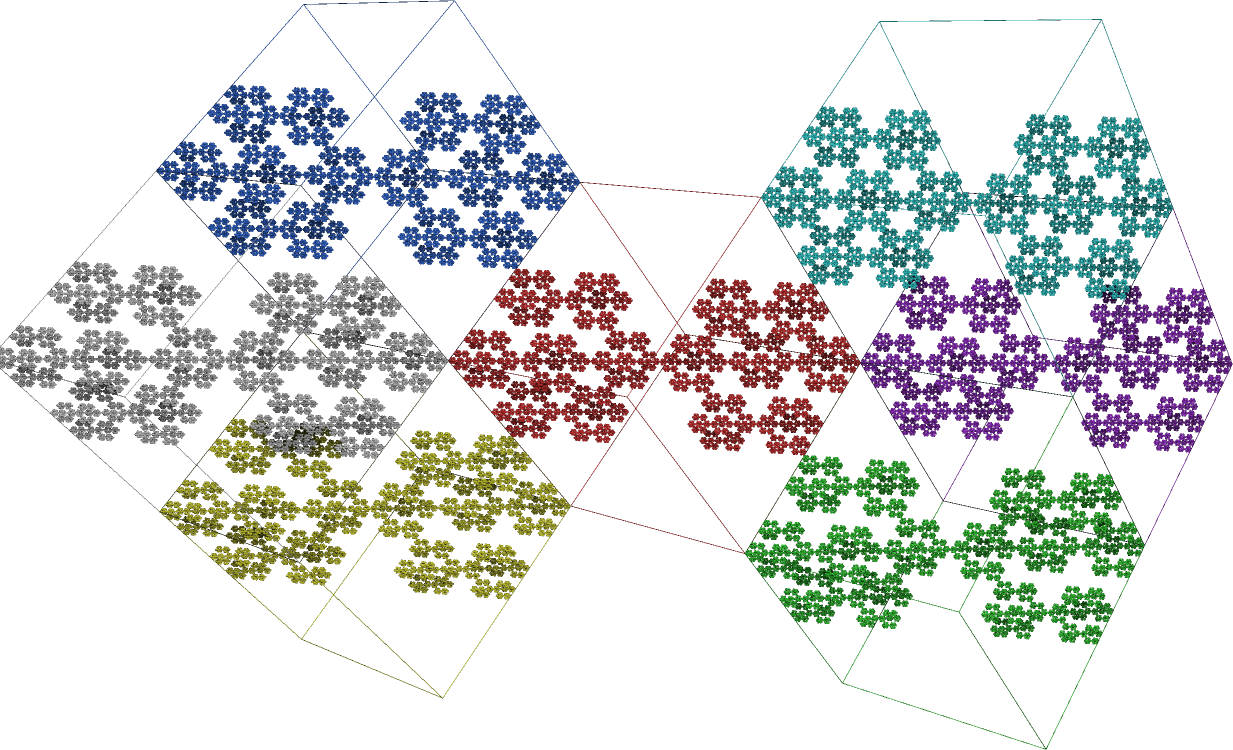} & 
    \includegraphics[width=0.23\textwidth] {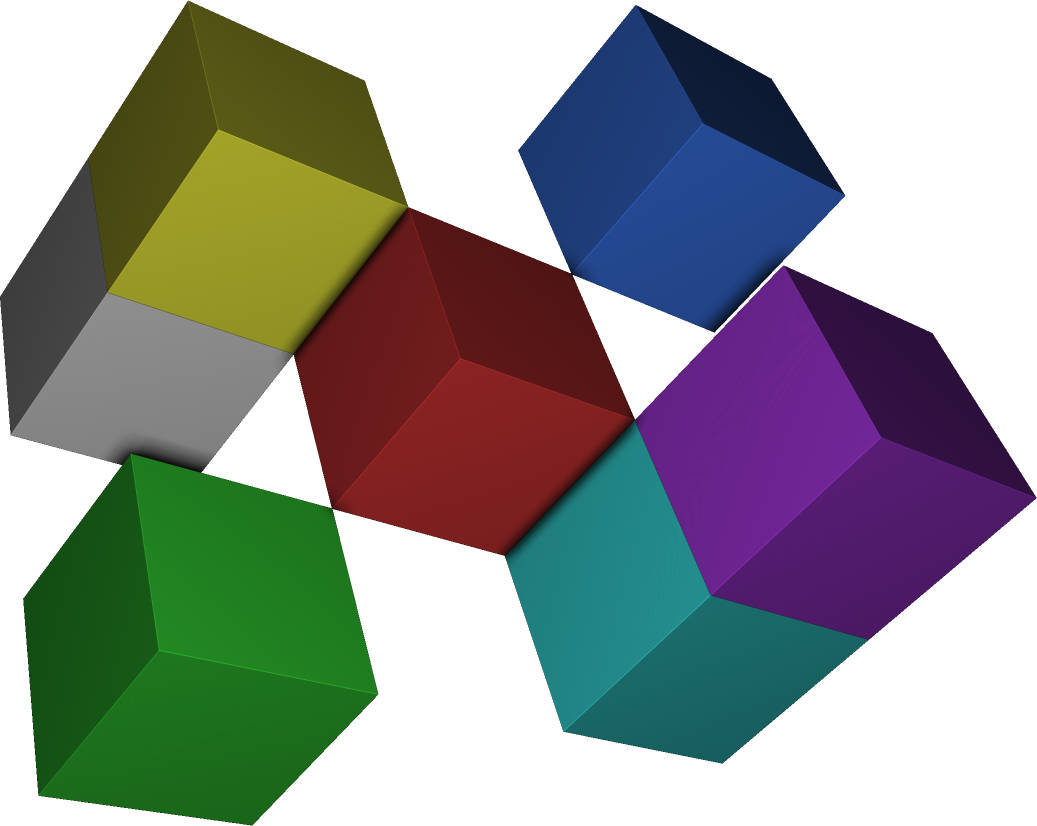} \includegraphics[width=0.23\textwidth] {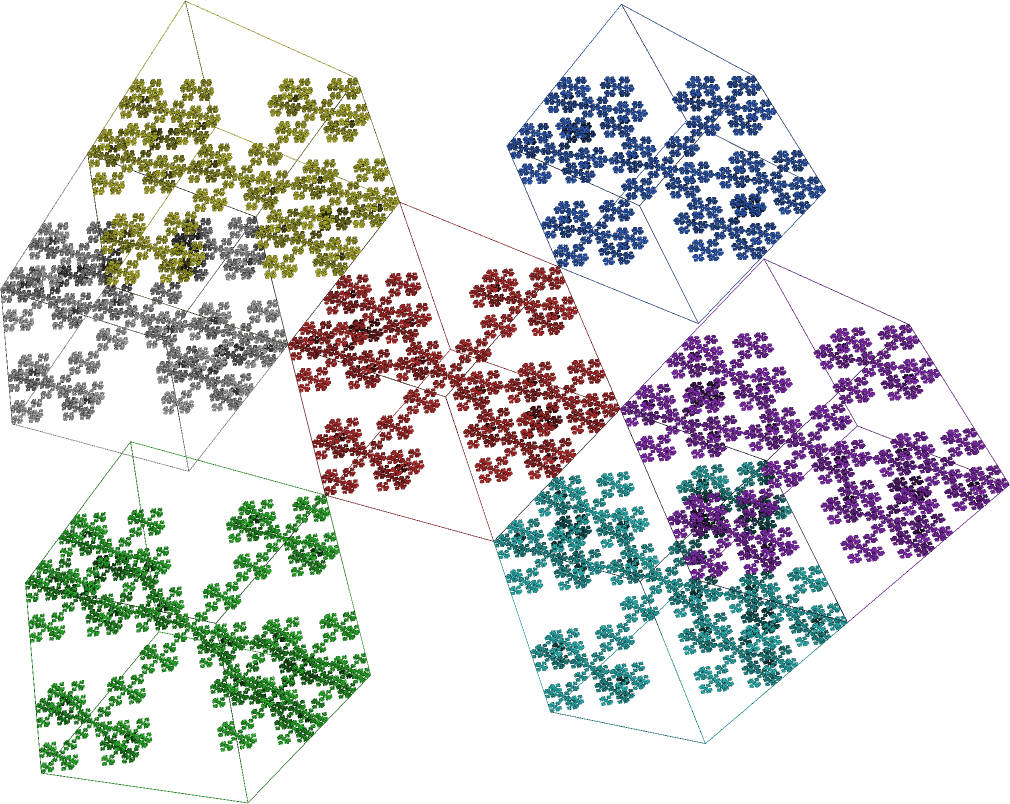}\\
    {\scriptsize$\mD=\{(0,0,0),\ (0,0,1),\ (1,0,0),\ (1,1,1),\ (1,2,2),\ (2,2,1),\ (2,2,2)\}$} &
    {\scriptsize$\mD=\{(0,0,0),\ (0,0,1),\ (0,2,2),\ (1,1,1),\ (2,0,0),\ (2,2,1),\ (2,2,2)\}$} \\
    \hline
    \includegraphics[width=0.23\textwidth] {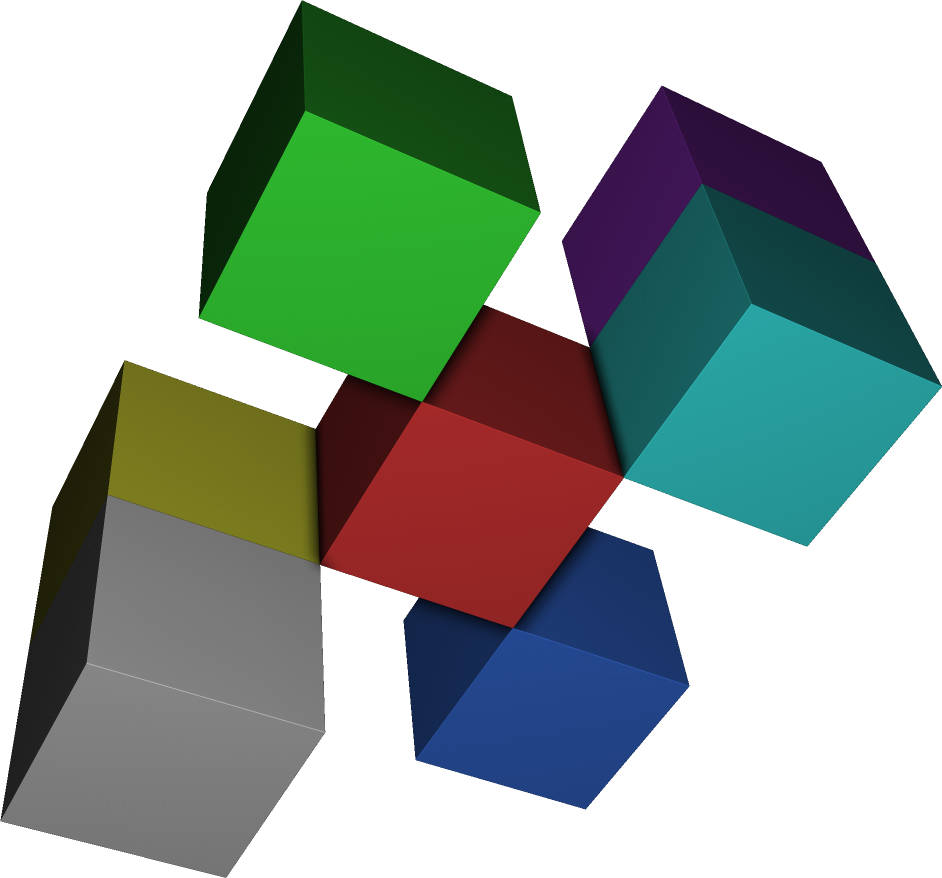}
    \includegraphics[width=0.23\textwidth] {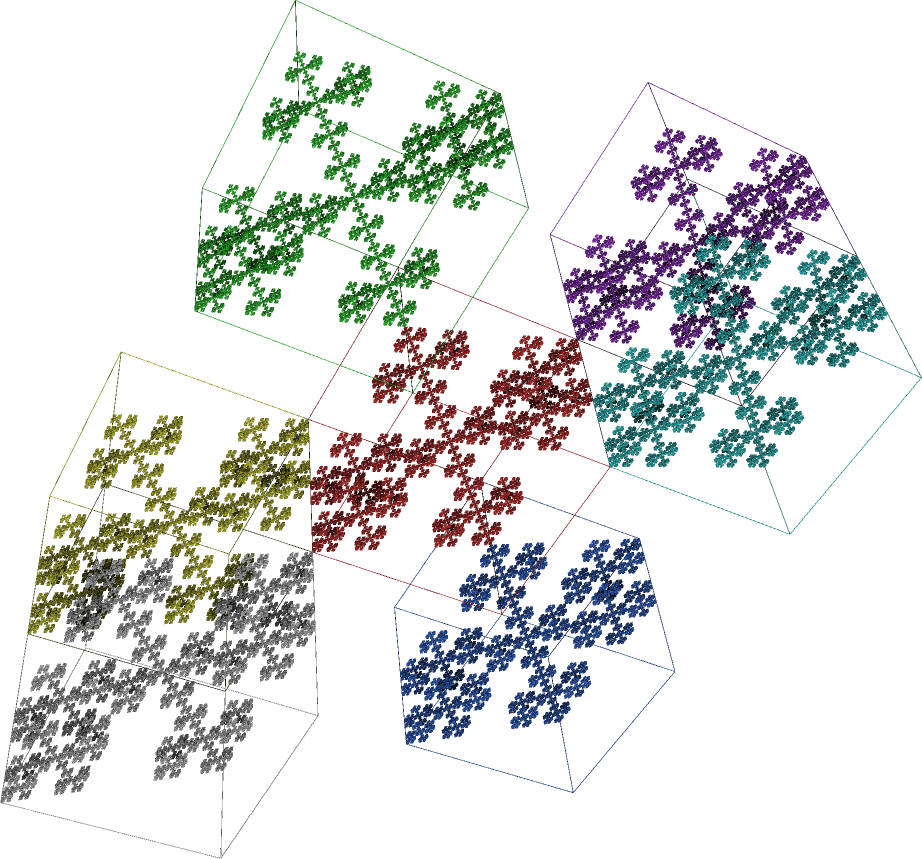} &\\
    {\scriptsize$\mD=\{(0,0,0),\ (0,0,1),\ (0,2,1),\ (1,1,1),\ (2,0,1),\ (2,2,1),\ (2,2,2)\}$} &\\
    \hline
\end{longtable}

\newpage
\begin{longtable}{|p{0.48\textwidth}|p{0.48\textwidth}|}
\caption{Dendrites corresponding to graph $7_{10}$ ($N=3$)}\label{tab:d2}\\
    \hline
    \multicolumn{2}{|c|}{\includegraphics{den2.pdf}} \\
    \hline
    \includegraphics[width=0.23\textwidth] {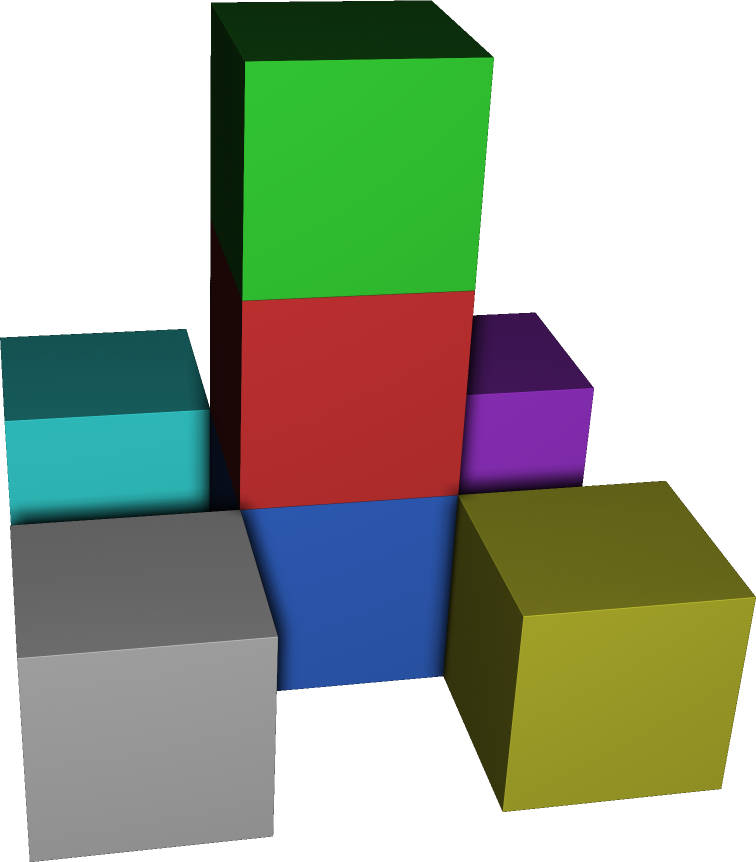}
    \includegraphics[width=0.23\textwidth] {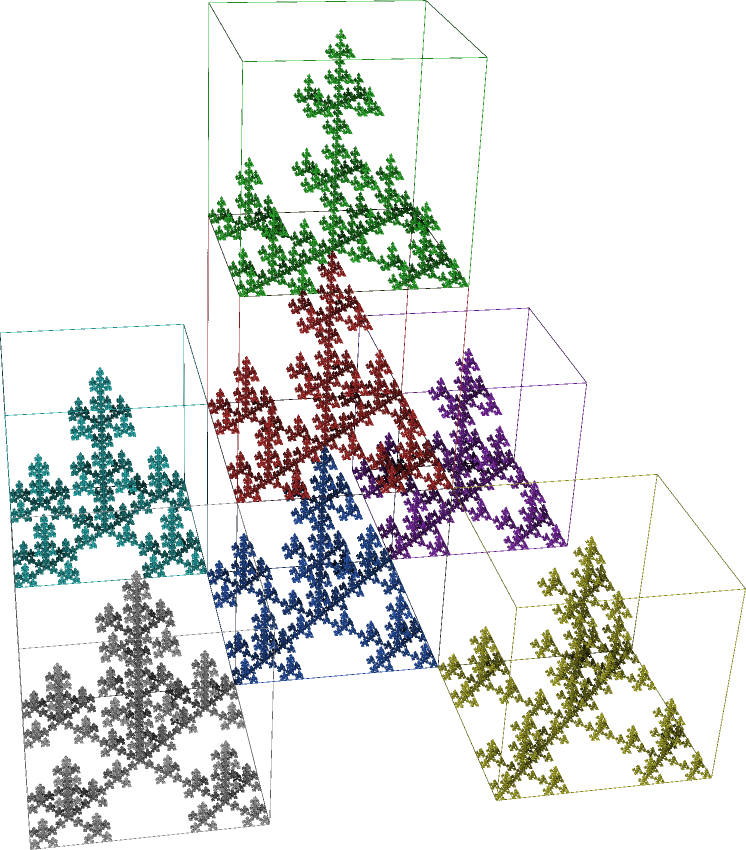} &
    \includegraphics[width=0.23\textwidth] {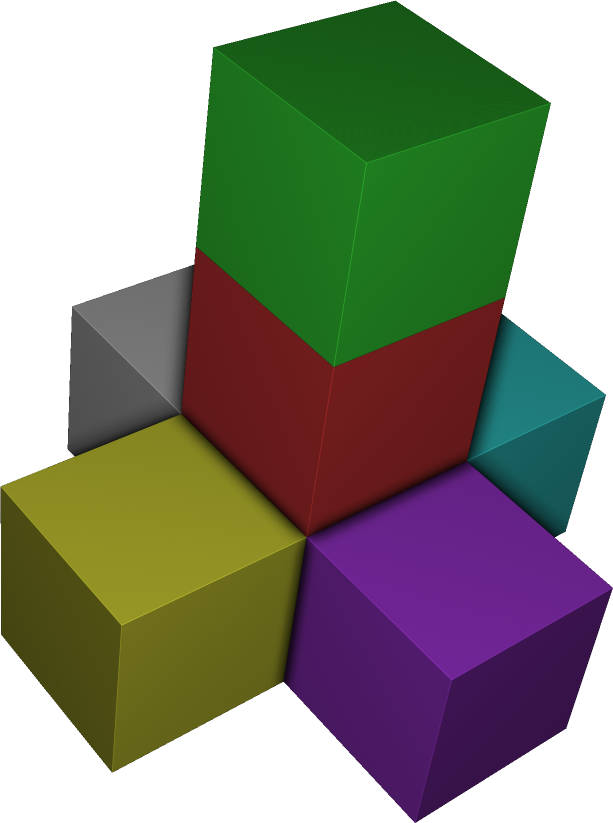}
    \includegraphics[width=0.23\textwidth] {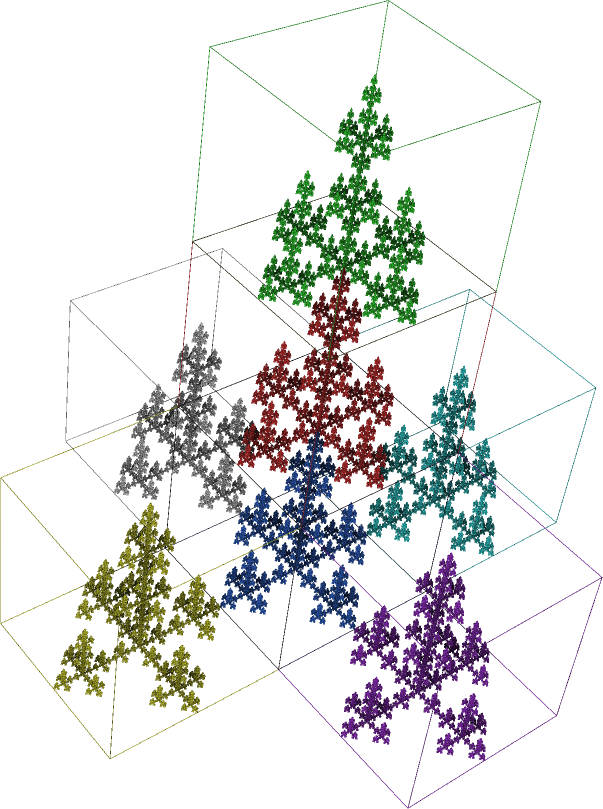} \\
    {\scriptsize$\mD=\{(0,0,0),\ (0,0,2),\ (1,0,1),\ (1,1,1),\ (1,2,1),\ (2,0,0),\ (2,0,2)\}$} & 
    {\scriptsize$\mD=\{(0,1,0),\ (1,0,0),\ (1,1,0),\ (1,1,1),\ (1,1,2),\ (1,2,0),\ (2,1,0)\}$} \\
    \hline
    \includegraphics[width=0.23\textwidth] {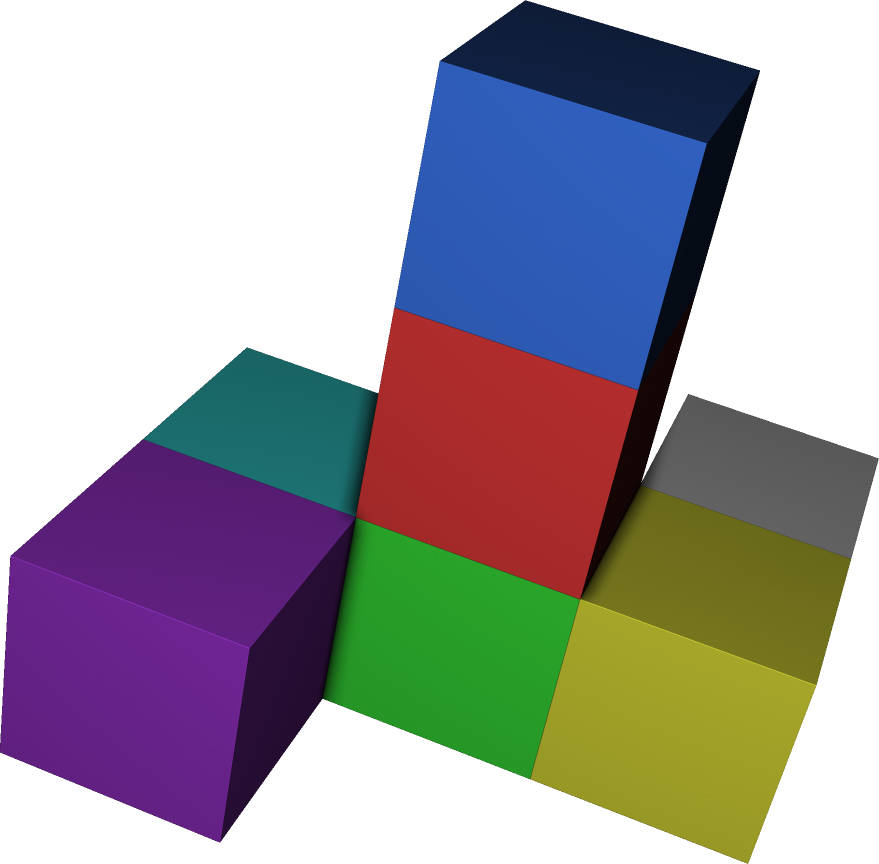}
    \includegraphics[width=0.23\textwidth] {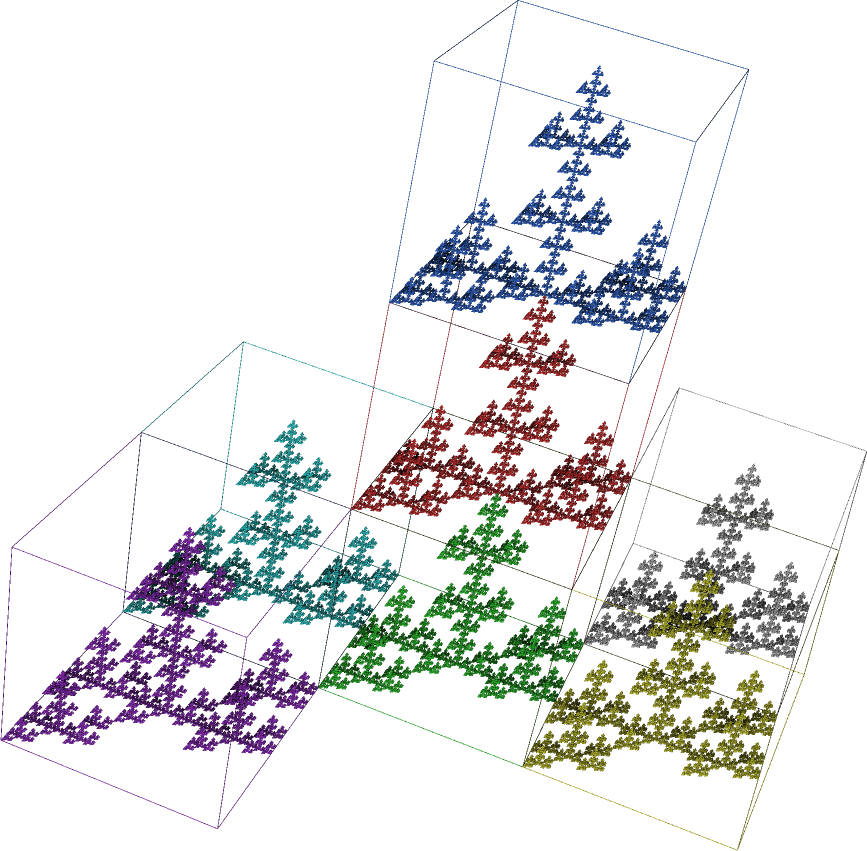} &\\
    {\scriptsize$\mD=\{(0,2,0),\ (0,2,1),\ (1,0,1),\ (1,1,1),\ (1,2,1),\ (2,2,1),\ (2,2,2)\}$} &\\
    \hline

\end{longtable}

\newpage
\begin{longtable}{|p{0.48\textwidth}|p{0.48\textwidth}|}
\caption{Dendrites corresponding to graph $7_{9}$ ($N=12$)}\label{tab:d3}\\
    \hline
    \multicolumn{2}{|c|}{\includegraphics{den3.pdf}} \\
    \hline
    \includegraphics[width=0.23\textwidth] {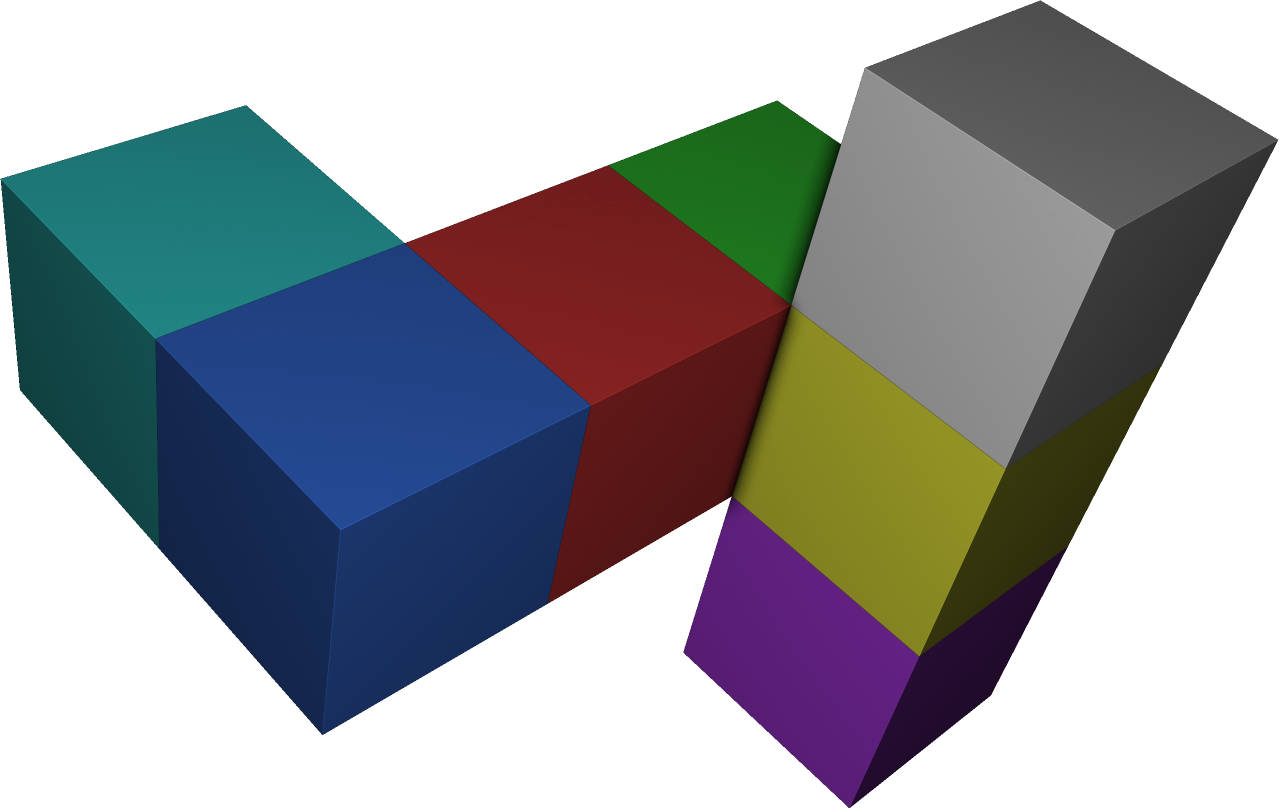}
    \includegraphics[width=0.23\textwidth] {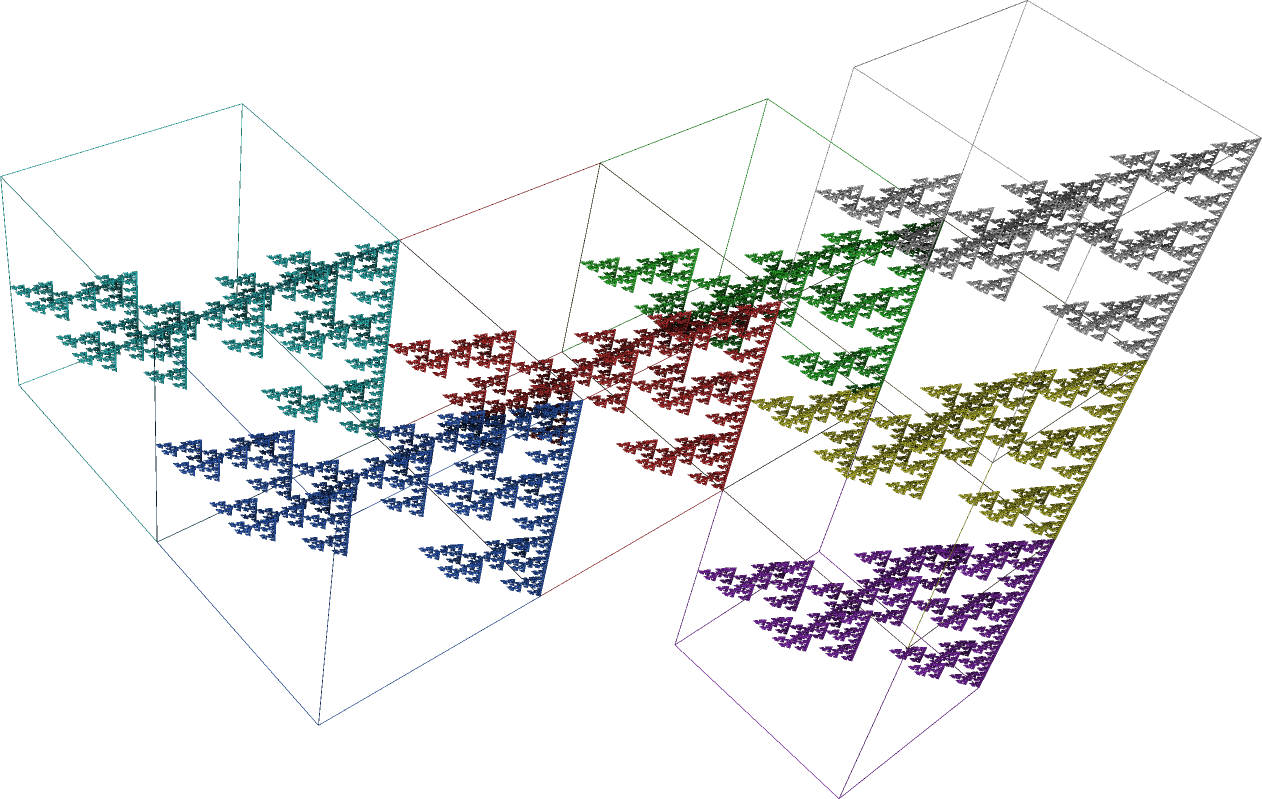} &
    \includegraphics[width=0.23\textwidth] {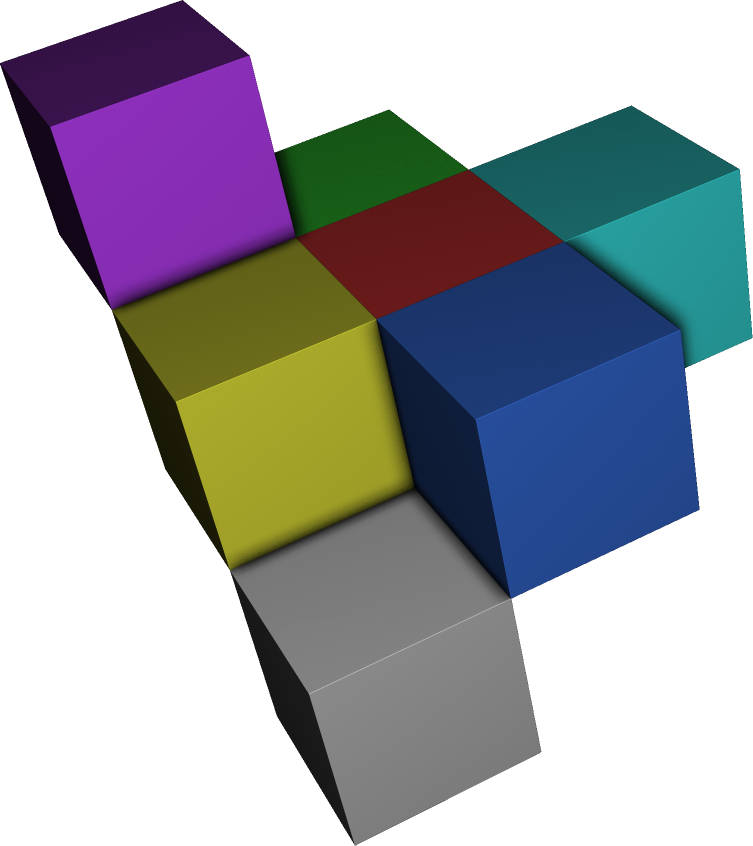}
    \includegraphics[width=0.23\textwidth] {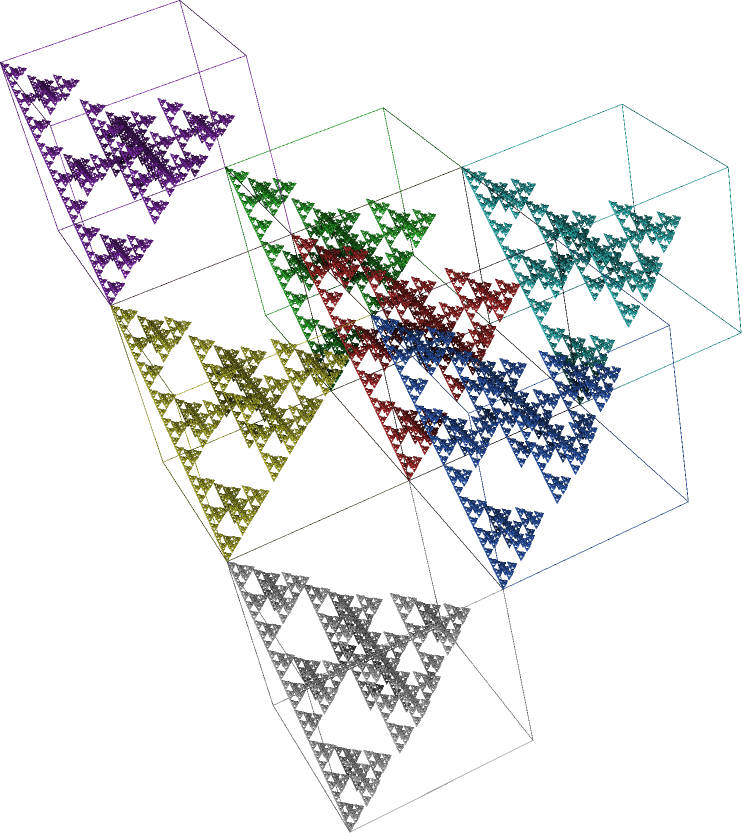} \\
    {\scriptsize$\mD=\{(0,0,2),\ (1,0,2),\ (1,1,0),\ (1,1,1),\ (1,1,2),\ (1,2,0),\ (2,0,2)\}$} & 
    {\scriptsize$\mD=\{(0,0,0),\ (1,0,1),\ (1,1,0),\ (1,1,1),\ (1,1,2),\ (1,2,1),\ (2,0,2)\}$} \\
    \hline
    \includegraphics[width=0.23\textwidth] {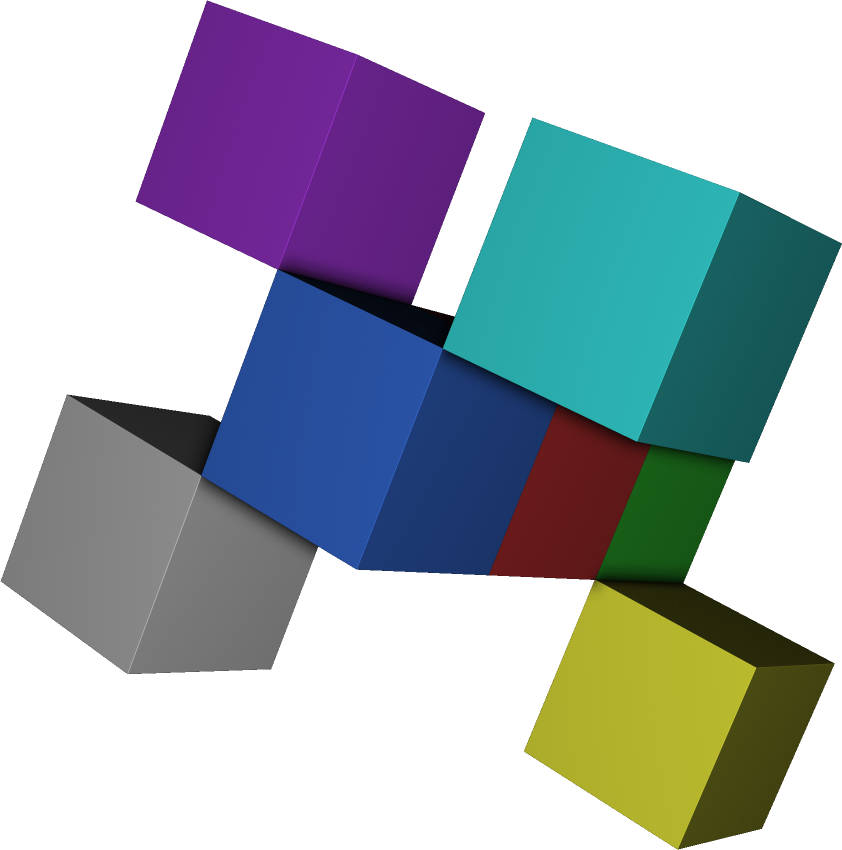}
    \includegraphics[width=0.23\textwidth] {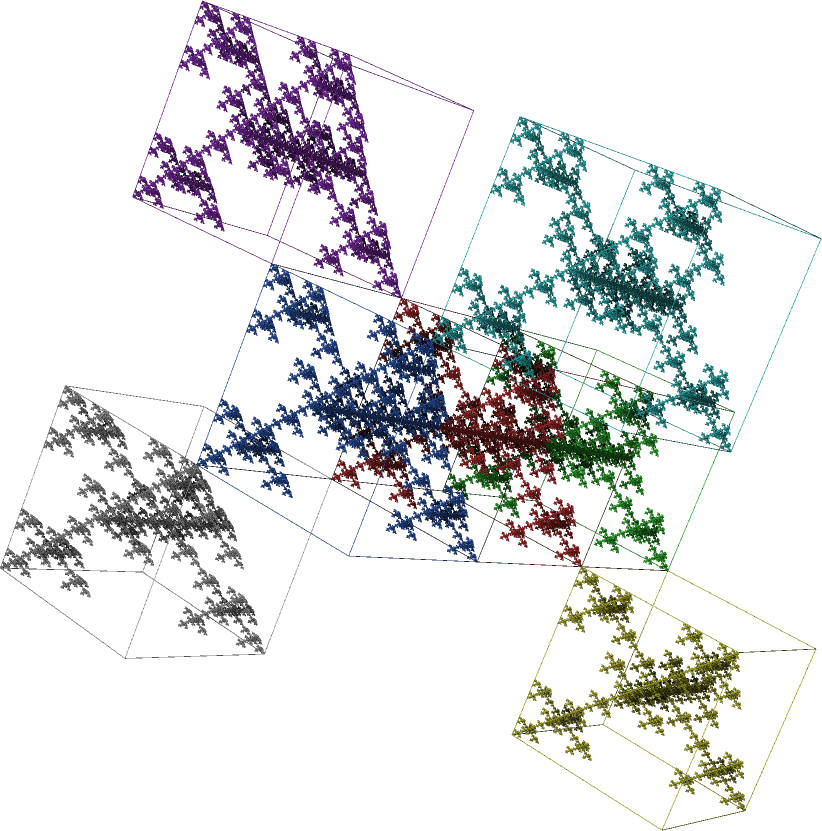} &
    \includegraphics[width=0.23\textwidth] {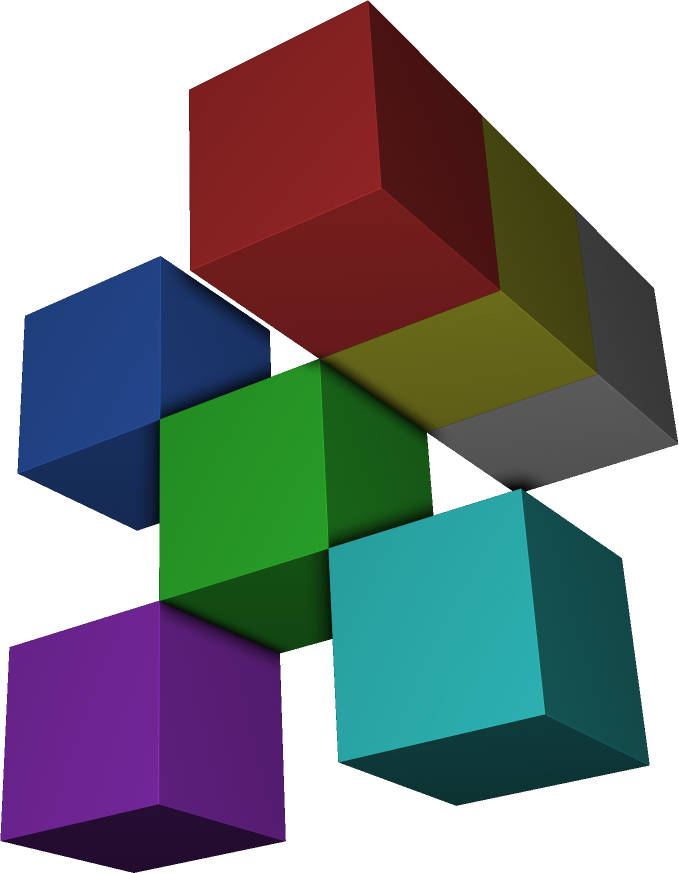}
    \includegraphics[width=0.23\textwidth] {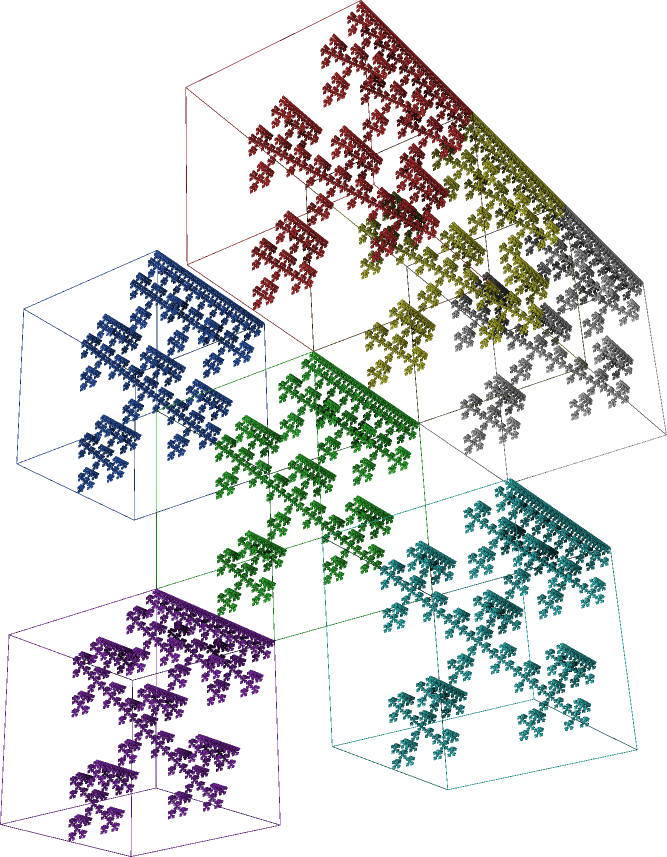} \\
    {\scriptsize$\mD=\{(0,0,2),\ (0,2,0),\ (1,0,1),\ (1,1,1),\ (1,2,1),\ (2,0,0),\ (2,0,2)\}$} & 
    {\scriptsize$\mD=\{(0,0,0),\ (0,1,0),\ (0,1,2),\ (0,2,0),\ (1,1,1),\ (2,1,0),\ (2,1,2)\}$} \\
    \hline
    \includegraphics[width=0.23\textwidth] {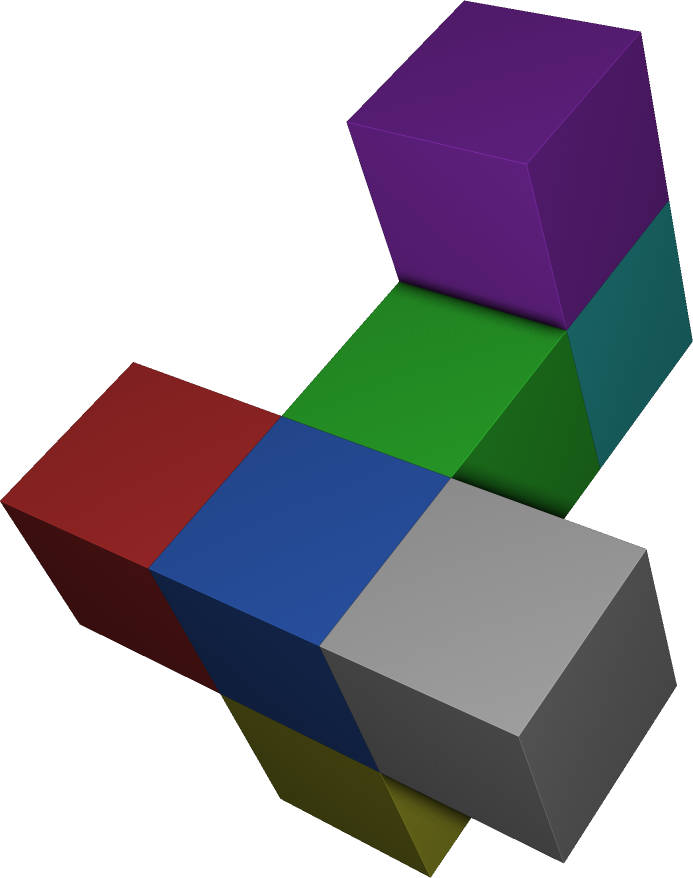}
    \includegraphics[width=0.23\textwidth] {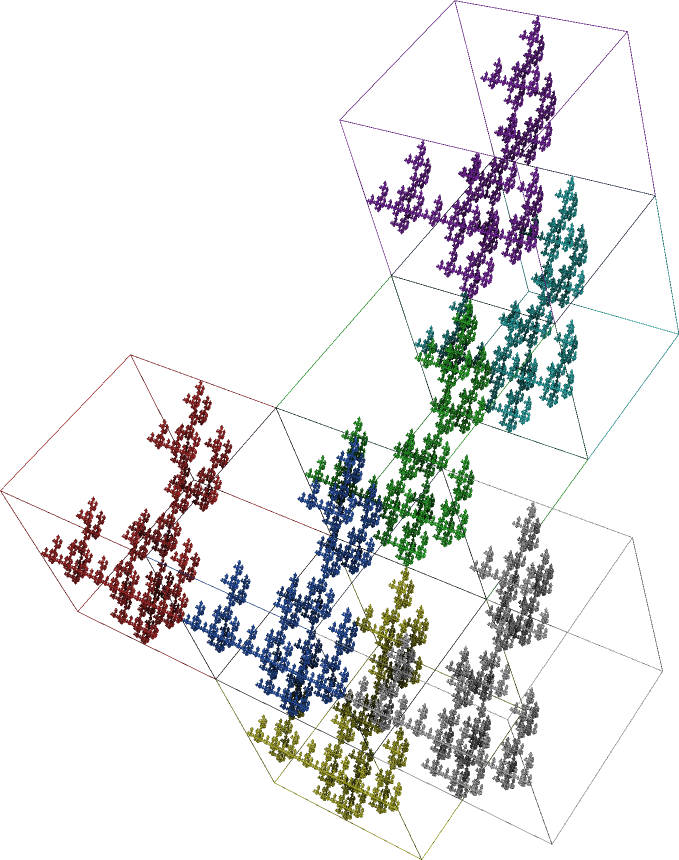} &
    \includegraphics[width=0.23\textwidth] {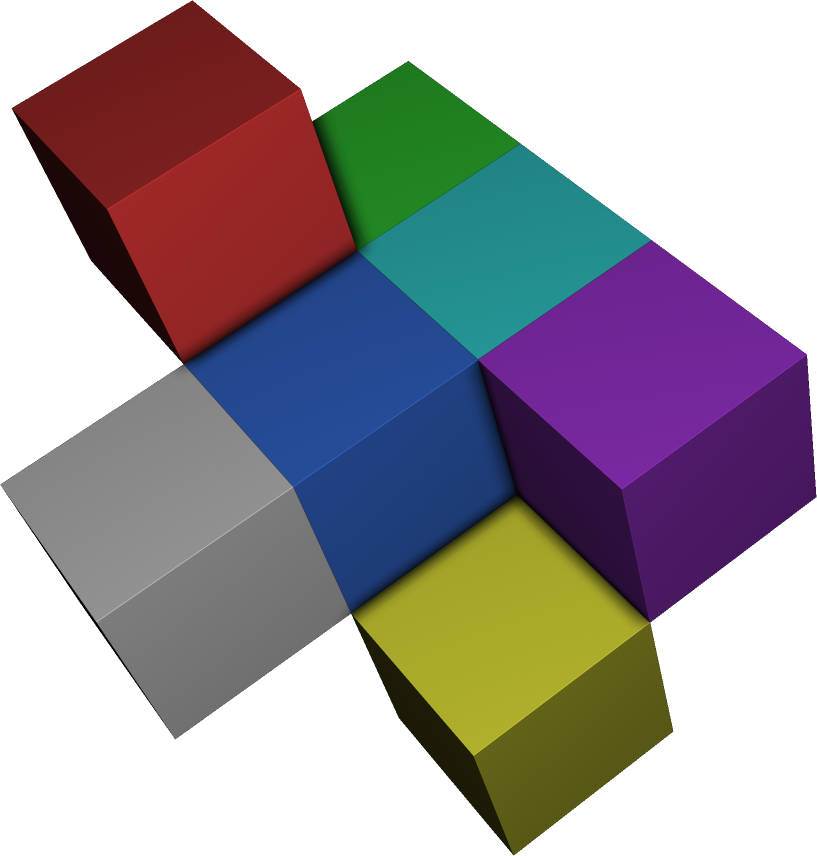}
    \includegraphics[width=0.23\textwidth] {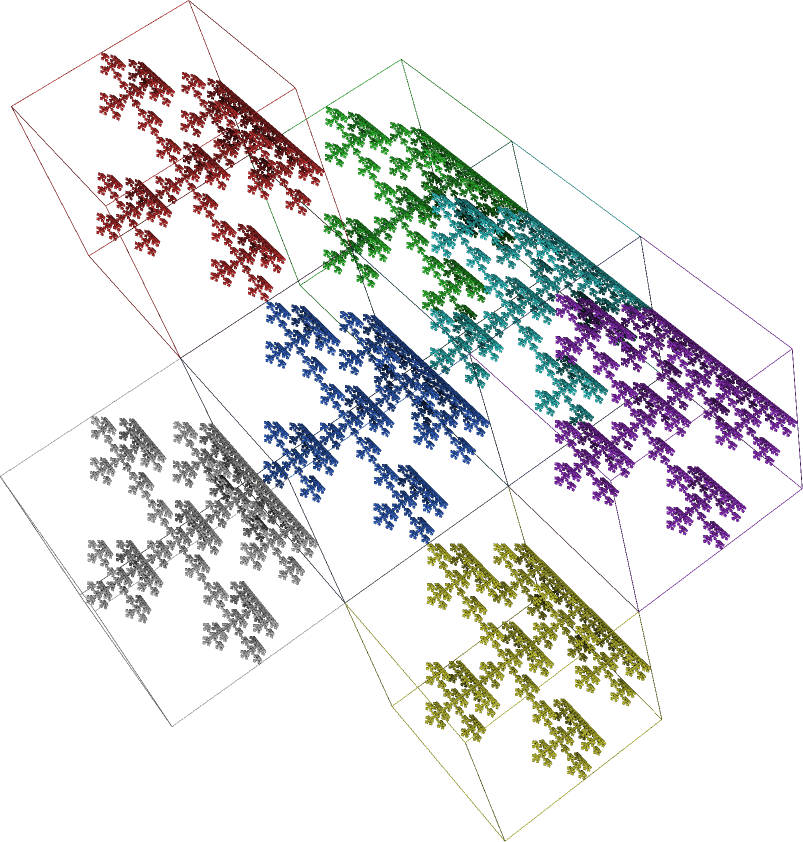} \\
    {\scriptsize$\mD=\{(0,0,1),\ (0,1,0),\ (0,1,1),\ (0,2,1),\ (1,1,1),\ (2,1,1),\ (2,1,2)\}$} & 
    {\scriptsize$\mD=\{(0,1,1),\ (1,0,2),\ (1,1,1),\ (1,2,0),\ (2,1,0),\ (2,1,1),\ (2,1,2)\}$} \\
    \hline
    \includegraphics[width=0.23\textwidth] {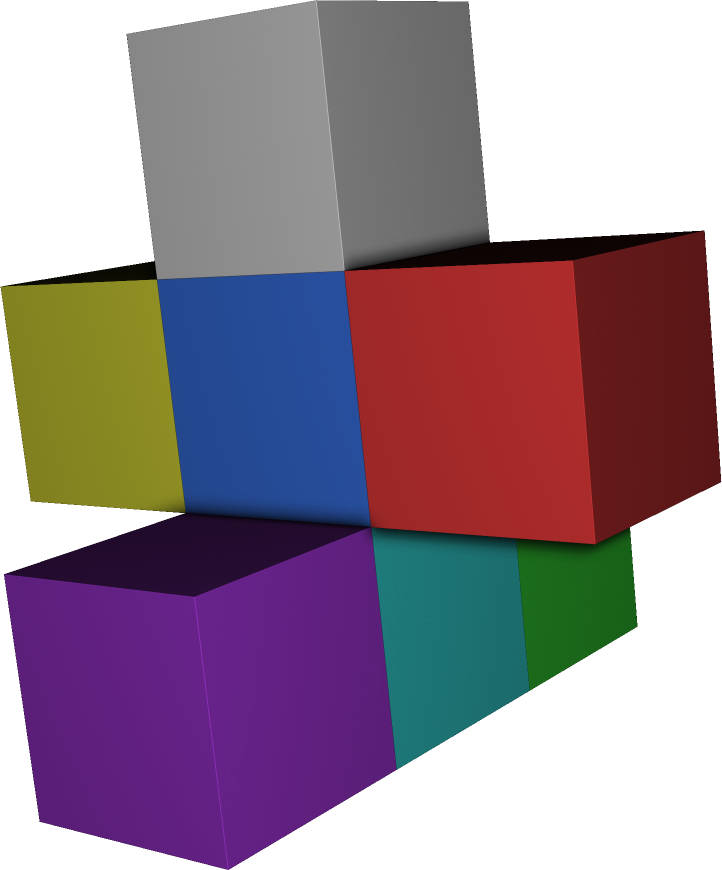}
    \includegraphics[width=0.23\textwidth] {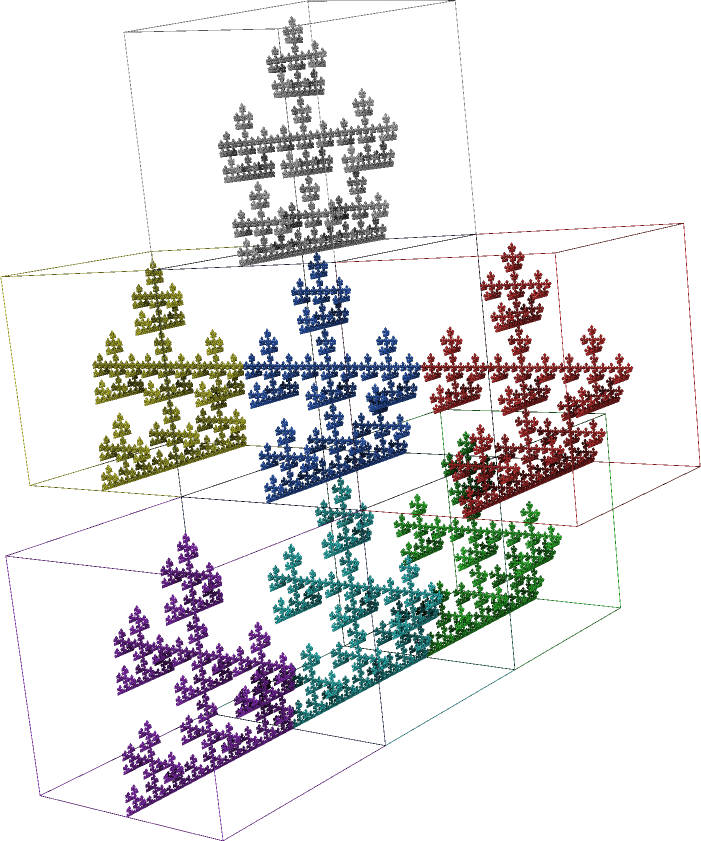} &
    \includegraphics[width=0.23\textwidth] {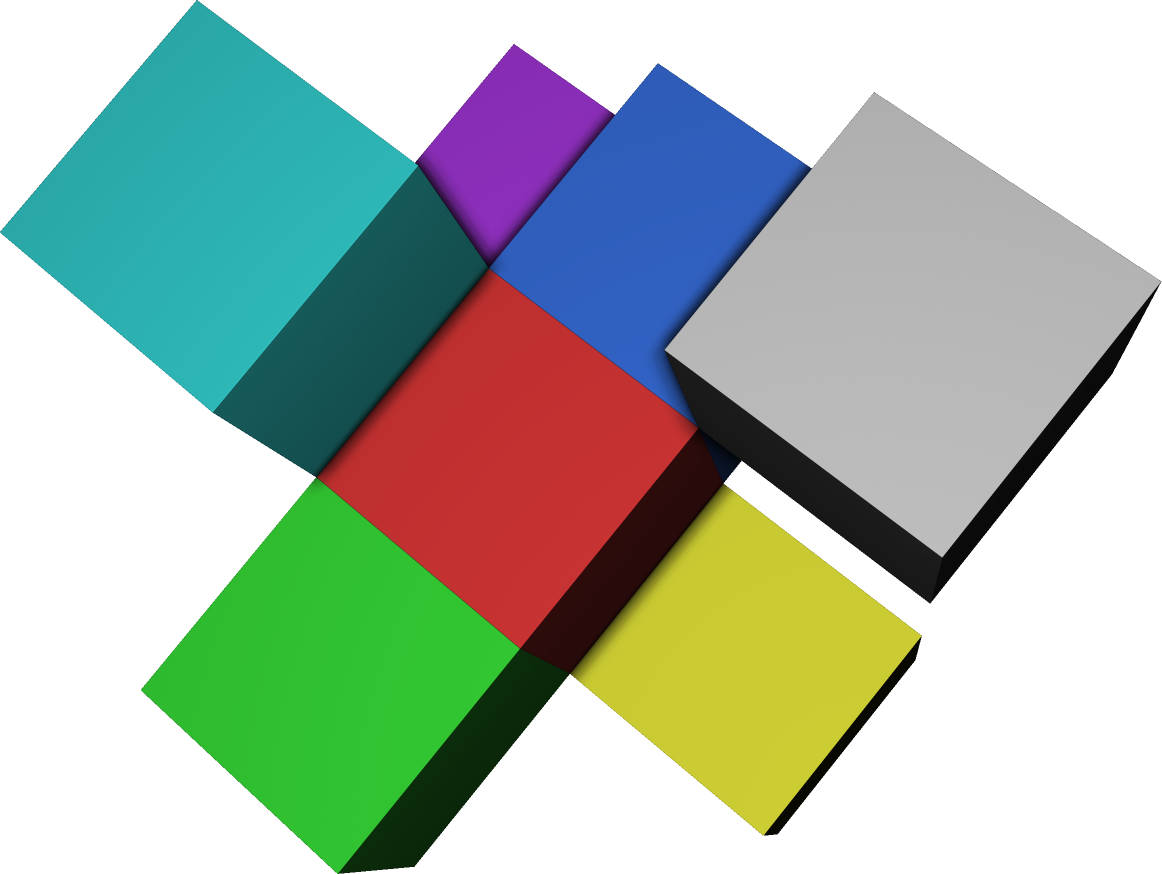}
    \includegraphics[width=0.23\textwidth] {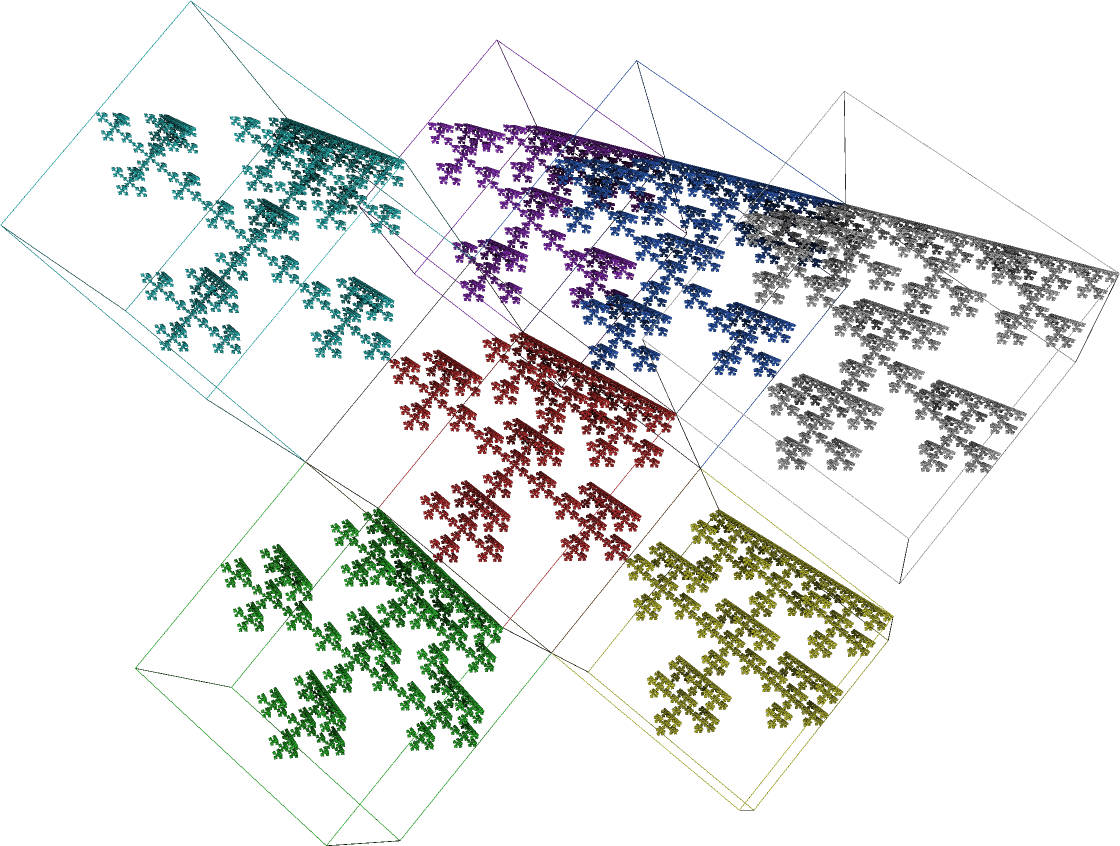} \\
    {\scriptsize$\mD=\{(0,1,1),\ (1,0,1),\ (1,1,1),\ (1,2,1),\ (2,1,0),\ (2,1,1),\ (2,1,2)\}$} & 
    {\scriptsize$\mD=\{(0,0,0),\ (0,2,1),\ (1,1,0),\ (1,1,1),\ (1,1,2),\ (2,0,1),\ (2,2,0)\}$} \\
    \hline
    \includegraphics[width=0.23\textwidth] {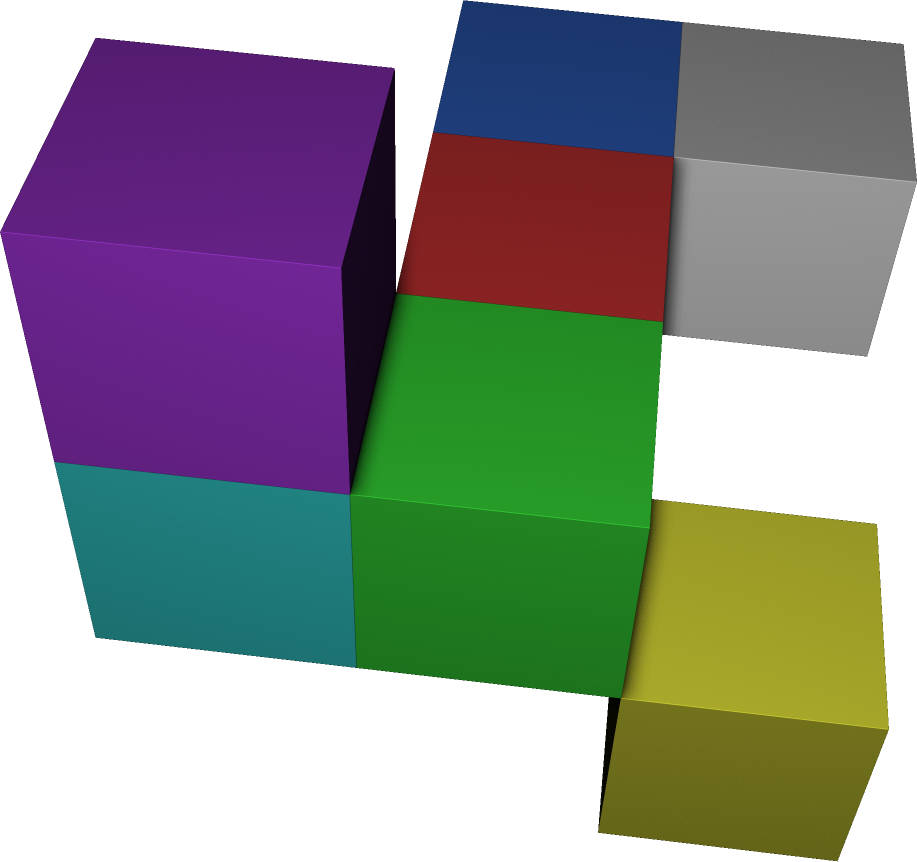}
    \includegraphics[width=0.23\textwidth] {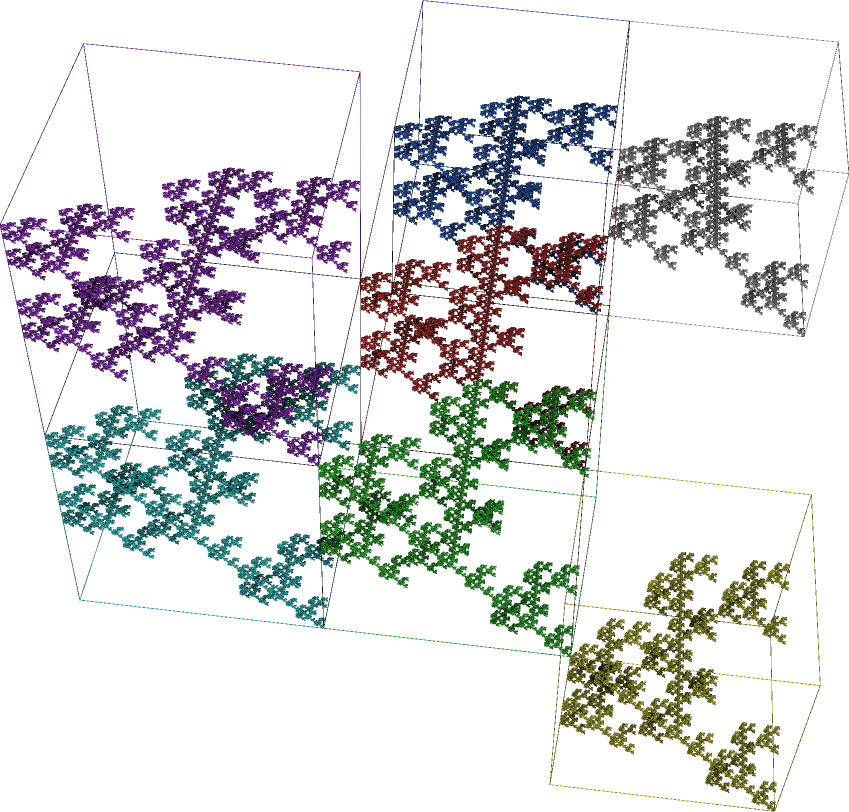} &
    \includegraphics[width=0.23\textwidth] {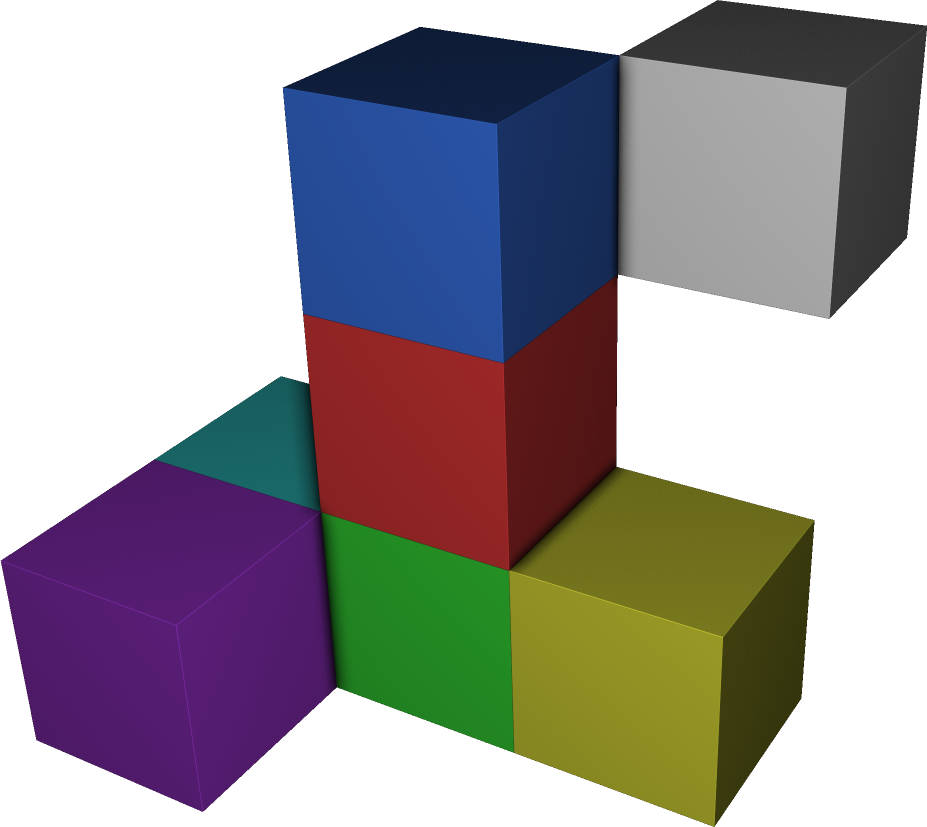}
    \includegraphics[width=0.23\textwidth] {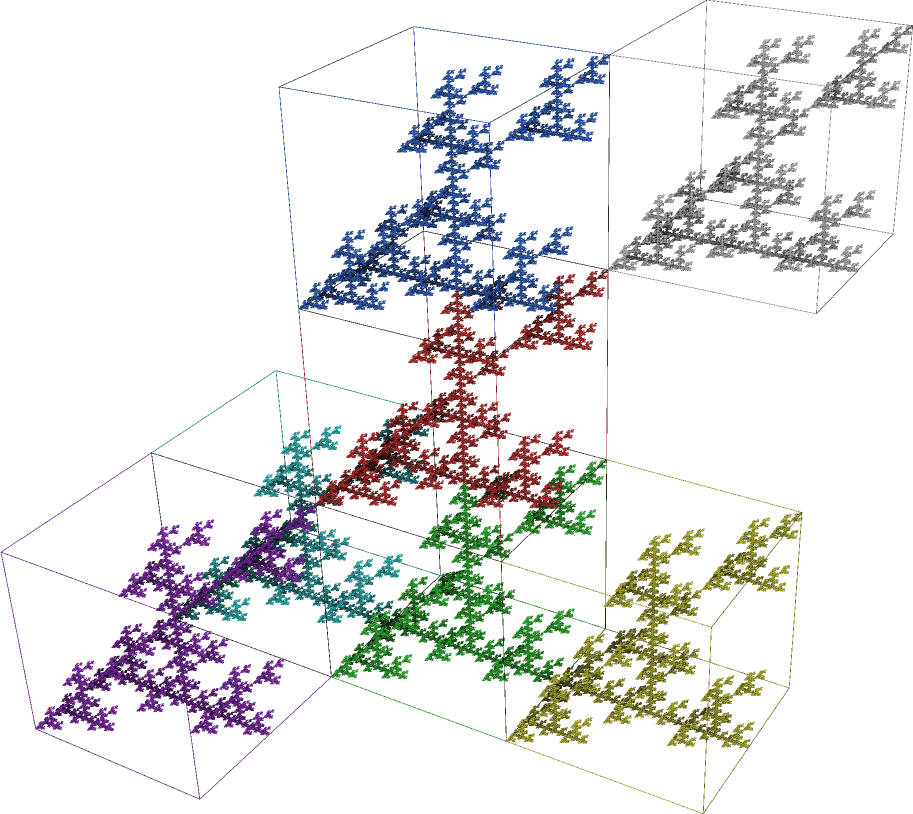} \\
    {\scriptsize$\mD=\{(0,0,1),\ (0,2,0),\ (1,0,1),\ (1,1,1),\ (1,2,1),\ (2,2,1),\ (2,2,2)\}$} & 
    {\scriptsize$\mD=\{(0,0,0),\ (0,2,1),\ (1,0,1),\ (1,1,1),\ (1,2,1),\ (2,2,1),\ (2,2,2)\}$} \\
    \hline
    \includegraphics[width=0.23\textwidth] {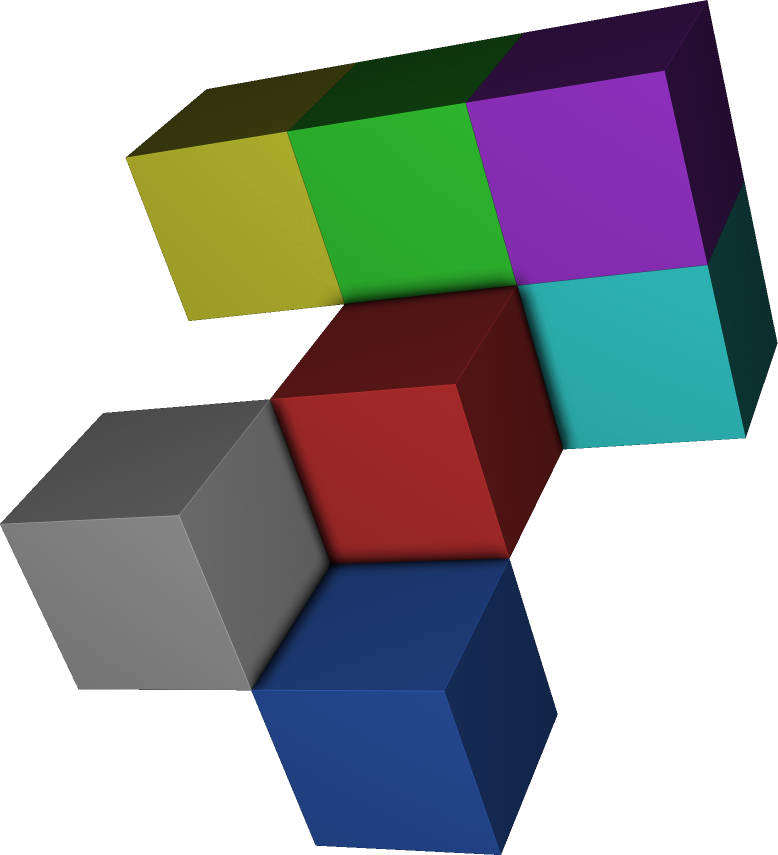}
    \includegraphics[width=0.23\textwidth] {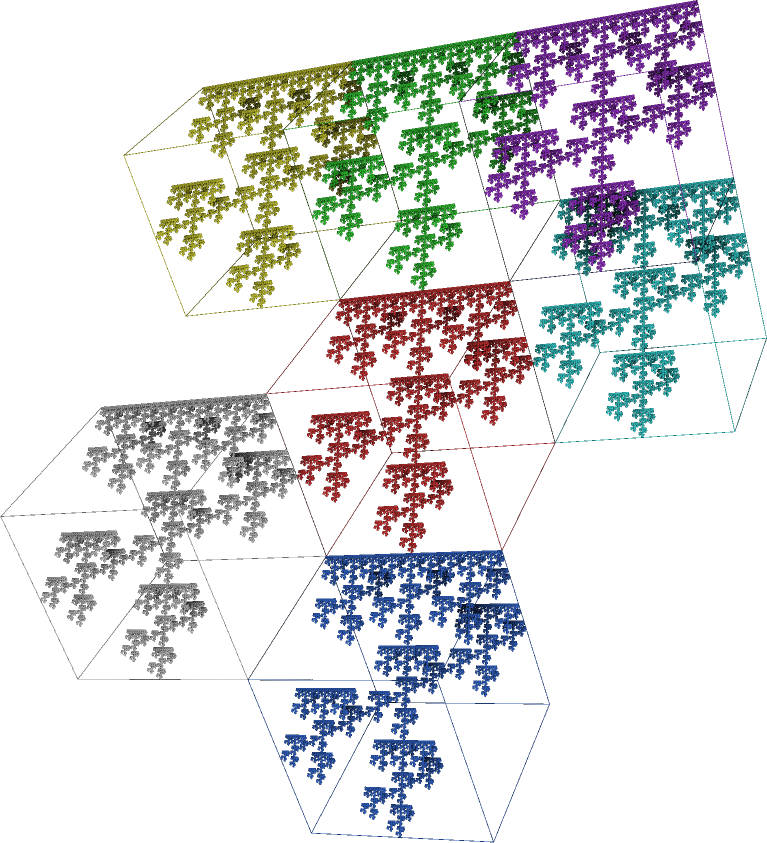} &
    \includegraphics[width=0.23\textwidth] {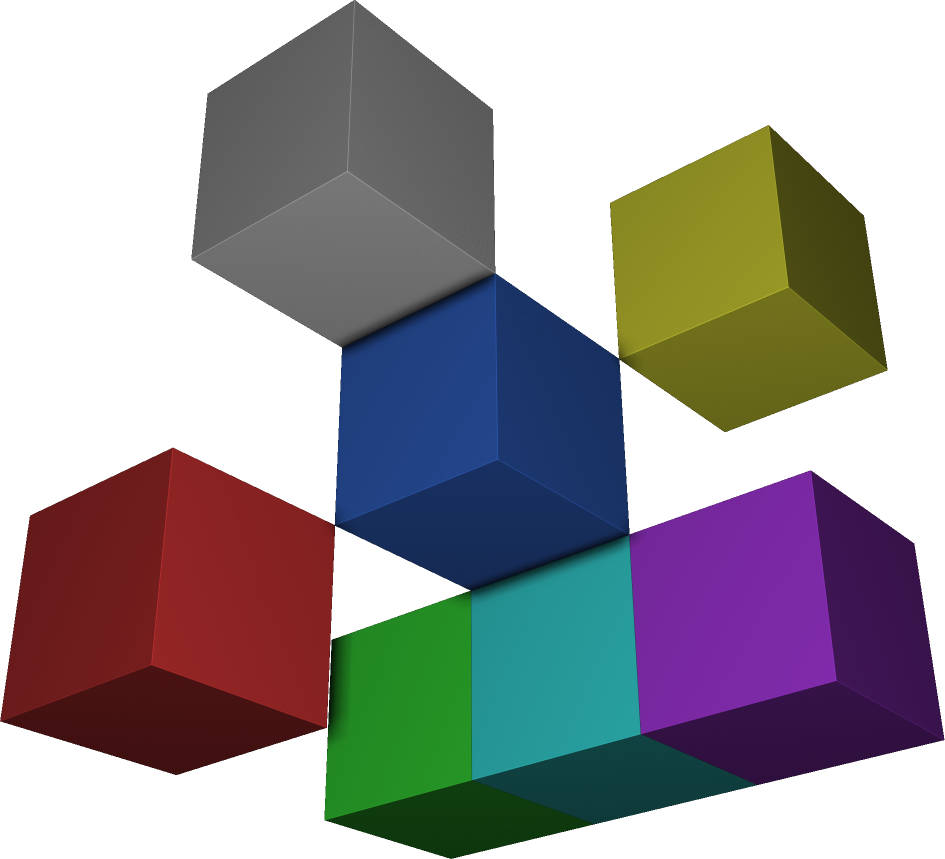}
    \includegraphics[width=0.23\textwidth] {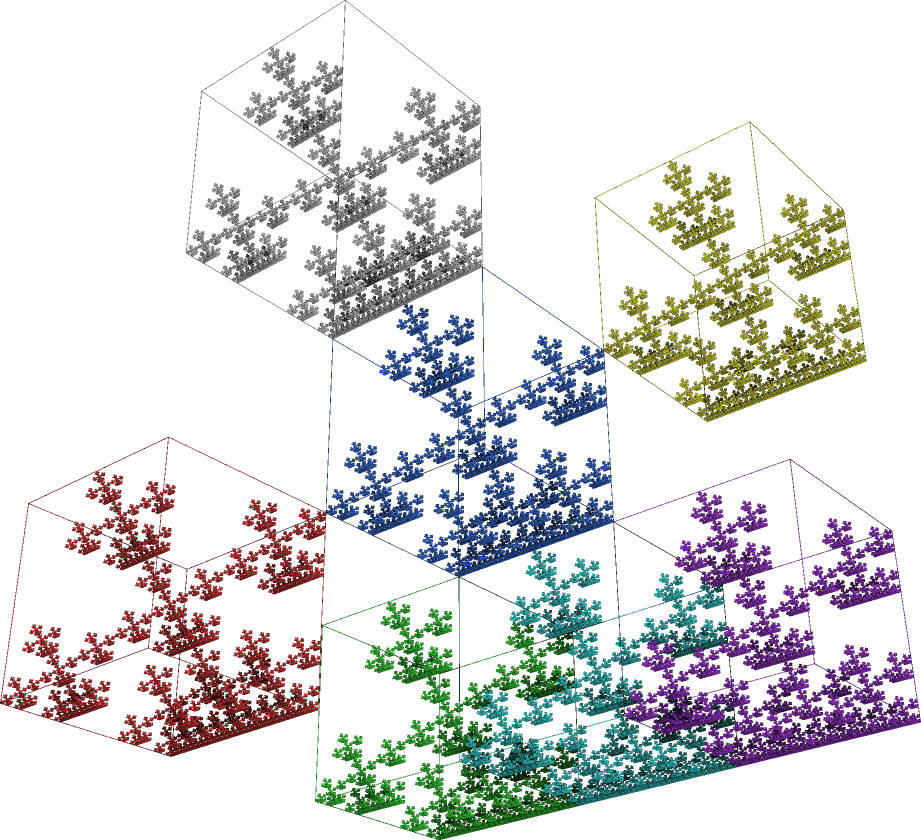} \\
    {\scriptsize$\mD=\{(0,0,1),\ (0,2,2),\ (1,0,0),\ (1,1,1),\ (1,2,2),\ (2,2,1),\ (2,2,2)\}$} & 
    {\scriptsize$\mD=\{(0,0,1),\ (0,2,2),\ (1,1,1),\ (2,0,0),\ (2,2,0),\ (2,2,1),\ (2,2,2)\}$} \\
    \hline
\end{longtable}

\newpage
\begin{longtable}{|p{0.48\textwidth}|p{0.48\textwidth}|}
    \caption{Dendrites corresponding to graph $7_{5}$ ($N=3$)}\label{tab:d4}\\
    \hline
    \multicolumn{2}{|c|}{\includegraphics{den4.pdf}} \\
    \hline
    \includegraphics[width=0.23\textwidth] {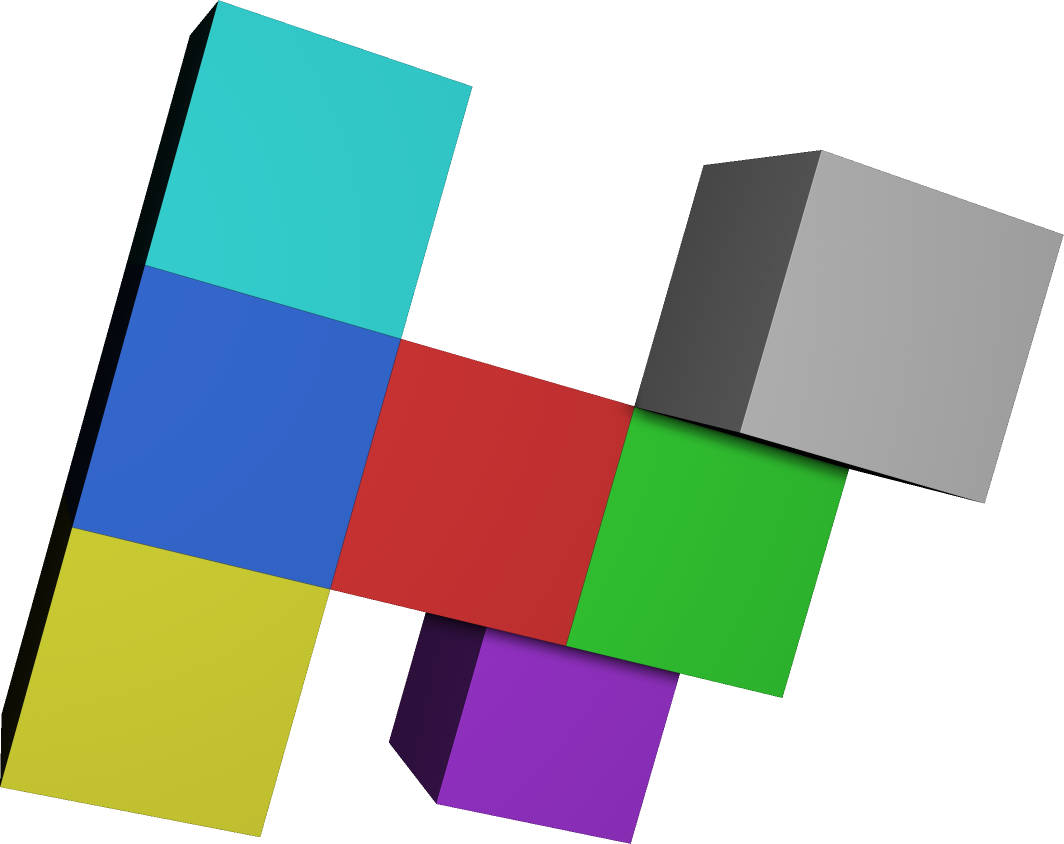}
    \includegraphics[width=0.23\textwidth] {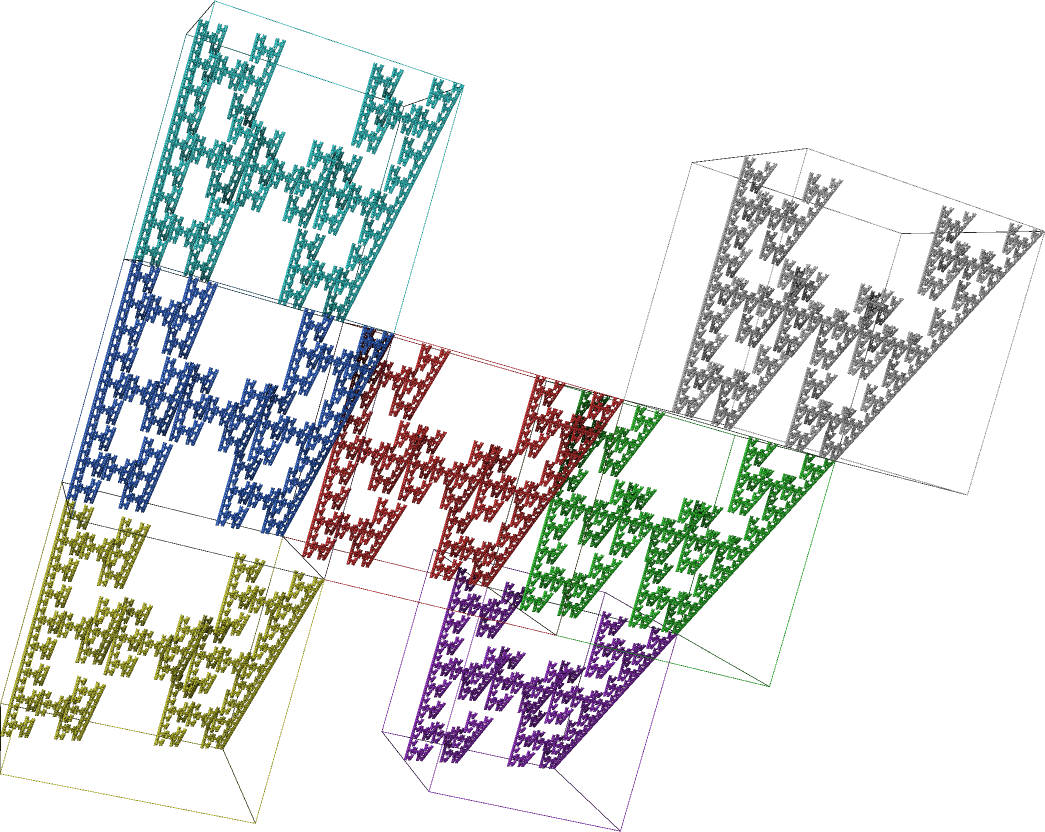} &
    \includegraphics[width=0.23\textwidth] {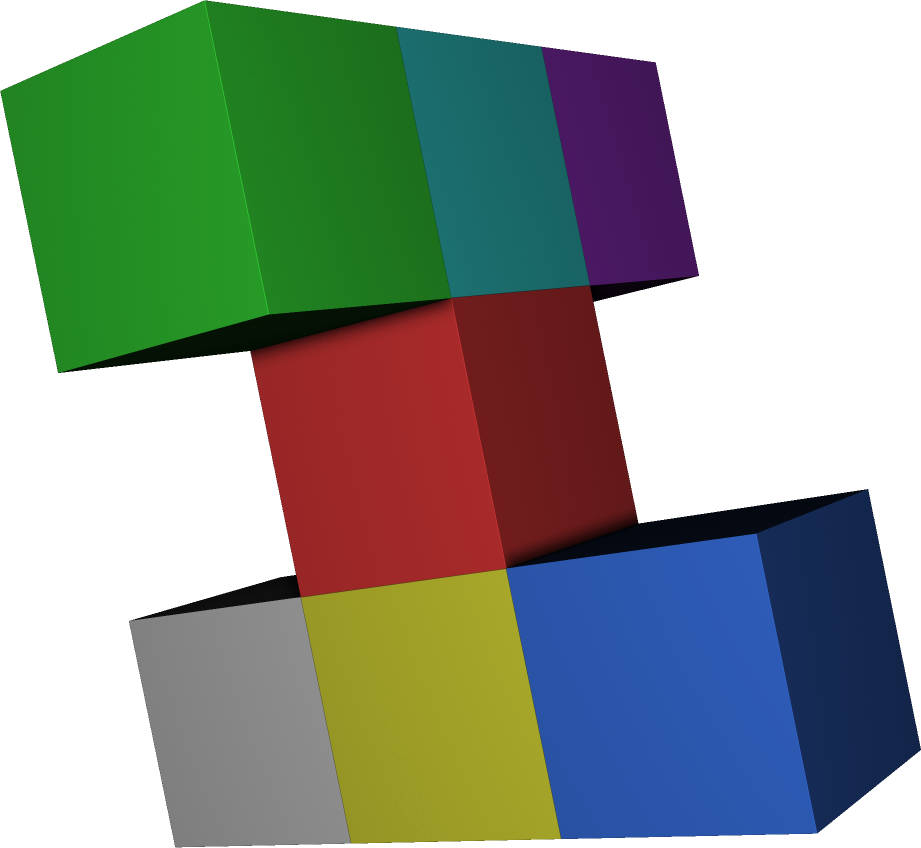}
    \includegraphics[width=0.23\textwidth] {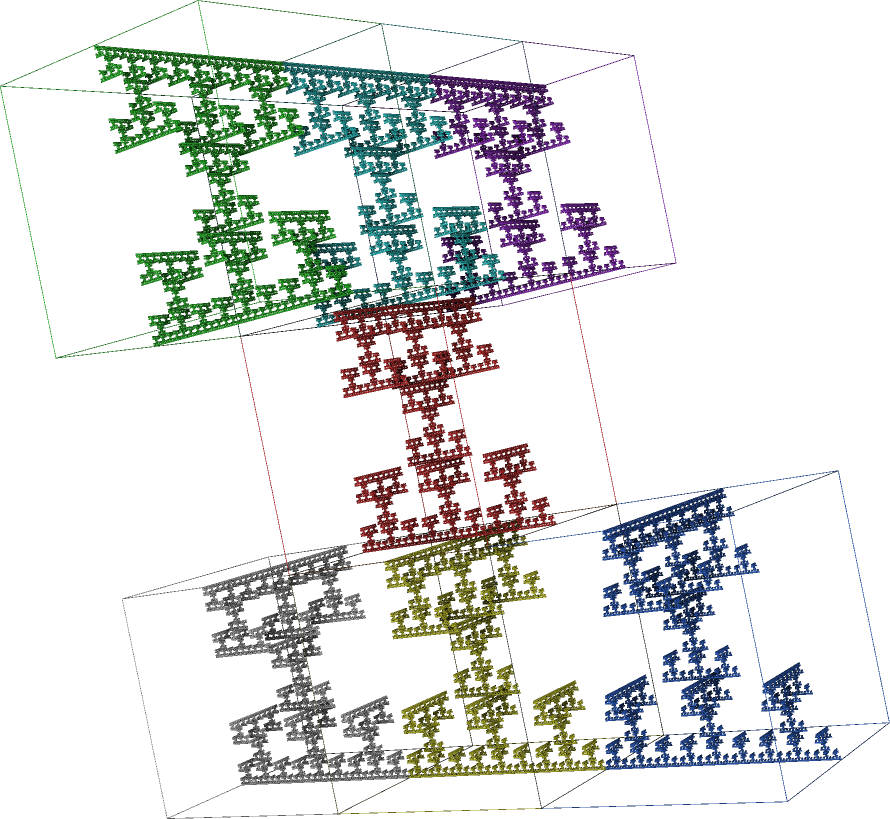} \\
    {\scriptsize$\mD=\{(0,2,2),\ (1,0,0),\ (1,1,0),\ (1,1,1),\ (1,1,2),\ (1,2,0),\ (2,0,2)\}$} & 
    {\scriptsize$\mD=\{(0,0,1),\ (0,1,1),\ (0,2,1),\ (1,1,1),\ (2,1,0),\ (2,1,1),\ (2,1,2)\}$} \\
    \hline
    \includegraphics[width=0.23\textwidth] {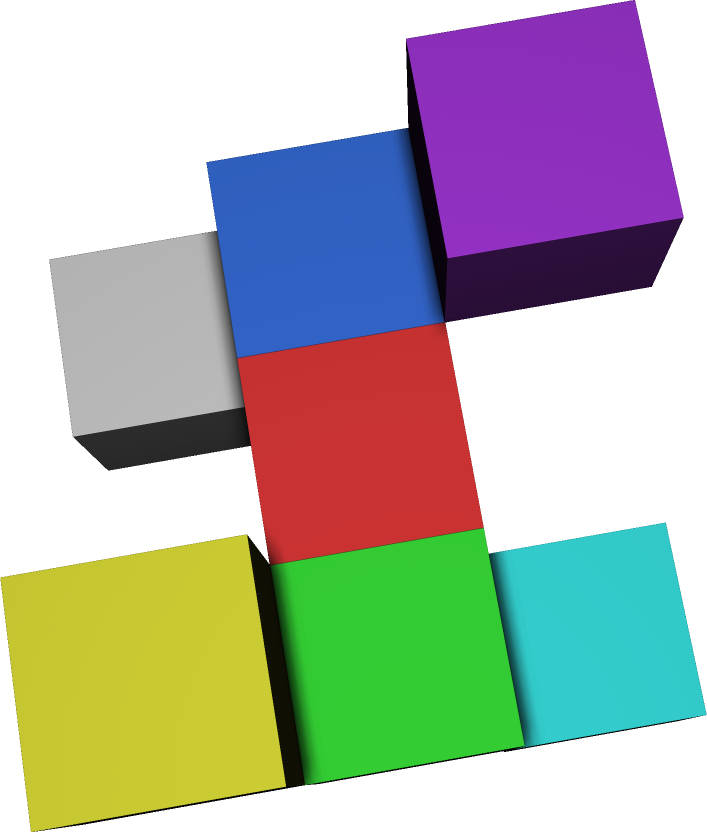}
    \includegraphics[width=0.23\textwidth] {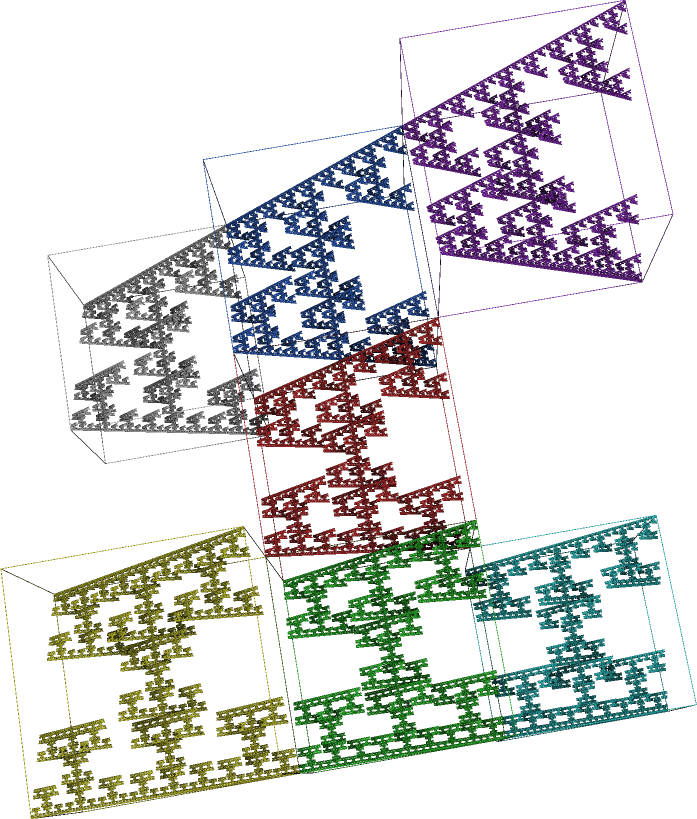} &\\
    {\scriptsize$\mD=\{(0,0,0),\ (0,2,2),\ (1,1,0),\ (1,1,1),\ (1,1,2),\ (2,0,2),\ (2,2,0)\}$} &\\
    \hline
\end{longtable}

\begin{longtable}{|p{0.48\textwidth}|p{0.48\textwidth}|}
    \caption{Dendrite corresponding to graph $7_{6}$ ($N=1$)}\label{tab:d5}\\
    \hline
    \multicolumn{2}{|c|}{\includegraphics{den5.pdf}} \\
    \hline
    \includegraphics[width=0.23\textwidth] {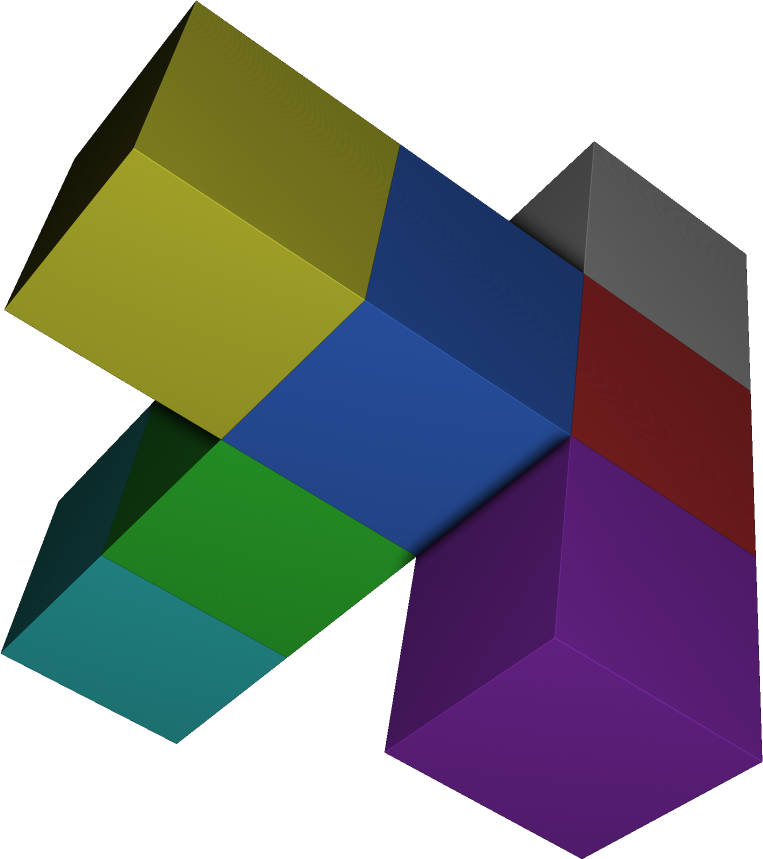}
    \includegraphics[width=0.23\textwidth] {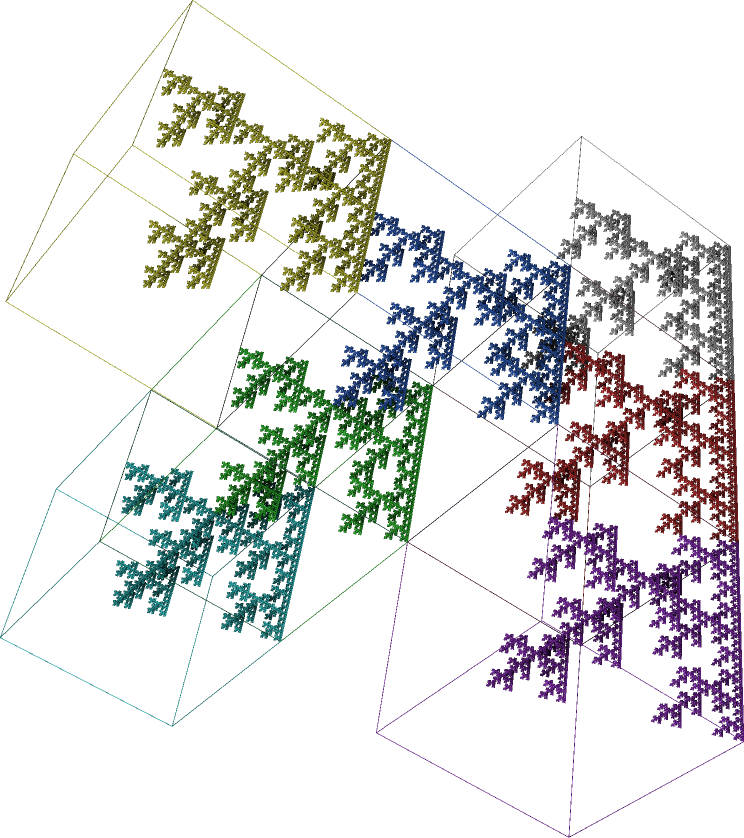} &\\
    {\scriptsize$\mD=\{(0,0,2),\ (1,0,0),\ (1,0,1),\ (1,0,2),\ (1,1,1),\ (1,2,1),\ (2,0,2)\}$} &\\
    \hline
\end{longtable}
\newpage

\begin{longtable}{|p{0.48\textwidth}|p{0.48\textwidth}|}
\caption{Non-dendrites of   type 1 ($N=3$)}\label{tab:n1}\\
    \hline
    \multicolumn{2}{|c|}{\includegraphics{nonden1.pdf}} \\
    \hline
    \includegraphics[width=0.23\textwidth] {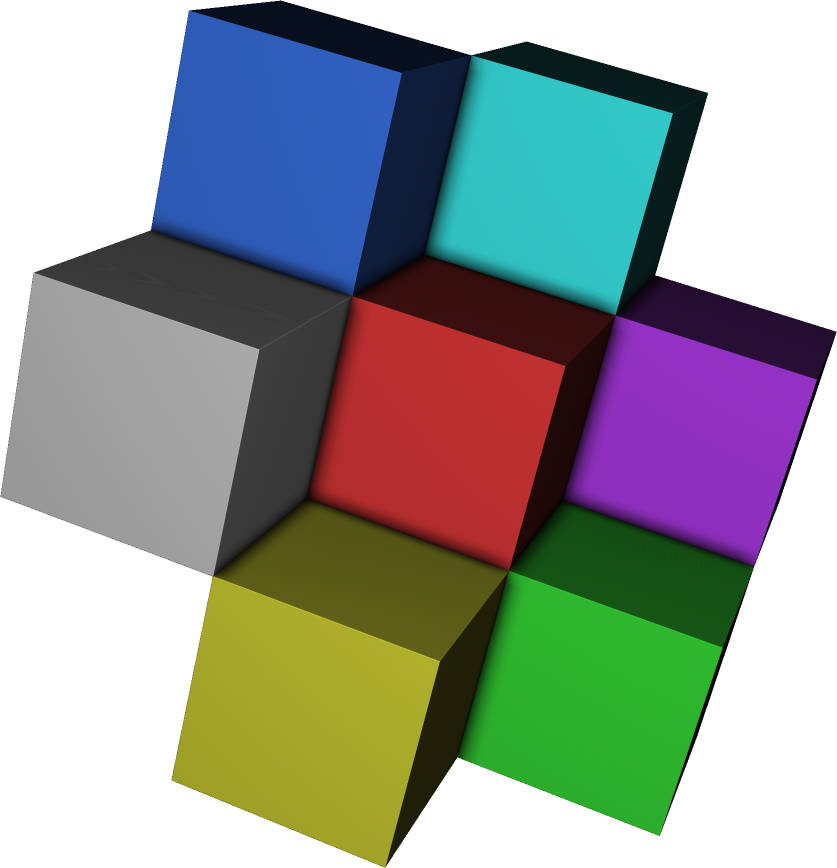}
    \includegraphics[width=0.23\textwidth] {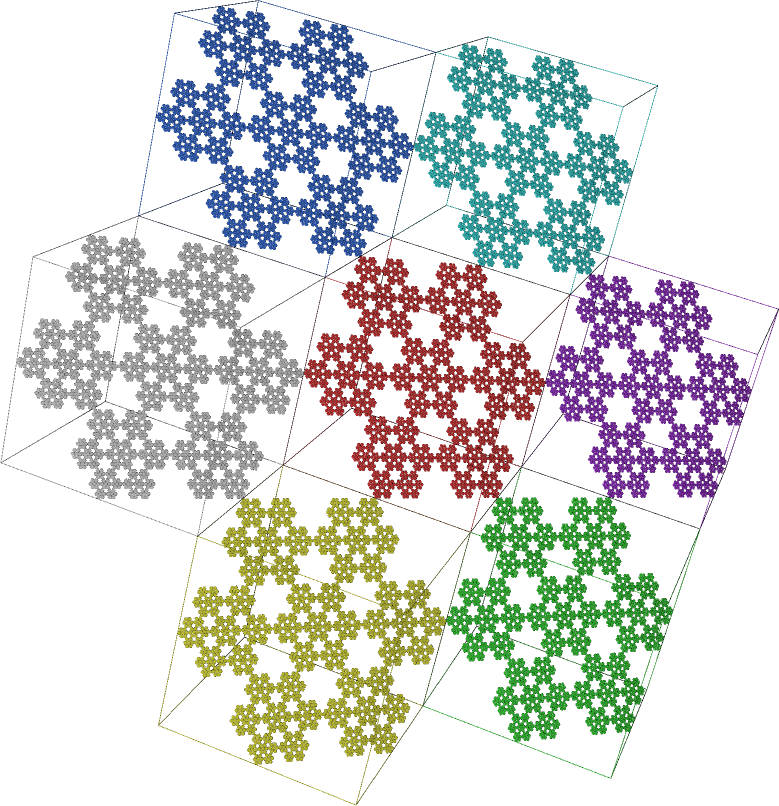} & \includegraphics[width=0.23\textwidth] {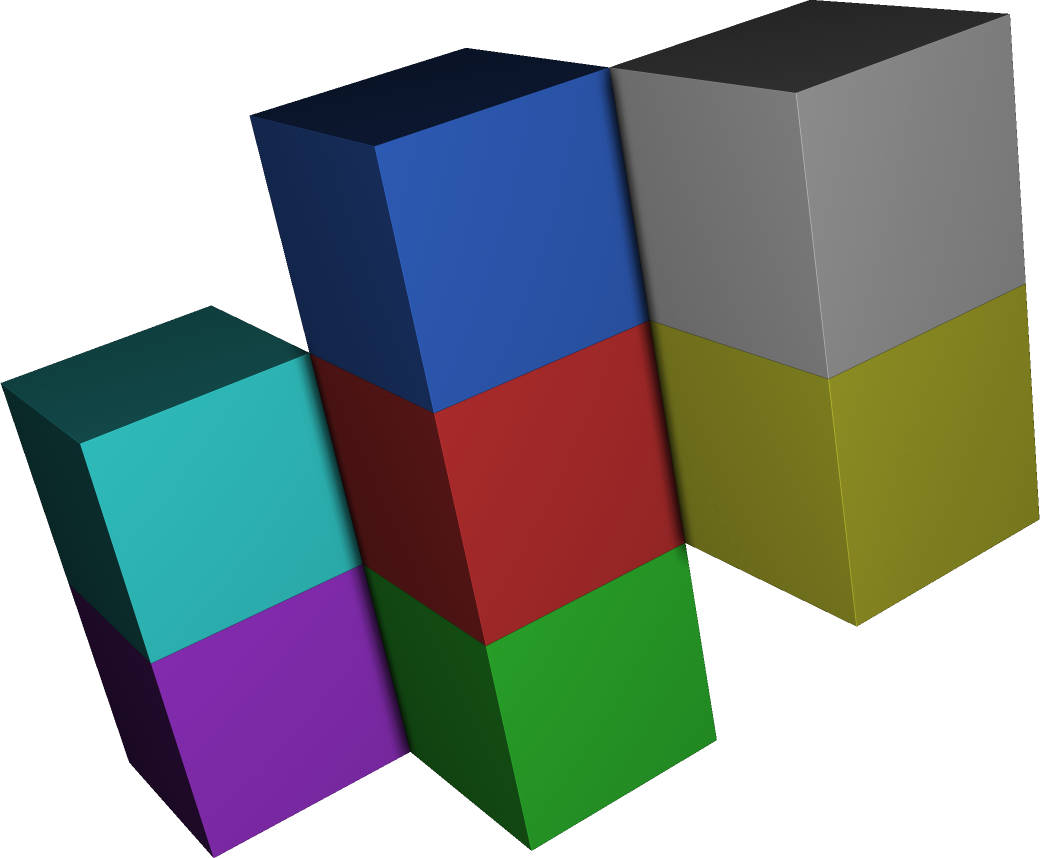}
    \includegraphics[width=0.23\textwidth] {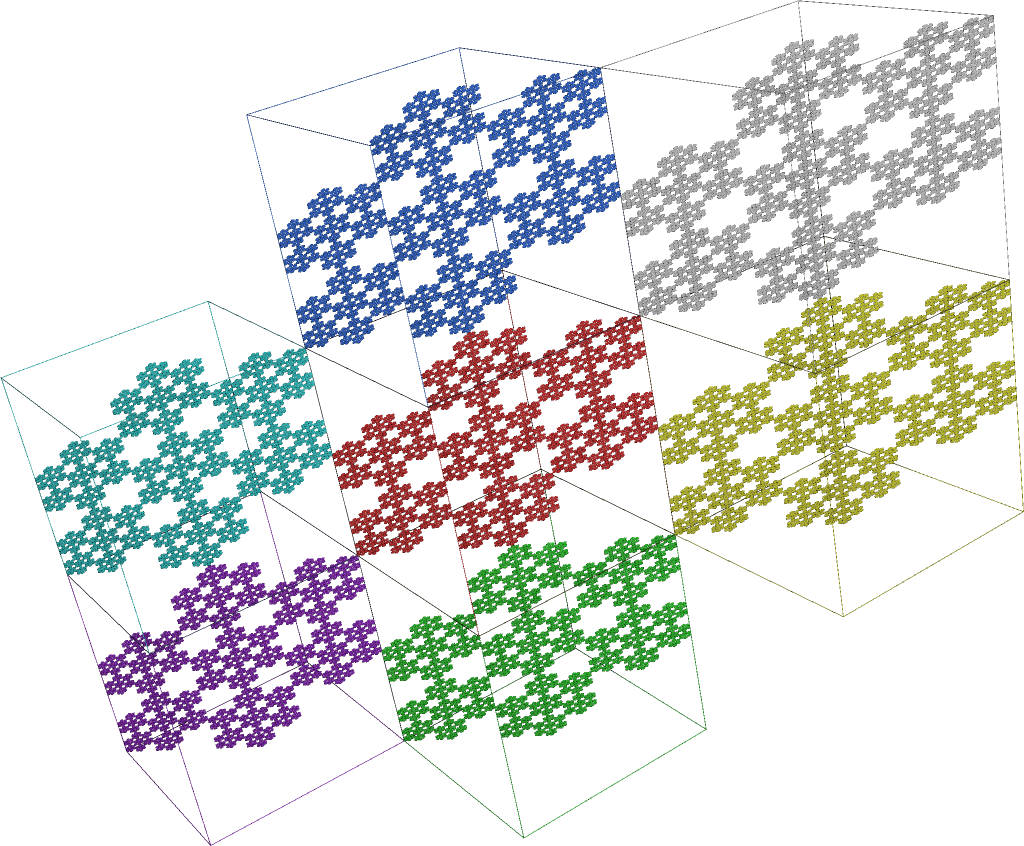} \\
    {\scriptsize$\mD=\{(0,1,2),\ (0,2,1),\ (1,0,2),\ (1,1,1),\ (1,2,0),\ (2,0,1),\ (2,1,0)\}$} & 
    {\scriptsize$\mD=\{(0,0,0),\ (0,0,1),\ (1,1,0),\ (1,1,1),\ (1,1,2),\ (2,2,1),\ (2,2,2)\}$} \\
    \hline
    \includegraphics[width=0.23\textwidth] {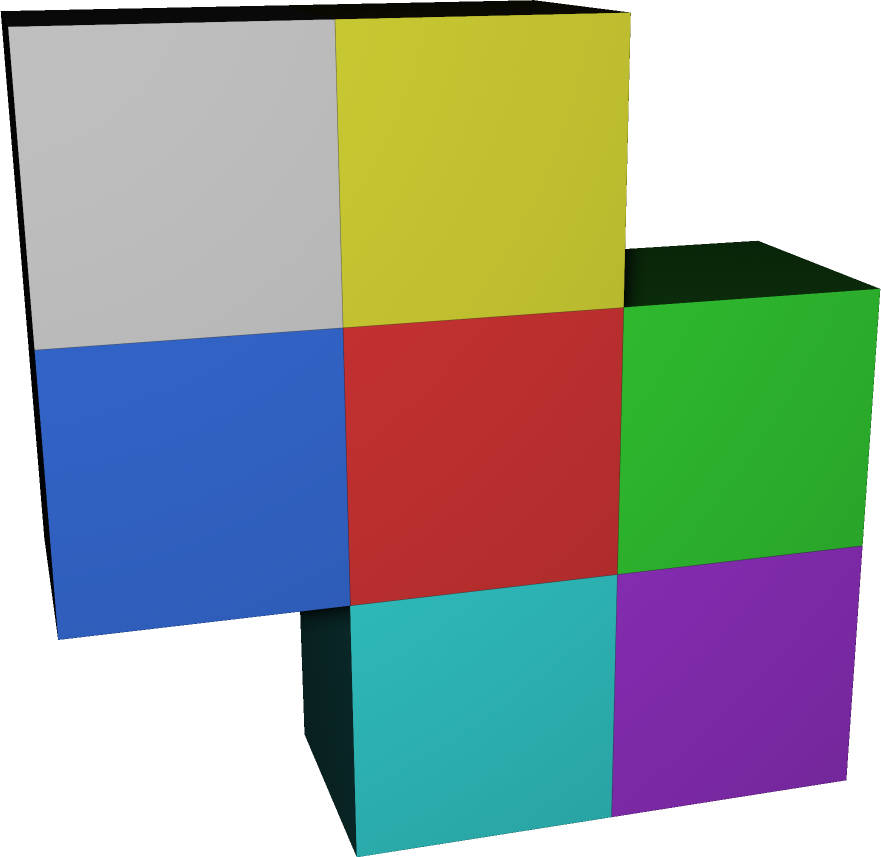}
    \includegraphics[width=0.23\textwidth] {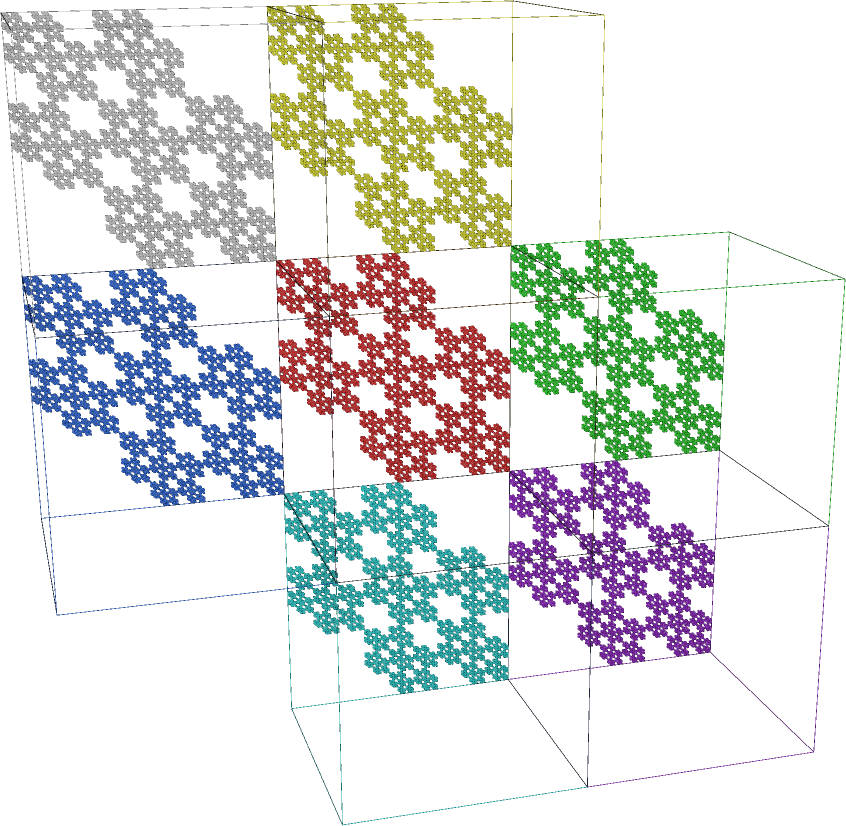} &\\
    {\scriptsize$\mD=\{(0,2,0),\ (0,2,1),\ (1,2,0),\ (1,2,1),\ (1,2,2),\ (2,2,1),\ (2,2,2)\}$} &\\
    \hline

\end{longtable}

\newpage
\begin{longtable}{|p{0.48\textwidth}|p{0.48\textwidth}|}
\caption{Non-dendrites of   type 2 ($N=8$) }\label{tab:n1}\\
    \hline
    \multicolumn{2}{|c|}{\includegraphics{nonden2.pdf}} \\
    \hline
    \includegraphics[width=0.23\textwidth] {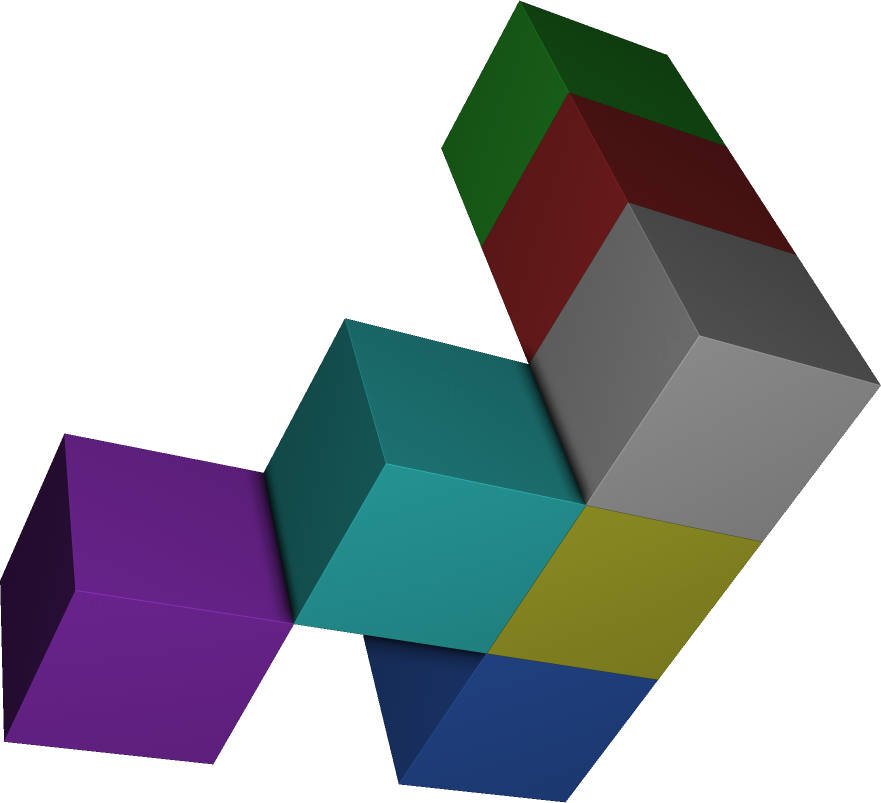}
    \includegraphics[width=0.23\textwidth] {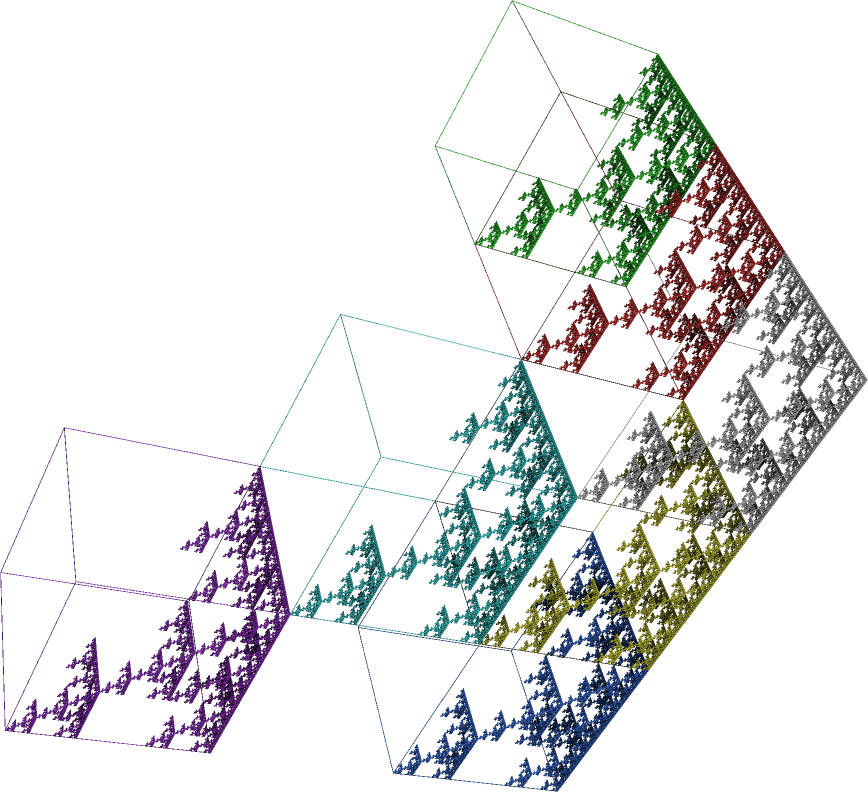} & \includegraphics[width=0.23\textwidth] {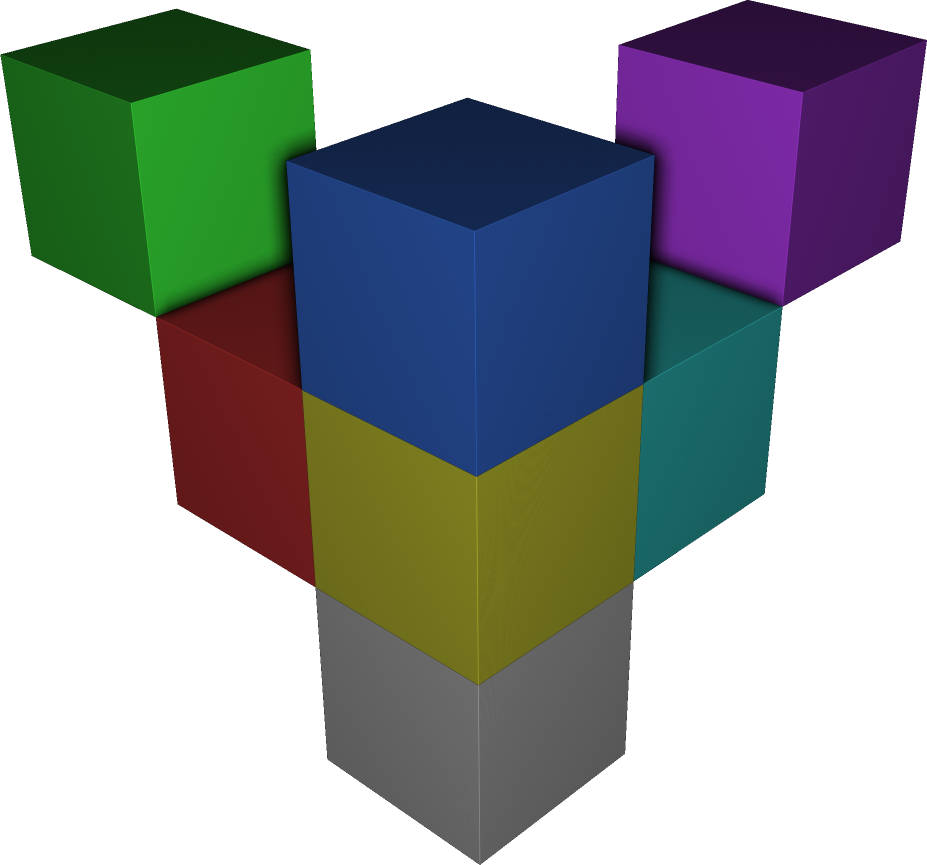}
    \includegraphics[width=0.23\textwidth] {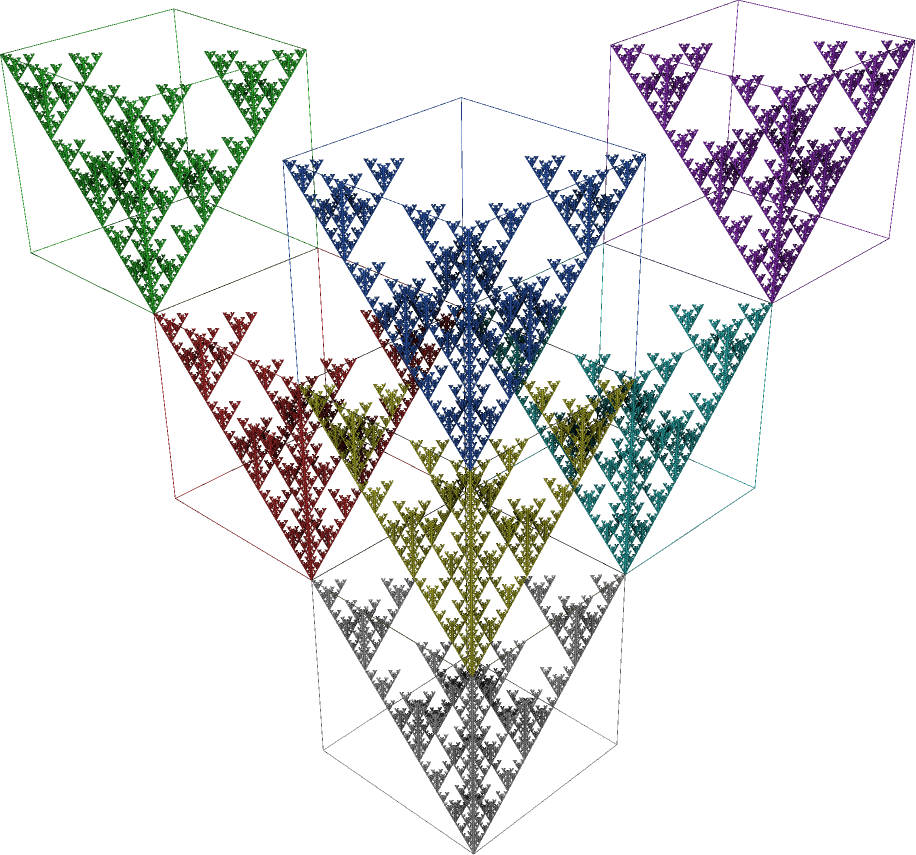} \\
    {\scriptsize$\mD=\{(0,0,0),\ (0,0,1),\ (0,0,2),\ (0,1,0),\ (0,2,0),\ (1,0,1),\ (2,0,2)\}$} & 
    {\scriptsize$\mD=\{(0,0,0),\ (0,0,1),\ (0,0,2),\ (0,1,1),\ (0,2,2),\ (1,0,1),\ (2,0,2)\}$} \\
    \hline
    \includegraphics[width=0.23\textwidth] {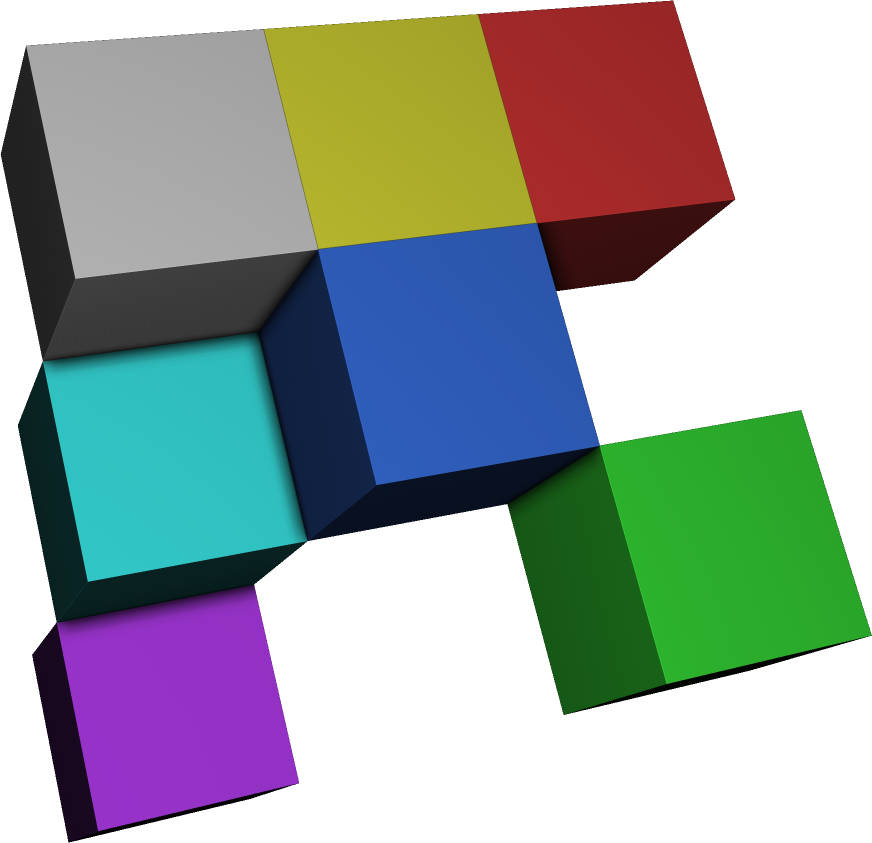}
    \includegraphics[width=0.23\textwidth] {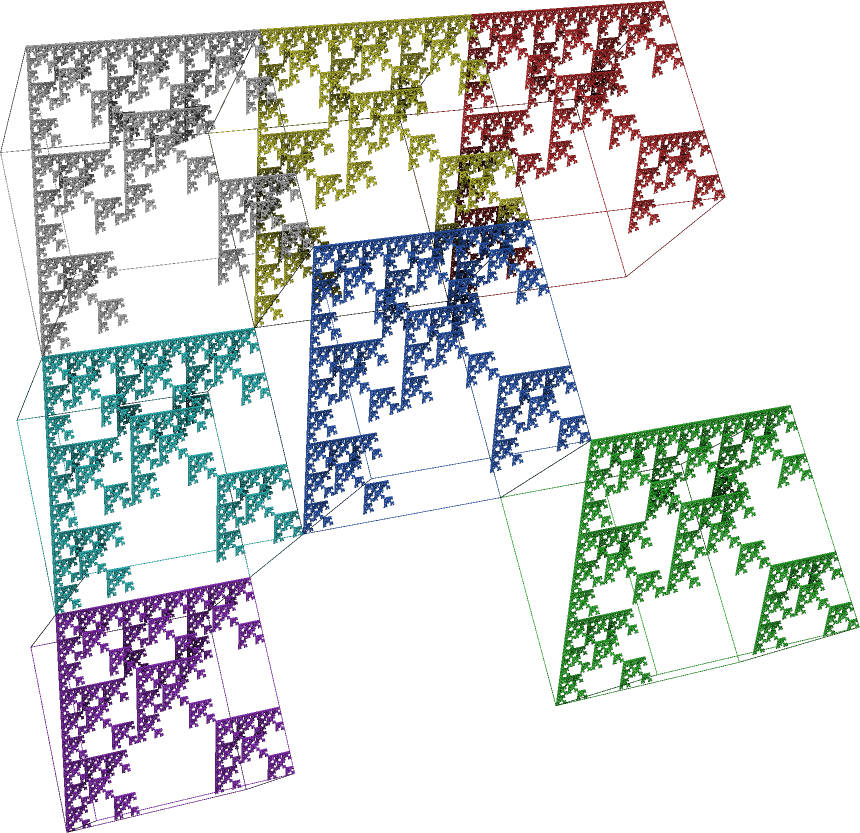} & \includegraphics[width=0.23\textwidth] {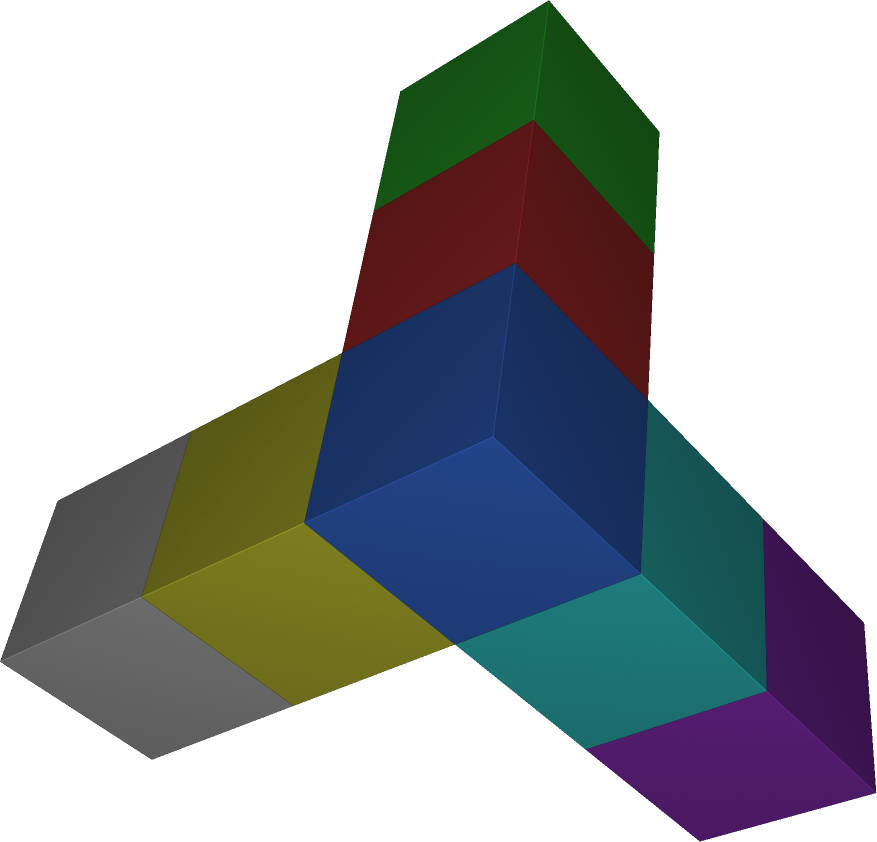}
    \includegraphics[width=0.23\textwidth] {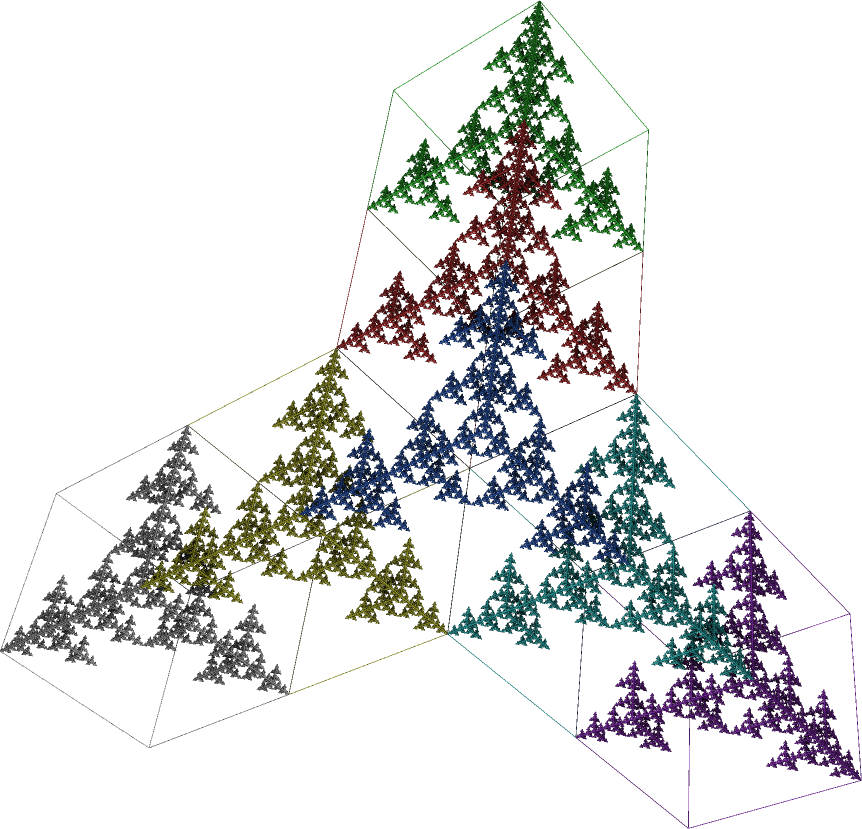} \\
    {\scriptsize$\mD=\{(0,0,0),\ (0,1,0),\ (0,1,1),\ (0,2,0),\ (0,2,2),\ (1,0,1),\ (2,0,2)\}$} & 
    {\scriptsize$\mD=\{(0,0,0),\ (0,0,1),\ (0,0,2),\ (0,1,2),\ (0,2,2),\ (1,0,2),\ (2,0,2)\}$} \\
    \hline
    \includegraphics[width=0.23\textwidth] {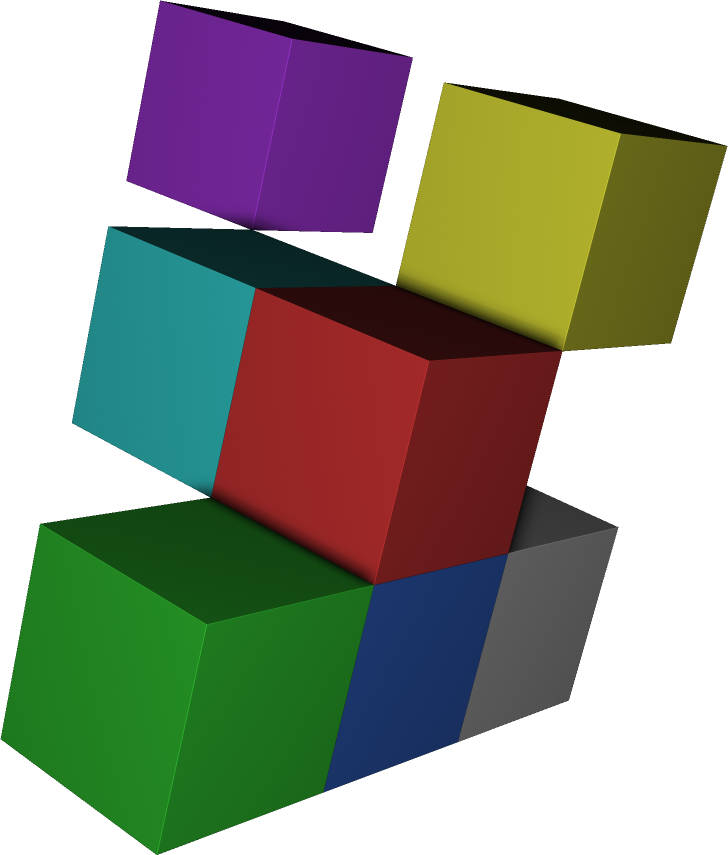}
    \includegraphics[width=0.23\textwidth] {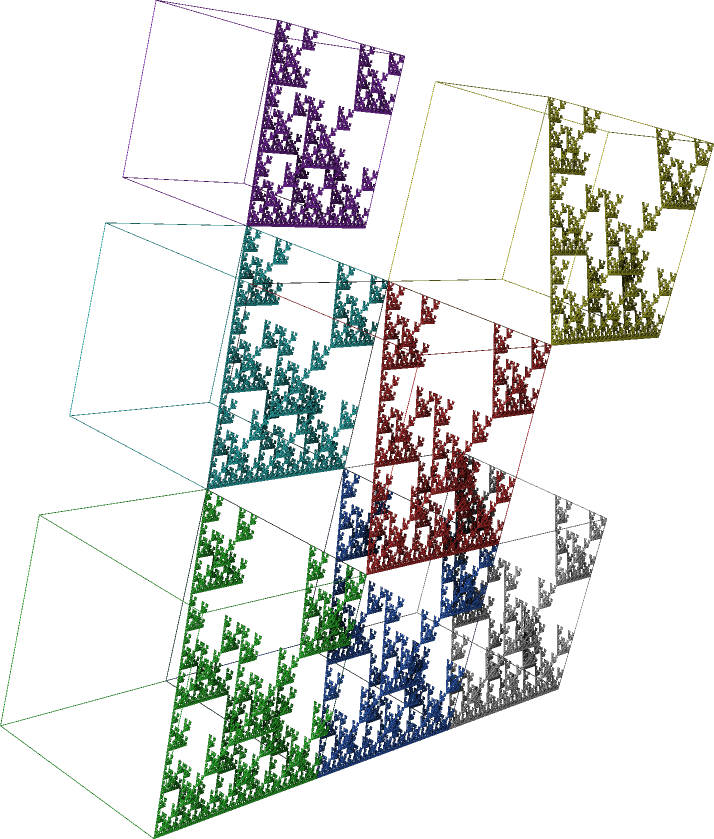} & \includegraphics[width=0.23\textwidth] {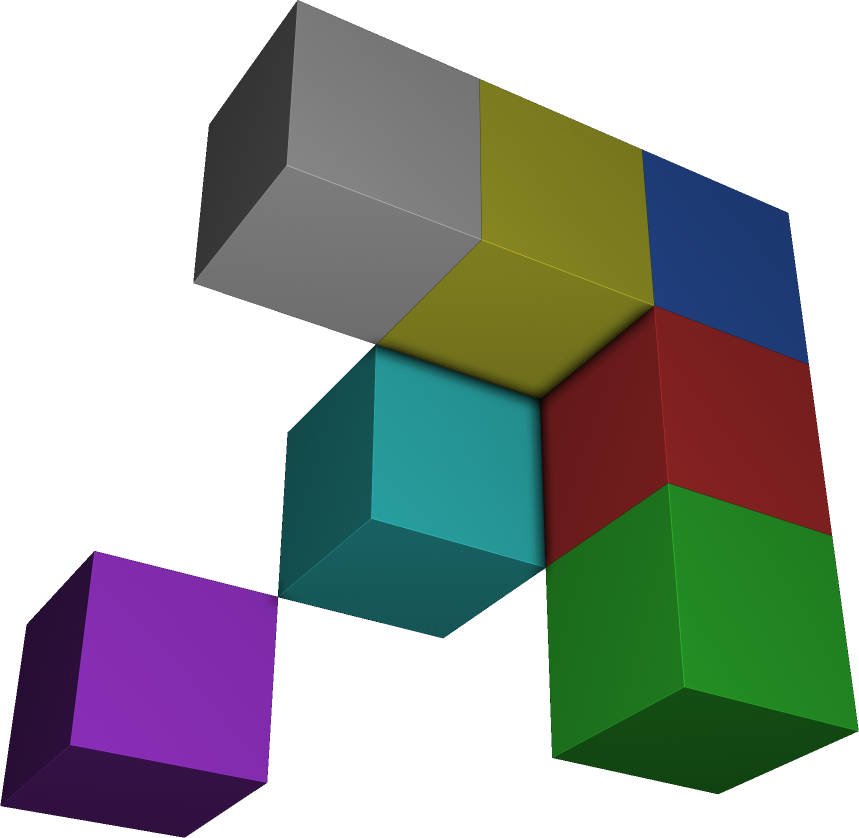}
    \includegraphics[width=0.23\textwidth] {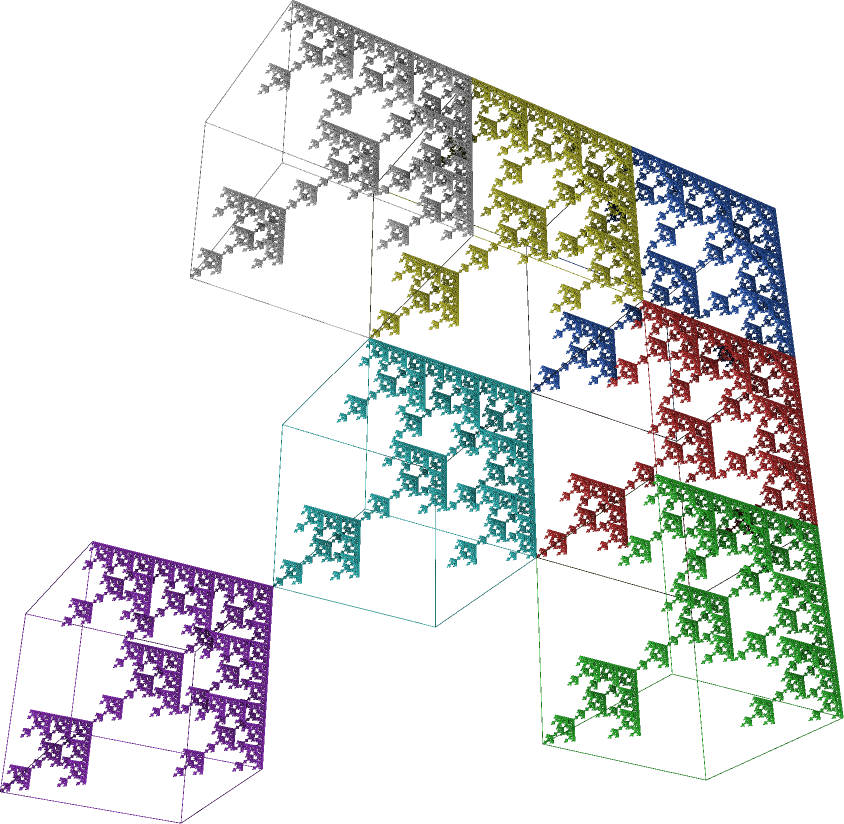} \\
    {\scriptsize$\mD=\{(0,0,0),\ (0,0,2),\ (0,1,0),\ (0,1,1),\ (0,2,0),\ (1,1,1),\ (2,0,2)\}$} & 
    {\scriptsize$\mD=\{(0,0,0),\ (0,1,0),\ (0,2,0),\ (0,2,1),\ (0,2,2),\ (1,1,1),\ (2,0,2)\}$} \\
    \hline
    \includegraphics[width=0.23\textwidth] {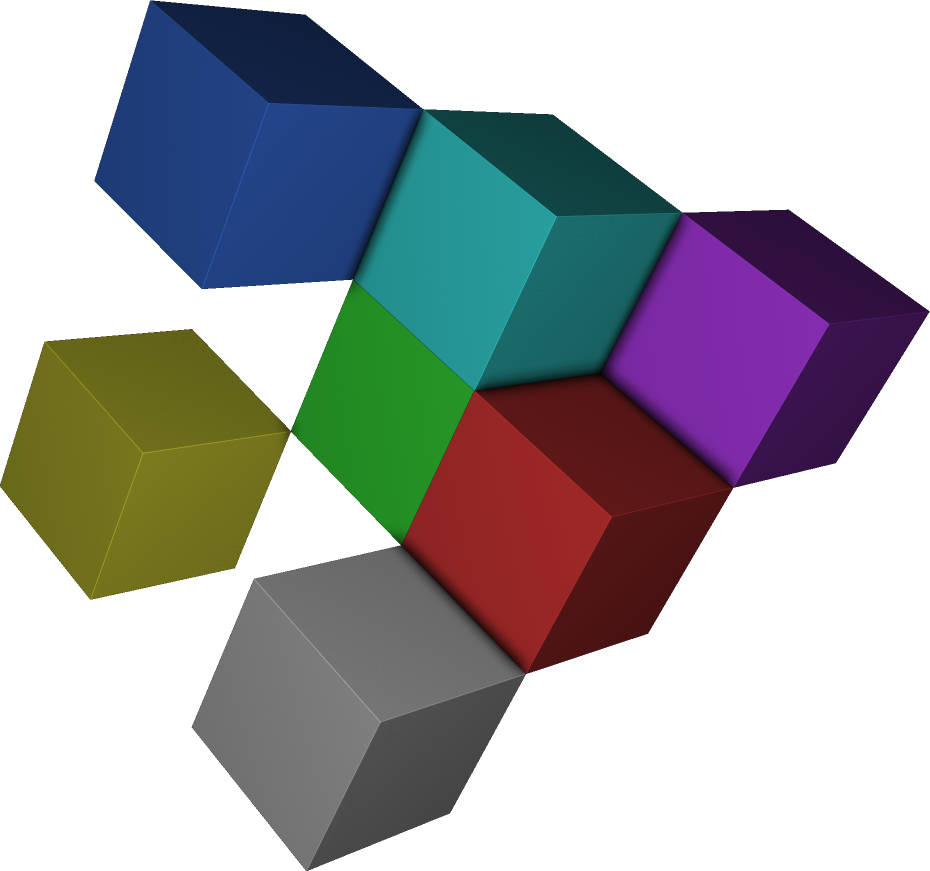}
    \includegraphics[width=0.23\textwidth] {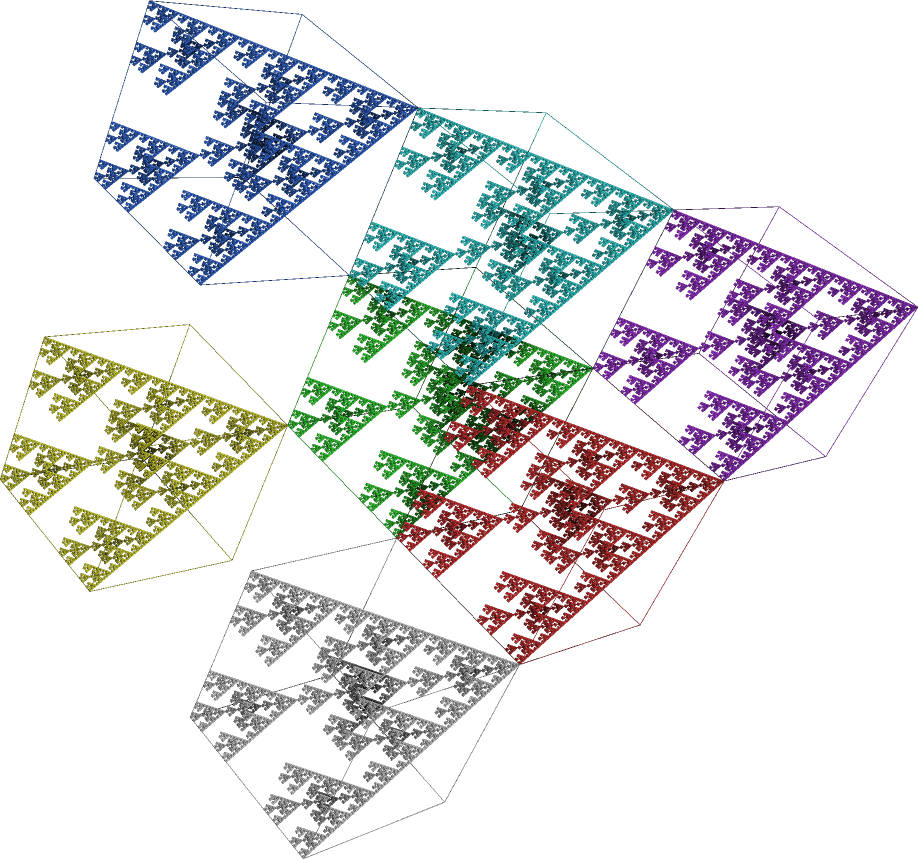} & \includegraphics[width=0.23\textwidth] {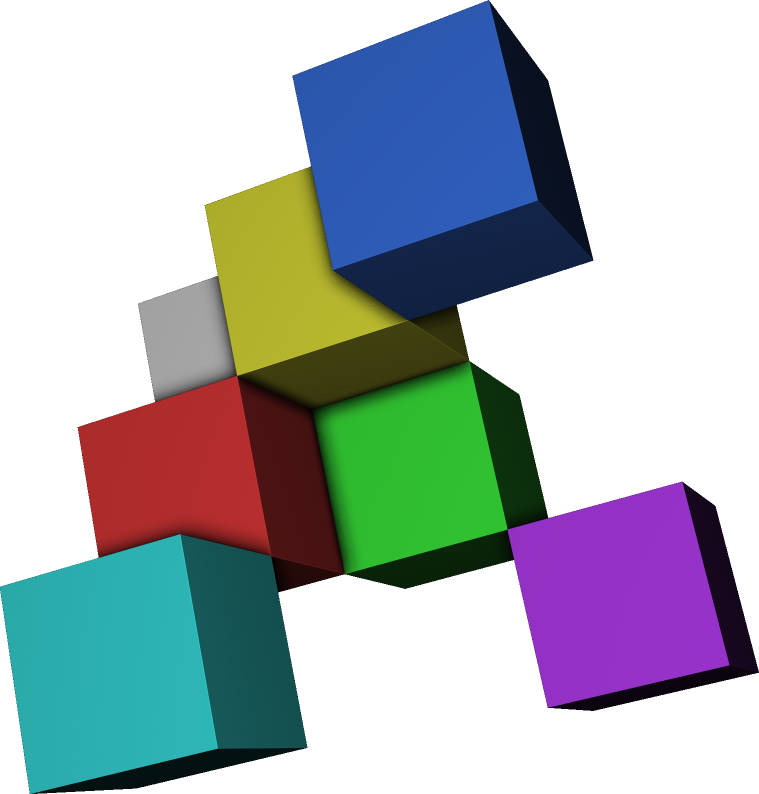}
    \includegraphics[width=0.23\textwidth] {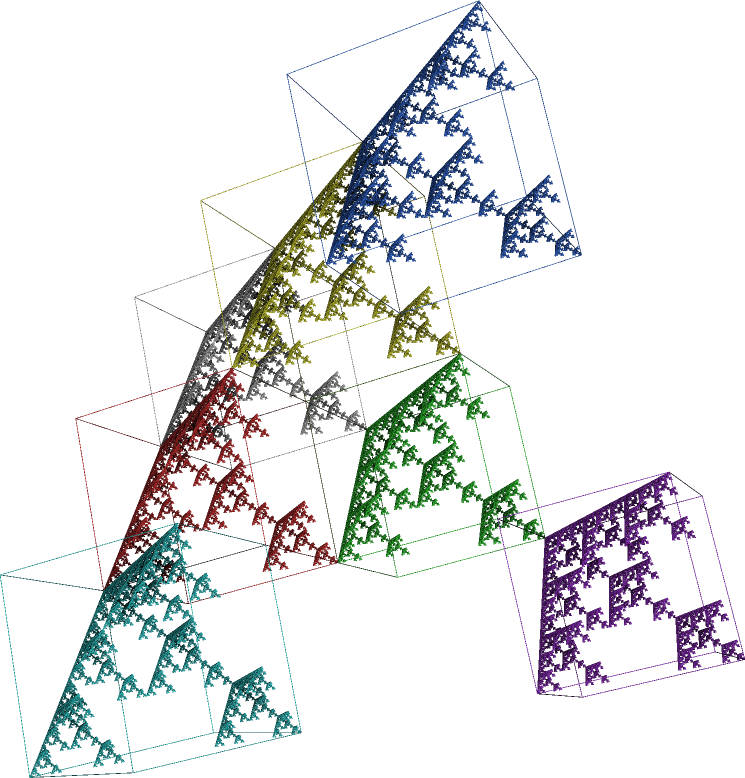} \\
    {\scriptsize$\mD=\{(0,0,0),\ (0,2,0),\ (0,2,2),\ (1,0,1),\ (1,1,1),\ (1,1,2),\ (2,0,2)\}$} & 
    {\scriptsize$\mD=\{(0,0,0),\ (0,1,1),\ (0,2,2),\ (1,0,1),\ (1,1,0),\ (2,0,2),\ (2,2,0)\}$} \\
    \hline
\end{longtable}

\newpage
\begin{longtable}{|p{0.48\textwidth}|p{0.48\textwidth}|}
\caption{Non-dendrites of   type 3 ($N=17$) }\label{tab:n3}\\
    \hline
    \multicolumn{2}{|c|}{\includegraphics{nonden3.pdf}} \\
    \hline
    \includegraphics[width=0.23\textwidth] {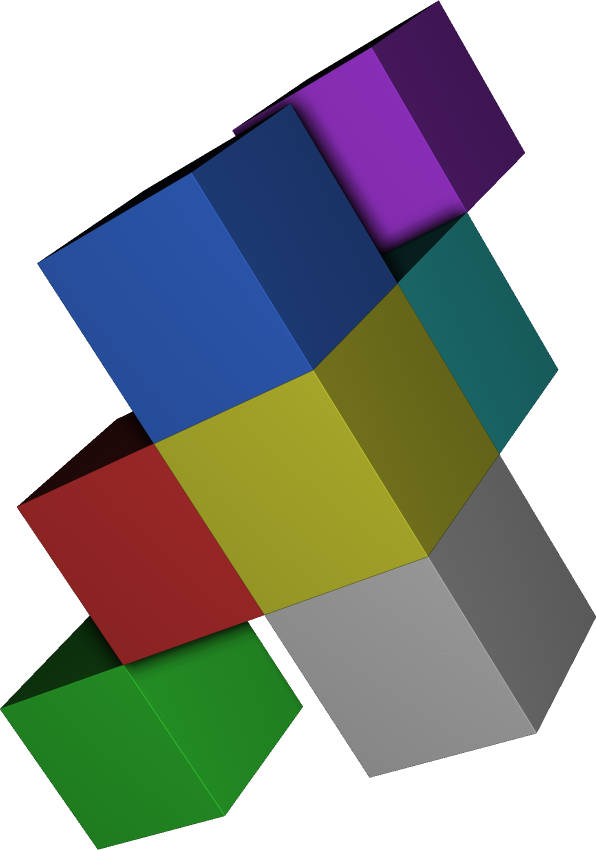}
    \includegraphics[width=0.23\textwidth] {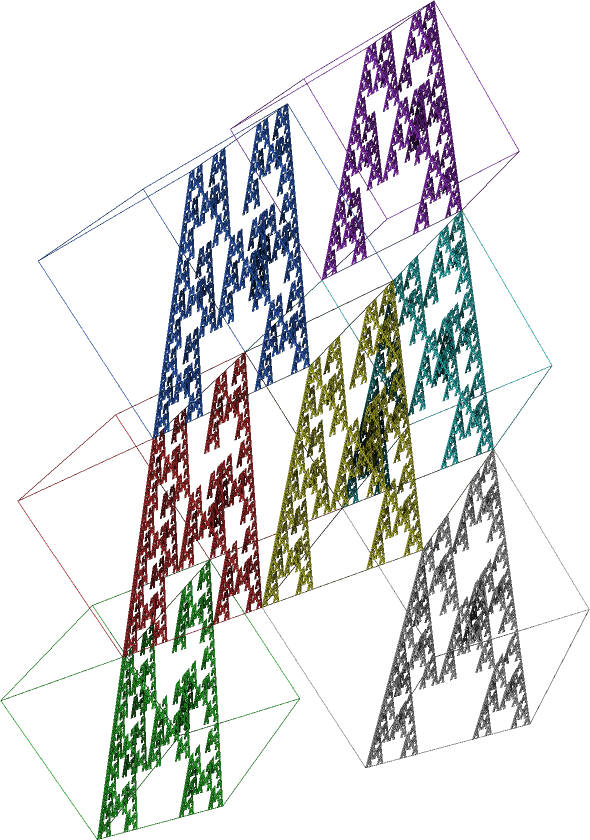} & \includegraphics[width=0.23\textwidth] {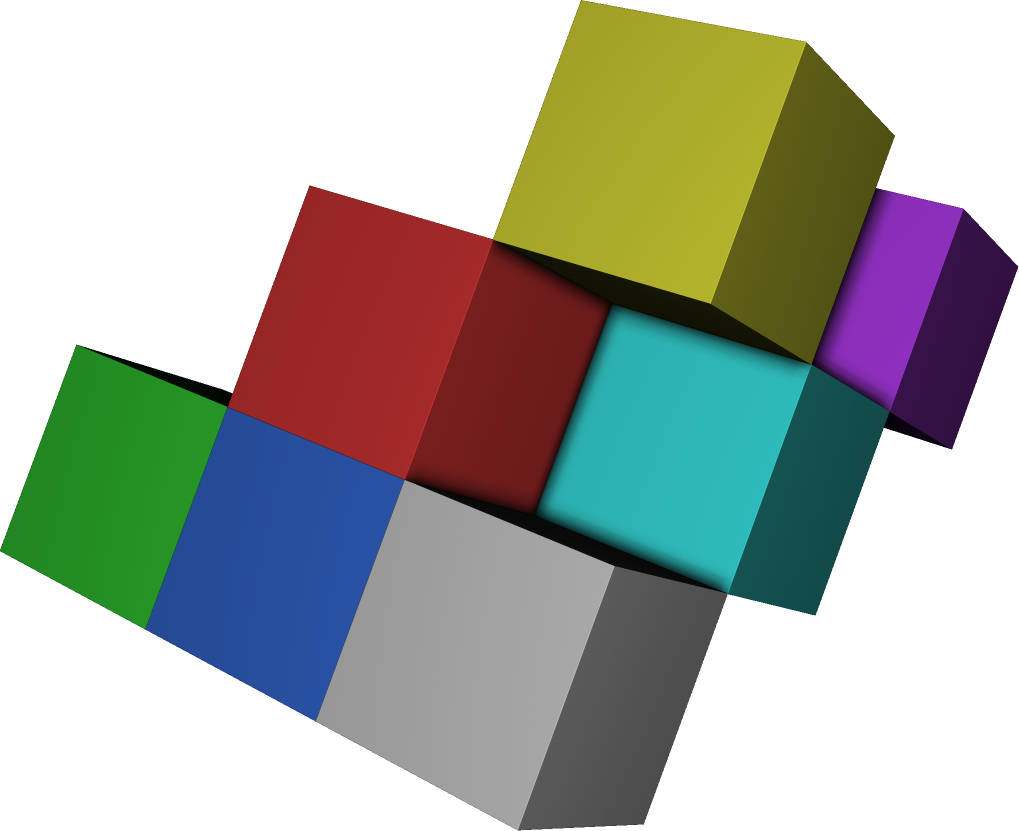}
    \includegraphics[width=0.23\textwidth] {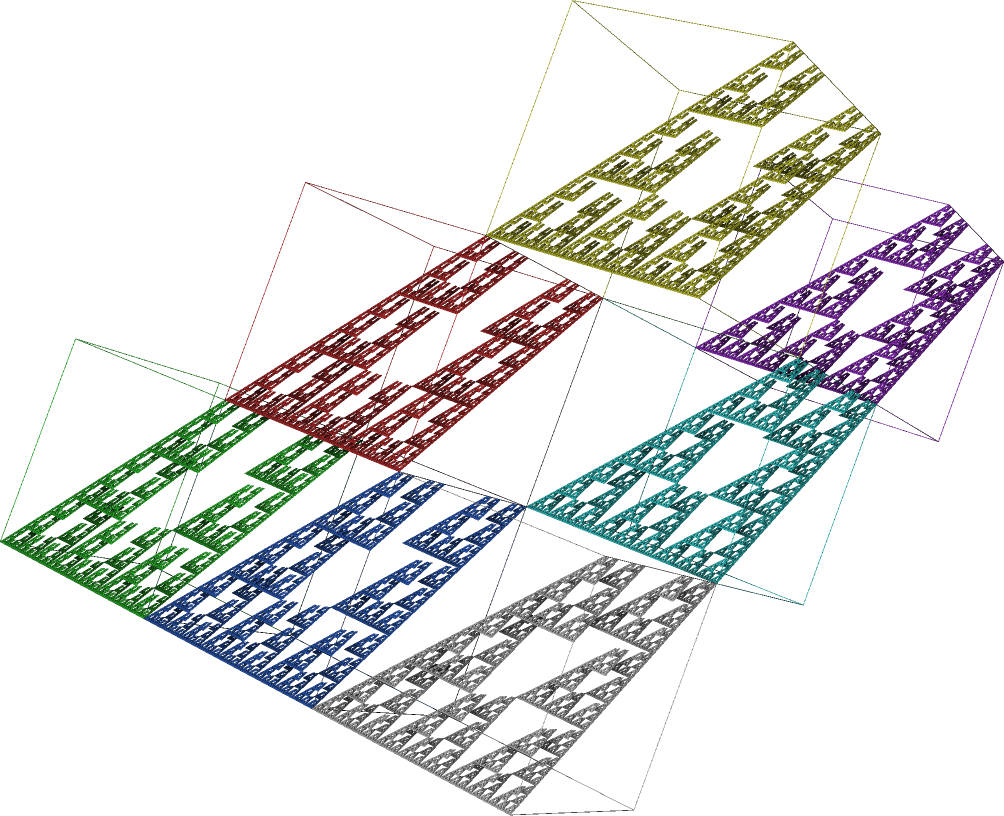} \\
    {\scriptsize$\mD=\{(0,0,0),\ (0,0,1),\ (0,0,2),\ (0,1,1),\ (0,2,0),\ (1,0,1),\ (2,0,2)\}$} & 
    {\scriptsize$\mD=\{(0,0,0),\ (0,0,2),\ (0,1,0),\ (0,1,1),\ (0,2,0),\ (1,0,1),\ (2,0,2)\}$} \\
    \hline
    \includegraphics[width=0.23\textwidth] {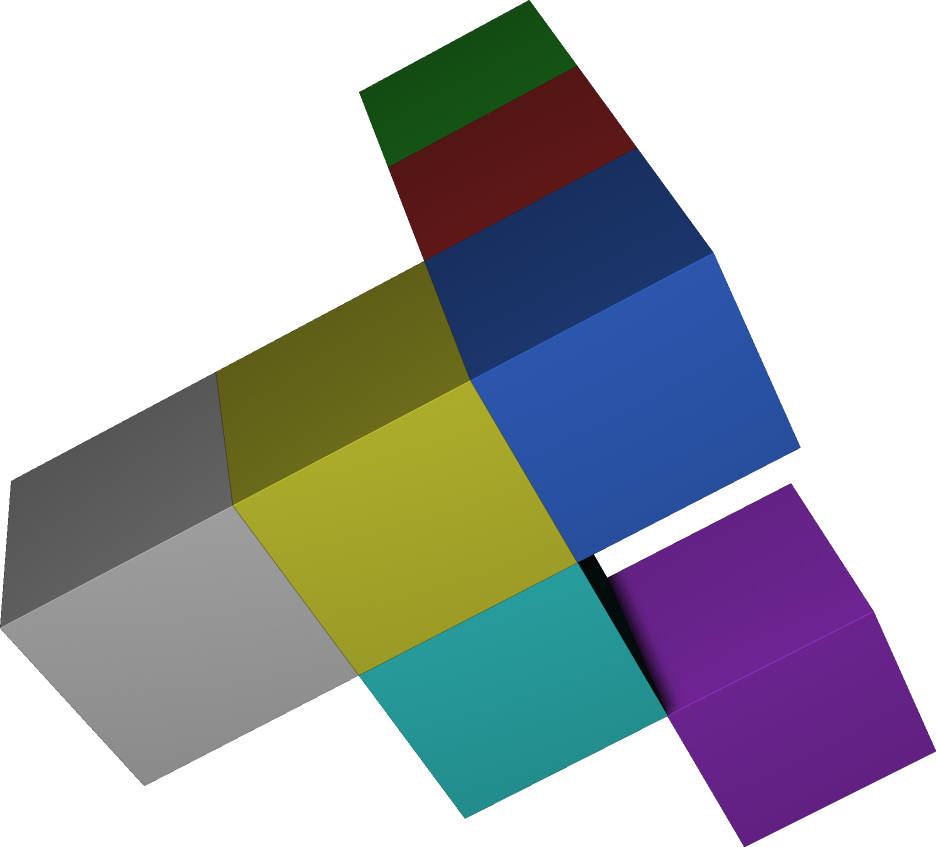}
    \includegraphics[width=0.23\textwidth] {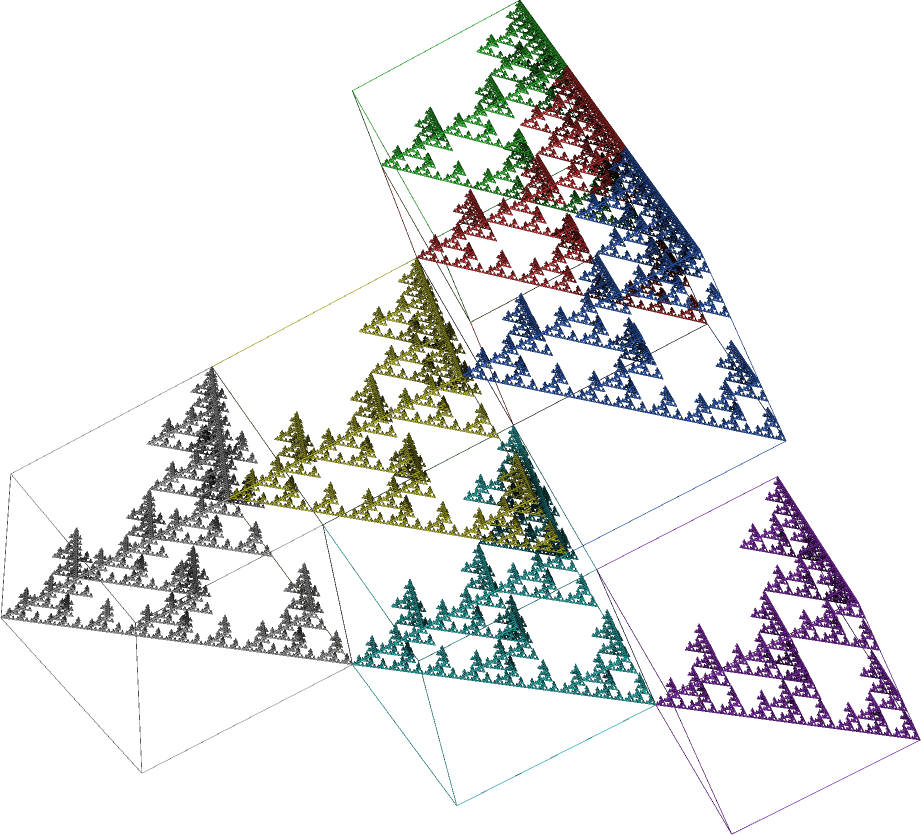} & \includegraphics[width=0.23\textwidth] {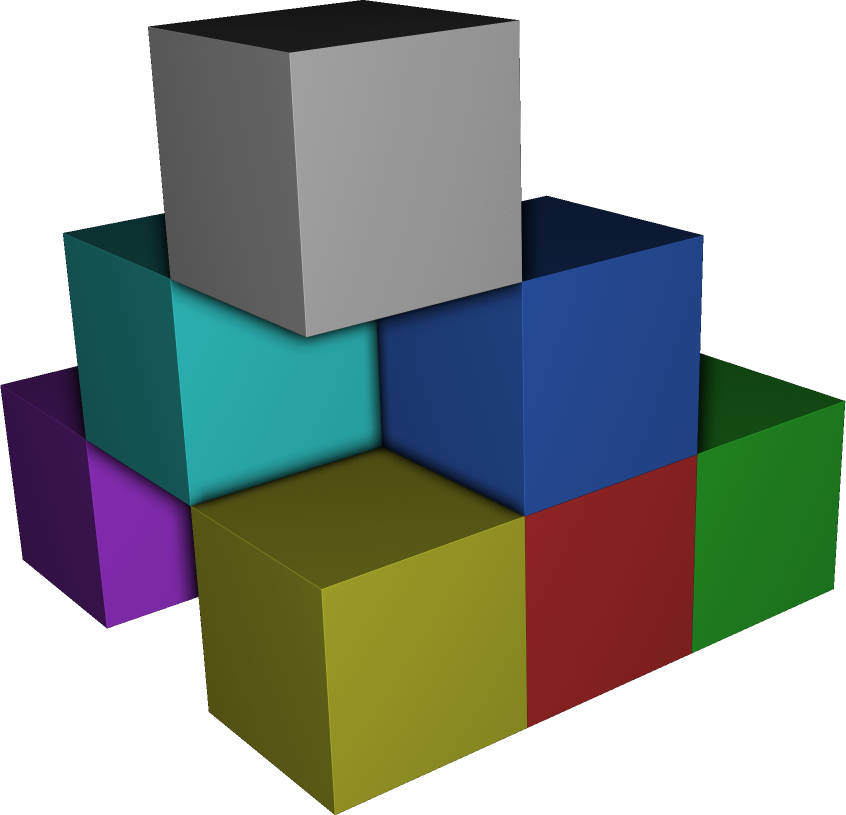}
    \includegraphics[width=0.23\textwidth] {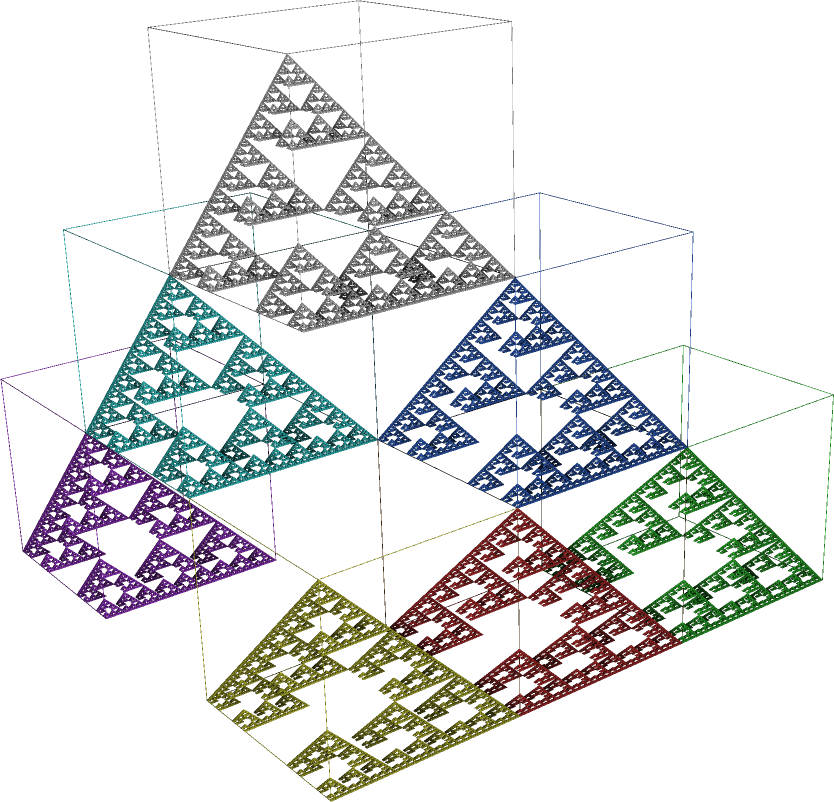} \\
    {\scriptsize$\mD=\{(0,0,0),\ (0,0,1),\ (0,0,2),\ (0,1,2),\ (0,2,2),\ (1,0,1),\ (2,0,2)\}$} & 
    {\scriptsize$\mD=\{(0,0,0),\ (0,0,2),\ (0,1,1),\ (0,1,2),\ (0,2,2),\ (1,0,1),\ (2,0,2)\}$} \\
    \hline
    \includegraphics[width=0.23\textwidth] {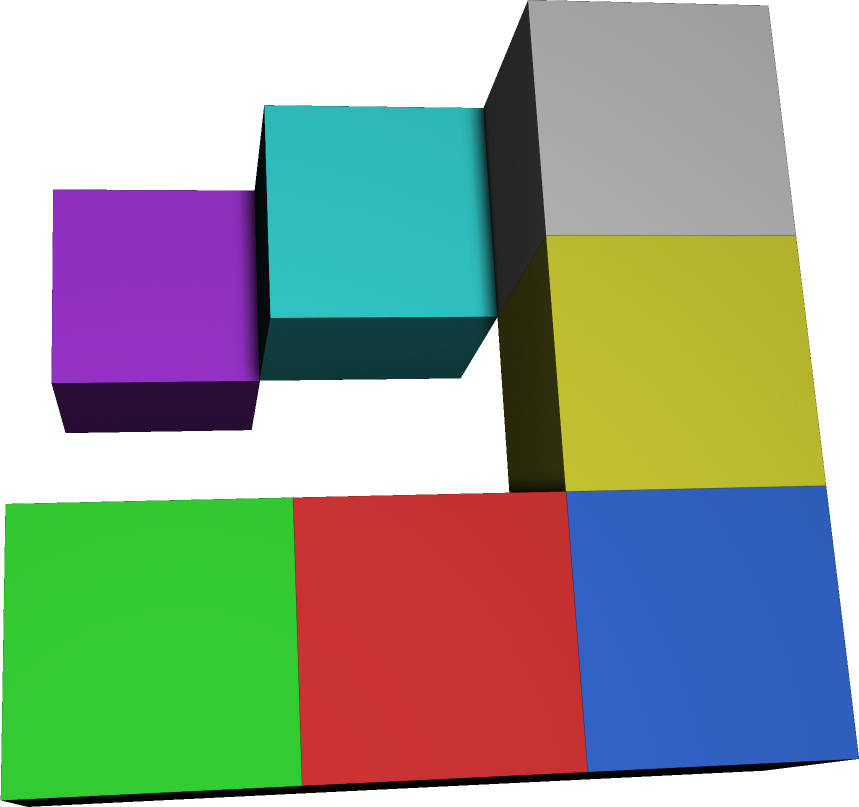}
    \includegraphics[width=0.23\textwidth] {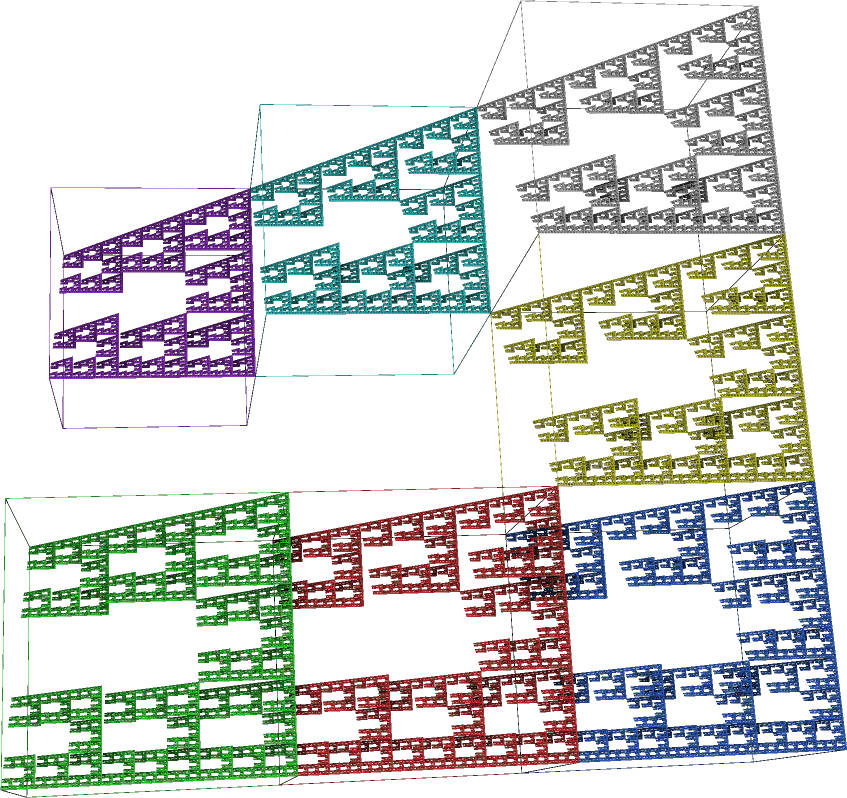} & \includegraphics[width=0.23\textwidth] {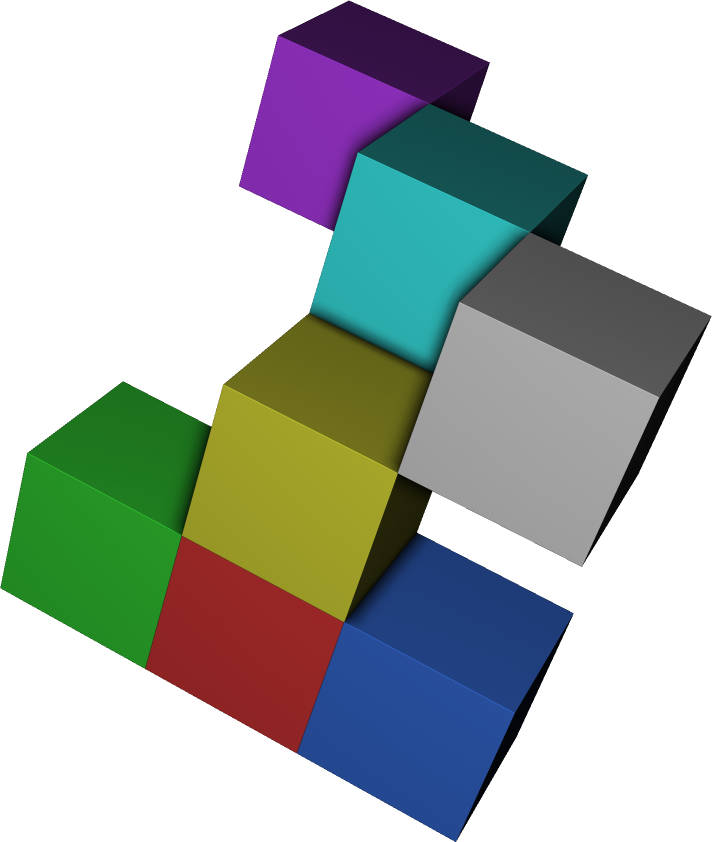}
    \includegraphics[width=0.23\textwidth] {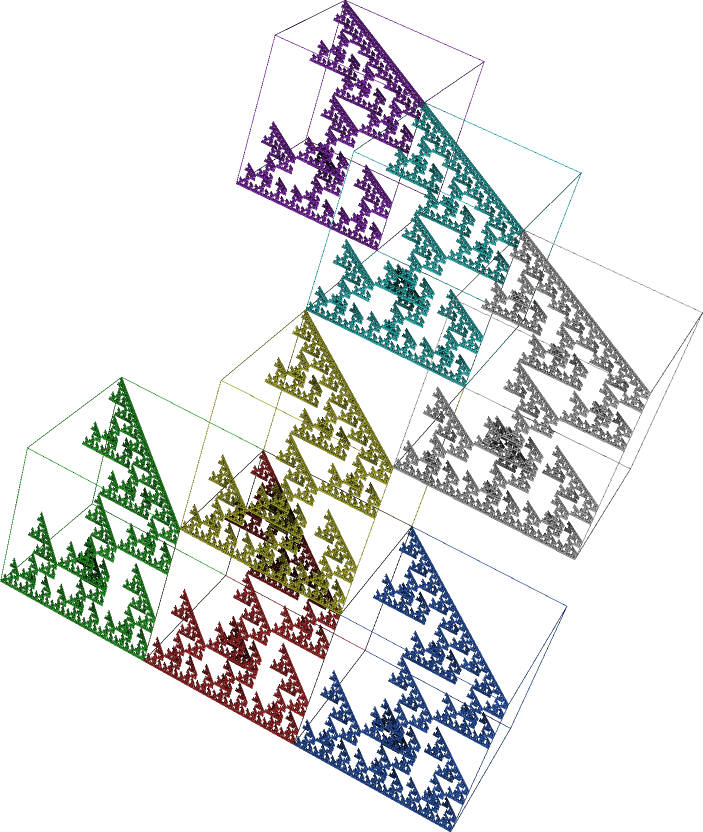} \\
    {\scriptsize$\mD=\{(0,0,0),\ (0,1,0),\ (0,2,0),\ (0,2,1),\ (0,2,2),\ (1,0,1),\ (2,0,2)\}$} & 
    {\scriptsize$\mD=\{(0,0,0),\ (0,1,1),\ (0,2,0),\ (0,2,1),\ (0,2,2),\ (1,0,1),\ (2,0,2)\}$} \\
    \hline
    \includegraphics[width=0.23\textwidth] {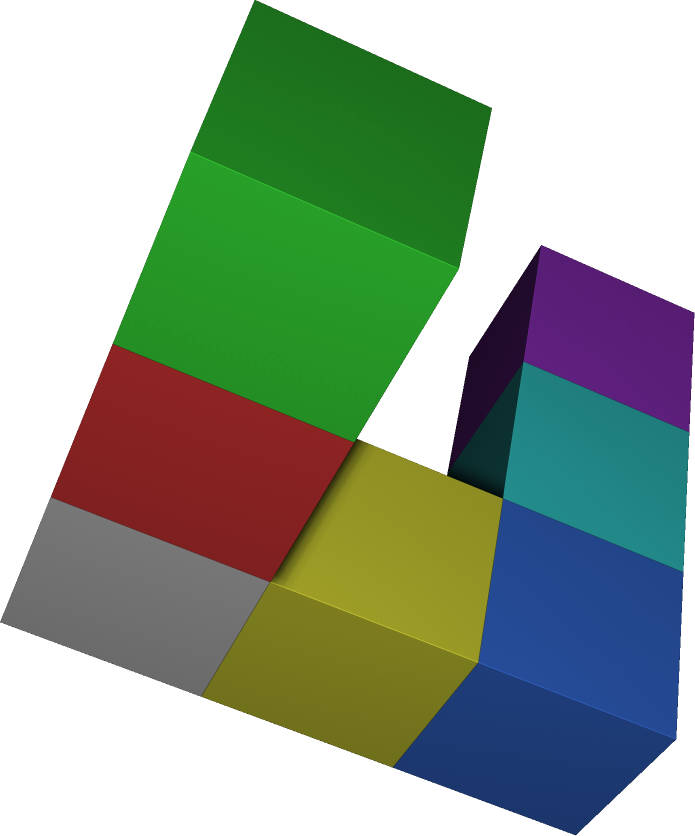}
    \includegraphics[width=0.23\textwidth] {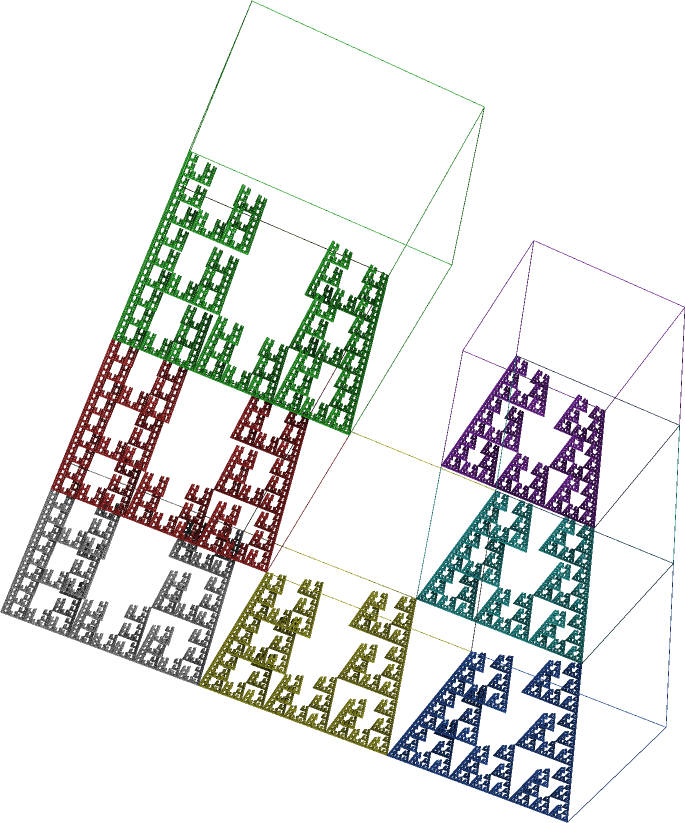} & \includegraphics[width=0.23\textwidth] {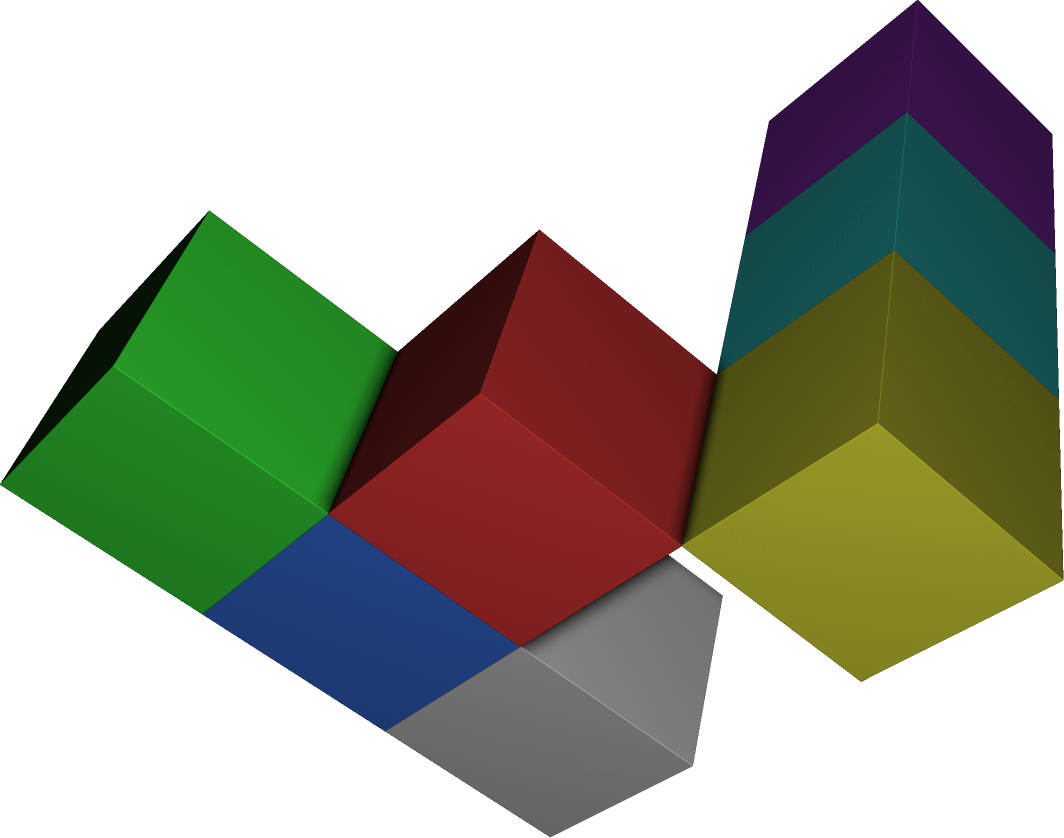}
    \includegraphics[width=0.23\textwidth] {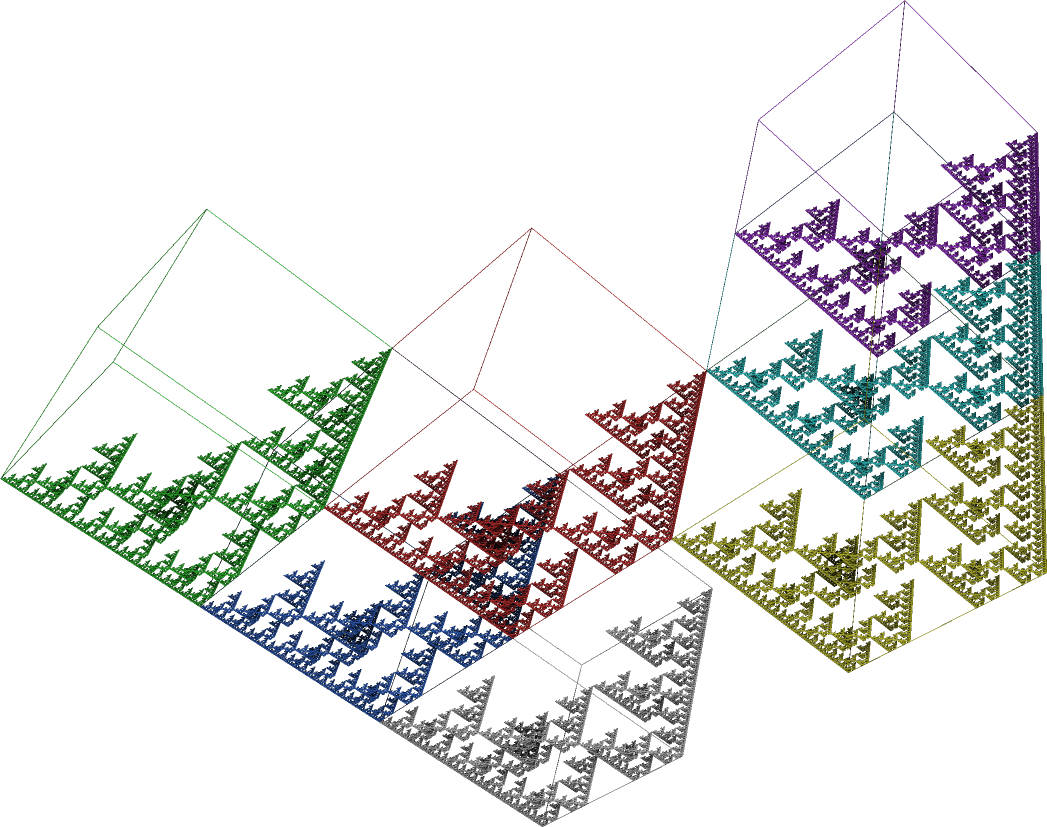} \\
    {\scriptsize$\mD=\{(0,0,0),\ (0,0,1),\ (0,0,2),\ (0,1,0),\ (0,2,0),\ (1,0,2),\ (2,0,2)\}$} & 
    {\scriptsize$\mD=\{(0,0,0),\ (0,0,2),\ (0,1,0),\ (0,1,1),\ (0,2,0),\ (1,0,2),\ (2,0,2)\}$} \\
    \hline
    \includegraphics[width=0.23\textwidth] {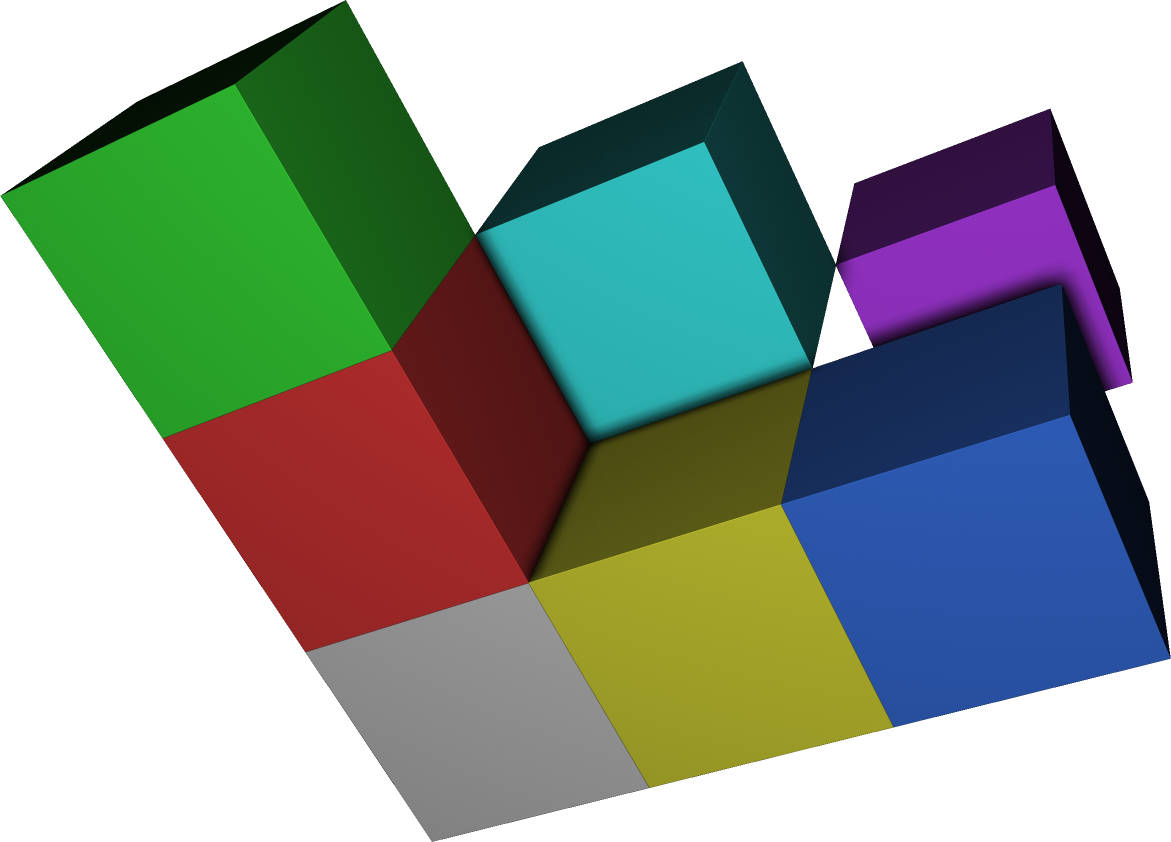}
    \includegraphics[width=0.23\textwidth] {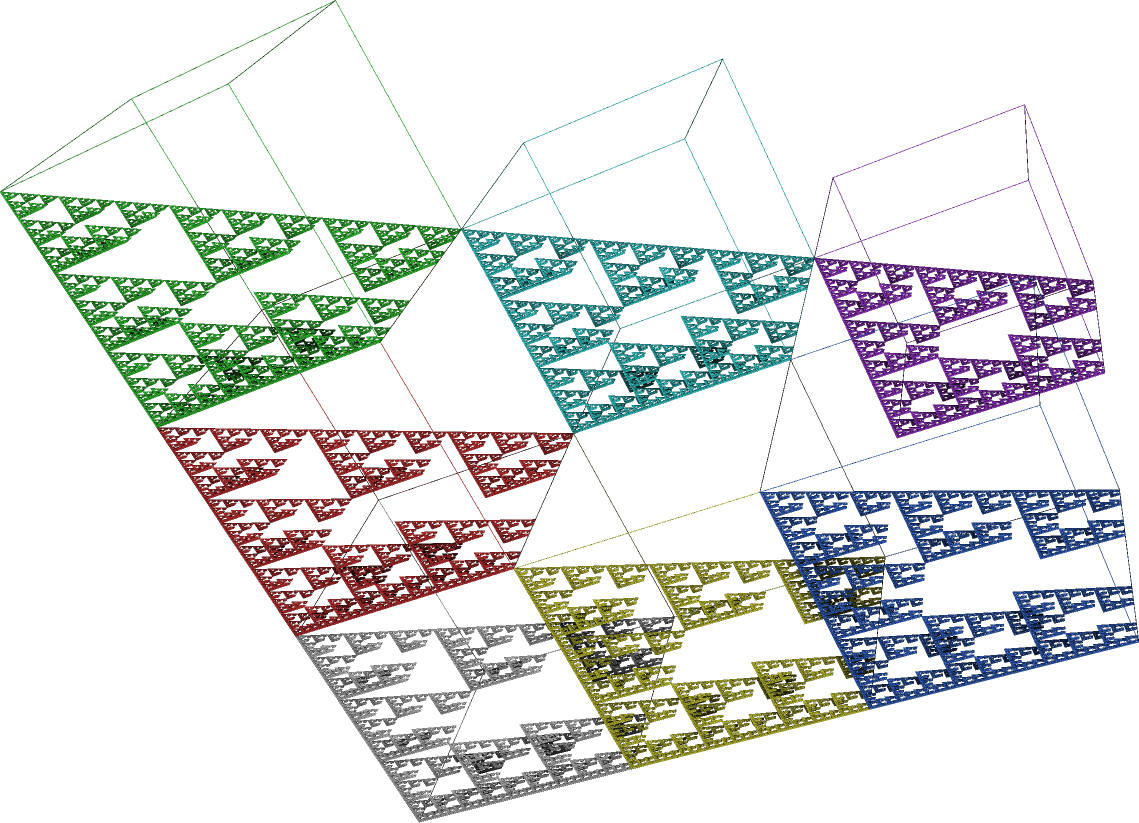} & \includegraphics[width=0.23\textwidth] {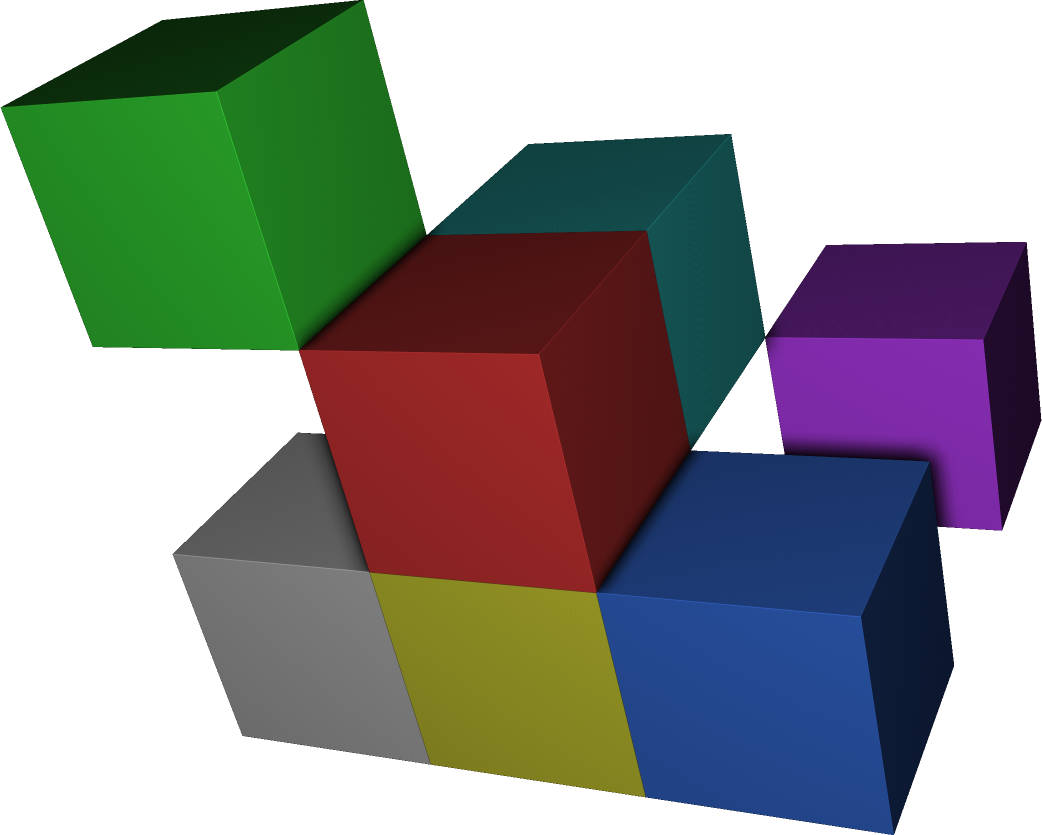}
    \includegraphics[width=0.23\textwidth] {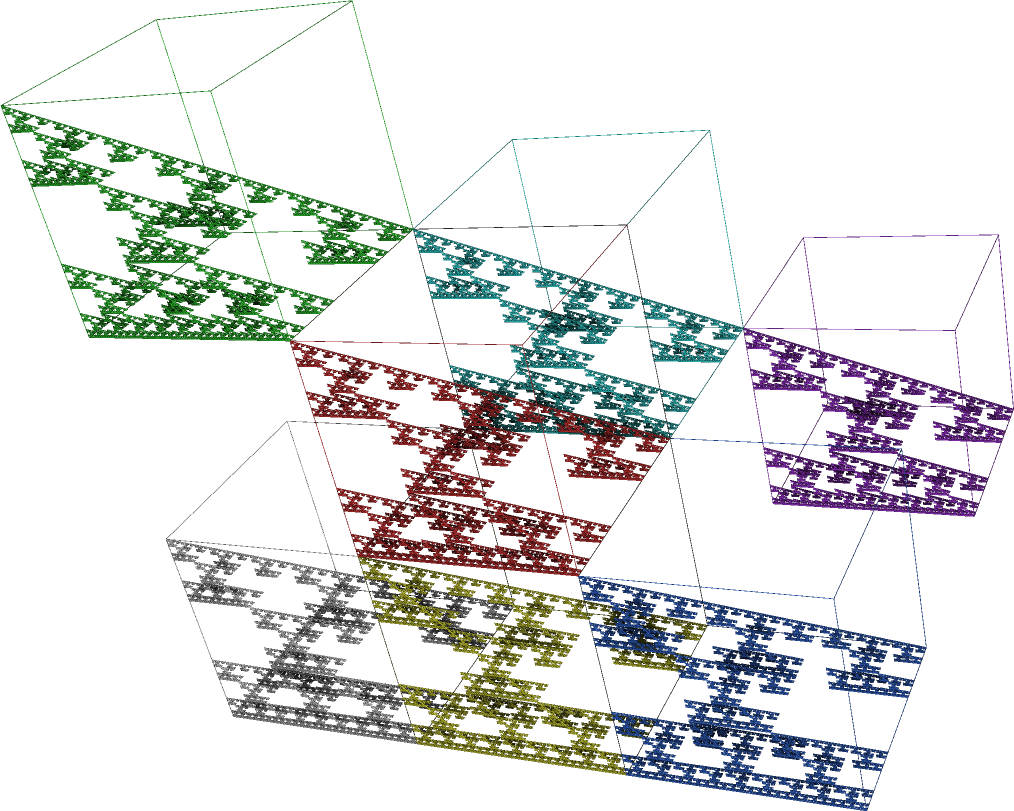} \\
    {\scriptsize$\mD=\{(0,0,0),\ (0,0,1),\ (0,0,2),\ (0,1,0),\ (0,2,0),\ (1,1,1),\ (2,0,2)\}$} & 
    {\scriptsize$\mD=\{(0,0,0),\ (0,0,1),\ (0,0,2),\ (0,1,1),\ (0,2,0),\ (1,1,1),\ (2,0,2)\}$} \\
    \hline
    \includegraphics[width=0.23\textwidth] {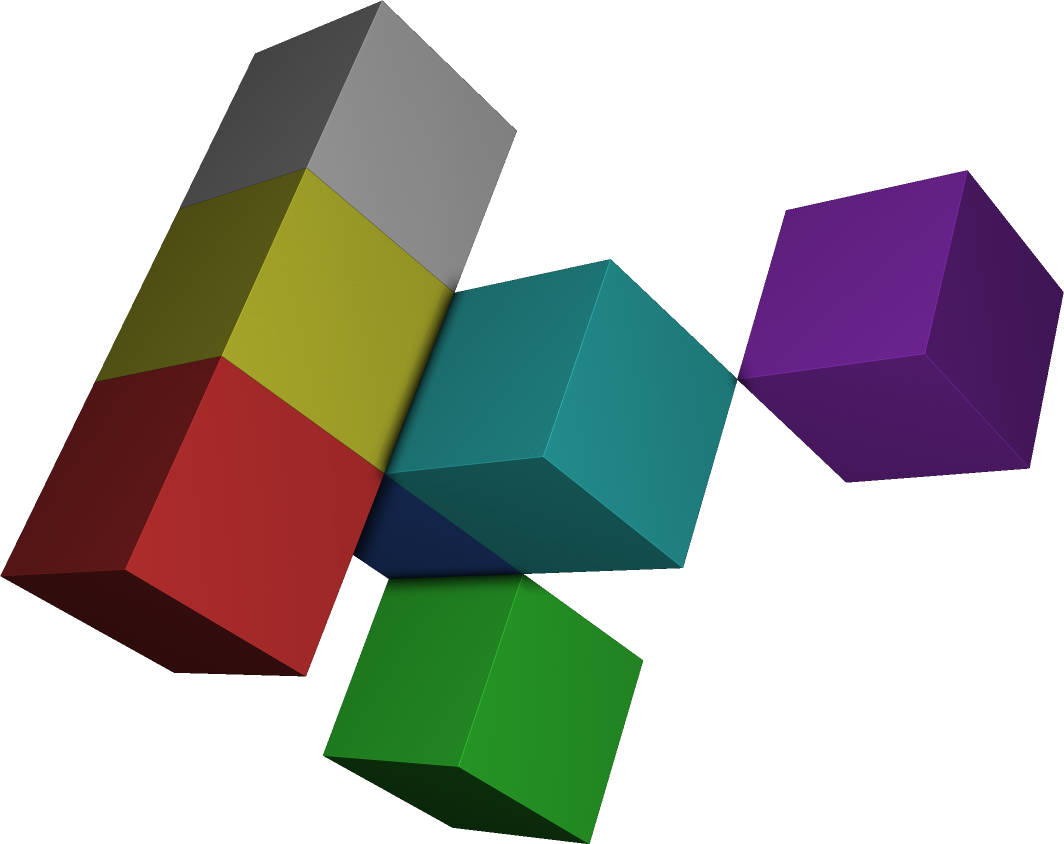}
    \includegraphics[width=0.23\textwidth] {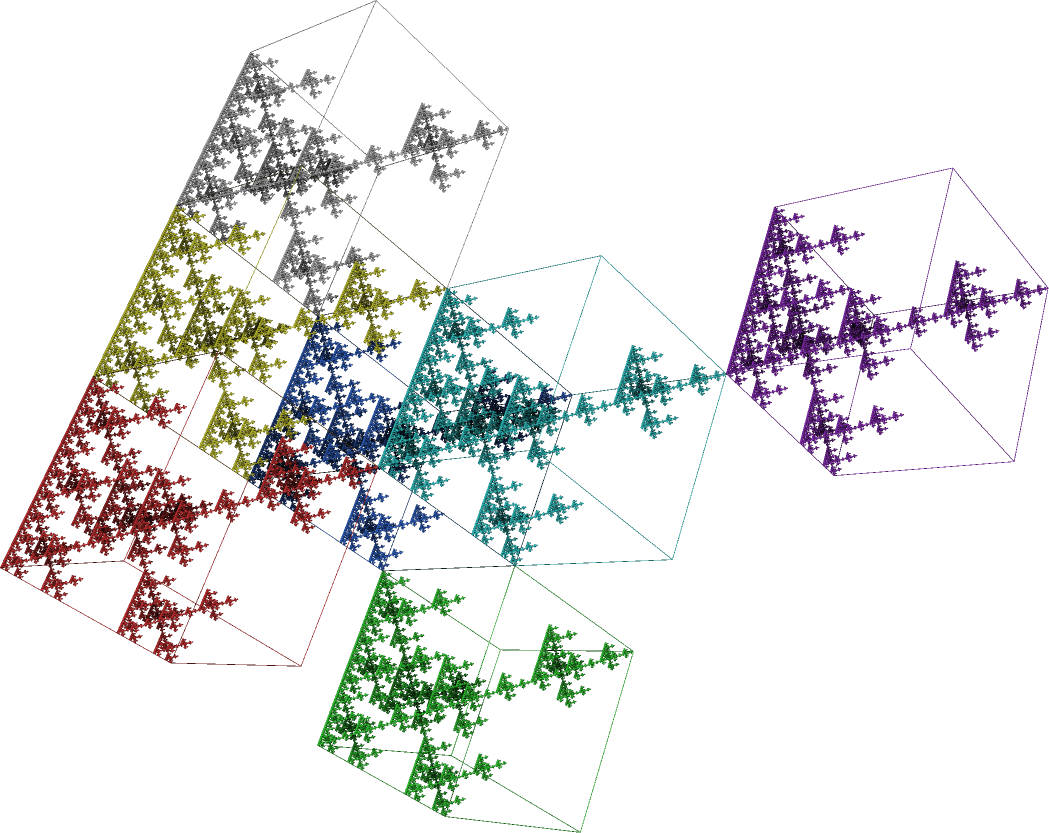} & \includegraphics[width=0.23\textwidth] {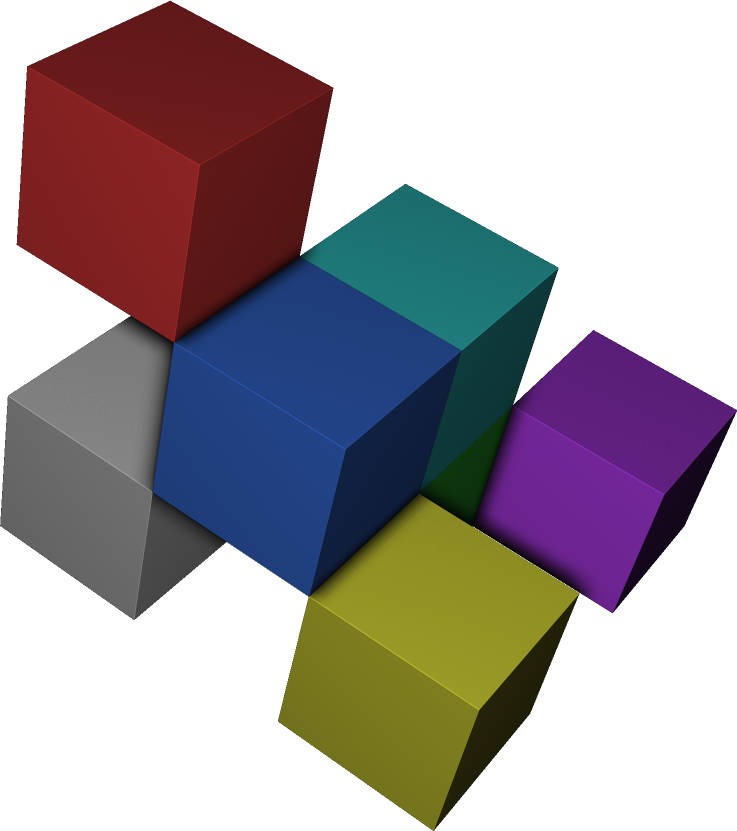}
    \includegraphics[width=0.23\textwidth] {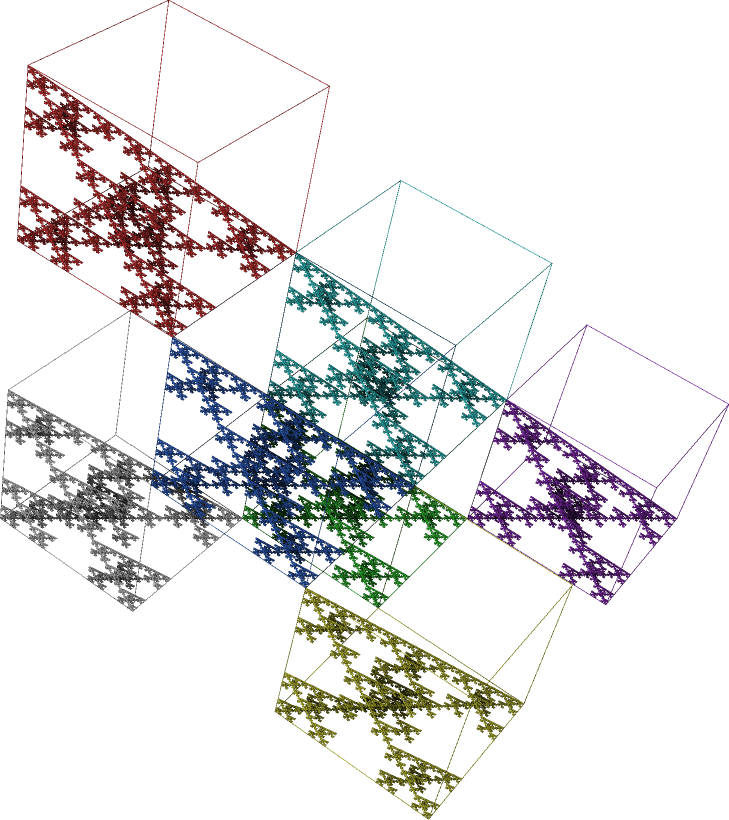} \\
    {\scriptsize$\mD=\{(0,0,0),\ (0,1,0),\ (0,1,1),\ (0,2,0),\ (0,2,2),\ (1,1,1),\ (2,0,2)\}$} & 
    {\scriptsize$\mD=\{(0,0,0),\ (0,0,2),\ (0,1,1),\ (0,2,0),\ (1,0,1),\ (1,1,1),\ (2,0,2)\}$} \\
    \hline
    \includegraphics[width=0.23\textwidth] {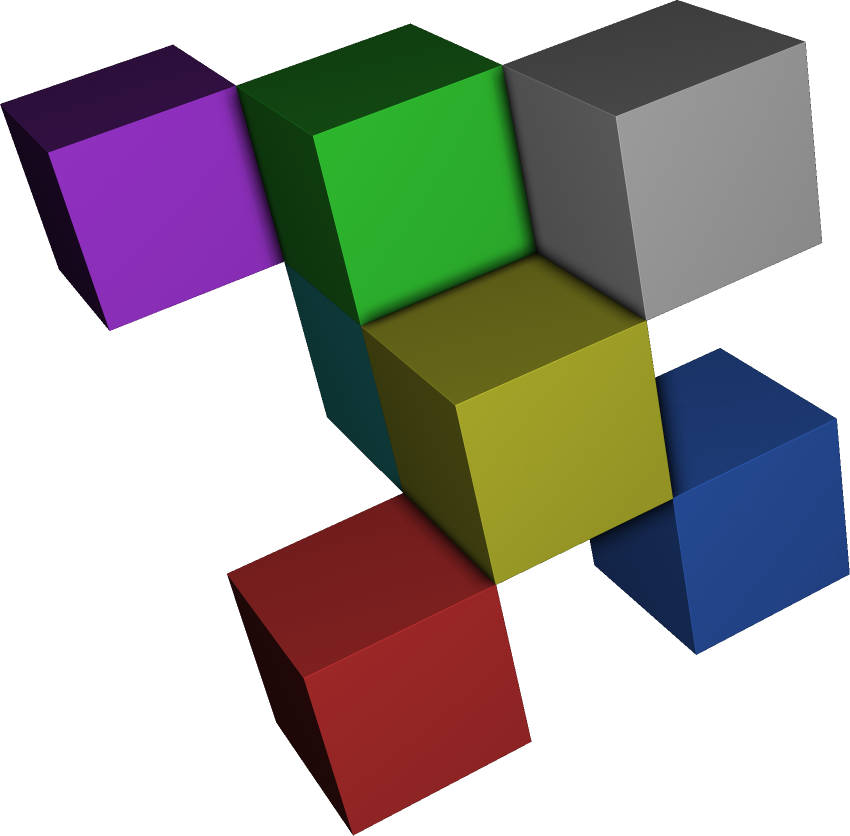}
    \includegraphics[width=0.23\textwidth] {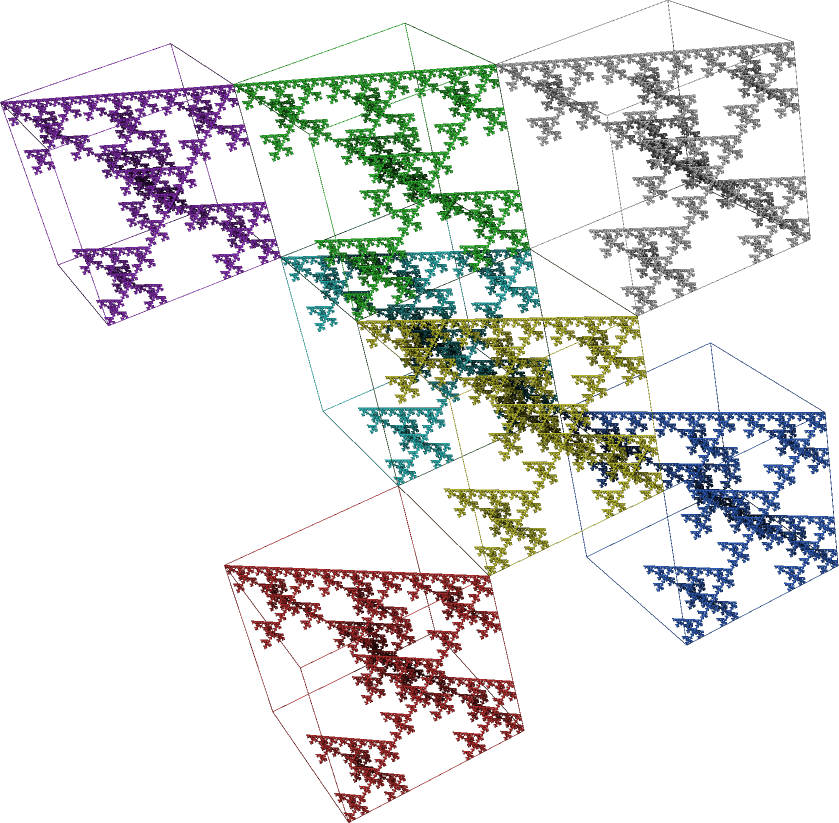} & \includegraphics[width=0.23\textwidth] {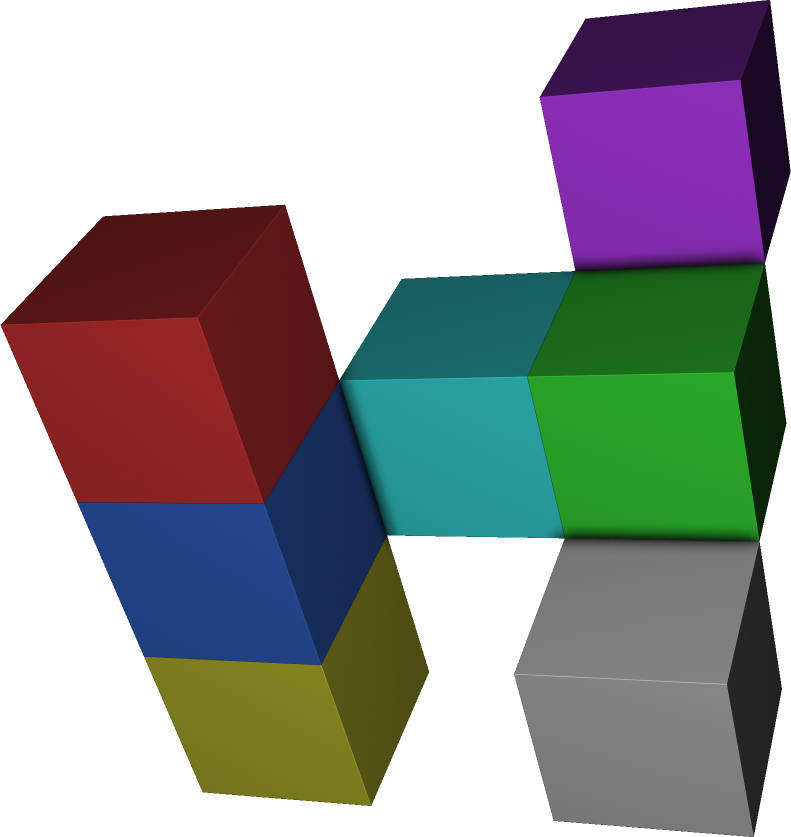}
    \includegraphics[width=0.23\textwidth] {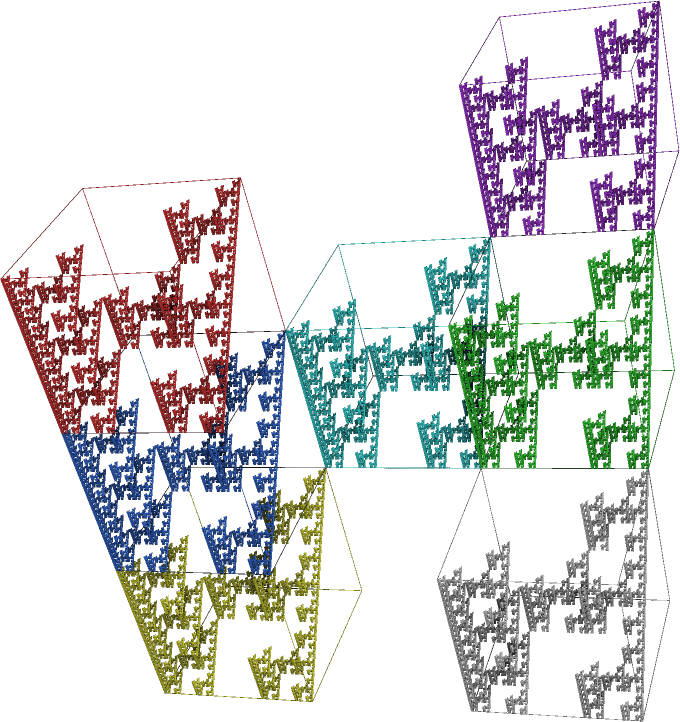} \\
    {\scriptsize$\mD=\{(0,0,0),\ (0,1,1),\ (0,2,0),\ (0,2,2),\ (1,0,1),\ (1,1,1),\ (2,0,2)\}$} & 
    {\scriptsize$\mD=\{(0,0,0),\ (0,2,0),\ (0,2,1),\ (0,2,2),\ (1,0,1),\ (1,1,1),\ (2,0,2)\}$} \\
    \hline
    \includegraphics[width=0.23\textwidth] {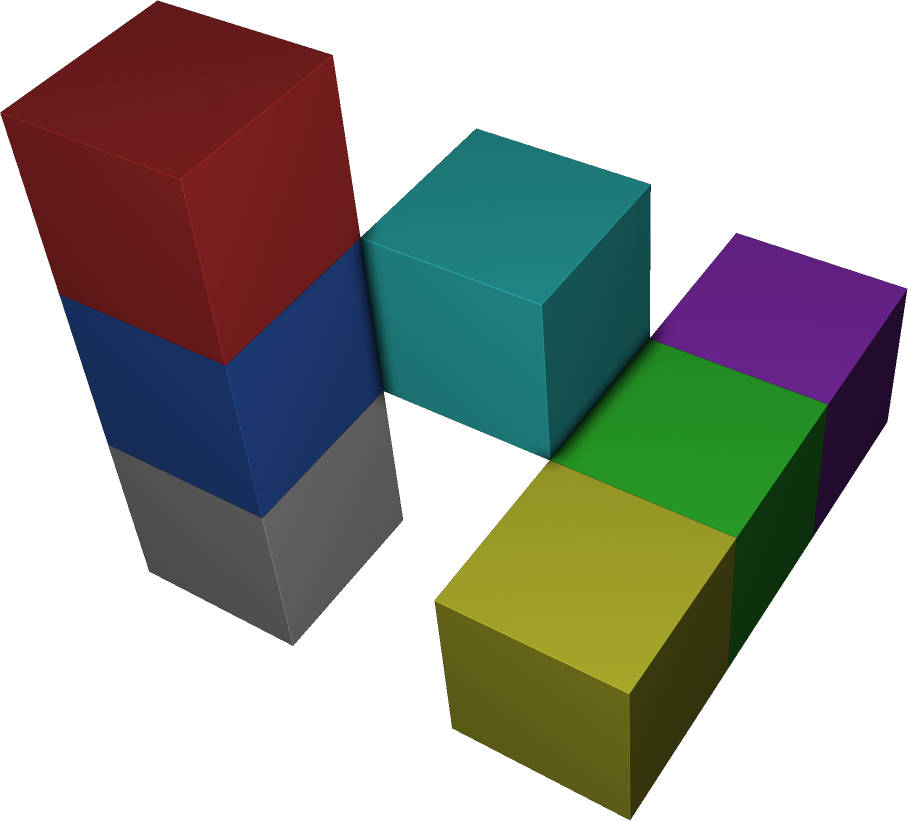}
    \includegraphics[width=0.23\textwidth] {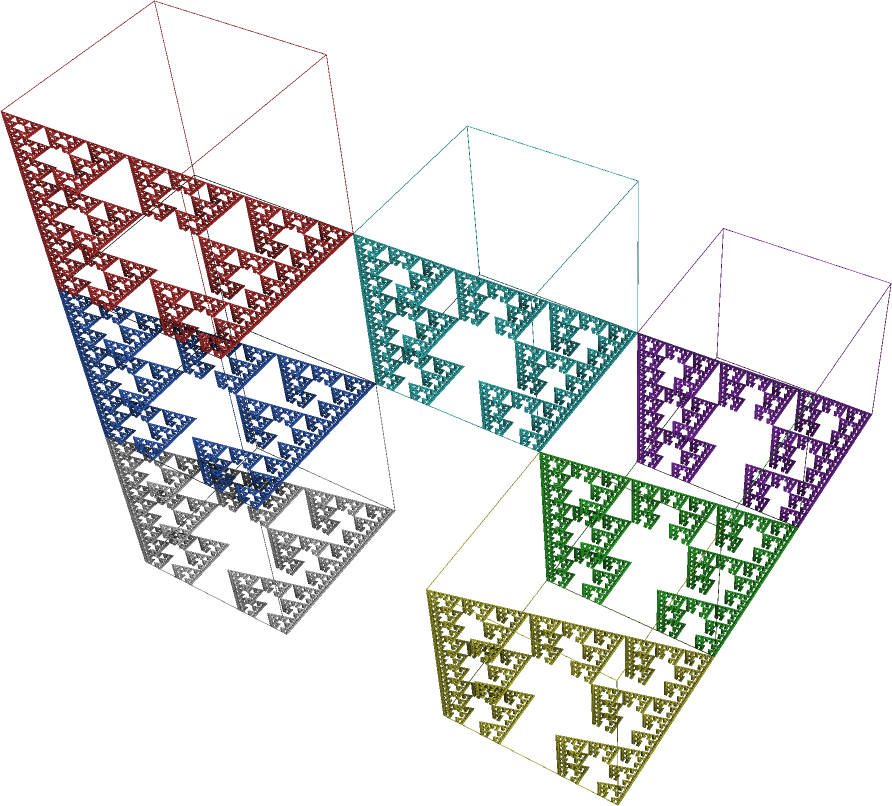} & \includegraphics[width=0.23\textwidth] {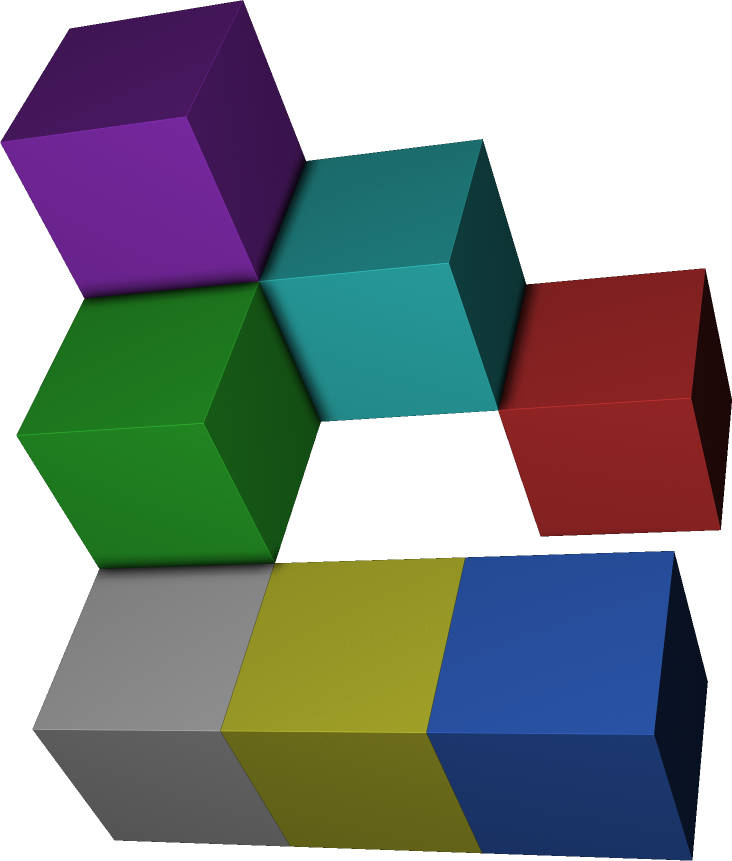}
    \includegraphics[width=0.23\textwidth] {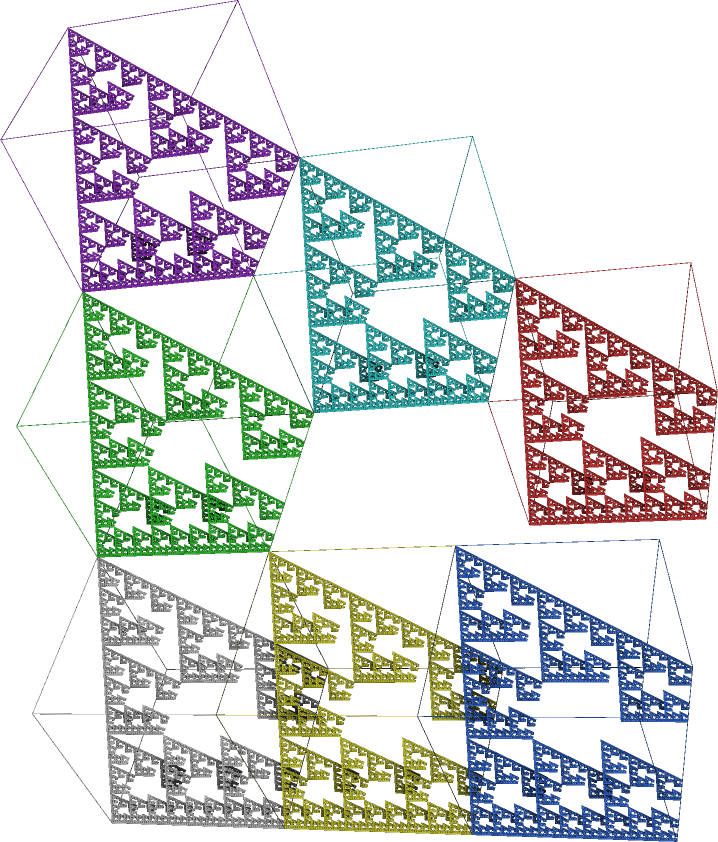} \\
    {\scriptsize$\mD=\{(0,0,0),\ (0,0,2),\ (0,1,0),\ (0,2,0),\ (1,0,2),\ (1,1,1),\ (2,0,2)\}$} & 
    {\scriptsize$\mD=\{(0,0,0),\ (0,1,0),\ (0,2,0),\ (0,2,2),\ (1,0,1),\ (1,1,2),\ (2,0,2)\}$} \\
    \hline
    \includegraphics[width=0.23\textwidth] {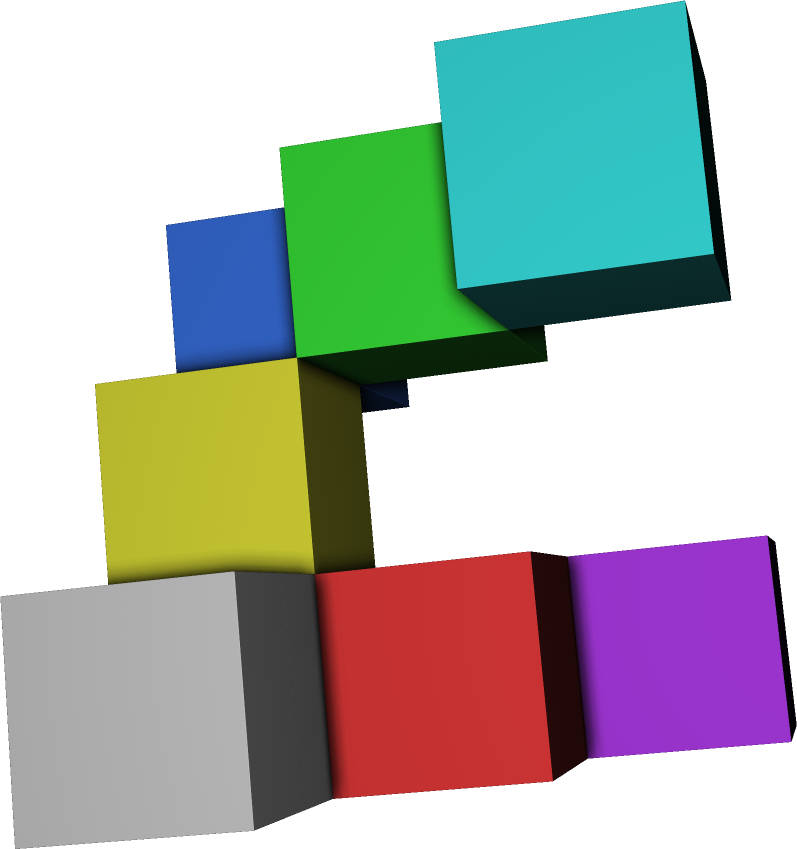}
    \includegraphics[width=0.23\textwidth] {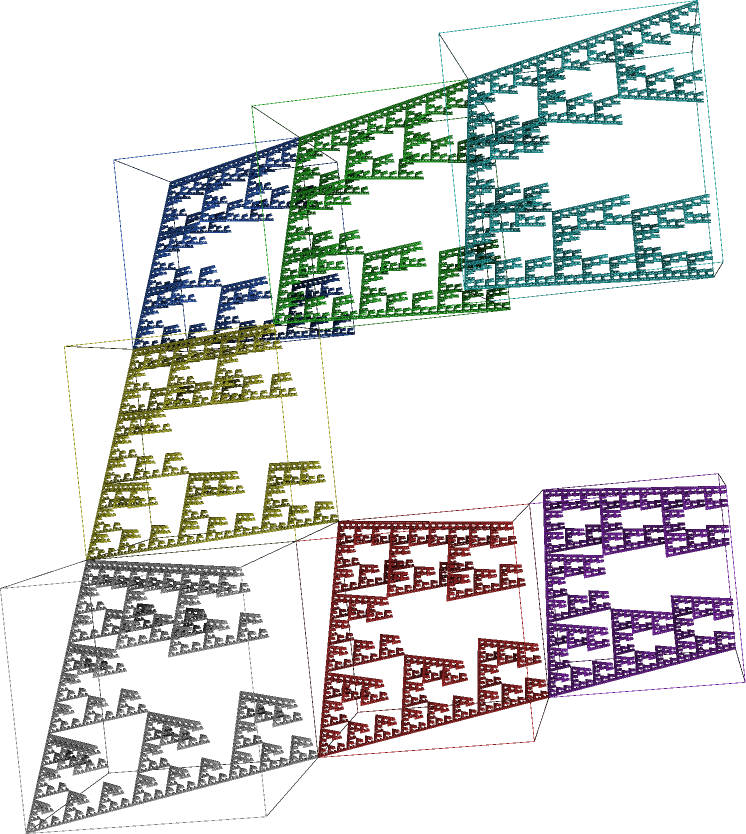} &\\
    {\scriptsize$\mD=\{(0,0,0),\ (0,1,1),\ (0,2,2),\ (1,1,0),\ (1,1,2),\ (2,0,2),\ (2,2,0)\}$} &\\
    \hline
\end{longtable}

\begin{longtable}{|p{0.48\textwidth}|p{0.48\textwidth}|}
\caption{Non-dendrites of   type 4 ($N=25$) }\label{tab:n4}\\
    \hline
    \multicolumn{2}{|c|}{\includegraphics{nonden4.pdf}} \\
    \hline
    \includegraphics[width=0.23\textwidth] {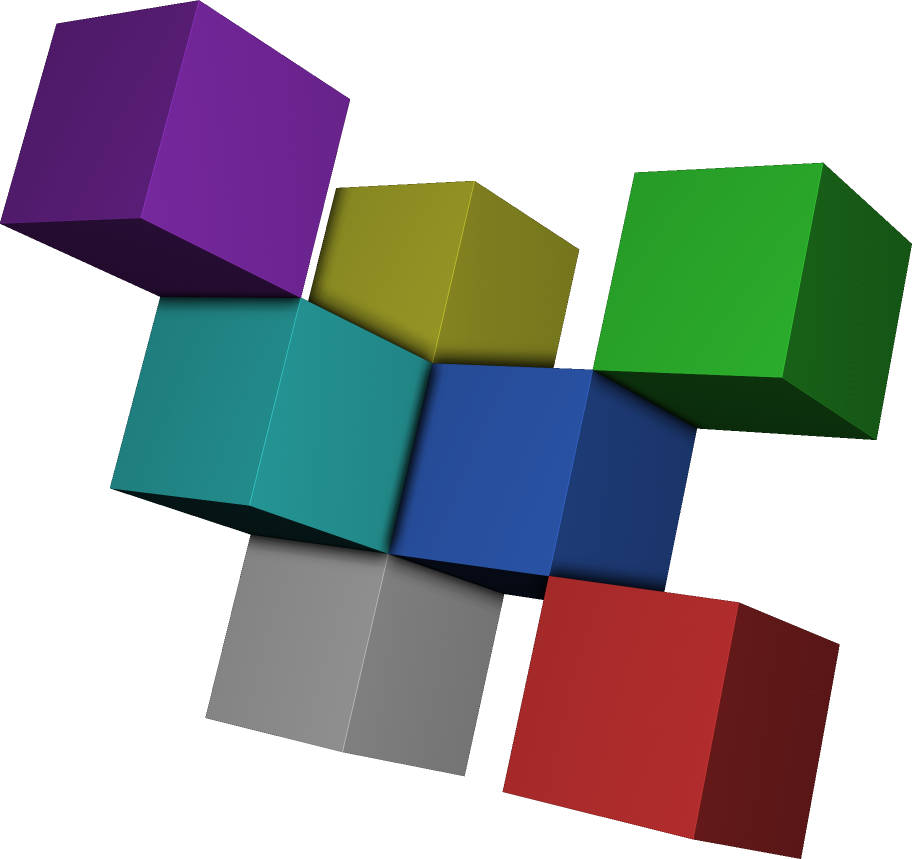}
    \includegraphics[width=0.23\textwidth] {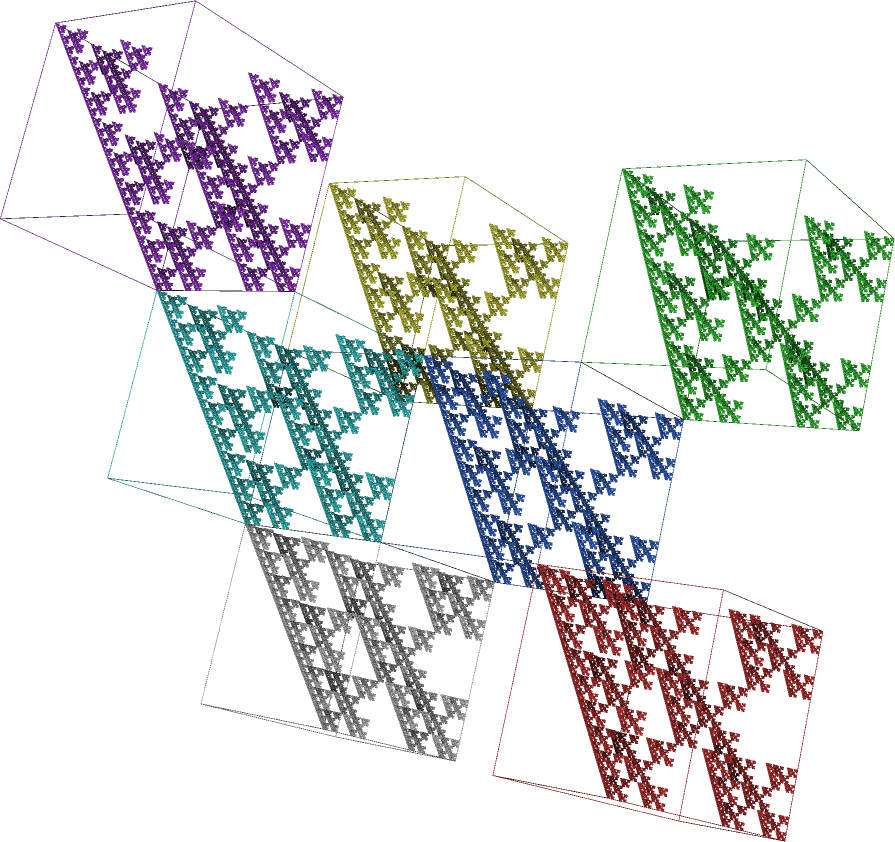} & 
    \includegraphics[width=0.23\textwidth] {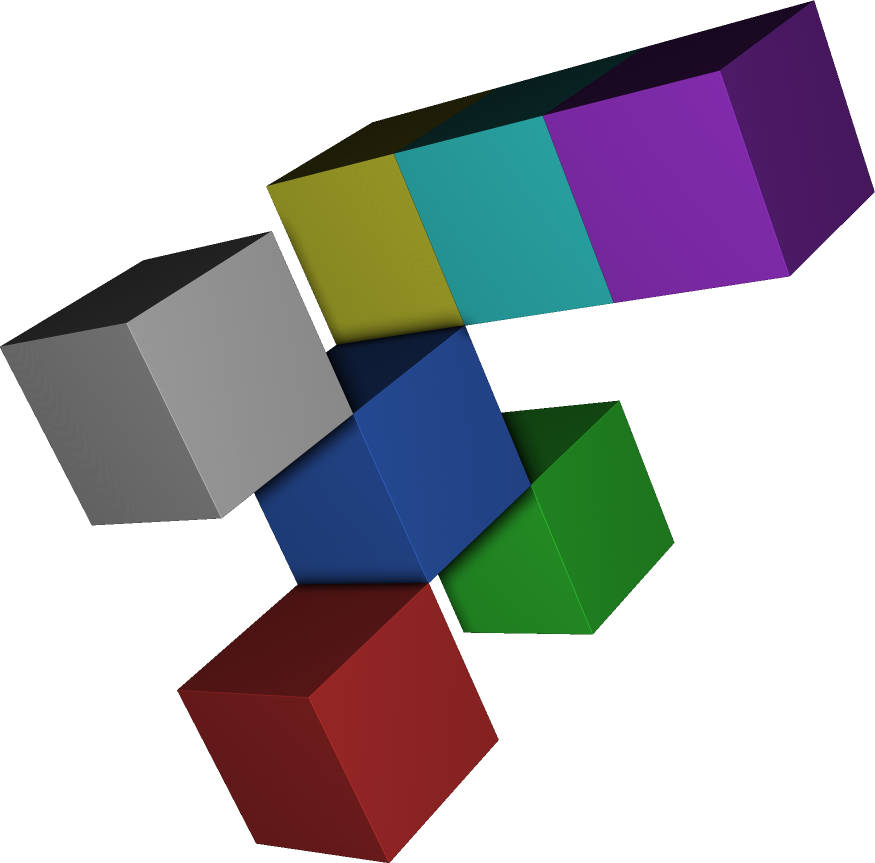}
    \includegraphics[width=0.23\textwidth] {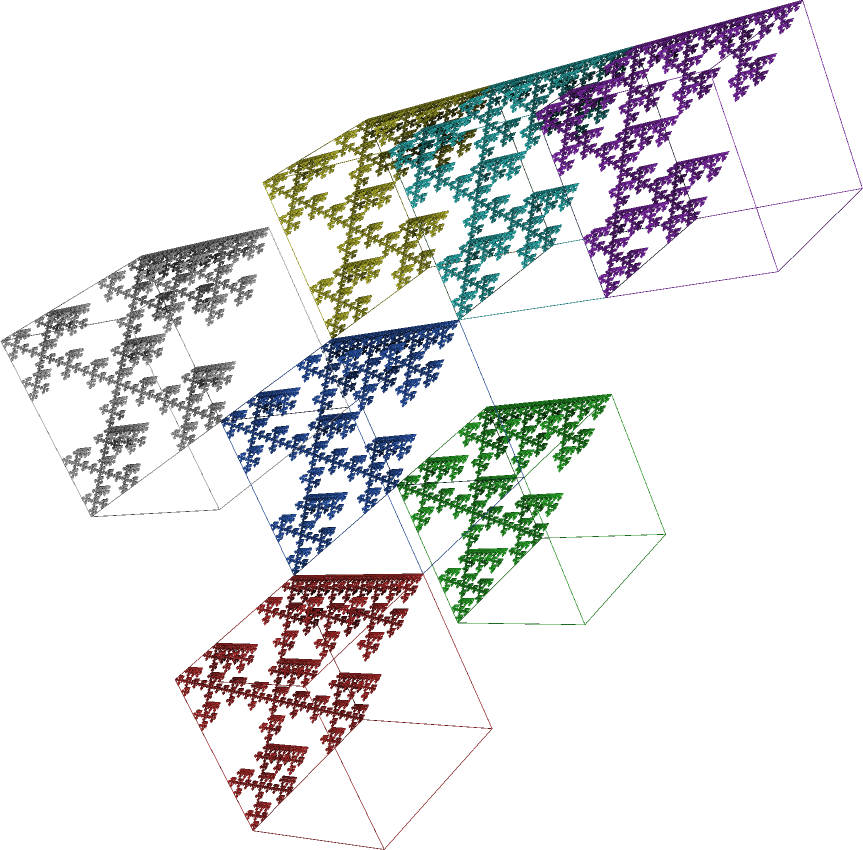} \\
    {\scriptsize$\mD=\{(0,0,0),\ (0,0,2),\ (0,1,1),\ (0,2,0),\ (0,2,2),\ (1,0,1),\ (2,0,2)\}$} & 
    {\scriptsize$\mD=\{(0,0,0),\ (0,0,2),\ (0,1,1),\ (0,2,0),\ (0,2,2),\ (1,0,2),\ (2,0,2)\}$} \\
    \hline
    \includegraphics[width=0.23\textwidth] {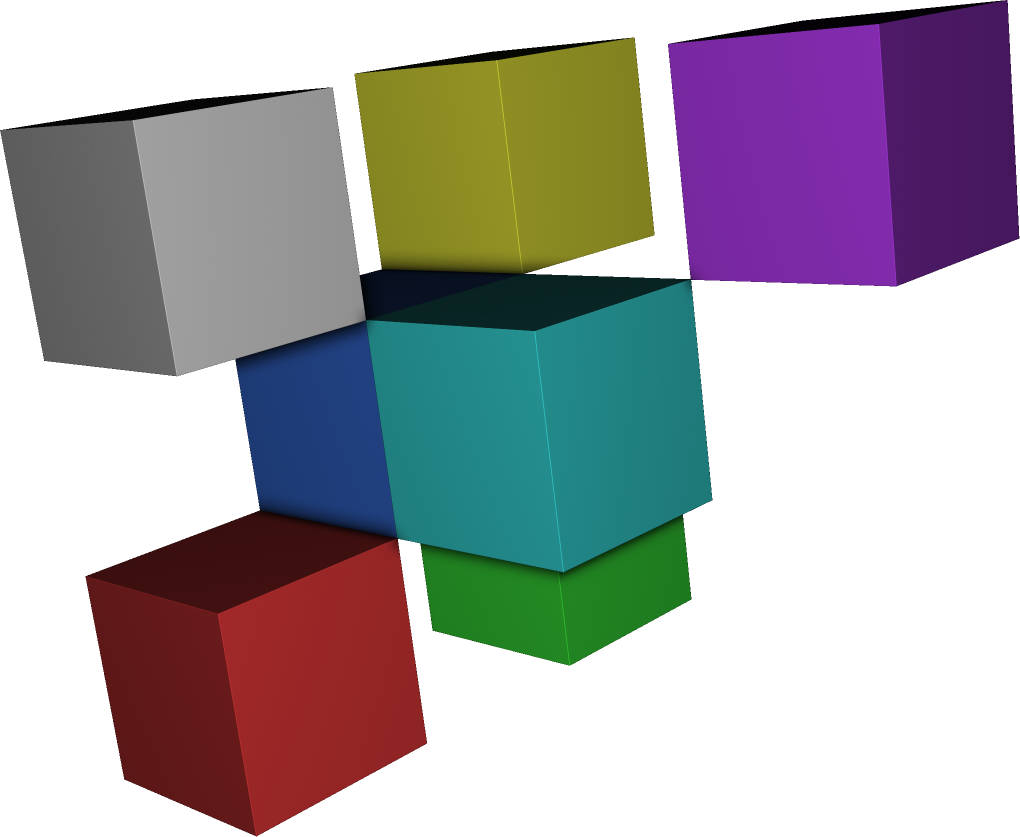}
    \includegraphics[width=0.23\textwidth] {n4_1057109__000_002_011_020_022_111_202__Q.jpg} & 
    \includegraphics[width=0.23\textwidth] {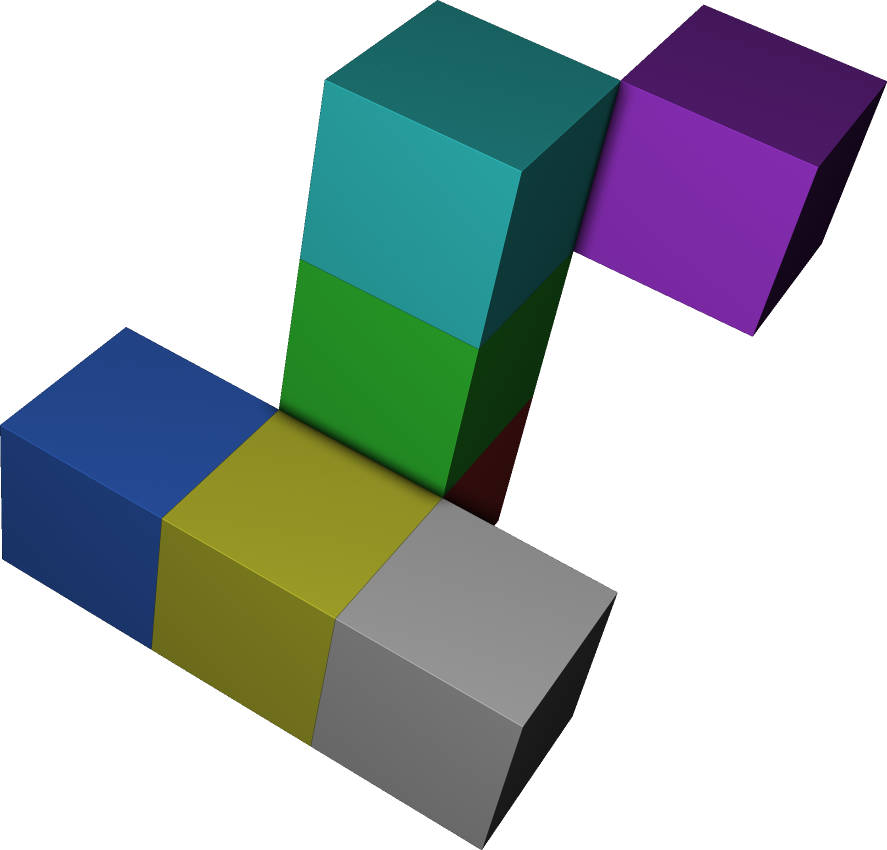}
    \includegraphics[width=0.23\textwidth] {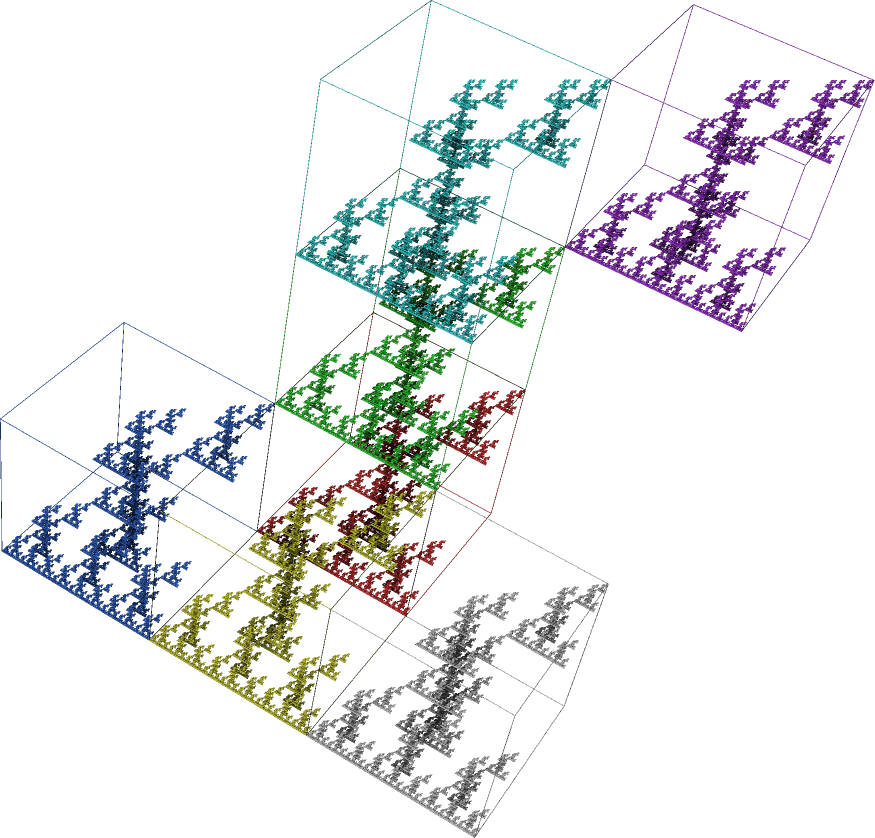} \\
    {\scriptsize$\mD=\{(0,0,0),\ (0,0,2),\ (0,1,1),\ (0,2,0),\ (0,2,2),\ (1,1,1),\ (2,0,2)\}$} & 
    {\scriptsize$\mD=\{(0,0,0),\ (0,1,0),\ (0,2,0),\ (1,1,0),\ (1,1,1),\ (1,1,2),\ (2,0,2)\}$} \\
    \hline
    \includegraphics[width=0.23\textwidth] {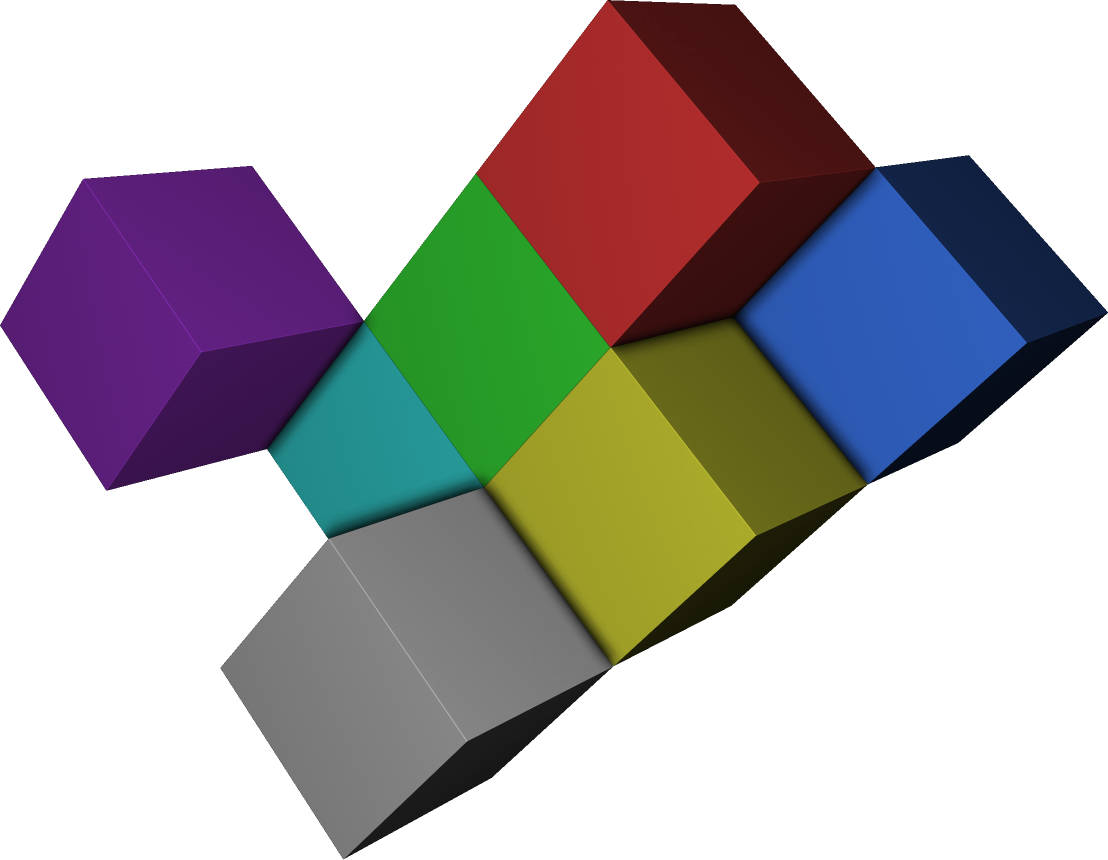}
    \includegraphics[width=0.23\textwidth] {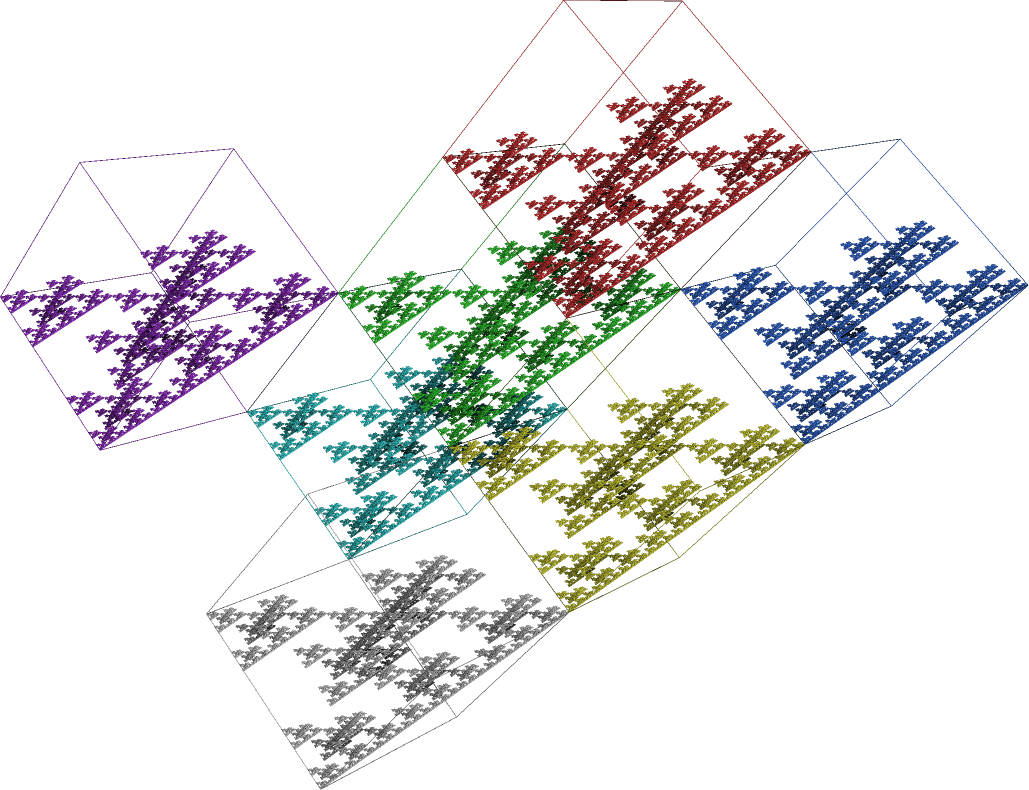} & 
    \includegraphics[width=0.23\textwidth] {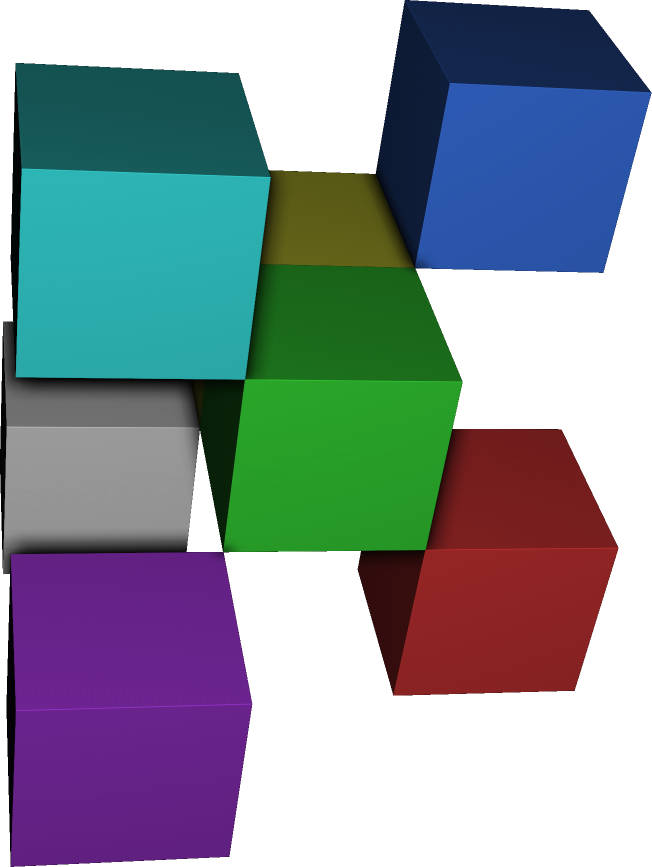}
    \includegraphics[width=0.23\textwidth] {n4_1188436__002_011_020_100_111_122_202__Q.jpg} \\
    {\scriptsize$\mD=\{(0,0,2),\ (0,1,1),\ (0,2,0),\ (1,1,0),\ (1,1,1),\ (1,1,2),\ (2,0,2)\}$} & 
    {\scriptsize$\mD=\{(0,0,2),\ (0,1,1),\ (0,2,0),\ (1,0,0),\ (1,1,1),\ (1,2,2),\ (2,0,2)\}$} \\
    \hline
    \includegraphics[width=0.23\textwidth] {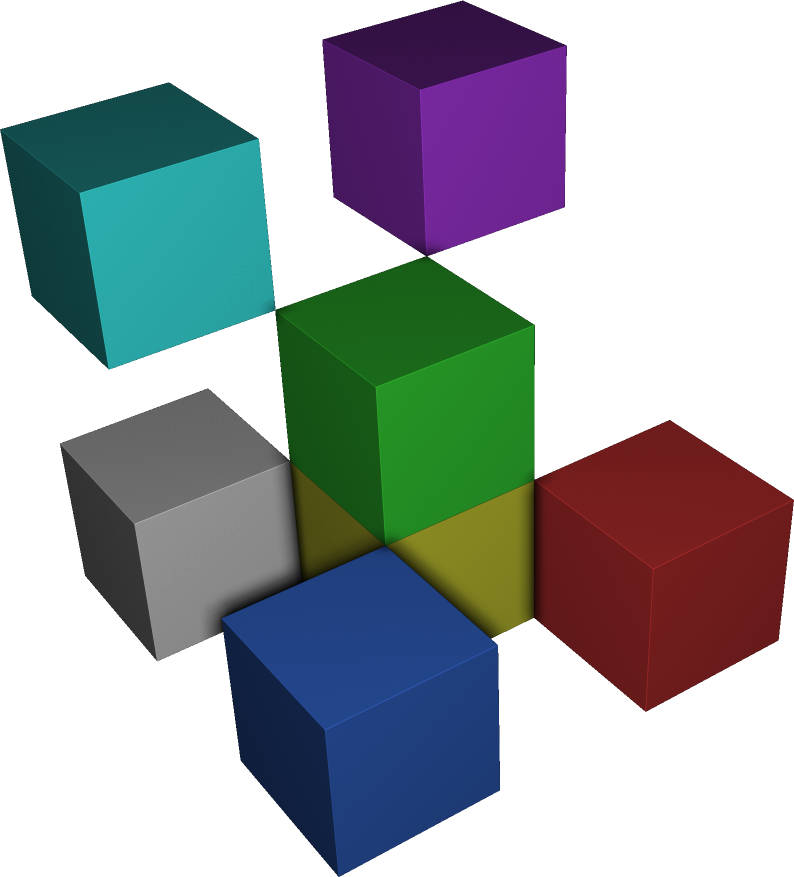}
    \includegraphics[width=0.23\textwidth] {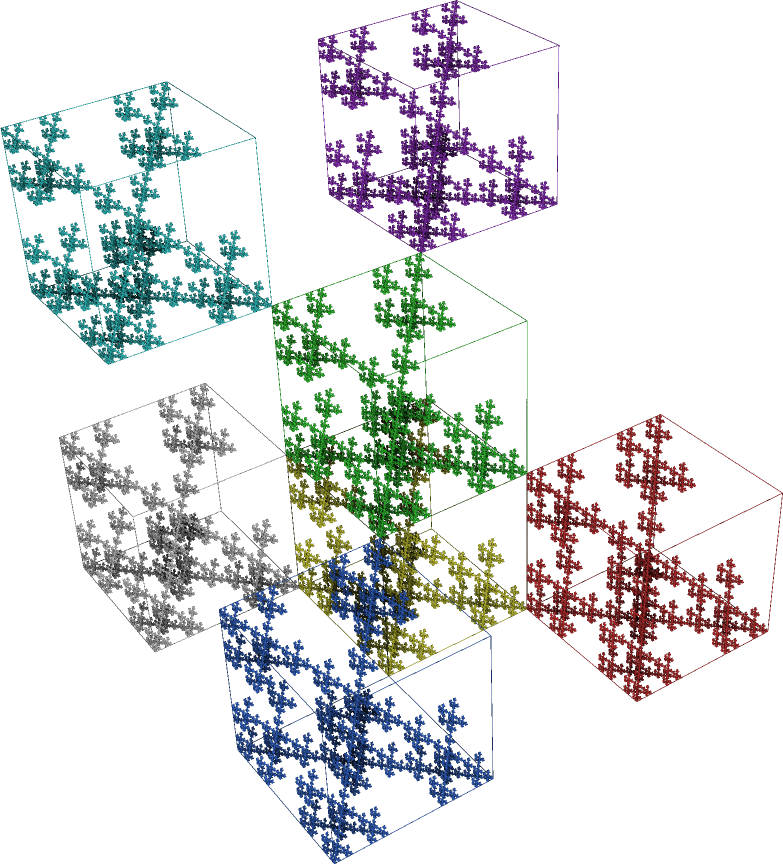} & 
    \includegraphics[width=0.23\textwidth] {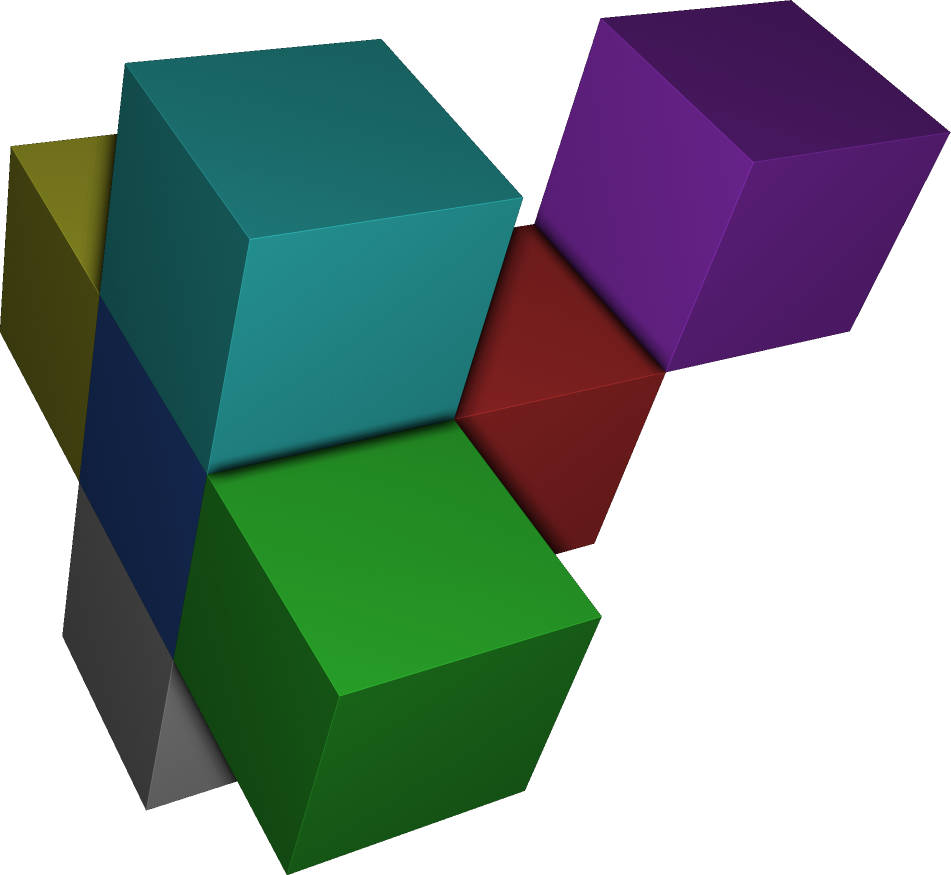}
    \includegraphics[width=0.23\textwidth] {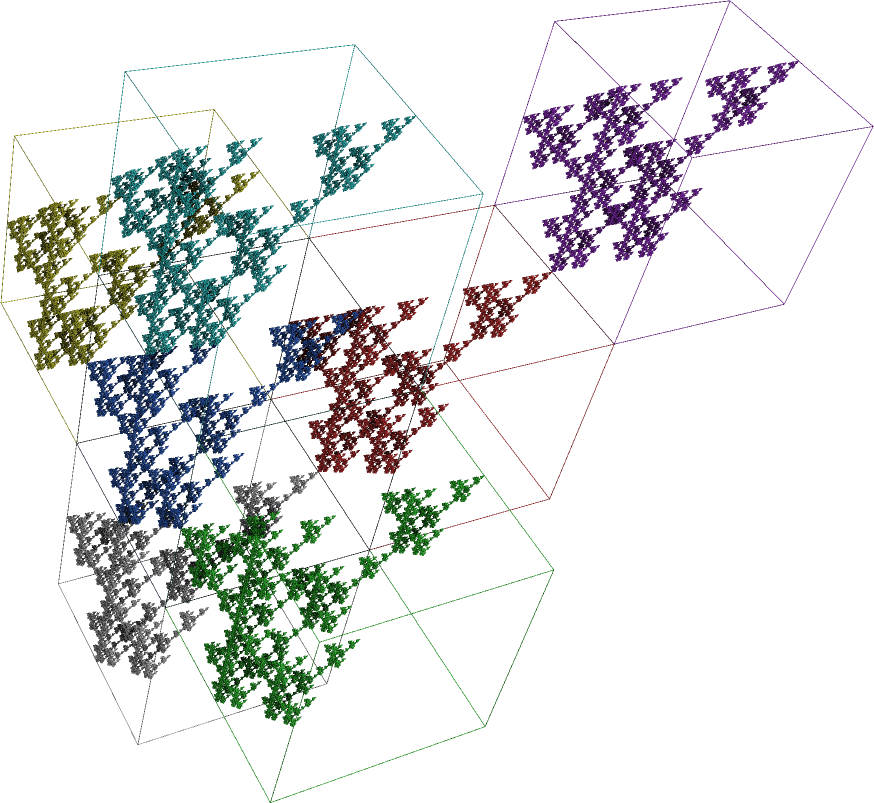} \\
    {\scriptsize$\mD=\{(0,0,0),\ (0,1,1),\ (0,2,0),\ (0,2,2),\ (1,1,1),\ (2,0,0),\ (2,0,2)\}$} & 
    {\scriptsize$\mD=\{(0,1,0),\ (1,0,0),\ (1,1,0),\ (1,1,1),\ (1,2,0),\ (2,1,0),\ (2,1,2)\}$} \\
    \hline
    \includegraphics[width=0.23\textwidth] {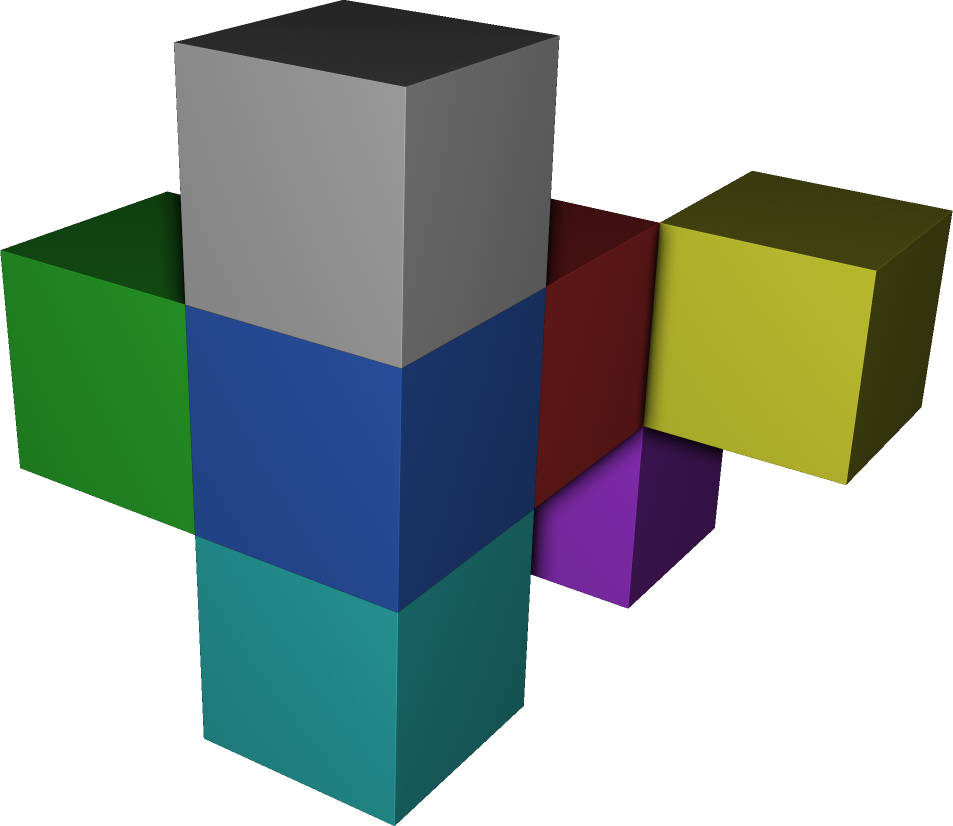}
    \includegraphics[width=0.23\textwidth] {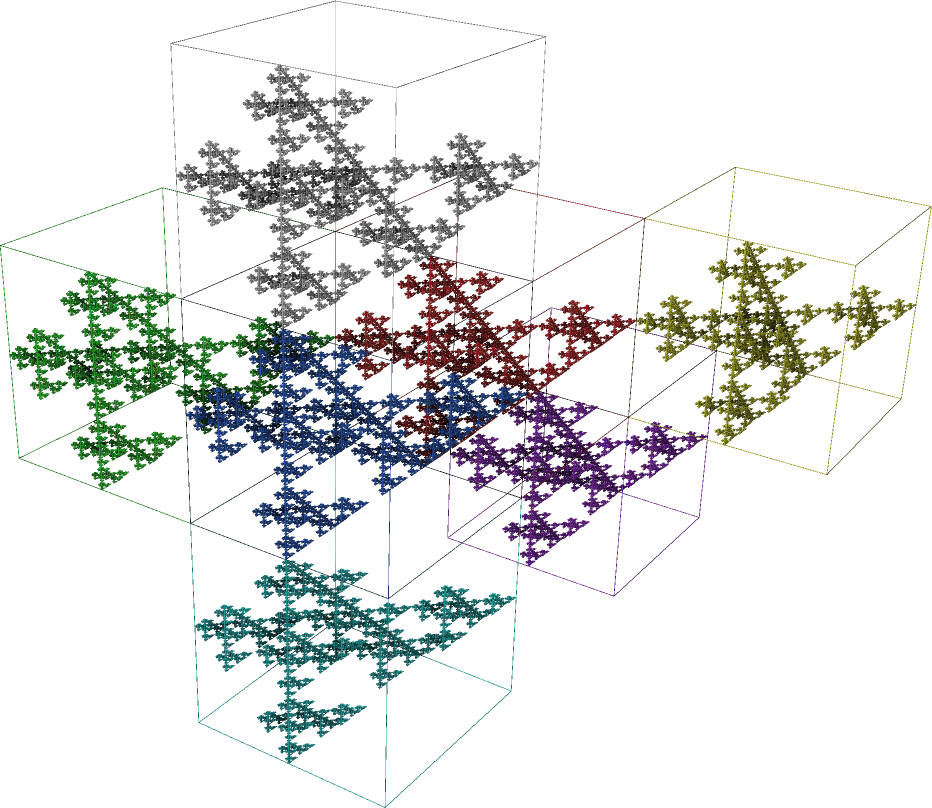} & 
    \includegraphics[width=0.23\textwidth] {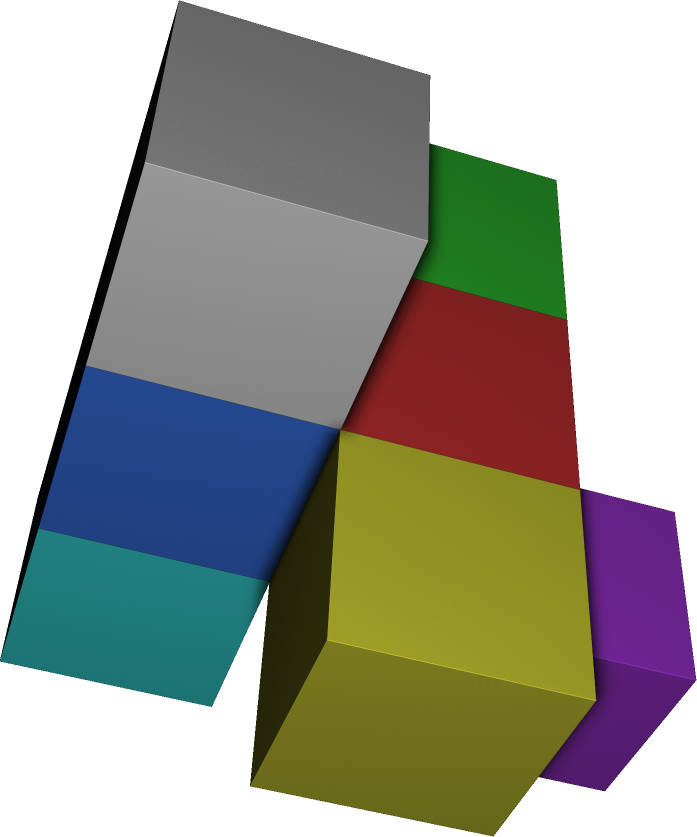}
    \includegraphics[width=0.23\textwidth] {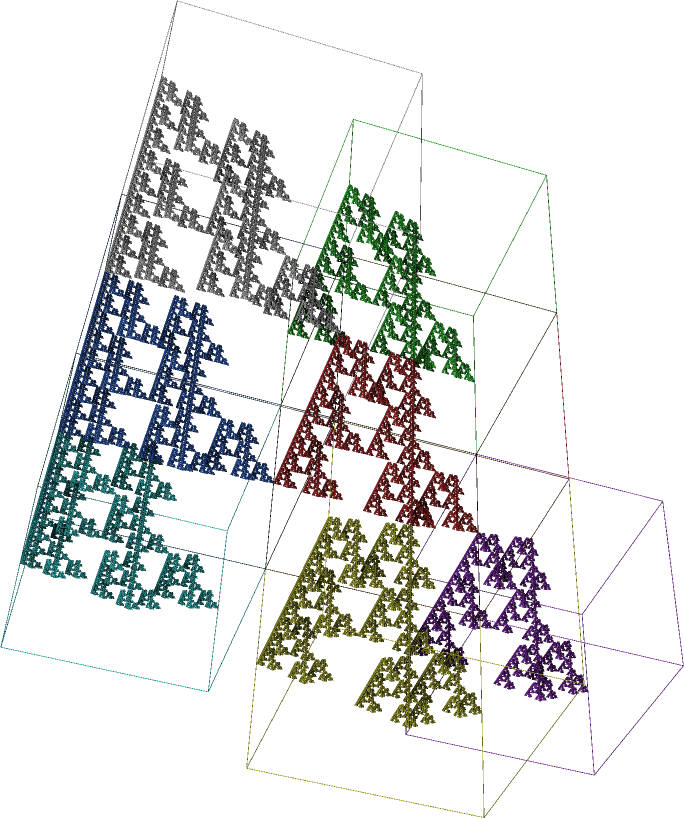} \\
    {\scriptsize$\mD=\{(0,1,0),\ (1,0,2),\ (1,1,0),\ (1,1,1),\ (1,2,0),\ (2,1,0),\ (2,1,2)\}$} & 
    {\scriptsize$\mD=\{(0,1,0),\ (1,0,1),\ (1,1,0),\ (1,1,1),\ (1,2,1),\ (2,1,0),\ (2,1,2)\}$} \\
    \hline
    \includegraphics[width=0.23\textwidth] {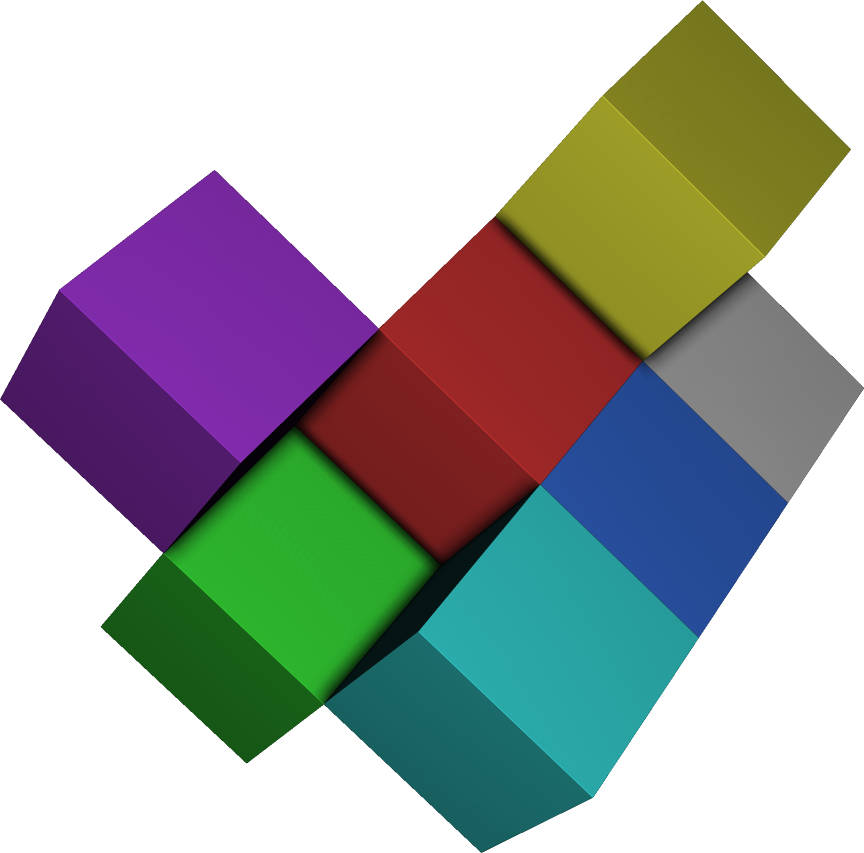}
    \includegraphics[width=0.23\textwidth] {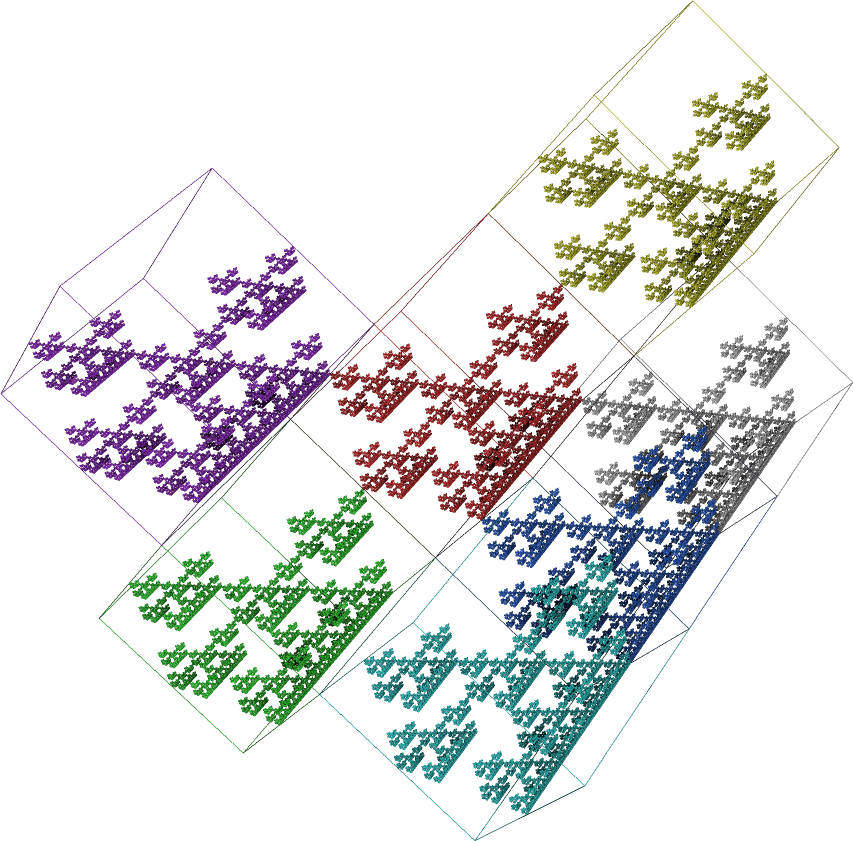} & 
    \includegraphics[width=0.23\textwidth] {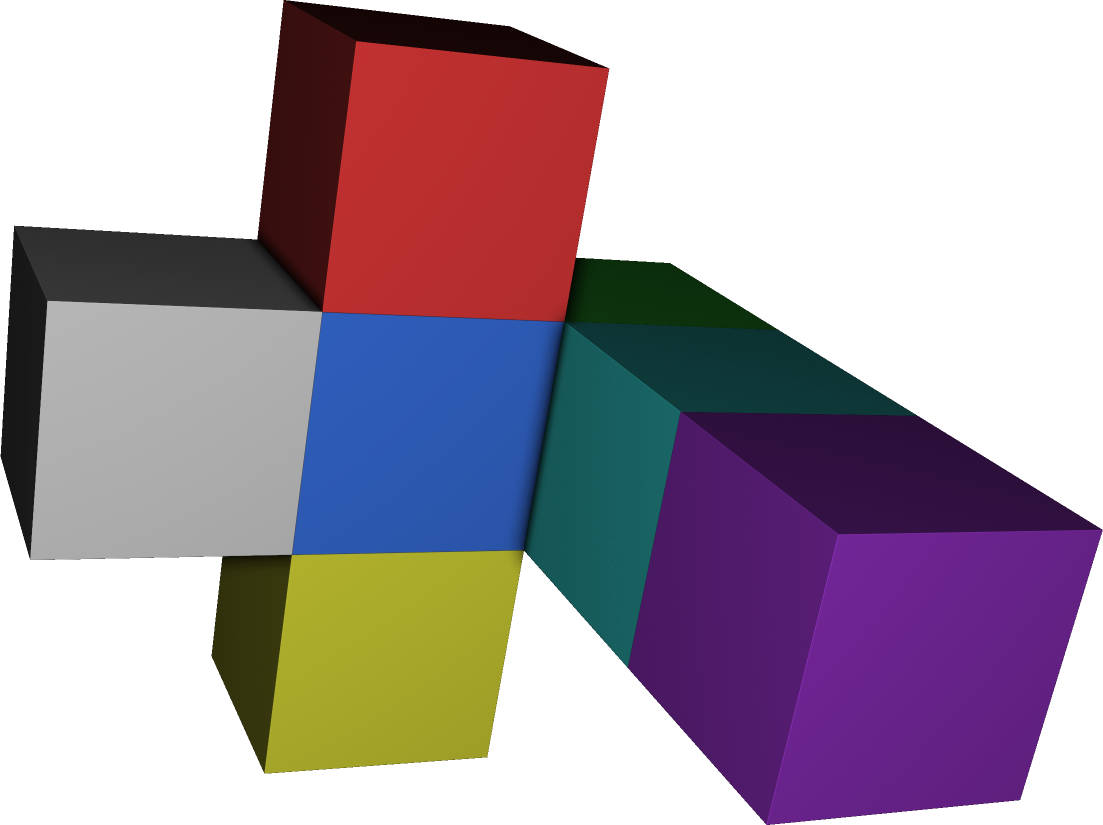}
    \includegraphics[width=0.23\textwidth] {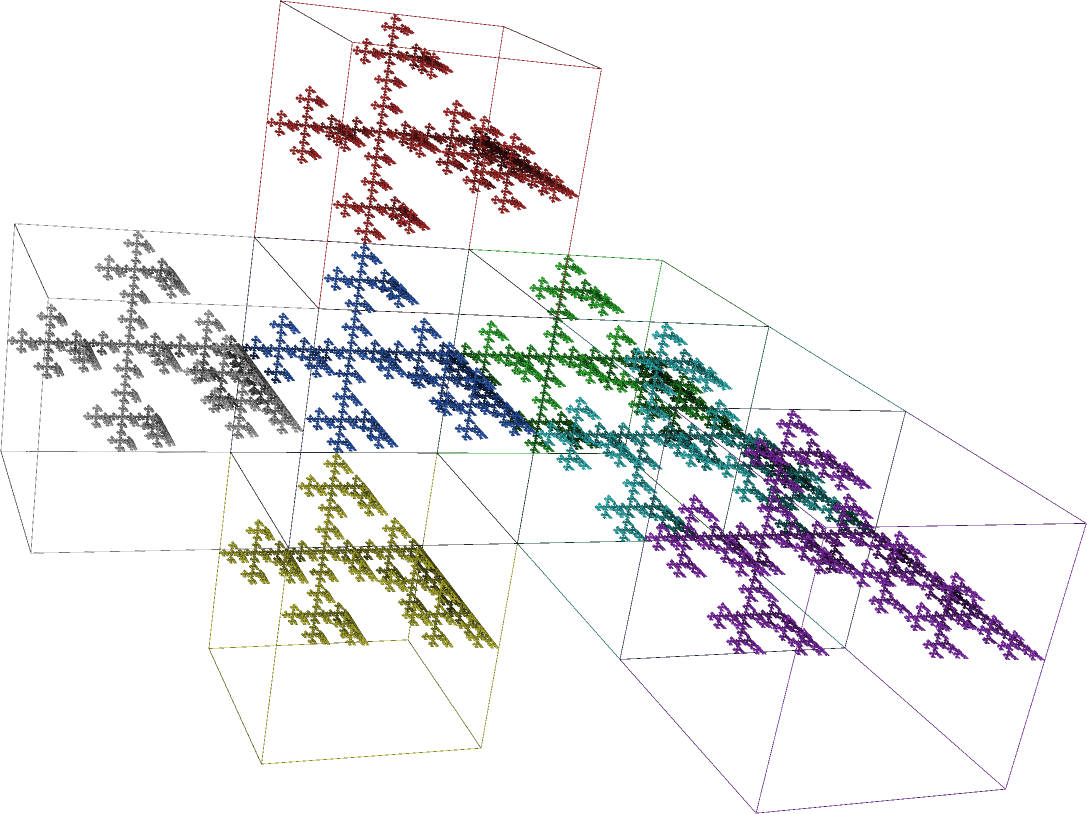} \\
    {\scriptsize$\mD=\{(0,1,0),\ (0,2,1),\ (1,1,0),\ (1,1,1),\ (2,0,1),\ (2,1,0),\ (2,1,2)\}$} & 
    {\scriptsize$\mD=\{(0,1,0),\ (1,0,0),\ (1,1,0),\ (1,2,0),\ (2,1,0),\ (2,1,1),\ (2,1,2)\}$} \\
    \hline
    \includegraphics[width=0.23\textwidth] {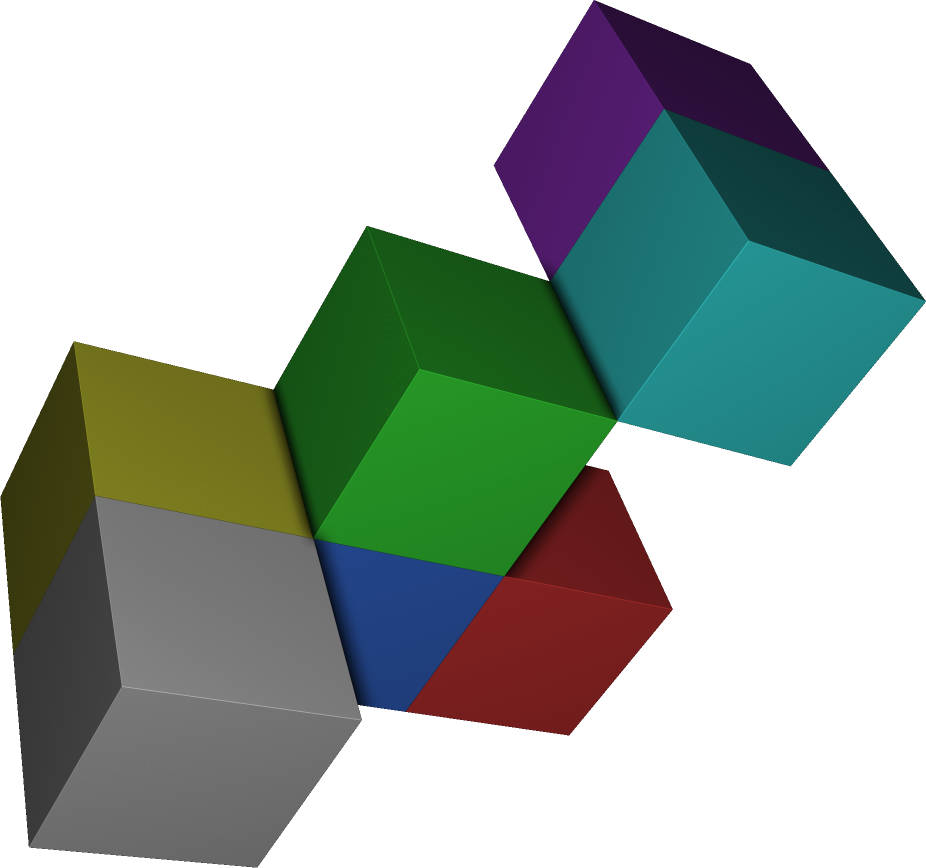}
    \includegraphics[width=0.23\textwidth] {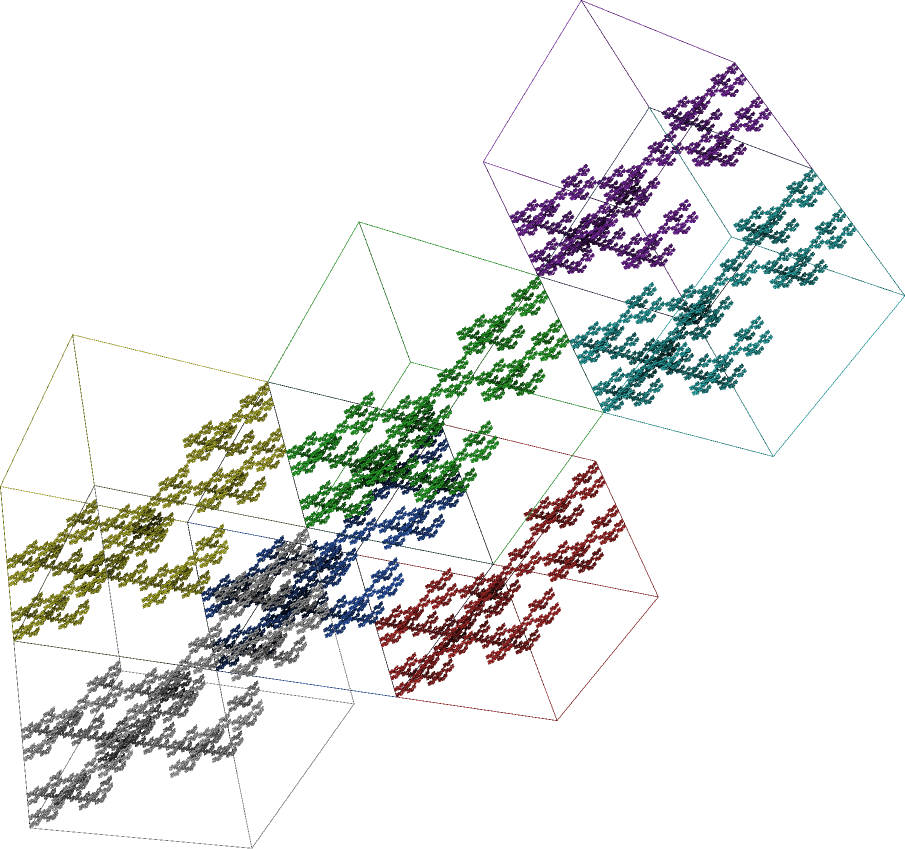} & 
    \includegraphics[width=0.23\textwidth] {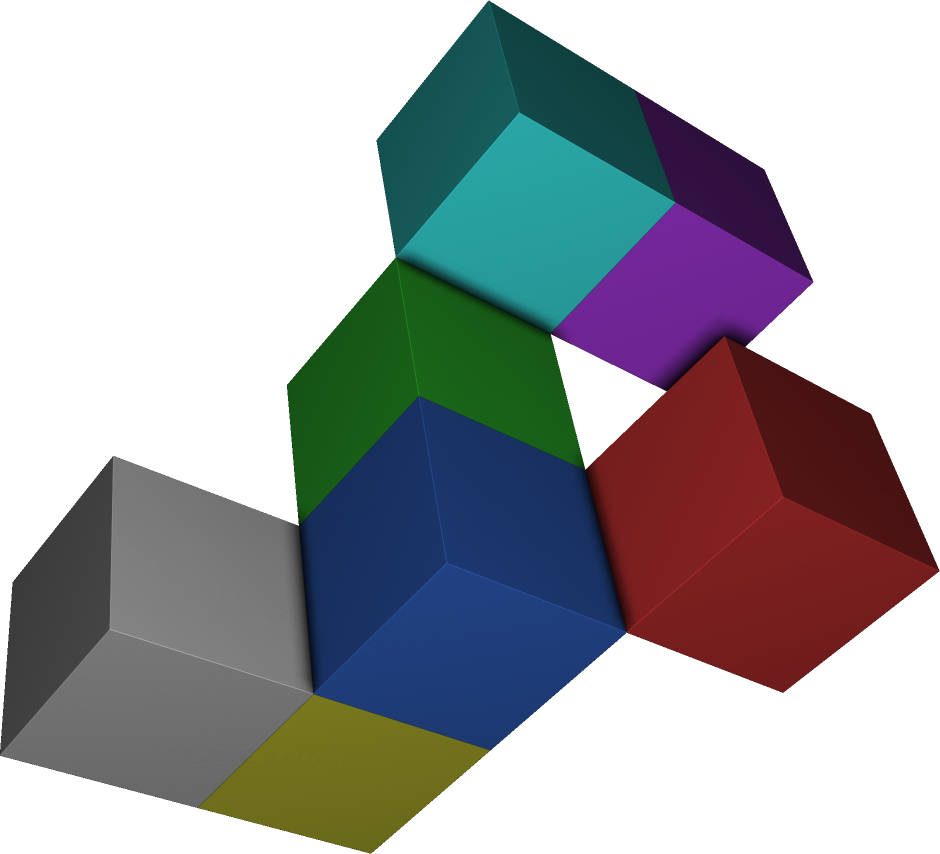}
    \includegraphics[width=0.23\textwidth] {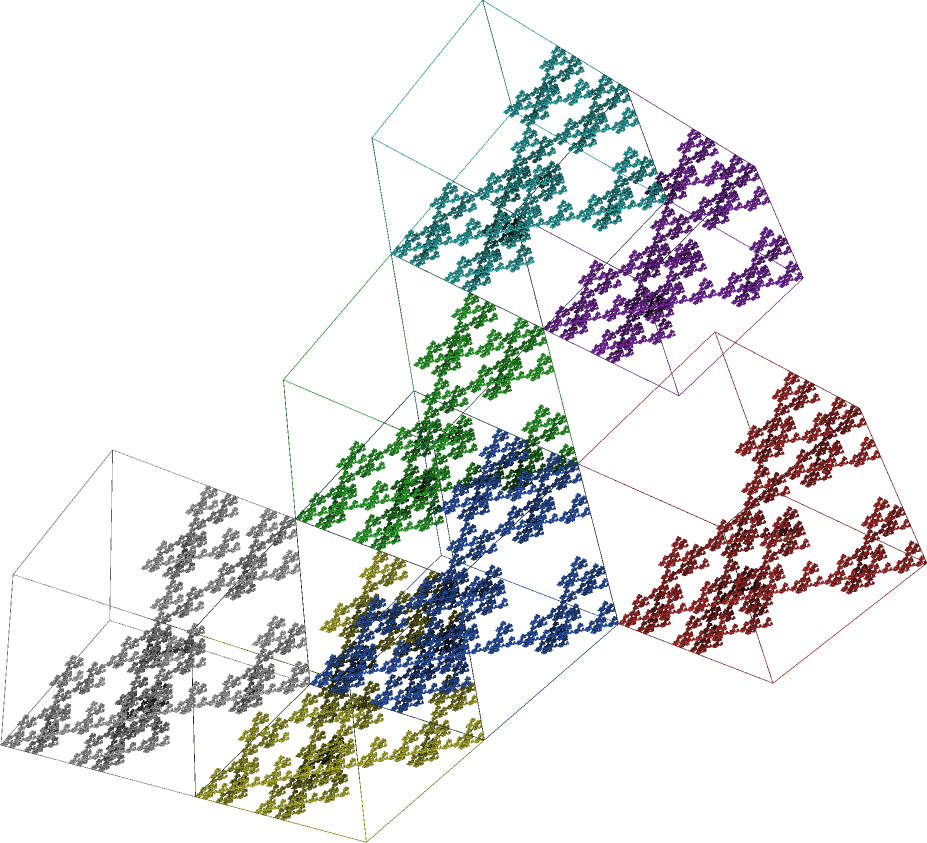} \\
    {\scriptsize$\mD=\{(0,0,0),\ (0,0,1),\ (0,1,1),\ (0,2,1),\ (1,1,1),\ (2,2,1),\ (2,2,2)\}$} & 
    {\scriptsize$\mD=\{(0,0,0),\ (0,0,1),\ (0,1,1),\ (0,2,2),\ (1,1,1),\ (2,2,1),\ (2,2,2)\}$} \\
    \hline
    \includegraphics[width=0.23\textwidth] {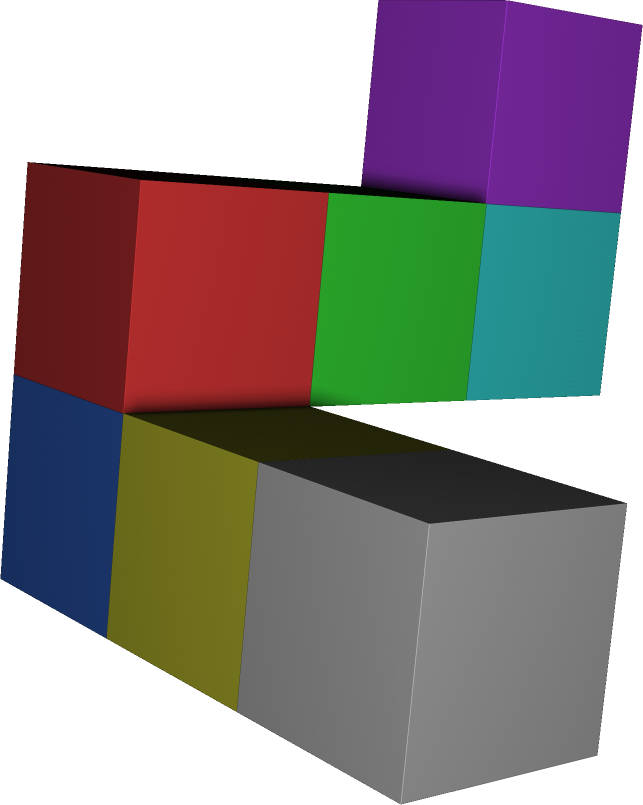}
    \includegraphics[width=0.23\textwidth] {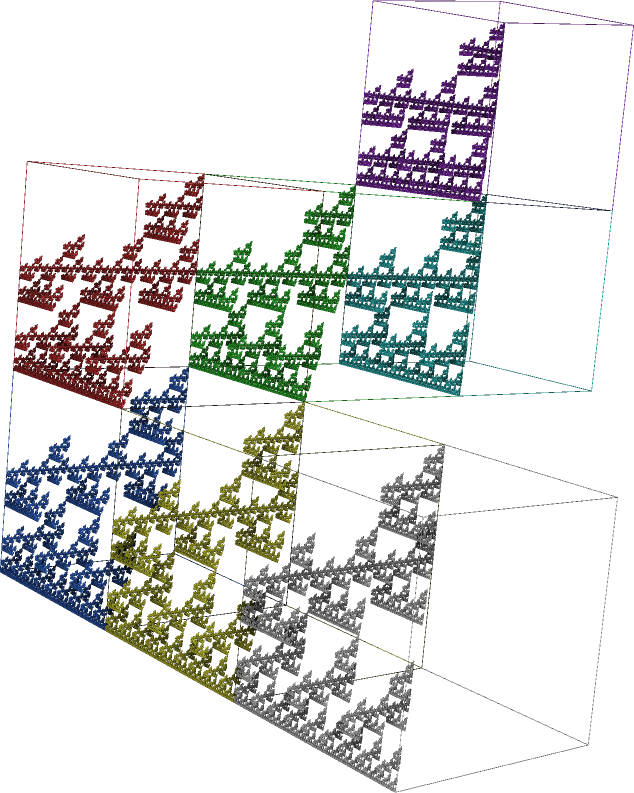} & 
    \includegraphics[width=0.23\textwidth] {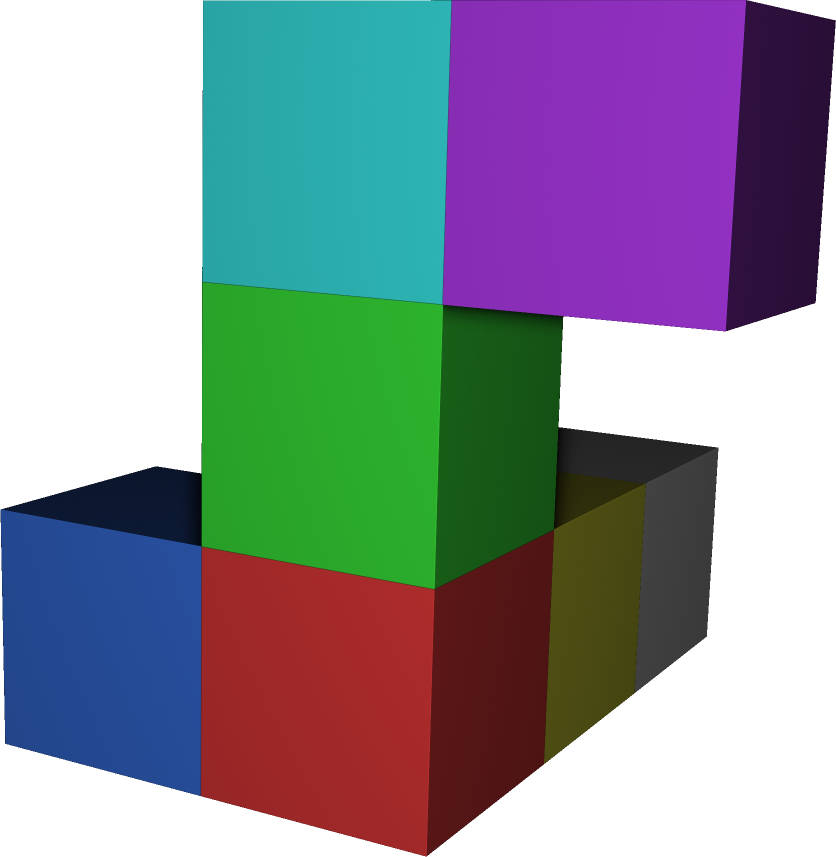}
    \includegraphics[width=0.23\textwidth] {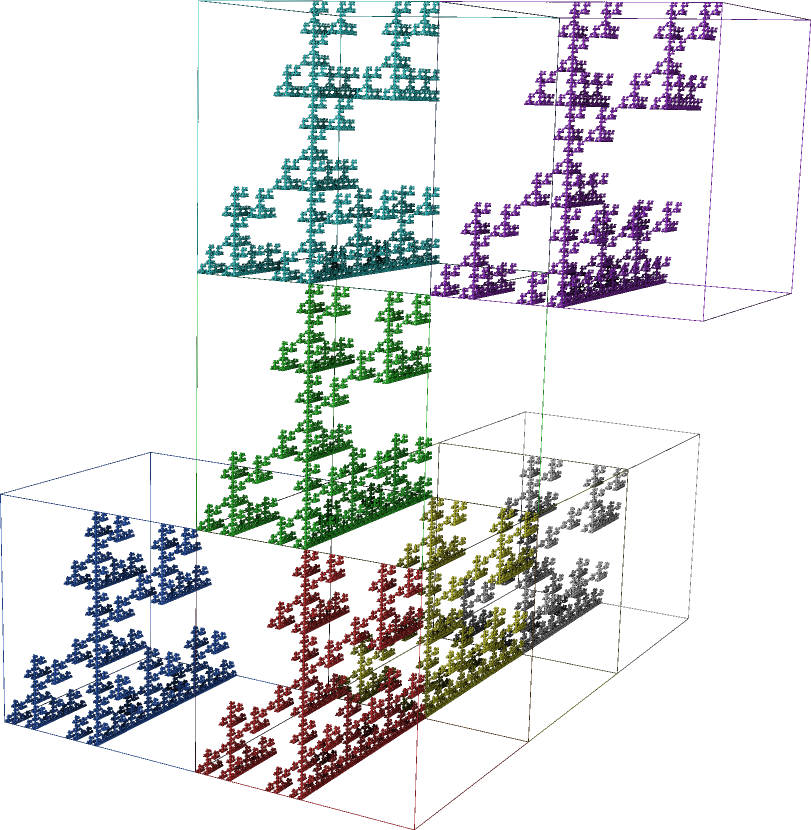} \\
    {\scriptsize$\mD=\{(0,0,0),\ (0,1,0),\ (0,2,0),\ (0,2,1),\ (1,2,1),\ (2,2,1),\ (2,2,2)\}$} & 
    {\scriptsize$\mD=\{(0,0,1),\ (0,1,1),\ (0,2,0),\ (0,2,1),\ (1,2,1),\ (2,2,1),\ (2,2,2)\}$} \\
    \hline
    \includegraphics[width=0.23\textwidth] {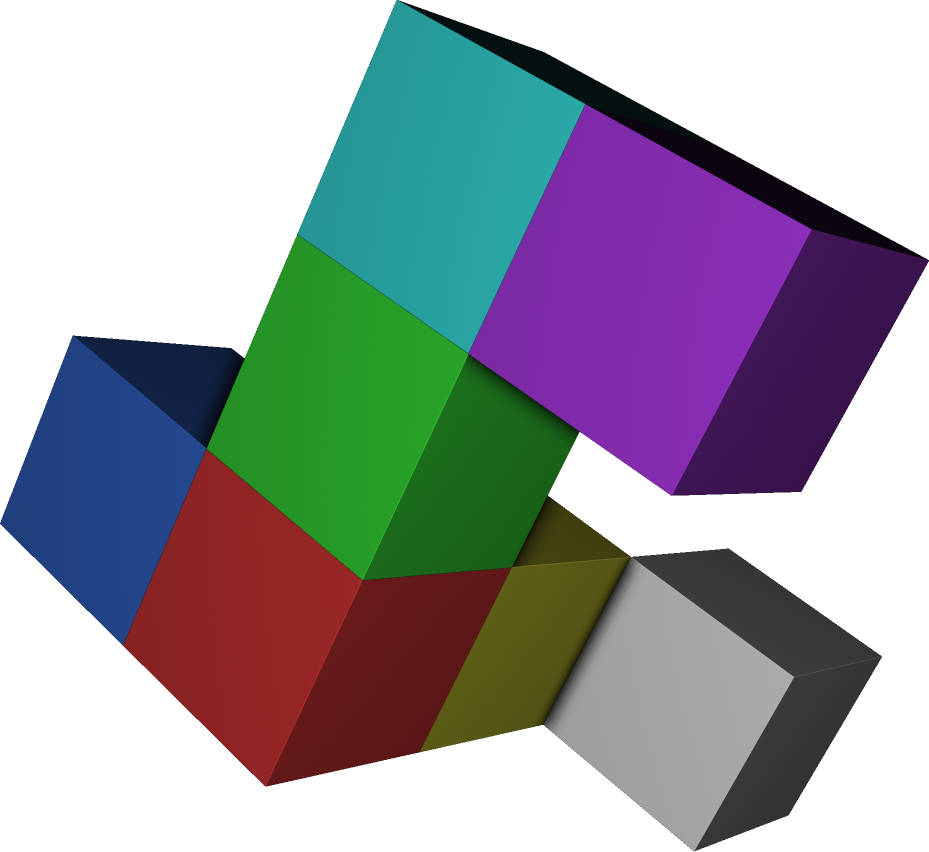}
    \includegraphics[width=0.23\textwidth] {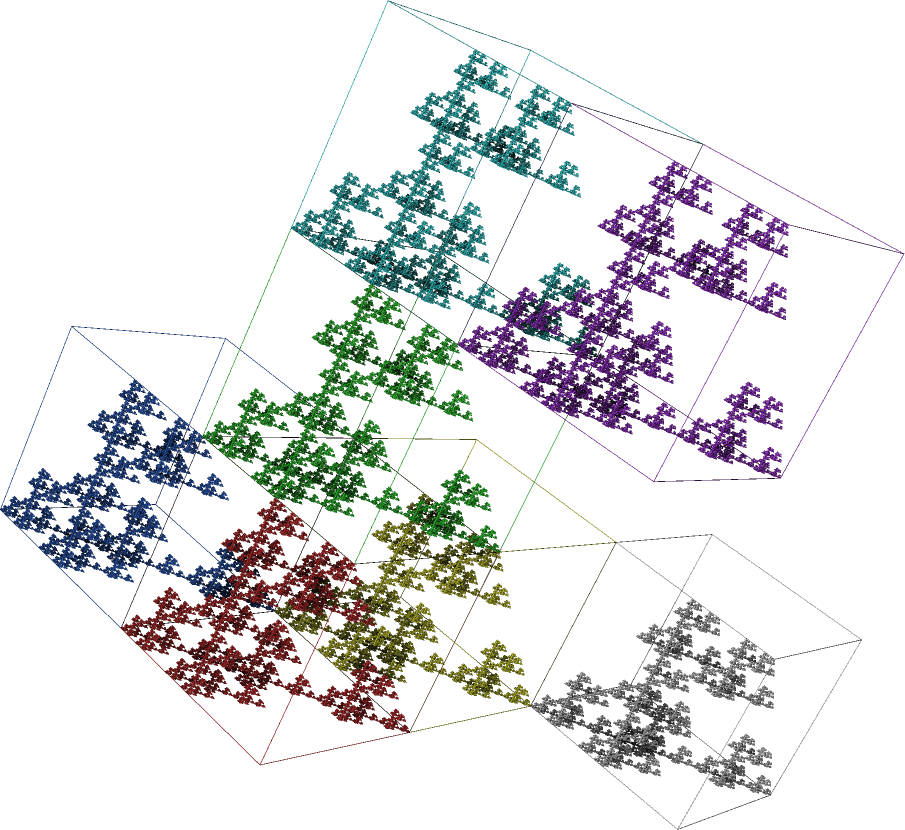} & 
    \includegraphics[width=0.23\textwidth] {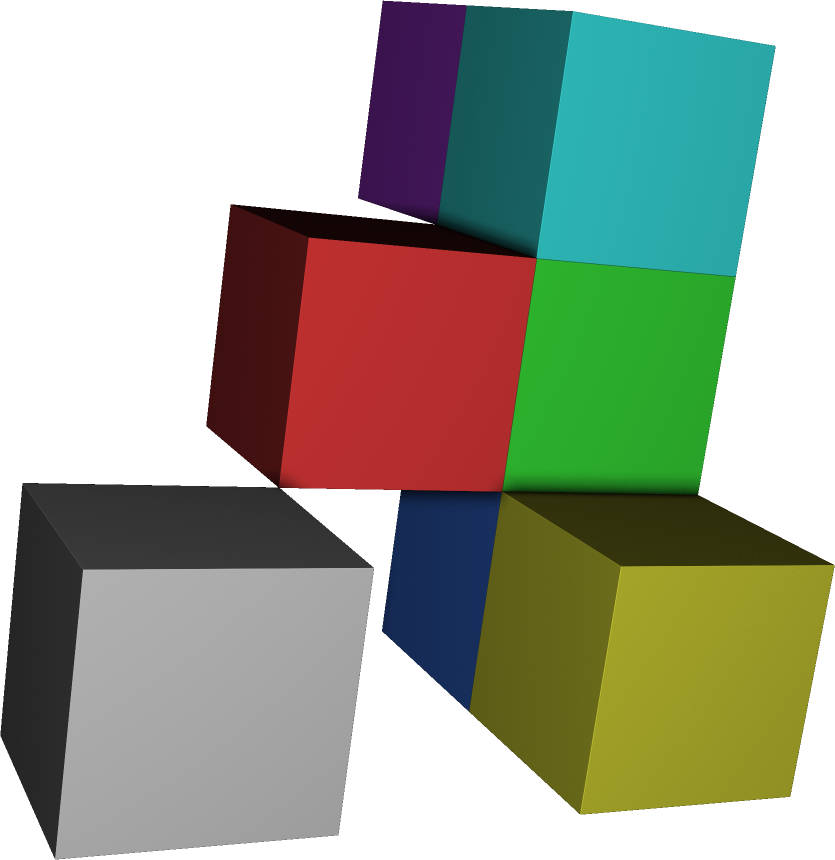}
    \includegraphics[width=0.23\textwidth] {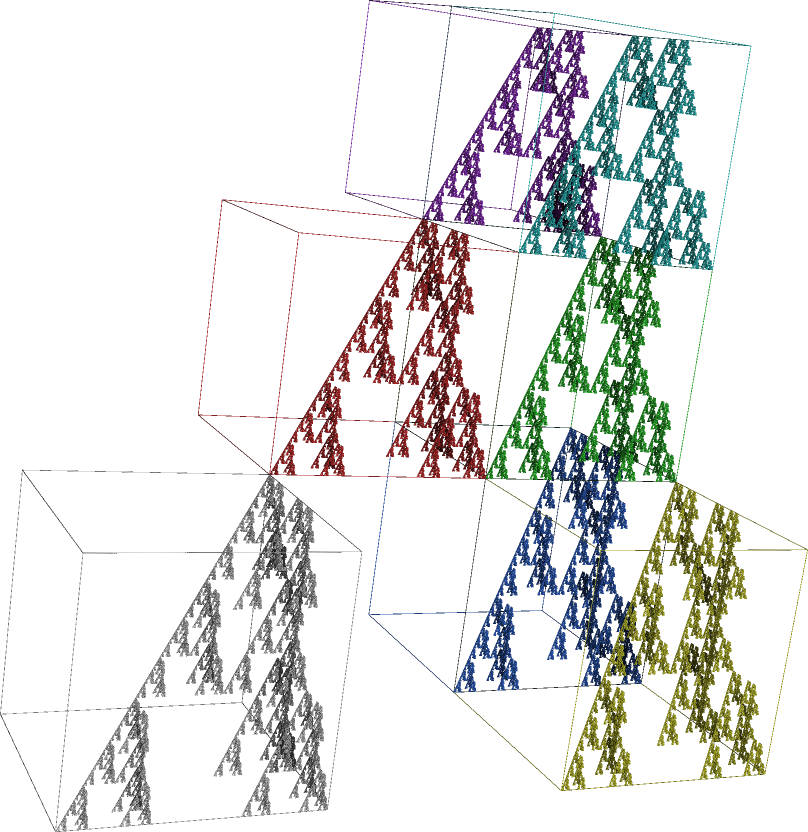} \\
    {\scriptsize$\mD=\{(0,0,2),\ (0,1,1),\ (0,2,0),\ (0,2,1),\ (1,2,1),\ (2,2,1),\ (2,2,2)\}$} & 
    {\scriptsize$\mD=\{(0,0,0),\ (0,2,0),\ (0,2,1),\ (1,1,1),\ (1,2,1),\ (2,2,1),\ (2,2,2)\}$} \\
    \hline
    \includegraphics[width=0.23\textwidth] {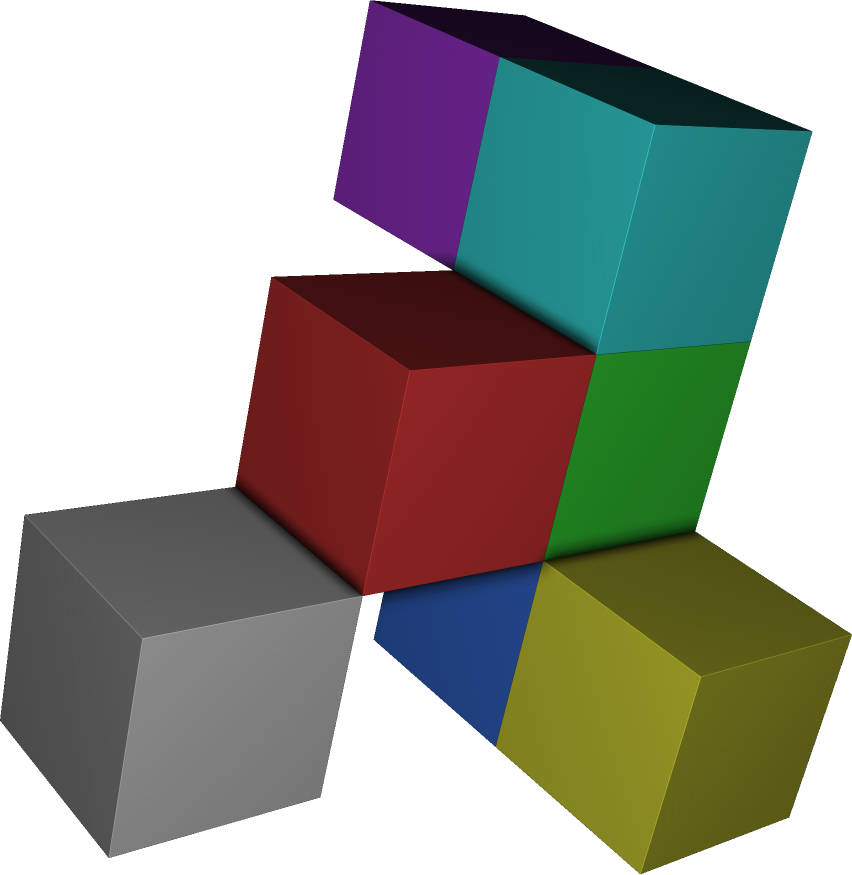}
    \includegraphics[width=0.23\textwidth] {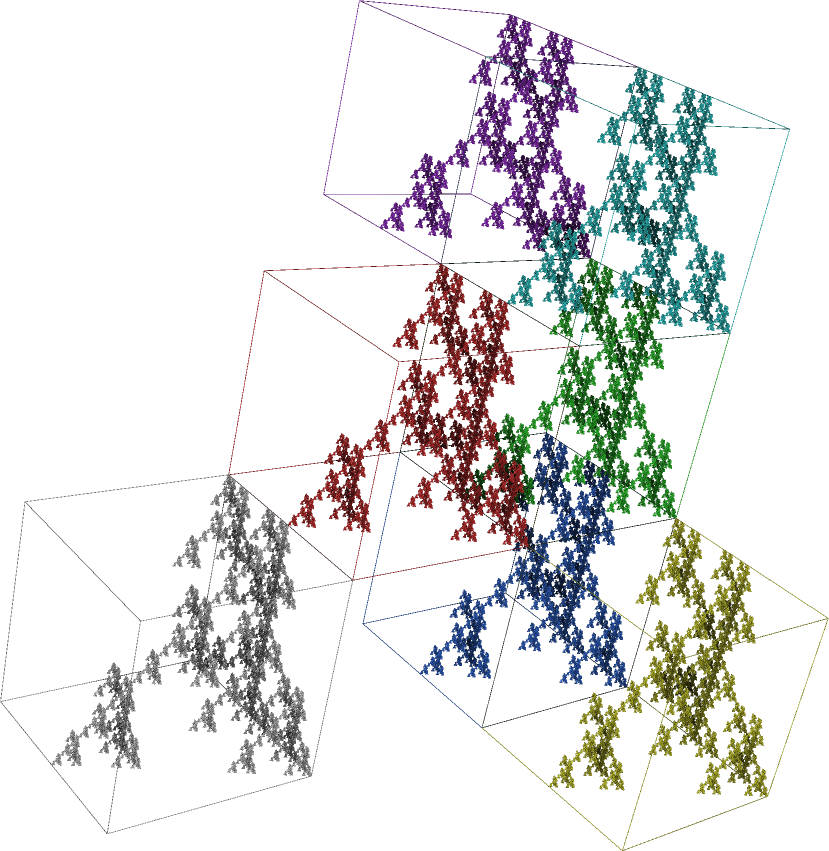} & 
    \includegraphics[width=0.23\textwidth] {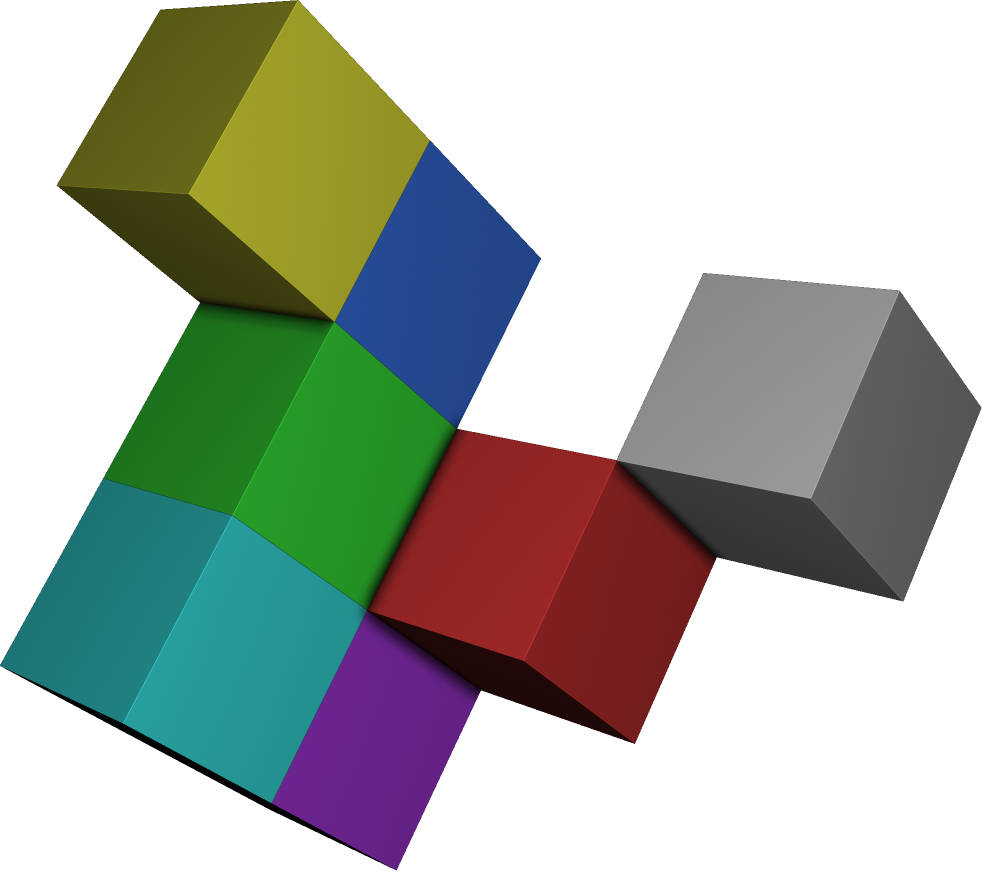}
    \includegraphics[width=0.23\textwidth] {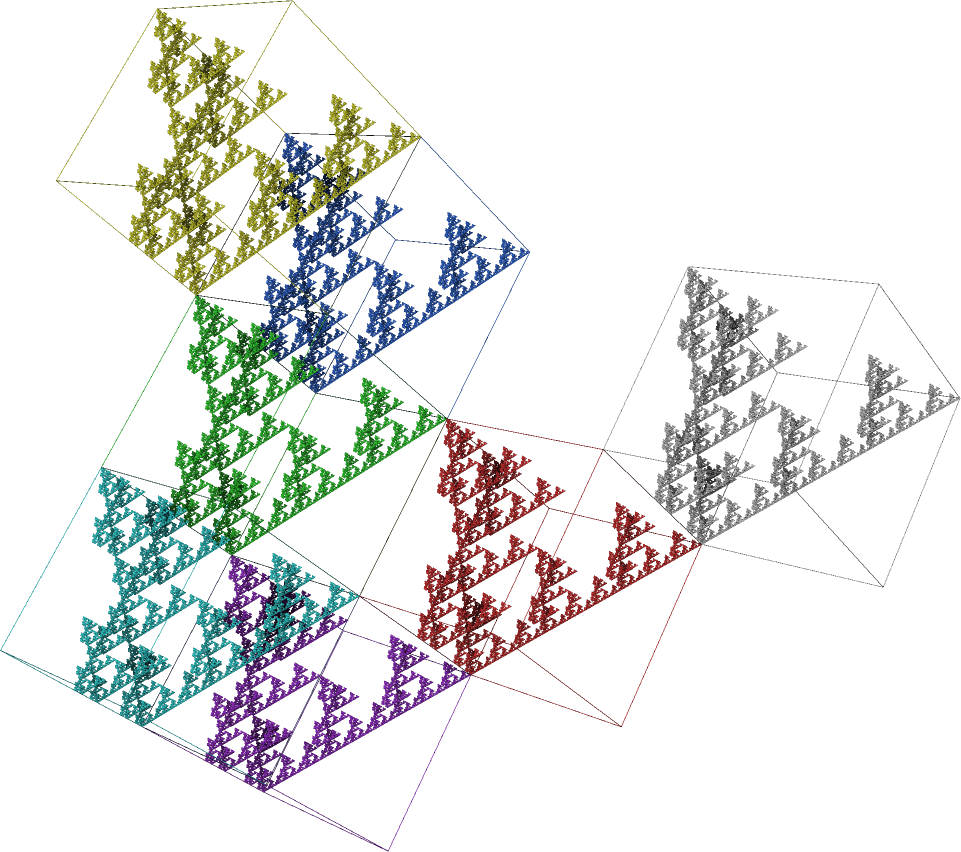} \\
    {\scriptsize$\mD=\{(0,0,1),\ (0,2,0),\ (0,2,1),\ (1,1,1),\ (1,2,1),\ (2,2,1),\ (2,2,2)\}$} & 
    {\scriptsize$\mD=\{(0,0,2),\ (0,2,0),\ (0,2,1),\ (1,1,2),\ (1,2,1),\ (2,2,1),\ (2,2,2)\}$} \\
    \hline
    \includegraphics[width=0.23\textwidth] {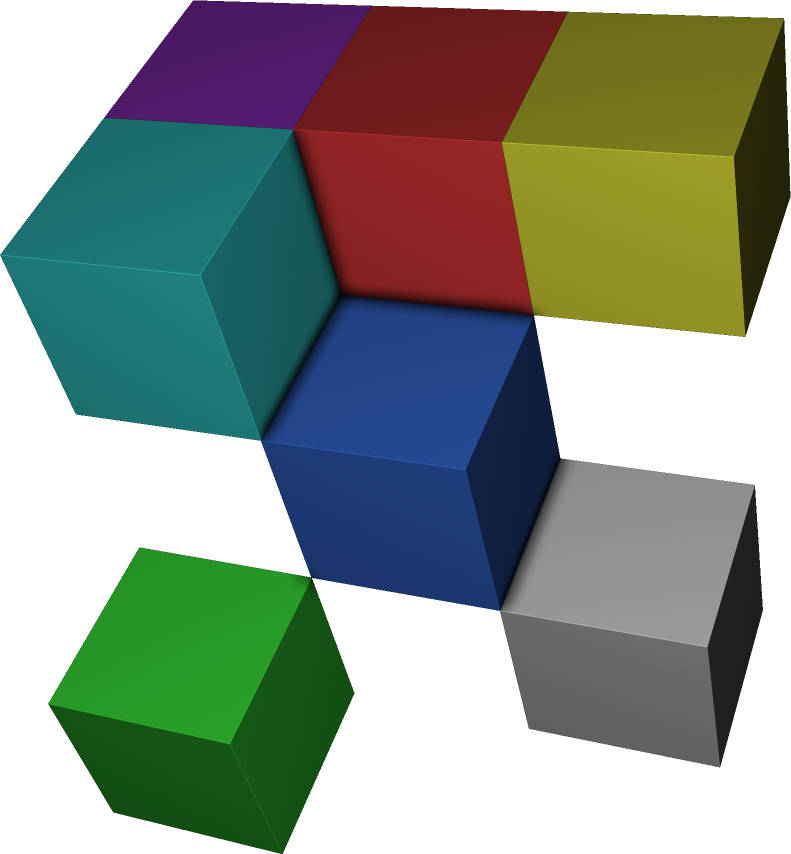}
    \includegraphics[width=0.23\textwidth] {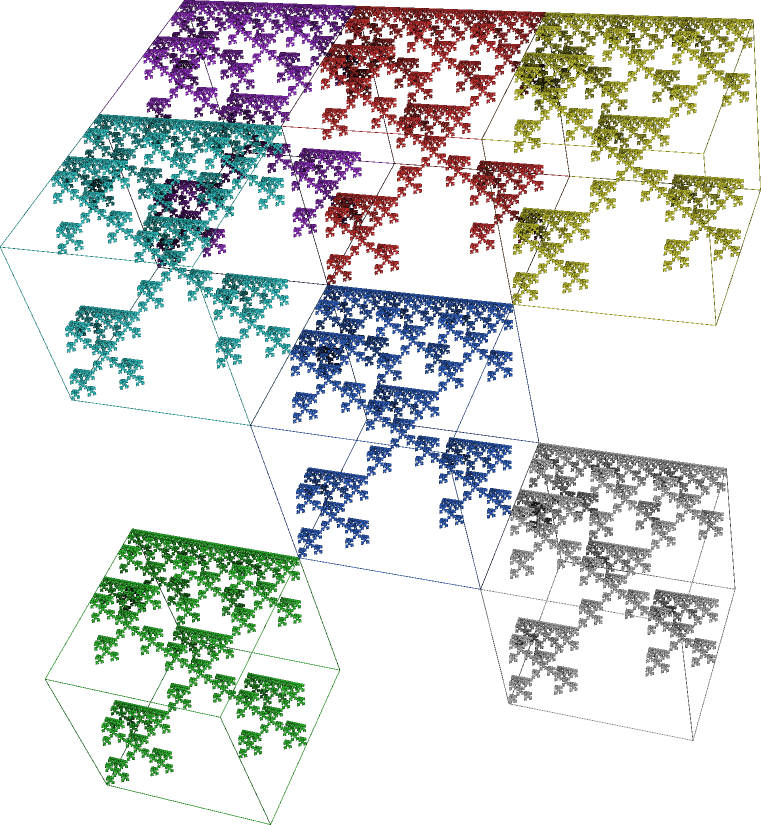} & 
    \includegraphics[width=0.23\textwidth] {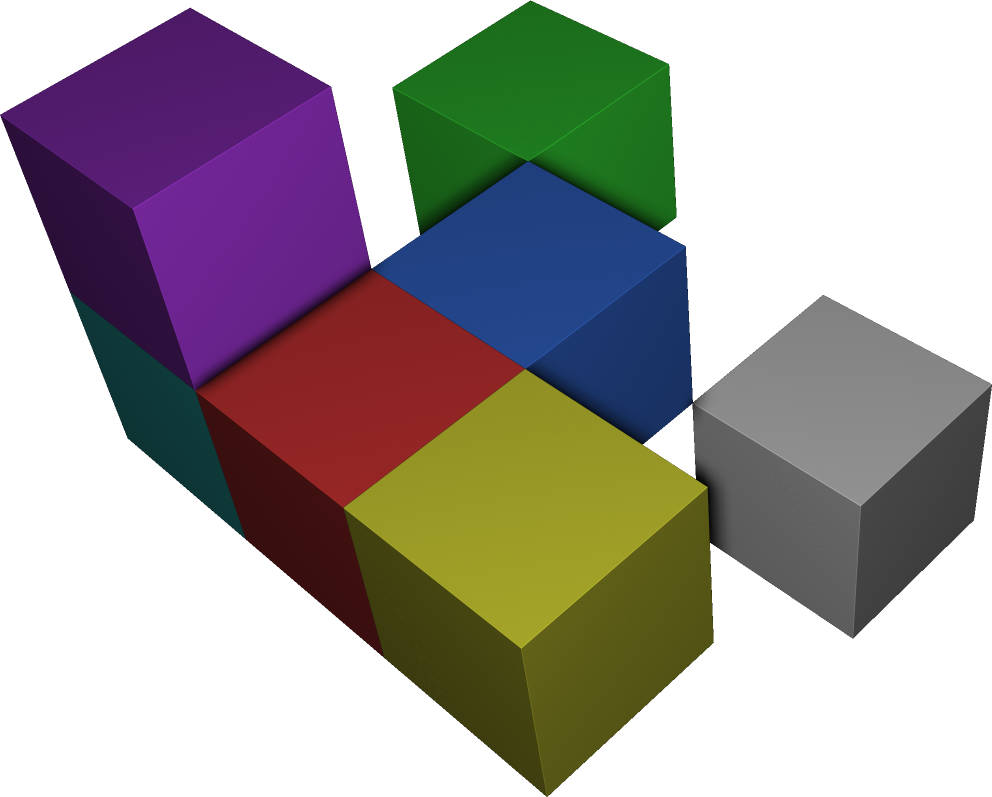}
    \includegraphics[width=0.23\textwidth] {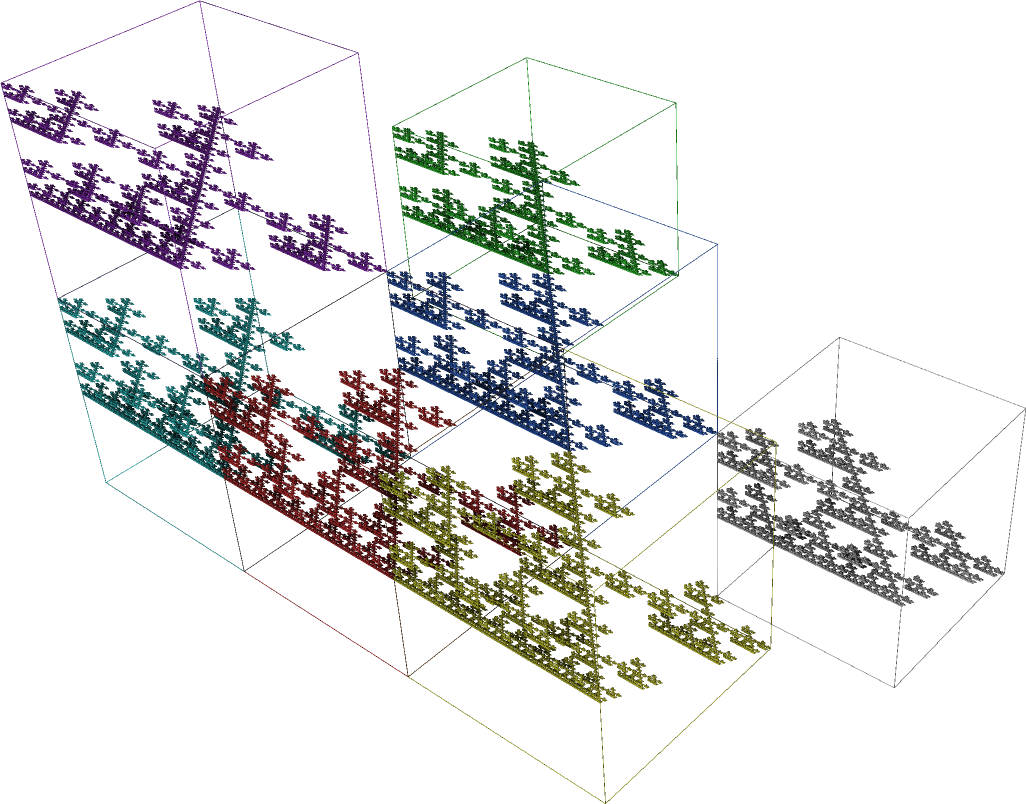} \\
    {\scriptsize$\mD=\{(0,0,1),\ (0,2,2),\ (1,1,1),\ (1,2,2),\ (2,0,0),\ (2,2,1),\ (2,2,2)\}$} & 
    {\scriptsize$\mD=\{(0,0,0),\ (0,2,1),\ (1,1,1),\ (1,2,1),\ (2,0,1),\ (2,2,1),\ (2,2,2)\}$} \\
    \hline
    \includegraphics[width=0.23\textwidth] {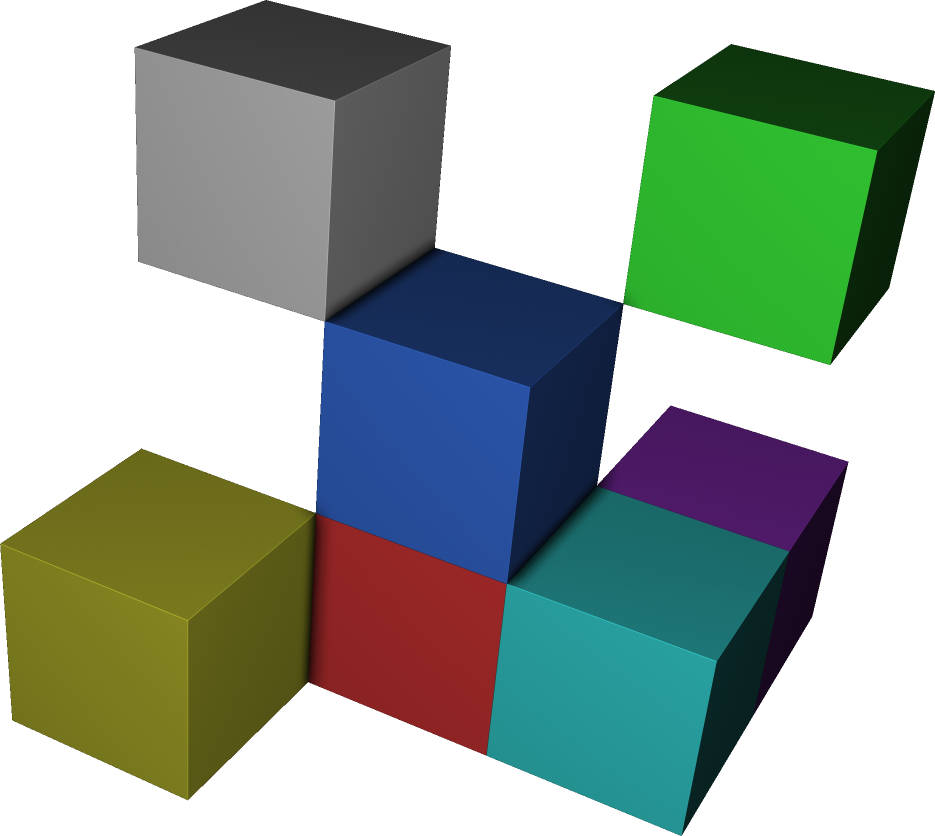}
    \includegraphics[width=0.23\textwidth] {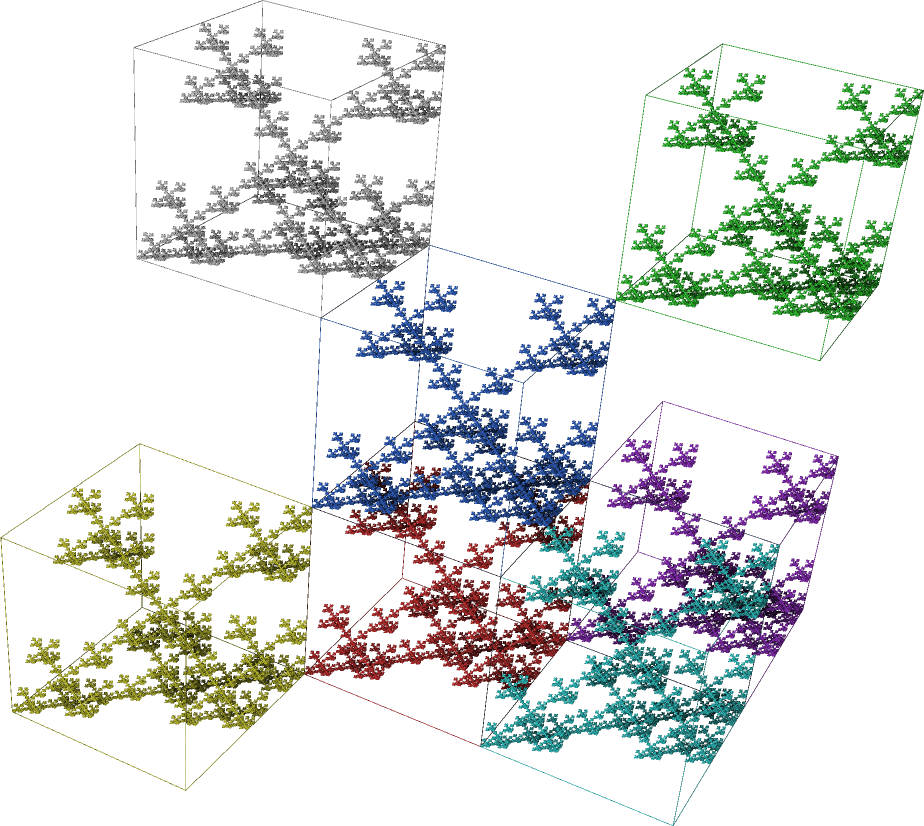} & 
    \includegraphics[width=0.23\textwidth] {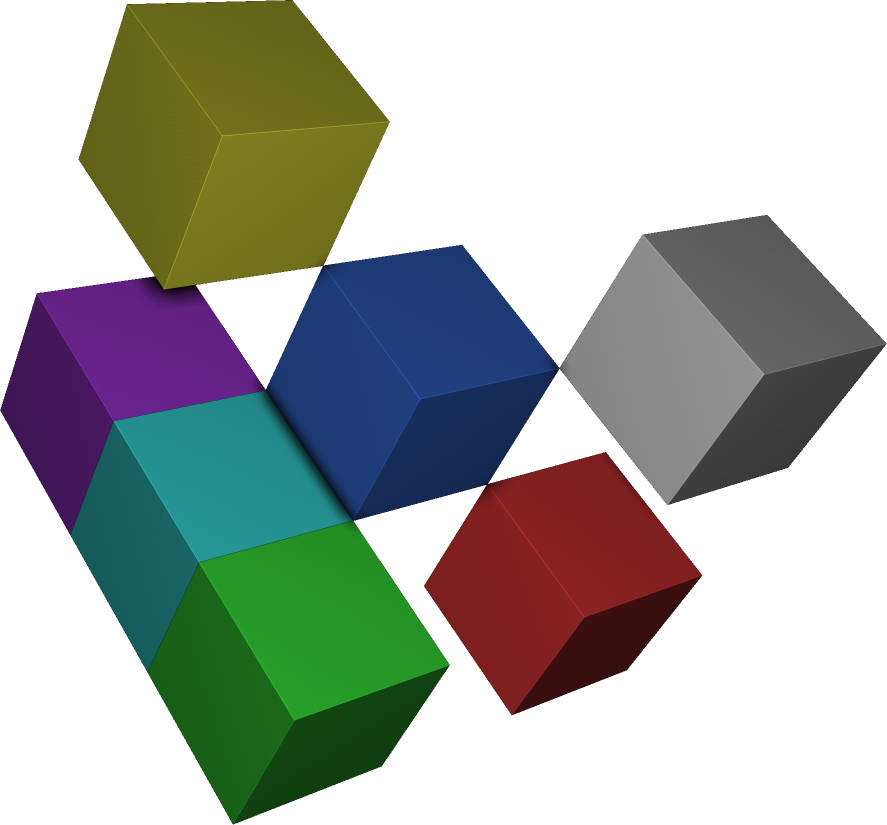}
    \includegraphics[width=0.23\textwidth] {n4_117711105__000_022_111_200_220_221_222__Q.jpg} \\
    {\scriptsize$\mD=\{(0,0,1),\ (0,2,0),\ (1,1,1),\ (1,2,1),\ (2,0,2),\ (2,2,1),\ (2,2,2)\}$} & 
    {\scriptsize$\mD=\{(0,0,0),\ (0,2,2),\ (1,1,1),\ (2,0,0),\ (2,2,0),\ (2,2,1),\ (2,2,2)\}$} \\
    \hline
    \includegraphics[width=0.23\textwidth] {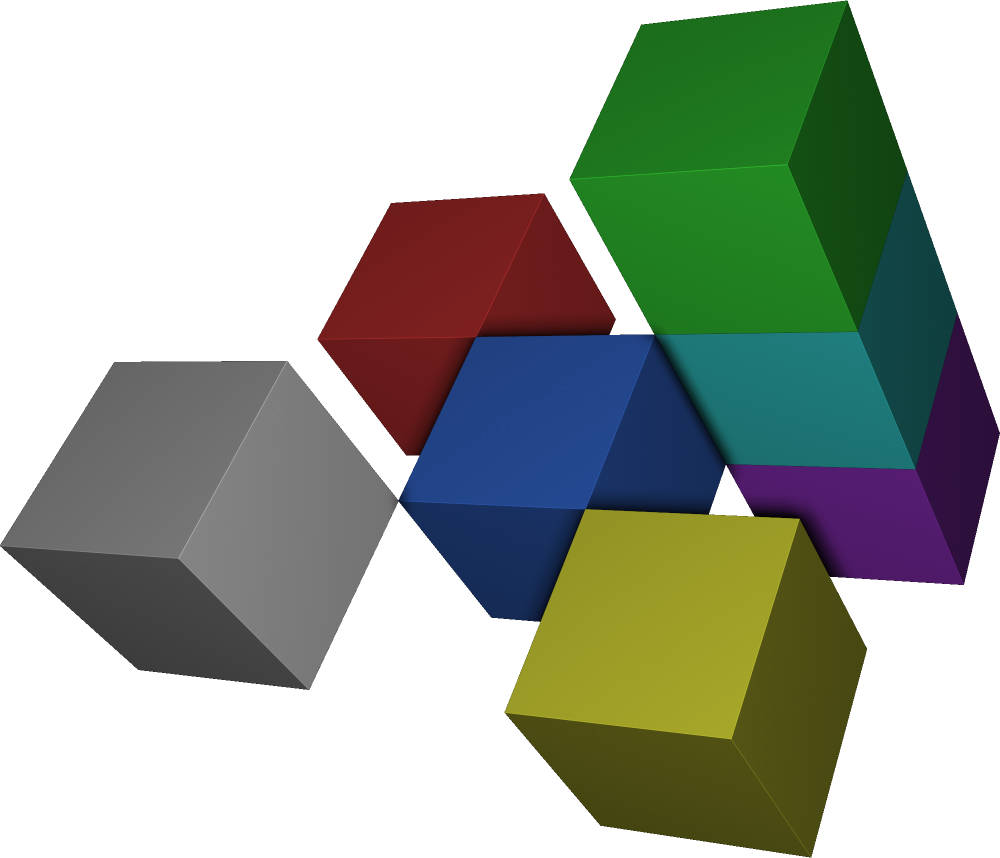}
    \includegraphics[width=0.23\textwidth] {n4_117973121__000_021_111_201_220_221_222__Q.jpg} & \\
    {\scriptsize$\mD=\{(0,0,0),\ (0,2,1),\ (1,1,1),\ (2,0,1),\ (2,2,0),\ (2,2,1),\ (2,2,2)\}$} &\\
    \hline
\end{longtable}

\begin{longtable}{|p{0.48\textwidth}|p{0.48\textwidth}|}
\caption{Non-dendrites of   type 5 ($N=2$) }\label{tab:n5}\\
    \hline
    \multicolumn{2}{|c|}{\includegraphics{nonden5.pdf}} \\
    \hline
    \includegraphics[width=0.23\textwidth] {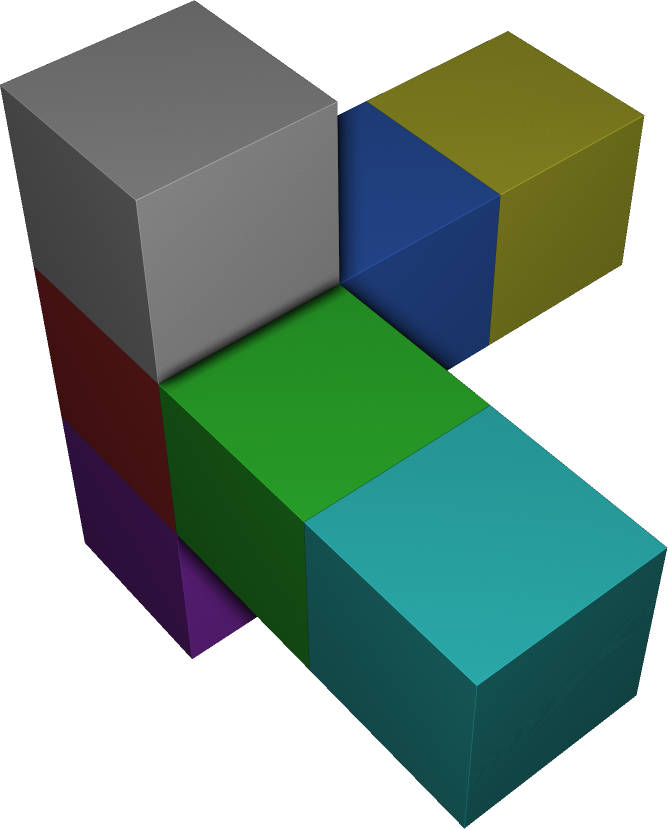}
    \includegraphics[width=0.23\textwidth] {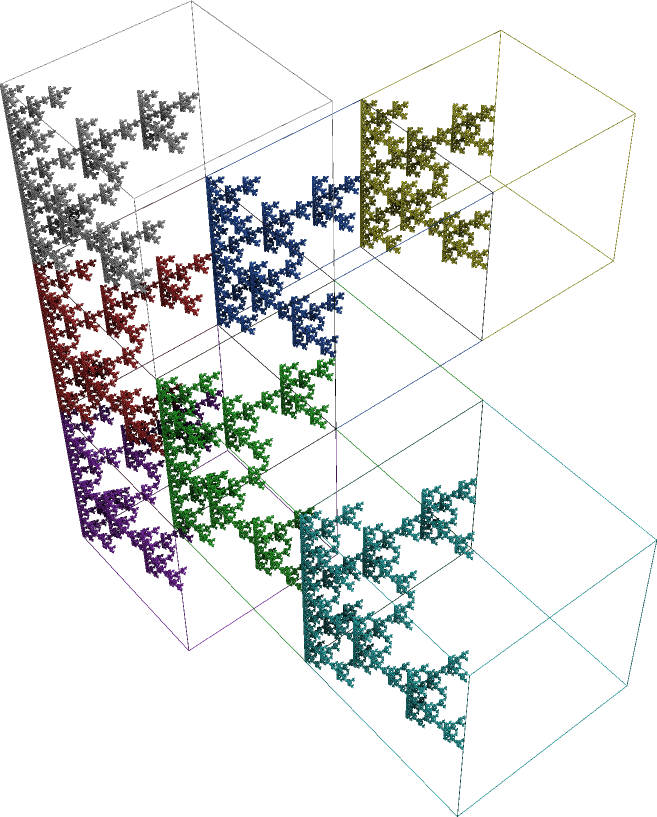} & 
    \includegraphics[width=0.23\textwidth] {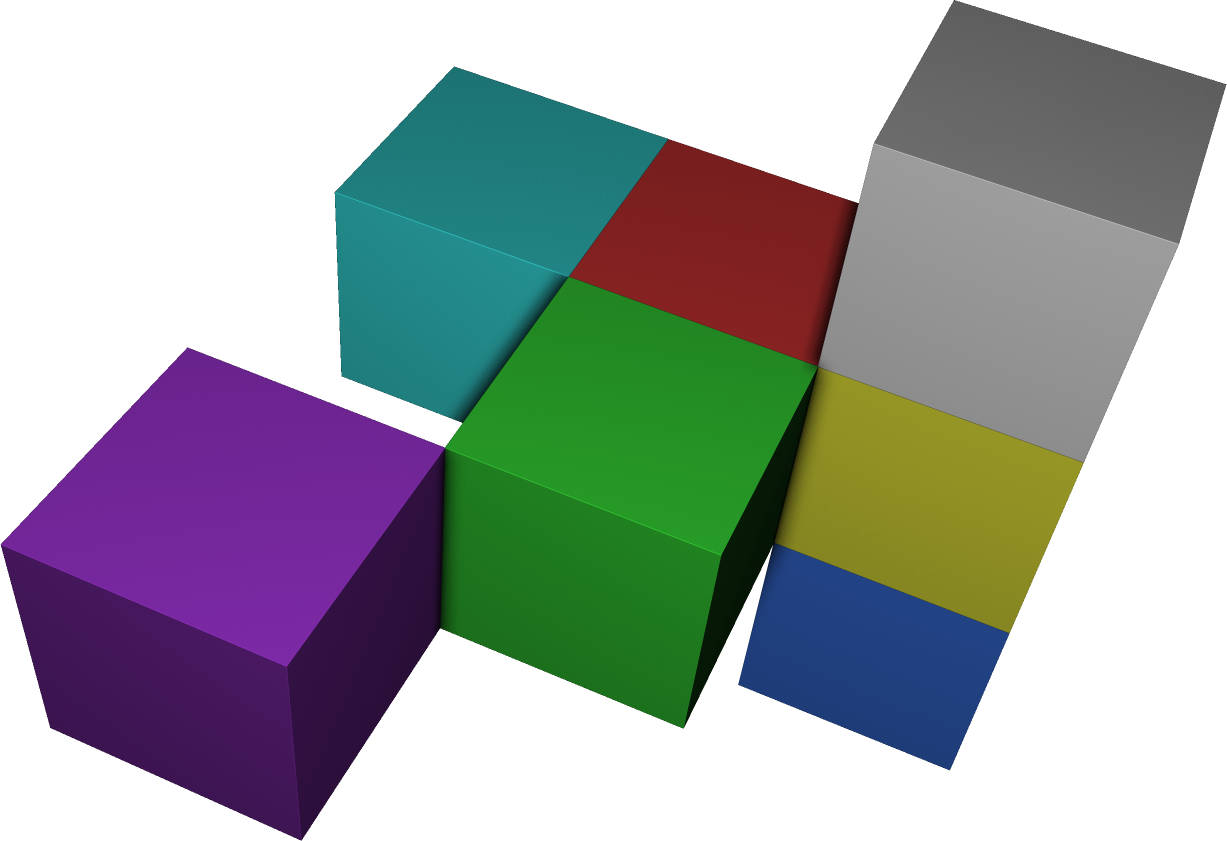}
    \includegraphics[width=0.23\textwidth] {n5_10498121__000_010_020_110_111_210_212__Q.jpg} \\
    {\scriptsize$\mD=\{(0,0,2),\ (1,0,0),\ (1,0,1),\ (1,0,2),\ (1,1,2),\ (1,2,2),\ (2,0,2)\}$} & 
    {\scriptsize$\mD=\{(0,0,0),\ (0,1,0),\ (0,2,0),\ (1,1,0),\ (1,1,1),\ (2,1,0),\ (2,1,2)\}$} \\
    \hline
\end{longtable}

\newpage

\begin{longtable}{|p{0.48\textwidth}|p{0.48\textwidth}|}
\caption{Non-dendrites of   type 6 ($N=9$) }\label{tab:n6}\\
    \hline
    \multicolumn{2}{|c|}{\includegraphics{nonden6.pdf}} \\
    \hline
    \includegraphics[width=0.23\textwidth] {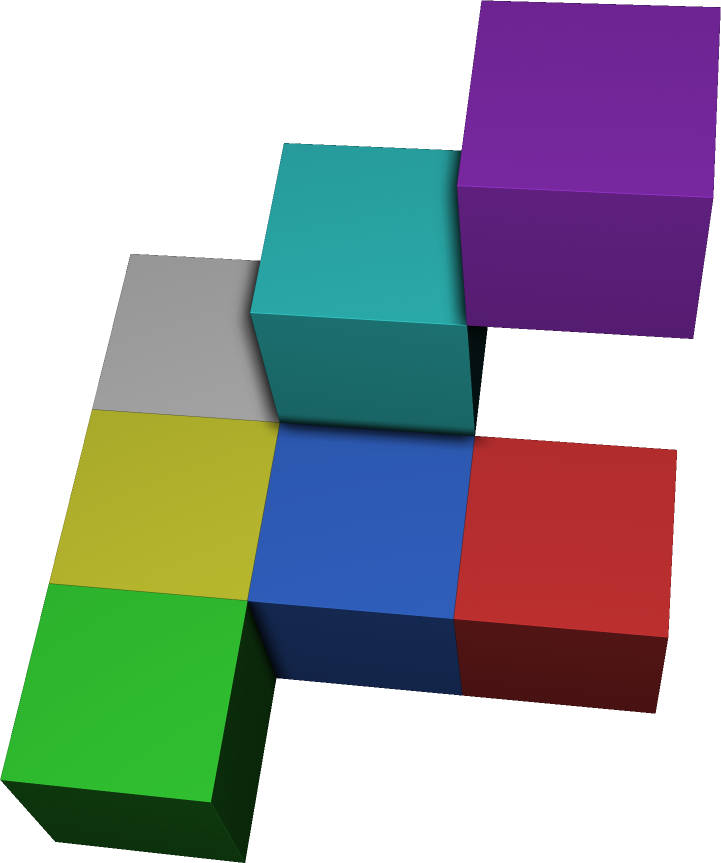}
    \includegraphics[width=0.23\textwidth] {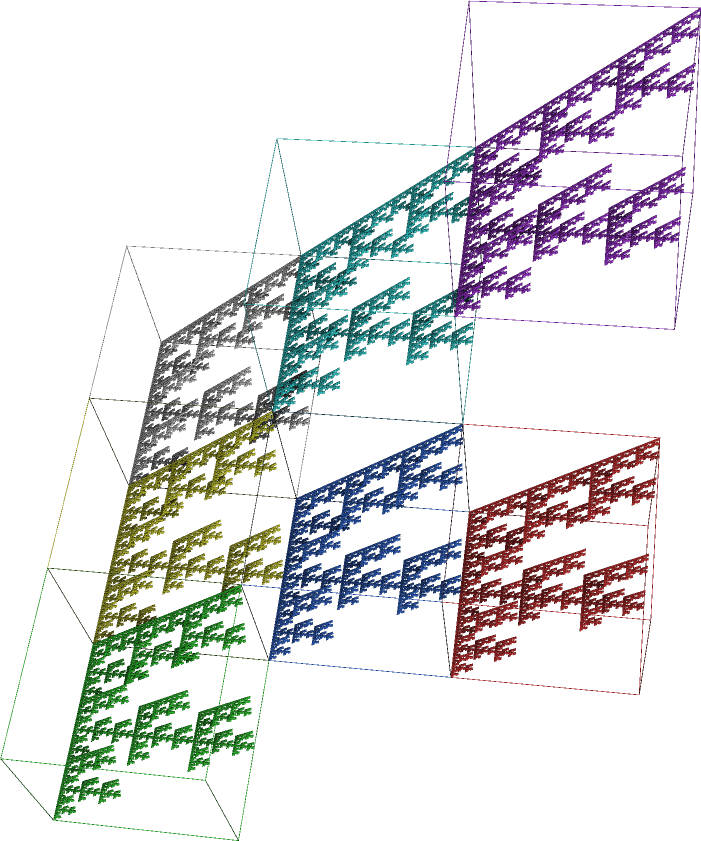} & 
    \includegraphics[width=0.23\textwidth] {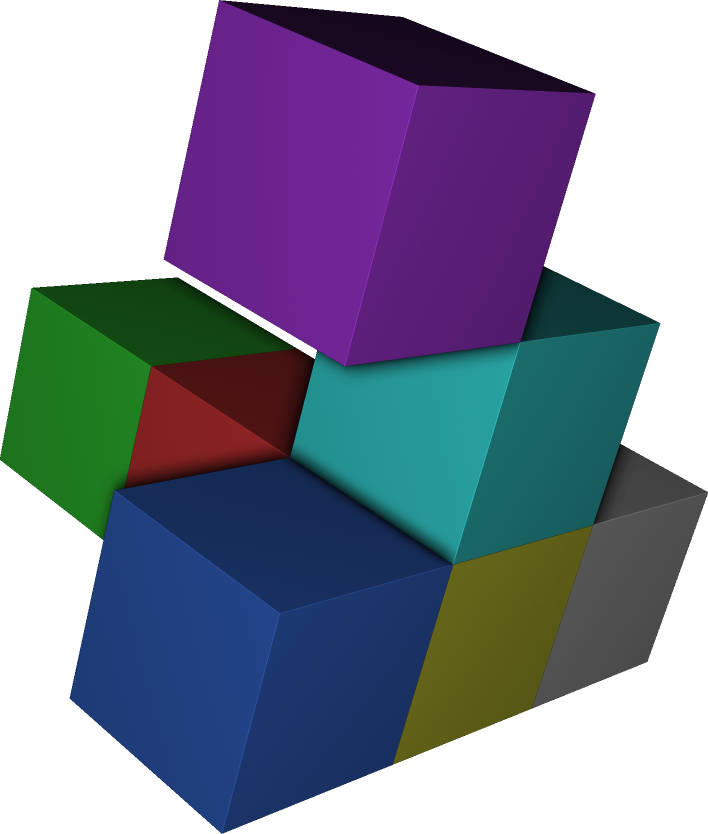}
    \includegraphics[width=0.23\textwidth] {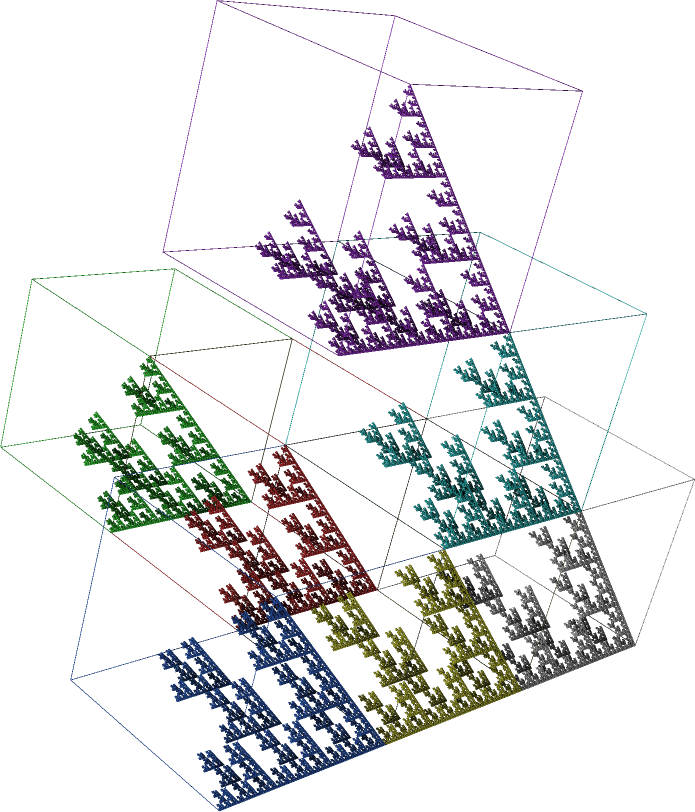} \\
    {\scriptsize$\mD=\{(0,0,0),\ (0,1,0),\ (0,1,1),\ (0,1,2),\ (0,2,0),\ (1,0,1),\ (2,0,2)\}$} & 
    {\scriptsize$\mD=\{(0,0,0),\ (0,0,1),\ (0,0,2),\ (0,1,1),\ (0,2,1),\ (1,0,1),\ (2,0,2)\}$} \\
    \hline
    \includegraphics[width=0.23\textwidth] {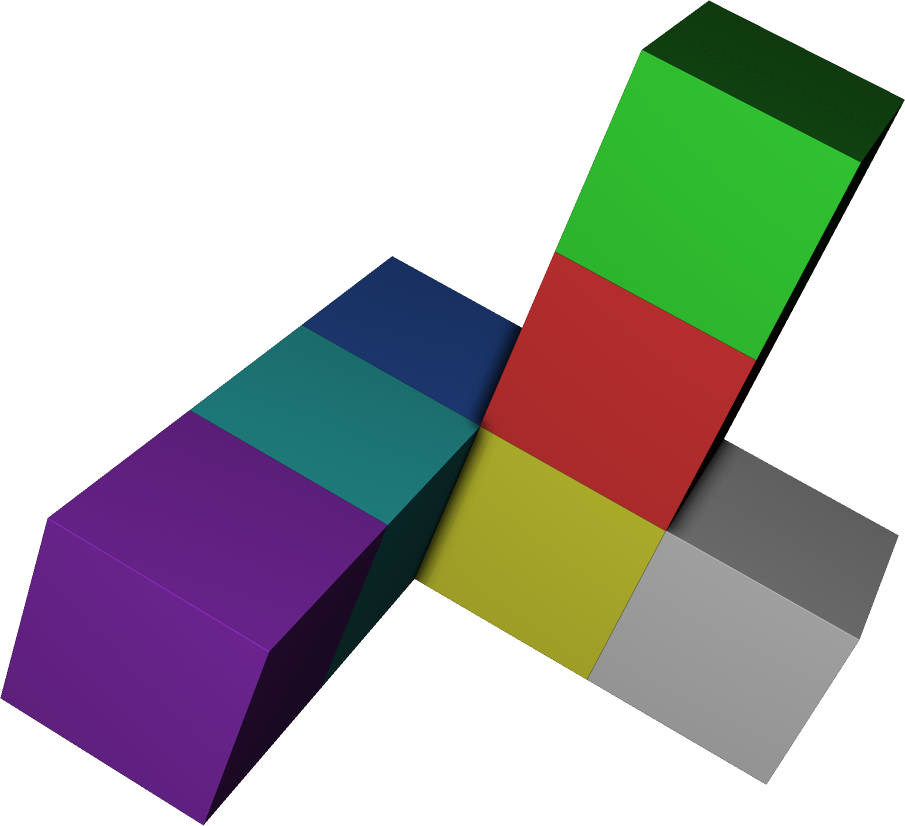}
    \includegraphics[width=0.23\textwidth] {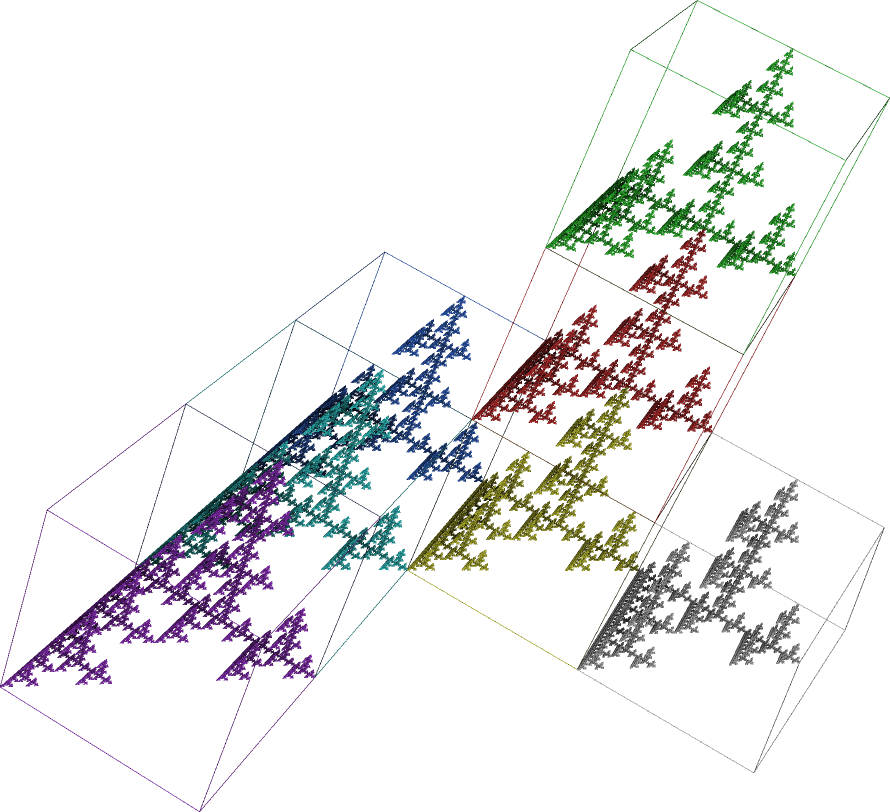} & 
    \includegraphics[width=0.23\textwidth] {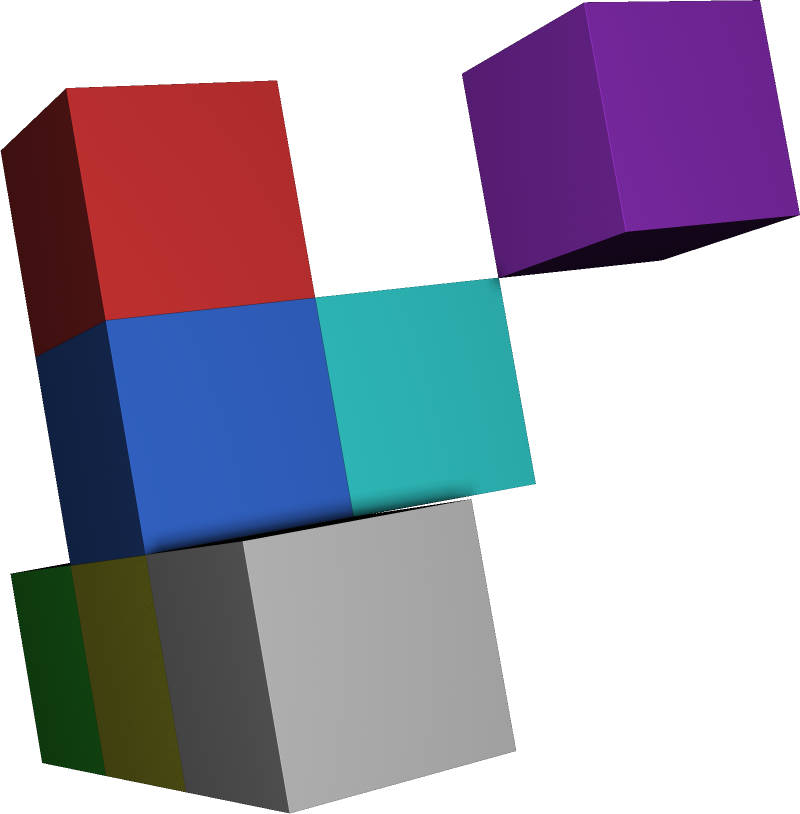}
    \includegraphics[width=0.23\textwidth] {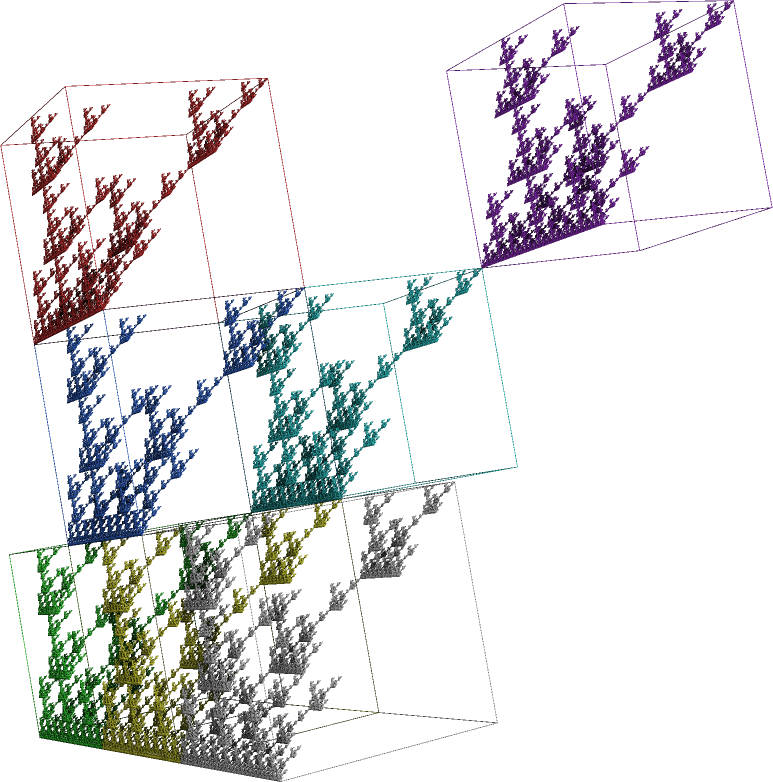} \\
    {\scriptsize$\mD=\{(0,0,0),\ (0,0,1),\ (0,0,2),\ (0,1,1),\ (0,2,1),\ (1,0,2),\ (2,0,2)\}$} & 
    {\scriptsize$\mD=\{(0,0,0),\ (0,1,0),\ (0,1,1),\ (0,1,2),\ (0,2,0),\ (1,1,1),\ (2,0,2)\}$} \\
    \hline
    \includegraphics[width=0.23\textwidth] {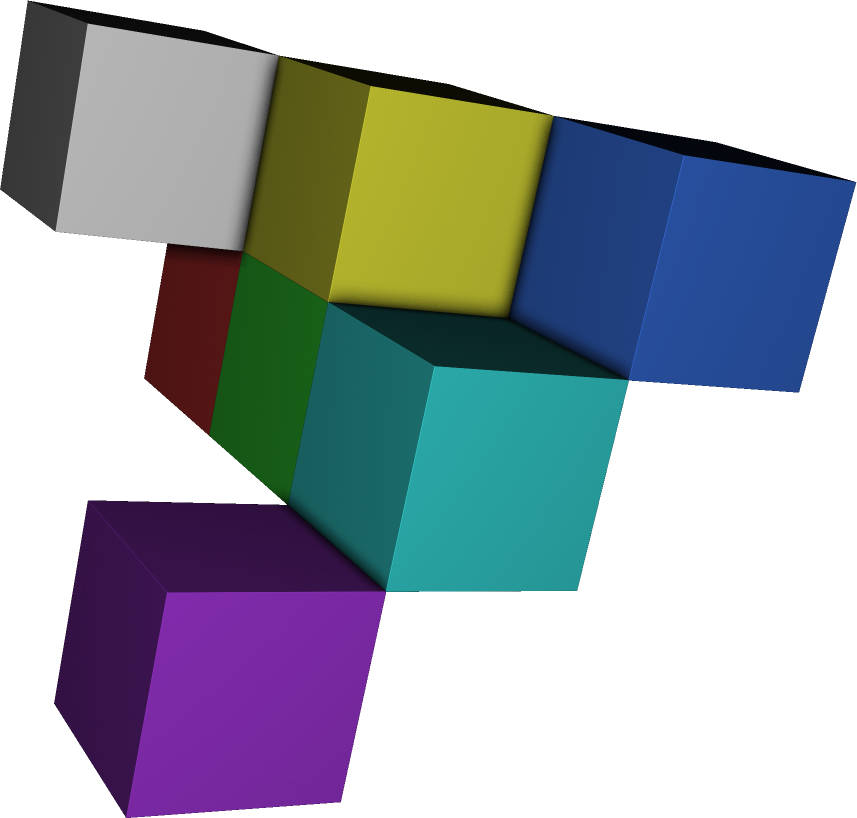}
    \includegraphics[width=0.23\textwidth] {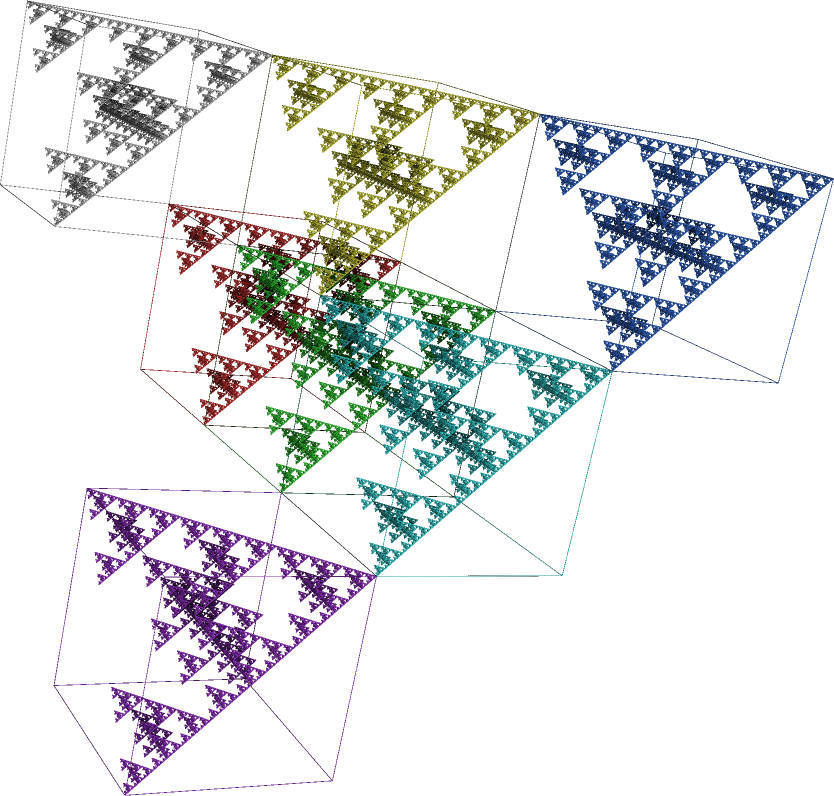} & 
    \includegraphics[width=0.23\textwidth] {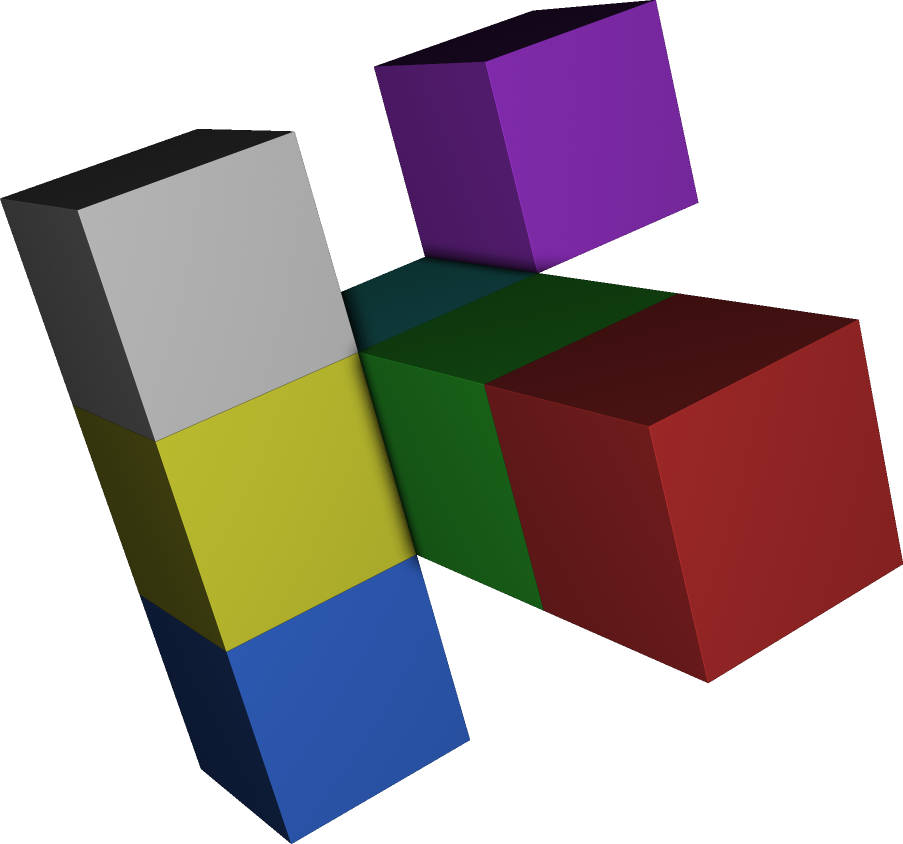}
    \includegraphics[width=0.23\textwidth] {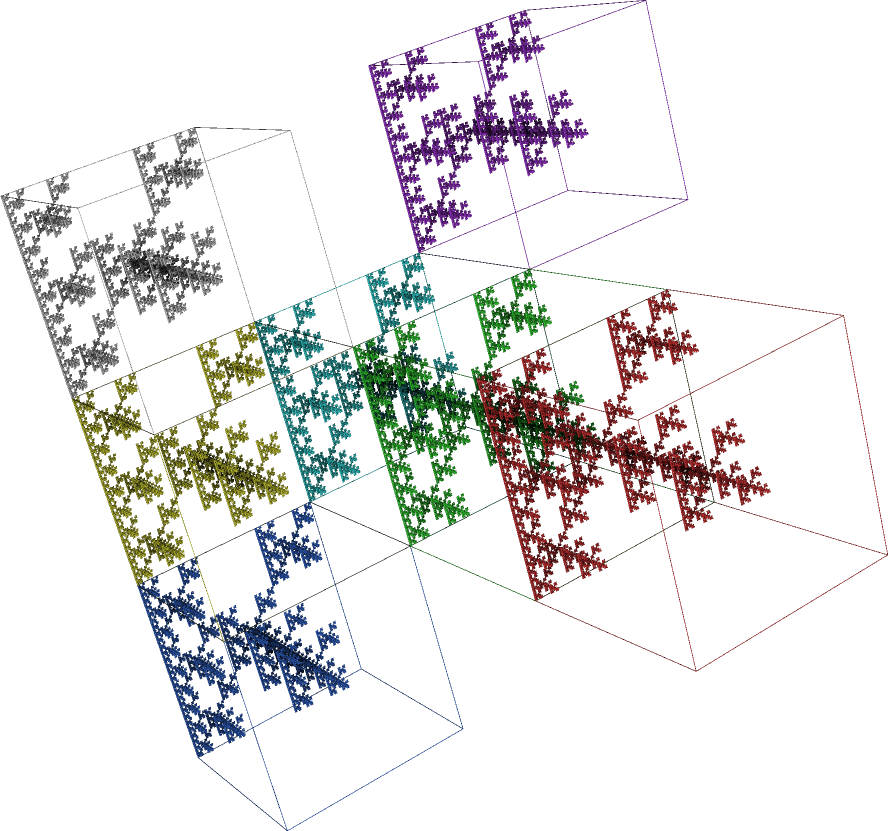} \\
    {\scriptsize$\mD=\{(0,0,0),\ (0,1,1),\ (0,2,2),\ (1,1,0),\ (1,1,1),\ (1,1,2),\ (2,0,2)\}$} & 
    {\scriptsize$\mD=\{(0,0,2),\ (0,1,2),\ (0,2,2),\ (1,1,0),\ (1,1,1),\ (1,1,2),\ (2,0,2)\}$} \\
    \hline
    \includegraphics[width=0.23\textwidth] {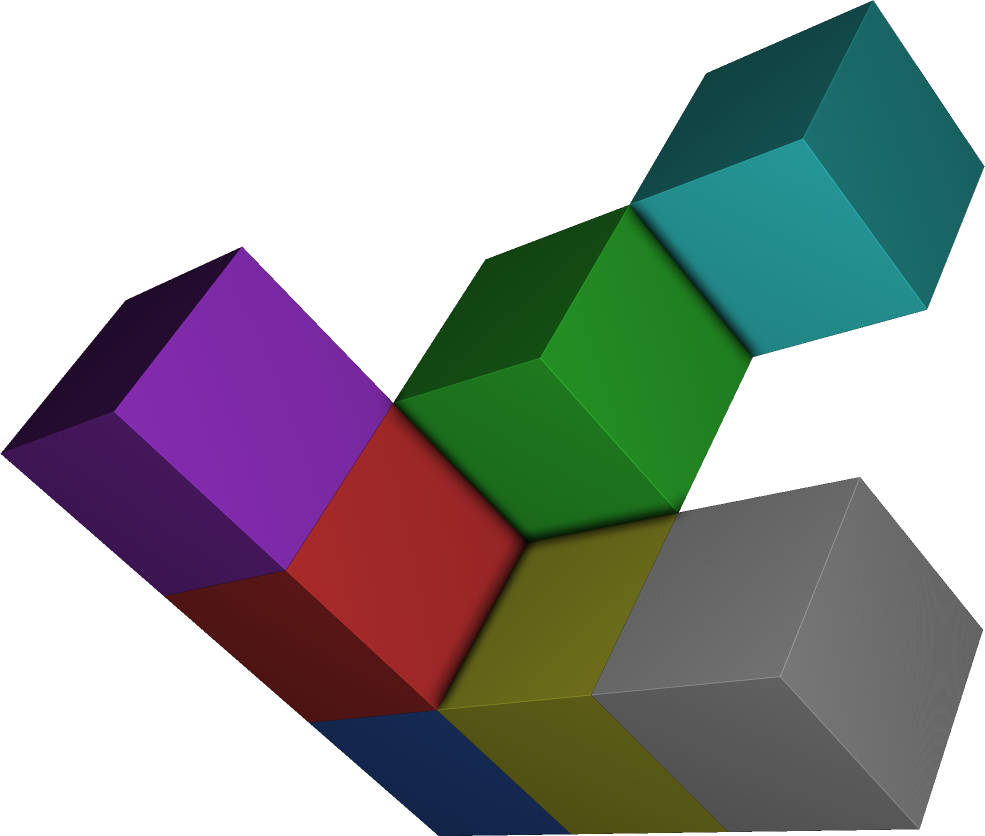}
    \includegraphics[width=0.23\textwidth] {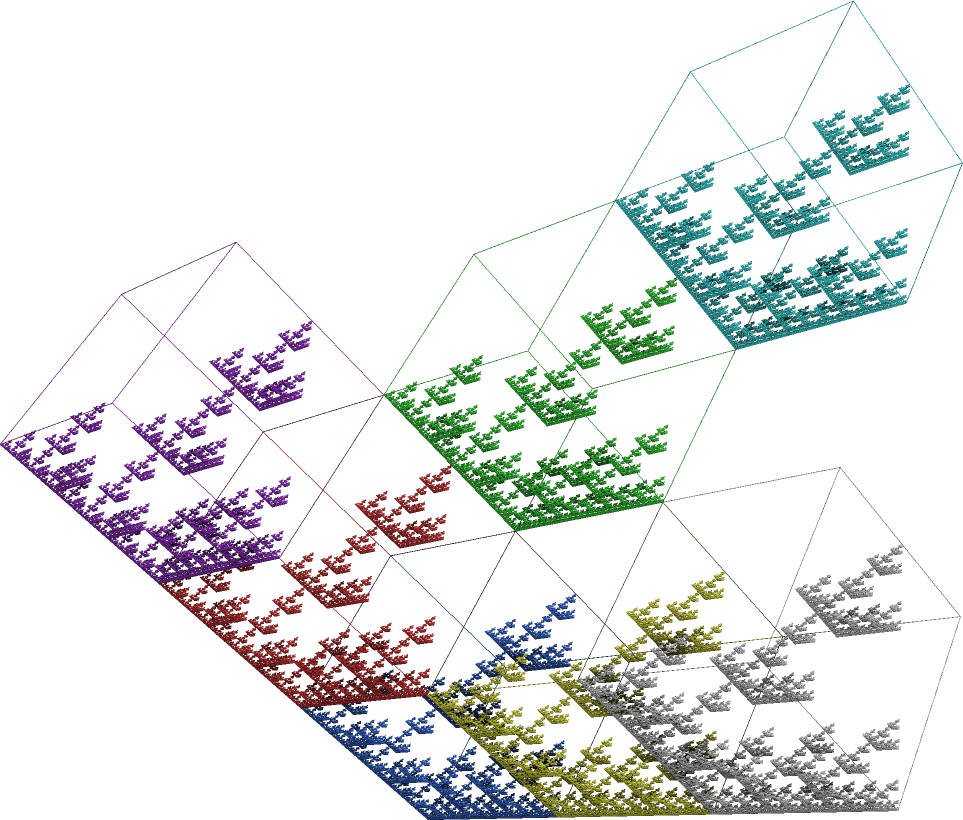} & 
    \includegraphics[width=0.23\textwidth] {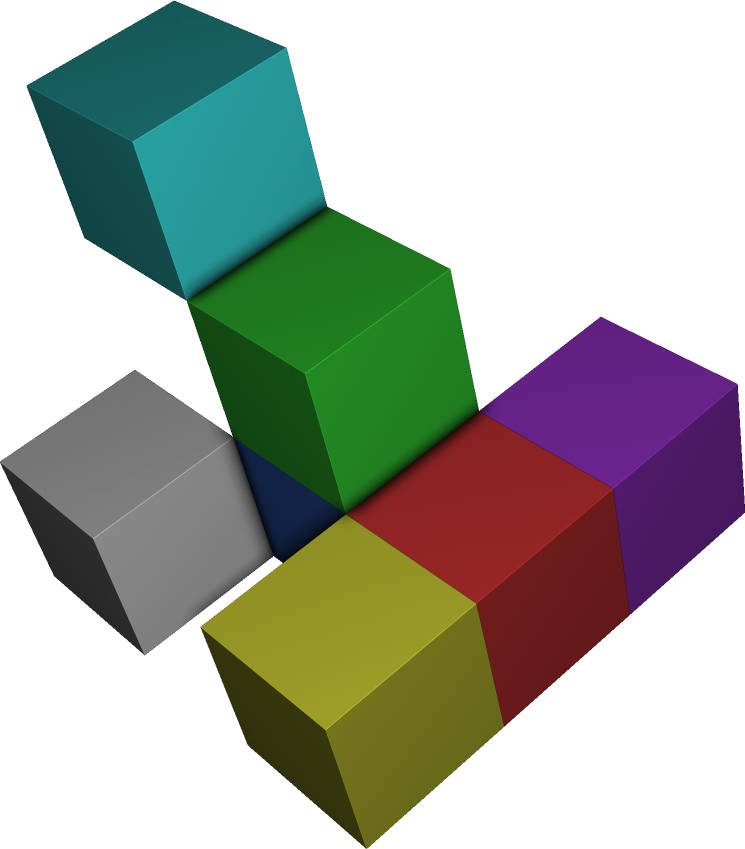}
    \includegraphics[width=0.23\textwidth] {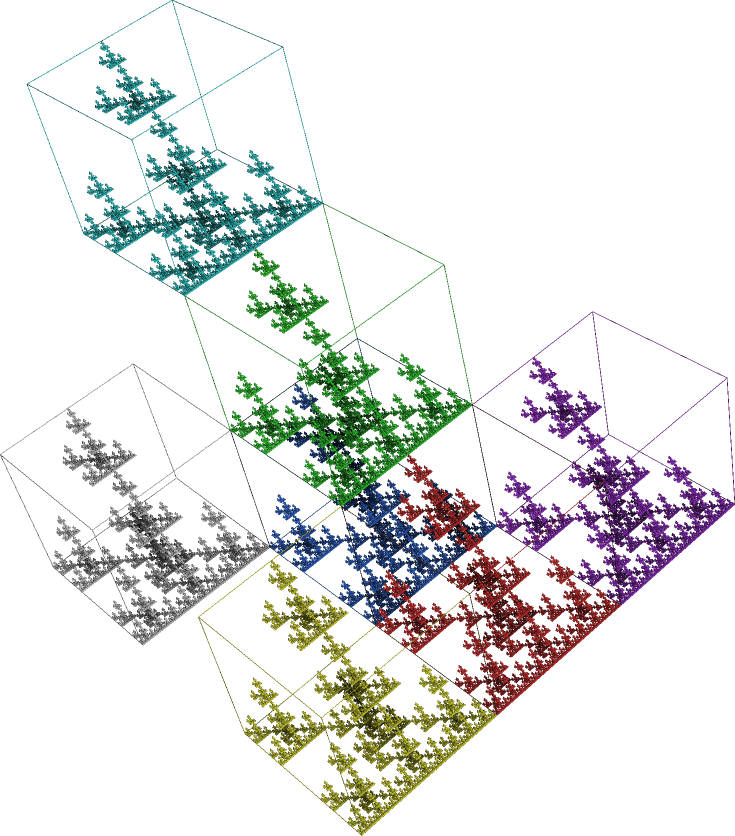} \\
    {\scriptsize$\mD=\{(0,0,0),\ (0,0,1),\ (0,0,2),\ (1,0,2),\ (1,1,1),\ (1,2,0),\ (2,0,2)\}$} & 
    {\scriptsize$\mD=\{(0,0,0),\ (0,0,2),\ (1,0,1),\ (1,0,2),\ (1,1,1),\ (1,2,0),\ (2,0,2)\}$} \\
    \hline
    \includegraphics[width=0.23\textwidth] {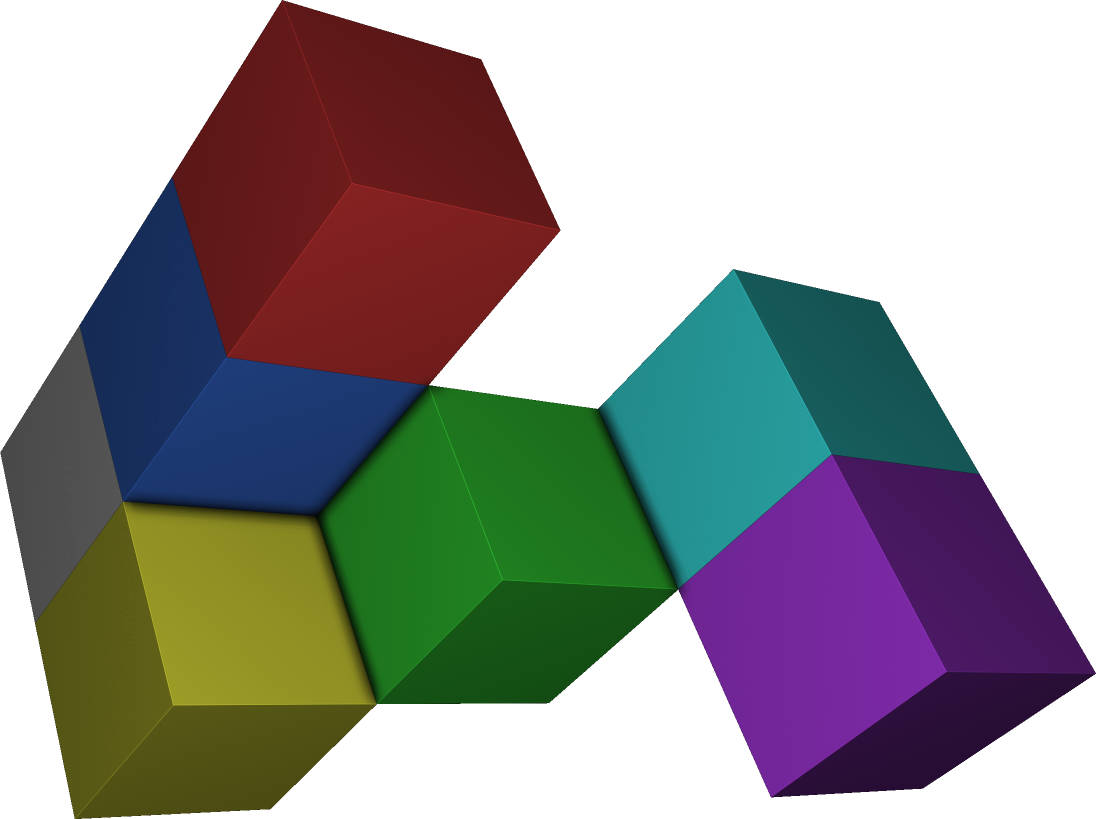}
    \includegraphics[width=0.23\textwidth] {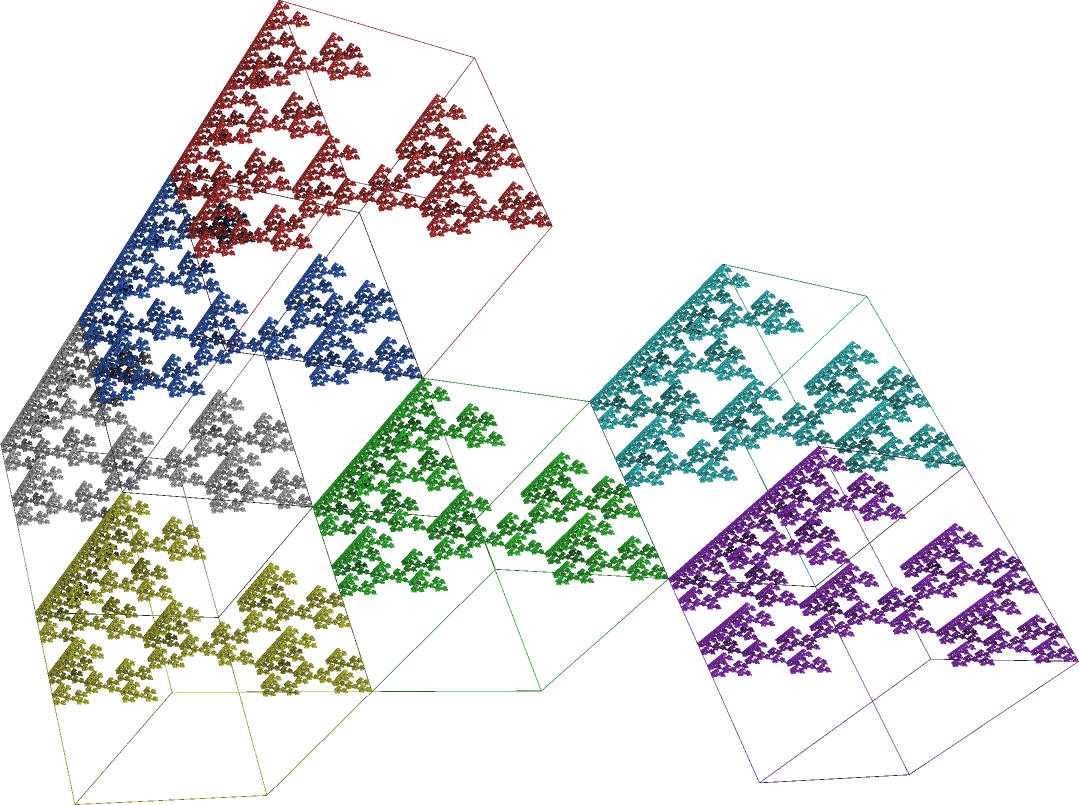} &\\
    {\scriptsize$\mD=\{(0,0,0),\ (0,0,1),\ (0,1,0),\ (0,2,0),\ (1,1,1),\ (2,2,1),\ (2,2,2)\}$} &\\
    \hline
\end{longtable}

\begin{longtable}{|p{0.48\textwidth}|p{0.48\textwidth}|}
\caption{Non-dendrites of   type 7 ($N=3$) }\label{tab:n7}\\
    \hline
    \multicolumn{2}{|c|}{\includegraphics{nonden7.pdf}} \\
    \hline
    \includegraphics[width=0.23\textwidth] {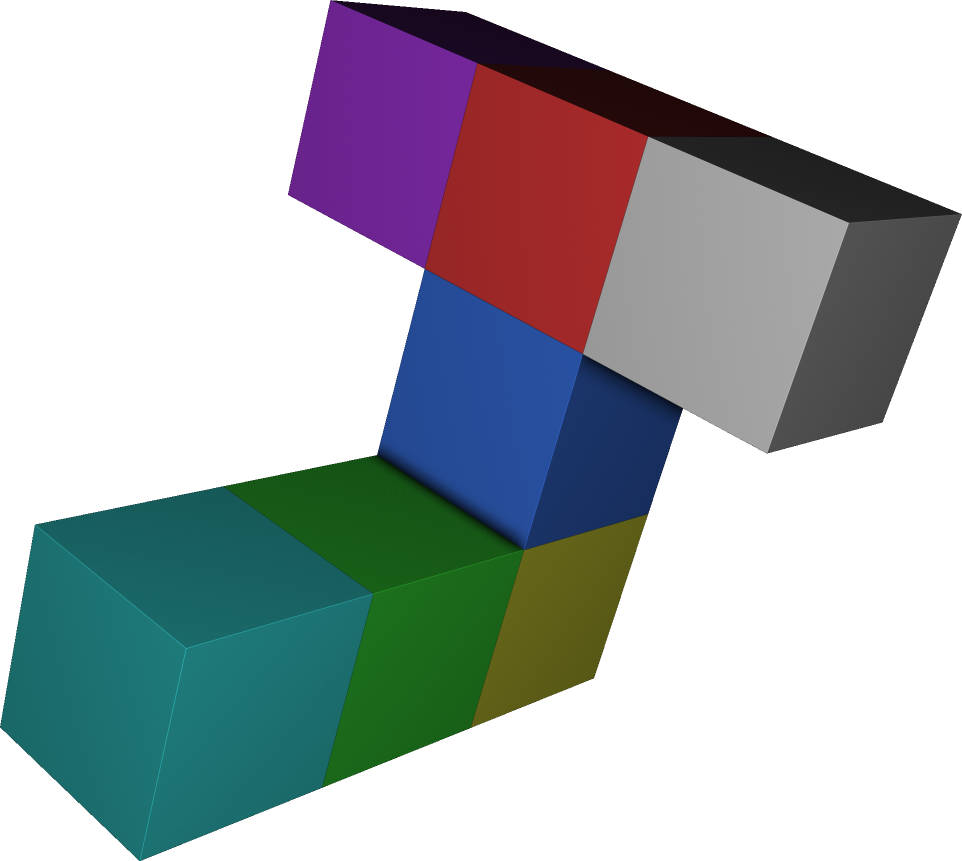}
    \includegraphics[width=0.23\textwidth] {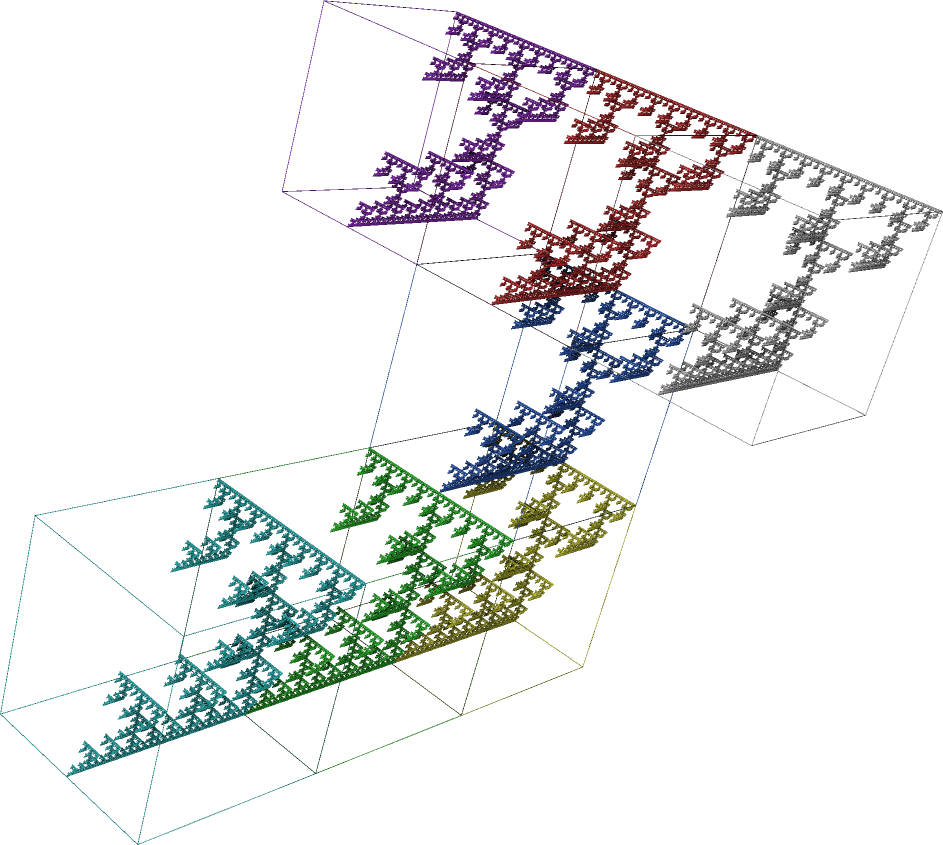} & 
    \includegraphics[width=0.23\textwidth] {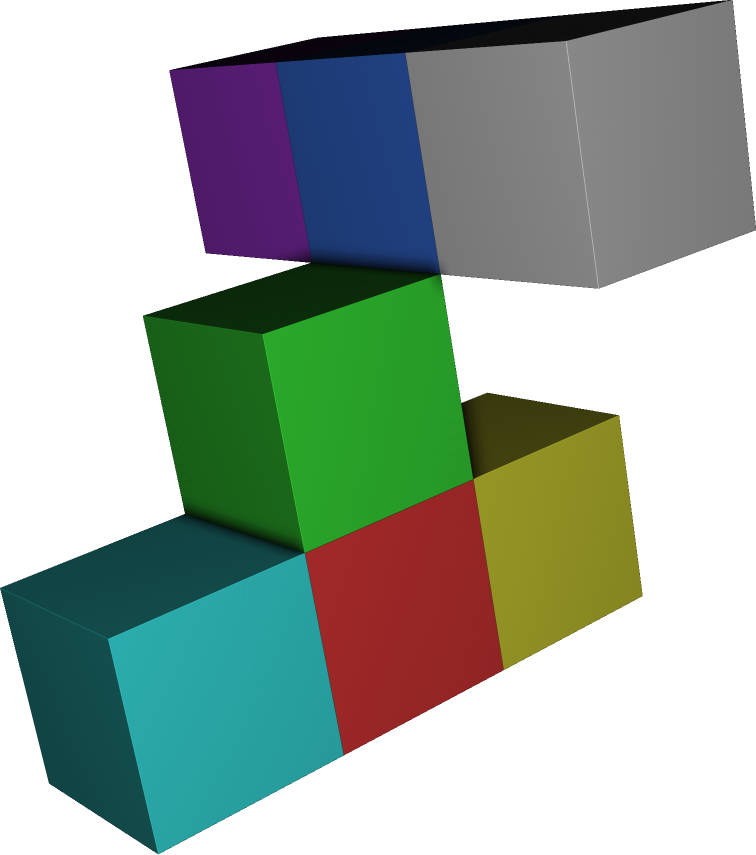}
    \includegraphics[width=0.23\textwidth] {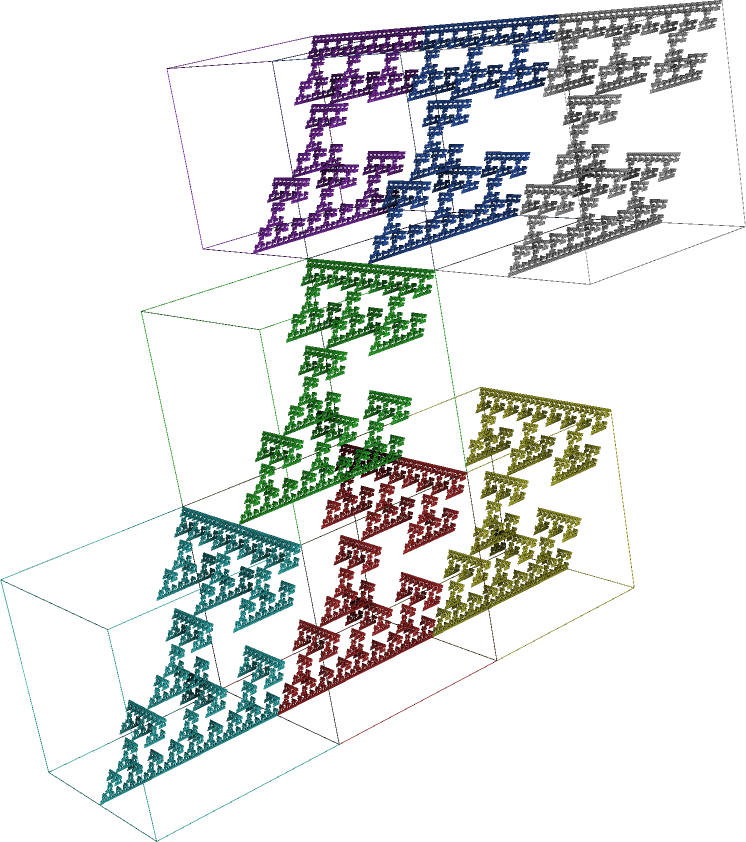} \\
    {\scriptsize$\mD=\{(0,0,2),\ (1,0,0),\ (1,0,1),\ (1,0,2),\ (1,1,0),\ (1,2,0),\ (2,0,2)\}$} & 
    {\scriptsize$\mD=\{(0,0,2),\ (1,0,0),\ (1,0,2),\ (1,1,0),\ (1,1,1),\ (1,2,0),\ (2,0,2)\}$} \\
    \hline
    \includegraphics[width=0.23\textwidth] {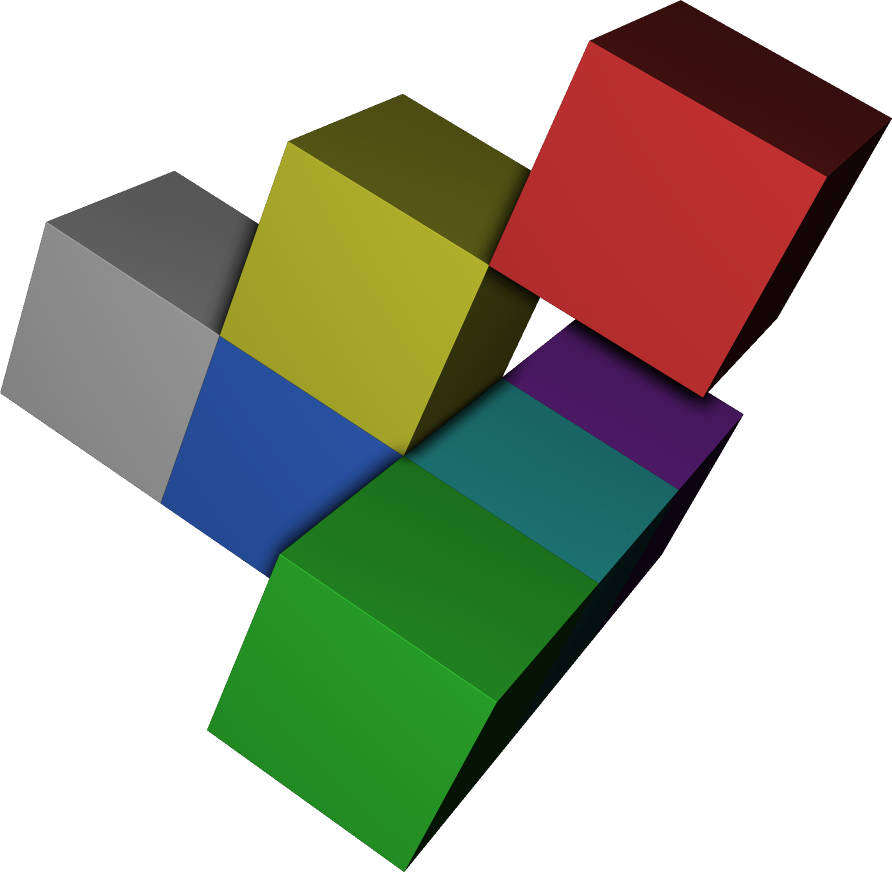}
    \includegraphics[width=0.23\textwidth] {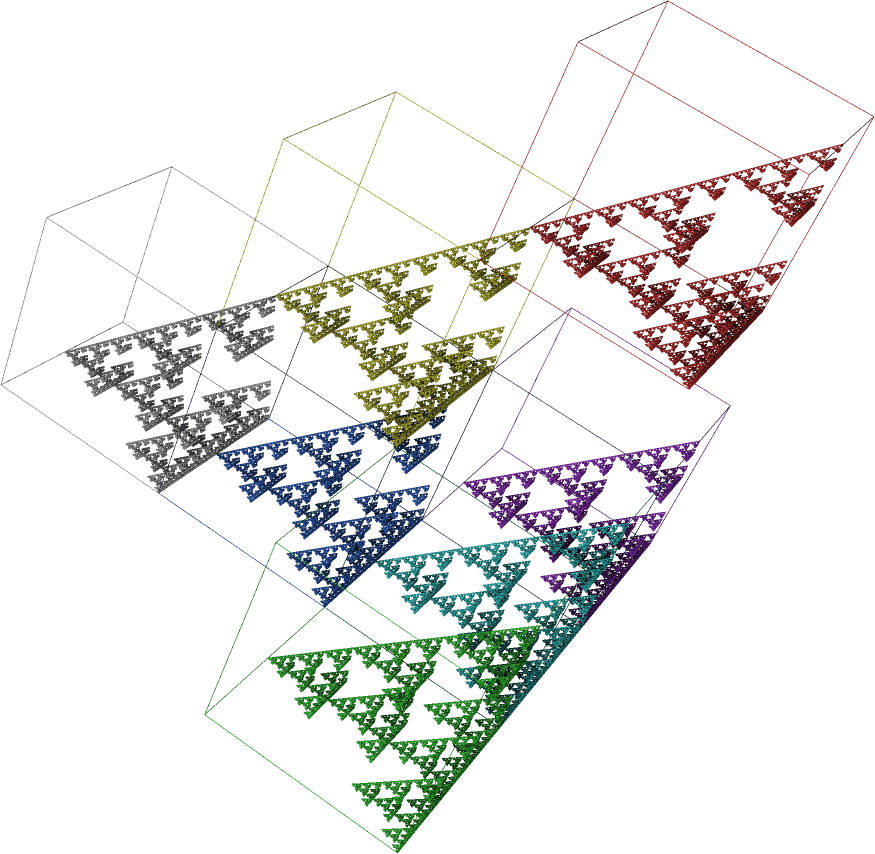} &\\
    {\scriptsize$\mD=\{(0,2,1),\ (1,1,1),\ (1,2,1),\ (2,0,1),\ (2,2,0),\ (2,2,1),\ (2,2,2)\}$} &\\
    \hline
\end{longtable}

\end{document}